\documentclass[12pt]{article}
\usepackage{amsmath,amssymb,amsfonts,amsthm}
\usepackage{fullpage}
\usepackage{enumitem}
\usepackage{caption}
\usepackage{shuffle}
\usepackage{subcaption}
\usepackage{graphicx}
\usepackage{slashed}
\usepackage{booktabs}
\usepackage{mciteplus}
\usepackage{xcolor}
\usepackage{bm}
\usepackage{fancyhdr, dsfont}
\usepackage{float}
\usepackage{url}
\usepackage{setspace}
\usepackage[nosort]{cite}
\usepackage{mathtools}
\allowdisplaybreaks
\input xypic 

\usepackage{color}

\usepackage{hyperref}
\definecolor{darkred}{rgb}{0.8,0.1,0.1}
\hypersetup{colorlinks=true, linkcolor=darkred, urlcolor=darkred, citecolor=blue, linktoc=page}

\usepackage{tikz}
\usetikzlibrary{calc} 
\usetikzlibrary{patterns} 
\usetikzlibrary{decorations.pathreplacing} 
\usetikzlibrary{decorations.markings} 
\usetikzlibrary{decorations.pathmorphing} 
\usetikzlibrary{positioning}
\usetikzlibrary{arrows.meta}

\newcommand{\Q}{{\mathbb{Q}}}
\newcommand{\I}{{\mathcal{I}}}
\newcommand{\F}{{\mathcal{F}}}
\newcommand{\Z}{{\mathcal{Z}}}
\newcommand{\MZ}{{\mathcal{MZ}}}
\newcommand{\FZ}{{\mathcal{FZ}}}
\newcommand{\U}{{\mathcal{U}}}

\newcommand{\ds}{{\mathfrak{ds}}}
\newcommand{\dmr}{{\mathfrak{dmr}}}
\newcommand{\grt}{{\mathfrak{grt}}}
\newcommand{\z}{{\mathfrak{z}}}
\newcommand{\fz}{{\mathfrak{fz}}}
\newcommand{\mz}{{\mathfrak{mz}}}
\newcommand{\beq}{\begin{equation}}
\newcommand{\eeq}{\end{equation}}

\def\nn{\nonumber}

\let\Im\relax

\DeclareMathOperator{\Im}{Im}

\newcommand{\nwc}{\newcommand}
\nwc{\ba}  {\begin{array}}
\nwc{\ea}  {\end{array}}
\nwc{\bdm} {\begin{displaymath}}
\nwc{\edm} {\end{displaymath}}
\nwc{\bda} {\bdm\ba{lcl}} 
\nwc{\eda} {\ea\edm}
\nwc{\bc}  {\begin{center}}
\nwc{\ec}  {\end{center}}
\nwc{\bmat}{\left(\ba}
\nwc{\emat}{\ea\right)}
\nwc{\nnn} {\nonumber \vspace{.2cm} \\ }
\nwc{\ra}  {\rightarrow}
\nwc{\lra} {\longrightarrow}
\nwc{\p} {\partial}
\nwc{\rcr} {\nabla_{\rm alt}}
\nwc{\barrcr} {\overline{\nabla_{\rm alt}}}
\nwc{\ep} {\epsilon}

\def\ad{\mathrm{ad}}
\def\Ad{\mathrm{Ad}}

\def\BF{{\rm BF}}

\newcommand{\eee}{\mathrm{e}}
\newcommand{\hhh}{\mathrm{h}}

\newcommand{\candh}{h}
\newcommand{\sigkey}{\sigma^{\rm key}}

\newcommand{\sigalg}{{\cal S}}

\usepackage{stackengine}

\definecolor{cadmiumgreen}{RGB}{60,142,23}
\definecolor{kleinblue}{RGB}{0,46,167}
\definecolor{venetianred}{RGB}{192,6,21}
\definecolor{forestgreen}{RGB}{34,139,34}

\usepackage{titlesec}
\titleformat*{\section}{\large \bfseries}

\newcommand{\stuffle}{*}
\newcommand{\circp}{\diamond}

\usepackage{thmtools}

\declaretheorem[numberwithin=subsection,name=Theorem]{thm}
\declaretheorem[name=Lemma, sharenumber=thm]{lemma}
\declaretheorem[name=Proposition, sharenumber=thm]{prop}

\declaretheorem[name=Corollary, sharenumber=thm]{cor}

\theoremstyle{definition}
\declaretheorem[name=Definition, sharenumber=thm]{dfn}

\declaretheorem[name=Remark, sharenumber=thm]{rmk}

\renewcommand{\qed}{\hfill{$\square$}\\[3mm]}

\hypersetup{pageanchor=false}
\begin{document}

 {\flushright  
 UUITP-16/24\\}

\begin{center}

\vspace{-2mm}

{\bf {\Large \sc Canonicalizing zeta generators: \\[2mm] genus zero and genus one }}

\vspace{2.5mm}

\normalsize
{\large  Daniele Dorigoni$^1$,
Mehregan Doroudiani$^2$,
Joshua Drewitt$^3$, \\ \vskip 0.45 em
Martijn Hidding$^{4,5}$, 
Axel Kleinschmidt$^{2,7}$,
Oliver Schlotterer$^{4,6}$, \\ \vskip 0.45 em
Leila Schneps$^8$
and Bram Verbeek$^{4}$}

\vspace{3mm}
${}^1${\it Centre for Particle Theory \& Department of Mathematical Sciences\\
Durham University, Lower Mountjoy, Stockton Road, Durham DH1 3LE, UK}
\vskip 0.4 em
${}^2${\it Max-Planck-Institut f\"{u}r Gravitationsphysik (Albert-Einstein-Institut)\\
Am M\"{u}hlenberg 1, 14476 Potsdam, Germany}
\vskip 0.4 em
${}^3${\it School of Mathematics, University of Bristol, Queens Road, Bristol, BS8 1QU, UK}
\vskip 0.4 em
${}^4${\it Department of Physics and Astronomy, Uppsala University, 75108 Uppsala, Sweden}
\vskip 0.4 em
${}^5${\it Institute for Theoretical Physics, ETH Zurich, 8093 Z\"urich, Switzerland}
\vskip 0.4 em
${}^6${\it Department of Mathematics, Centre for Geometry and Physics, \\Uppsala University, 75106 Uppsala, Sweden}
\vskip 0.4 em
${}^7${\it International Solvay Institutes, 
ULB-Campus Plaine CP231, 1050 Brussels, Belgium}
\vskip 0.4 em
${}^8${\it CNRS, Sorbonne Universit\'e, Campus Pierre et Marie Curie\\
4 place Jussieu, 75005 Paris, France}

\vspace{4mm}

\hrule

\vspace{4mm}

\begin{tabular}{p{150mm}}
Zeta generators are derivations associated with odd Riemann zeta values that act freely on the Lie algebra of the fundamental group of Riemann surfaces with marked points. The genus-zero incarnation of zeta generators are Ihara derivations of certain Lie polynomials in two generators that can be obtained from the Drinfeld associator. We characterize a canonical choice of these polynomials, together with their non-Lie counterparts at even degrees $w\geq 2$, through the action of the dual space of formal and motivic multizeta values.
Based on these canonical polynomials, we propose a canonical isomorphism that maps motivic multizeta values into the $f$-alphabet. 

The canonical Lie polynomials from the genus-zero setup determine canonical zeta generators in genus one that act on the two generators of Enriquez' elliptic associators. Up to a single contribution at fixed degree, the zeta generators in genus one are systematically expanded in terms of Tsunogai's geometric derivations dual to holomorphic Eisenstein series, leading to a wealth of explicit high-order computations. Earlier ambiguities in defining the non-geometric part of genus-one zeta generators are resolved by imposing a new representation-theoretic condition. The tight interplay between zeta generators in genus zero and genus one unravelled in this work connects the construction of single-valued multiple polylogarithms on the sphere with iterated-Eisenstein-integral representations of modular graph forms.
\end{tabular}

\end{center}

\thispagestyle{empty}

\newpage
\hypersetup{pageanchor=true}
\setcounter{page}{1}

\numberwithin{equation}{section}

\setcounter{tocdepth}{2}
\tableofcontents


\vskip.5cm
\noindent\rule{\textwidth}{.3mm}
\vskip1cm

\section{Introduction}\label{sec:Intro}

A wealth of recent interactions between mathematics and physics evolves around the appearance of multizeta values in period integrals and scattering amplitudes.
In the first place, multizeta values arise as real numbers defined by the infinite sums 
\beq
\zeta_{k_1,k_2,\ldots,k_r} \coloneqq \sum_{1\leq n_1<n_2<\ldots<n_r}^{\infty} n_1^{-k_1} n_2^{-k_2}\ldots n_r^{-k_r}\, ,
\label{MZV.real}
\eeq
where $k_1,\ldots,k_r \in \mathbb N$ and $k_r>1$ in order to ensure convergence of the sum. While many of their number-theoretic properties, including transcendentality of Riemann zeta values, remain conjectural at the level of real numbers, there are rigorous results on their motivic versions.\footnote{The definition of motivic multizeta values originates from algebro-geometric considerations as explained for example in~\cite{Goncharov:2005sla,Brown:2011mot,Brown2014MotivicPA}.} In particular, motivic multizeta values admit a Hopf algebra structure which was observed to universally govern quantum-field-theory and string-theory amplitudes. These structures also facilitate and guide many concrete investigations of period computations. In this article we set forth some canonical features of motivic multizeta values at the level of their Hopf algebra comodule $\MZ$ and its graded dual.

Both practical computations and structural understanding of multizetas benefit from organizing period integrals into generating series. The non-commuting bookkeeping variables of these generating series can be identified with generators of certain Lie algebras associated with punctured Riemann surfaces of different genus. A central ingredient of generating-series approaches to multizeta values are {\it zeta generators} that have been studied before from different angles in the literature but whose definition depended on certain ad-hoc choices.

In this work, we present canonical zeta generators in genus zero and genus one.
Furthermore, our results also imply a canonical map from motivic multizeta values to the so-called $f$-alphabet~\cite{Brown:2011ik, Brown:2011mot}, a representation of $\MZ$ that is widely used but has eluded a canonical form until this work. The methods we present in this work are fully constructive.

One of the main motivations for seeking concrete expressions for zeta generators stems from the key role they play for the construction of single-valued multiple polylogarithms on the sphere \cite{svpolylog, Broedel:2016kls, DelDuca:2016lad, Frost:2023stm} and of modular equivariant iterated integrals of Eisenstein series~\cite{Brown:2017qwo, Brown:2017qwo2, Dorigoni:2022npe, Dorigoni:2024oft}.
In a physics context, our results on zeta generators facilitate and organize the low-energy expansion of perturbative string-theory scattering amplitudes~\cite{DHoker:2015wxz, DHoker:2016mwo, Gerken:review, Berkovits:2022ivl, DHoker:2022dxx}. In particular, the contributions from string world-sheets of genus zero and one are intertwined through the connection between the associated zeta generators, bringing out universal structures of importance also for string dualities.

\subsection{The canonical zeta generators in genus zero} 

Our first main contribution is the definition of a canonical set of generators for the graded dual
$\MZ^\vee$ of the Hopf algebra comodule $\MZ$ of motivic multizeta values. These canonical generators will be encoded in a family of homogeneous degree-$w$ polynomials 
\beq
g_w(x,y)\in \Q \langle x,y\rangle\, ,\ \ \ \ w\ge 2
\eeq
in two non-commutative variables $x,y$. The polynomials $g_w$ satisfy three natural conditions presented in Theorem \ref{thmintro.1} below and related to the intrinsic structure of $\MZ$, see Remark \ref{firstexs} and \eqref{explgws} for examples. For odd values of $w$, the polynomials $g_w$ are Lie polynomials which provide a set of canonical generators for the genus zero motivic Lie algebra which is well-known to be a free Lie algebra with one generator in each odd degree $w\ge 3$ (a result  established in \cite{Levine}).\footnote{More formally, this is the Lie algebra of the pro-unipotent radical of the fundamental group of the Tannakian category of mixed Tate motives unramified over $\mathbb{Z}$ \cite{Deligne:2005}.} This freeness property is essential throughout this work but only established for the Lie-algebra structure underlying {\it motivic} multizetas as opposed to the incarnation (\ref{MZV.real}) of multizetas as real numbers. Hence, our main results are stated only for {\it motivic} multizetas, although they are expected to apply to real multizetas (\ref{MZV.real}) in identical form by the conjectural isomorphism between real and motivic multizetas.

The key tool used to define the polynomials $g_w$ is the {\em Z-map}, first introduced in \cite{Schneps:2013} and explained here in section~\ref{sec:Zmap}, which is a canonical linear isomorphism from $\MZ^\vee$ to $\MZ$. The Z-map comes from the canonical isomorphism of vector spaces
\beq
\Q\langle x,y\rangle\rightarrow \Q[Z(w)]\,,
\eeq
where the space on the right-hand side is the $\Q$-vector space on symbols $Z(w)$ indexed by all monomials $w$ in the letters $x,y$, and the isomorphism is given simply by mapping $w\mapsto Z(w)$. Identifying $\Q\langle x,y\rangle$ with the dual space of $\Q[Z(w)]$ and considering the bases of monomials $w$ and of symbols $Z(w)$ as dual bases makes this into an isomorphism of dual vector spaces.
We will also use variants of the Z-map associated with different quotients of $\Q[Z(w)]$ or subspaces of $\Q\langle x,y\rangle$. 
One such quotient is given by $\MZ$, obtained by imposing the linear relations between motivic multizeta values on the symbols $Z(w)$.

Let $\mz$ denote the quotient of $\MZ$ modulo the linear subspace spanned by constants, non-trivial products of motivic multizeta values, and the motivic single zeta value $\zeta^{\mathfrak{m}}_2$. Then, $\mz$ inherits the structure of a Lie coalgebra from the Hopf-algebra-comodule structure of $\MZ$ (cf.~section \ref{sec:Background.2}). Let $\mz^\vee\subset\MZ^\vee\subset \Q\langle x,y\rangle$ denote its dual space, which is a Lie algebra equipped with the Ihara bracket below (see also section \ref{sec:Zmap.1}). Like $\MZ$ and $\MZ^\vee$, the spaces $\mz$ and $\mz^\vee$ are graded by the (homogeneity) degree in $x,y$ or {\it weight}, with finite-dimensional graded parts for fixed weight. A major structure theorem by Brown \cite{Brown:2011mot} has shown that $\mz^\vee$ is freely generated by one Lie polynomial (of depth 1 in the sense of Definition \ref{gwcan}) in each odd homogeneous weight $w\ge 3$. 

The universal enveloping algebra ${\cal U}\mz^\vee$ is freely generated by the generators of $\mz^\vee$ under the Poincar\'e--Birkhoff--Witt multiplication, which we denote by $\circp$. In the case where $g\in \mz^\vee$ and $h\in {\cal U}\mz^\vee$, this multiplication rule has a simple form:
\beq\label{simplediamond}
g\circp h=gh+D_g(h)\, ,
\eeq
where $D_g$ is the Ihara derivation of $\Q\langle x,y\rangle$ defined by $D_g(x)=0$ and $D_g(y)=[y,g]$. The space $\MZ^\vee$ is a module over the Hopf algebra ${\cal U}\mz^\vee$.

Let us write $\mz^\vee_{\ge 2}$ for the subspace of $\mz^\vee$ spanned by Ihara brackets $\{ g,h\} \coloneqq g\circp h - h\circp g$ of the generators; this is a canonical subspace independent of any actual choice of generators. 
The spaces $\MZ$, $\MZ^\vee$, $\mz$, $\mz^\vee$ and $\mz^\vee_{\ge 2}$ are all weight-graded spaces; we write
$\MZ_w$, $\MZ^\vee_w$ etc.~to indicate their graded parts of weight $w$, all of which are finite-dimensional. 
Each graded piece $\MZ_w$ contains a canonical {\it reducible} subspace $\hat R_w$ spanned by all weight-$w$ products of lower-weight multizeta values. We write $R_w \coloneqq\hat R_w$ if $w$ is odd, and if $w$ is even we let $R_w$ denote the subspace of $\hat R_w$ spanned by all products except for $(\zeta_2^\mathfrak{m})^{w/2}$, so that
\beq\begin{cases}
    \hat R_w=R_w&\hbox{if $w$ is odd}\, ,\\
    \hat R_w=\Q\zeta^{\mathfrak{m}}_w\oplus R_w&\hbox{if $w$ is even}\, ,
\end{cases}
\eeq
where $\zeta^{\mathfrak{m}}_w$ denotes the single zeta value in weight $w$.
We then have $\mz_w=\MZ_w/\hat R_w$ for $w\ge 3$.
We further define canonical subspaces of {\it irreducible} multizeta values (resp.~non-single irreducible multizeta values) in $\MZ_w$ for each weight $w\ge 2$ by setting
\beq
    \hat I_w \coloneqq Z(\mz^\vee_w) \, ,\ \ \ \ 
    I_w \coloneqq Z\bigl((\mz^\vee_{\ge 2})_w\bigr) \, ,
\eeq
where we note that 
\beq\begin{cases}
    \hat I_w=I_w&\hbox{if $w$ is even} \, , \\
    \hat I_w=\Q\zeta^{\mathfrak{m}}_w\oplus I_w&\hbox{if $w$ is odd}\, .
\end{cases}
\eeq
In this way, we obtain a canonical decomposition of $\MZ_w$ into single, irreducible and reducible parts:
\beq\label{canondec0}\MZ_w=\Q\zeta^{\mathfrak{m}}_w\oplus I_w\oplus R_w
\hspace{10mm}\text{for all $w\ge 2$}\, .
\eeq

\vspace{.3cm}
\noindent {\bf Examples.} While for all $w\leq7$ the irreducible parts are trivial, e.g.
\begin{equation}
 \MZ_5=\Q\zeta_5^{\mathfrak{m}}\oplus R_5 = \langle\zeta_5^{\mathfrak{m}}\rangle\oplus \langle \zeta_2^{\mathfrak{m}}\zeta_3^{\mathfrak{m}}\rangle\,,
\end{equation}
for $w\geq 8$ we have that generically $I_w\neq \emptyset$.
The first non-trivial instance of this decomposition~is
\begin{align}
\MZ_8&=\langle\zeta_8^{\mathfrak{m}}\rangle\oplus \langle Z_{35}\rangle\oplus \langle \zeta_3^{\mathfrak{m}}\zeta_5^{\mathfrak{m}}\,,\ \zeta_2^{\mathfrak{m}}(\zeta_3^{\mathfrak{m}})^2\rangle =\Q\zeta_8^{\mathfrak{m}}\oplus I_8\oplus R_8\,,
\end{align}
where the canonical choice of irreducible mutizeta,
\begin{align}
Z_{35} &\coloneqq Z(\{g_3,g_5\})=
-\tfrac{1105181}{80}\zeta_8^{\mathfrak{m}}
+\tfrac{24453}{5}\zeta_{3,5}^{\mathfrak{m}}
+\tfrac{28743}{2}\zeta_3^{\mathfrak{m}}\zeta_5^{\mathfrak{m}}-1683\,\zeta_2^{\mathfrak{m}}(\zeta_3^{\mathfrak{m}})^2\,,\label{eq:Z35Intro}
\end{align}
is dictated by the procedure described above. We refer to section \ref{sec:Zmap.4a} for more details and examples at higher weight.

\medskip

The explicit form of the canonical polynomials $g_w$ can be obtained
from the motivic Drinfeld associator~\cite{Drinfeld:1989st, Drinfeld2} 
\beq
\Phi^{\mathfrak{m}}_{\rm KZ}(x,y)\in \Q\langle\langle x,y\rangle\rangle\otimes_\Q\MZ\,,
\eeq
which for our work can be thought of as a generating series of motivic multizeta values \cite{LeMura}. For convenience, we work with the motivic power series $\Phi^{\mathfrak{m}}(x,y) \coloneqq \Phi^{\mathfrak{m}}_{\rm KZ}(x,-y)$. 
Apart from the definition of the Z-map and the canonical decomposition \eqref{canondec0}, the main results of sections \ref{sec:Background} and \ref{sec:Zmap} are summarized by:

\begin{thm}\label{thmintro.1}  Write the expansion of $\Phi^{\mathfrak{m}}$ in $x,y$ in any basis of motivic multizetas adapted to the canonical decomposition \eqref{canondec0}, and for each $w\ge 2$, set 
\beq
g_w \coloneqq \Phi^{\mathfrak{m}}|_{\zeta^{\mathfrak{m}}_w} \, .
\eeq
Then the polynomials $g_w$ lie in $\MZ_w^\vee$. Equivalently, $g_w$ can be identified (with no reference to $\Phi^{\mathfrak{m}}$) as the unique polynomial in $\MZ_w^\vee$ satisfying the following three properties: 
\begin{itemize}
    \item[(i)] $\langle g_w,\zeta^{\mathfrak{m}}_w\rangle=1$, where $\langle \cdot, \cdot \rangle$ denotes the canonical action of the dual space $\MZ^\vee$ on $\MZ$ (see~\eqref{eq:Zmap.01} below),
    \item[(ii)] $g_w$ annihilates the reducible subspace $R_w\subset \MZ_w$,
    \item[(iii)] $Z(g_w)\in \Q\zeta^{\mathfrak{m}}_w\oplus R_w$, i.e.\ it does not contain any irreducible multizeta values in $I_w$.
\end{itemize}
The $g_w$ for odd $w\ge 3$ form a canonical set of generators for the Lie algebra $\mz^\vee$, and the $g_w$ for all $w\ge 2$ form a set of generators for the Hopf algebra module $\MZ^\vee$ over the Hopf algebra ${\cal U}\mz^\vee$. More precisely, every element of $\MZ^\vee$ can be written uniquely as a product
\beq
g_{w_1}\circp\cdots\circp g_{w_r}\circp g_k \, ,
\eeq
where the $w_i$ are all odd $\ge 3$ and $k\ge 2$, and the multiplication proceeds from right to left using the rule
\eqref{simplediamond}.
\end{thm}

\begin{rmk}
\label{firstexs}
For both even and odd $w\geq 2$, the polynomials $g_w$ are canonical since the subspaces $R_w, I_w$ in part (ii) and (iii) of Theorem \ref{thmintro.1} are. Their simplest instances are given by $g_2 = [x,y]$ and $g_3 = [x-y,[x,y]]$, with more examples in~\eqref{explgws}. 
The ancillary files of the arXiv submission of this work
contain the explicit form of all $g_w$ with $w\leq 12$.
\end{rmk}

\begin{rmk} \label{rmk:113}
In \cite{Ecalledim}, \'Ecalle gave an alternative method to specify a depth 1 polynomial $g_w$ in each odd weight $w$, by defining an inner product on monomials and then choosing $g_w$ to be the unique element orthogonal to all elements in the associated graded Lie algebra. Using the language of moulds reviewed in section \ref{sec:prop1}, he gave two different possibilities for inner products with good symmetry properties (called $kya$ and $kwa$, cf.~section 19 of \cite{Ecalledim}). 
\end{rmk}

\begin{rmk} \label{rmk:114}
In his Ph.D.\ thesis \cite{Keilthy}, Keilthy gave a similar method, but using the ``trivial'' inner product on $\Q\langle x,y\rangle$ (for which the inner product of two monomials $u$ and $v$ is $\delta_{u,v}$). This is analogous to our use of the action of the dual space $\Q\langle x,y\rangle^\vee$ on $\Q\langle x,y\rangle$ described in section \ref{sec:Zmap.1} below (cf.~(\ref{eq:Zmap.01})). As we explain in more detail in section~\ref{sec:Zmap.4}, this construction agrees with ours for odd $w$; however our definition applies uniformly to produce canonical elements for both odd and also for even~$w$.
\end{rmk}

\subsection{\texorpdfstring{The canonical $f$-alphabet isomorphism}{The canonical f-alphabet isomorphism}}

Brown proved in \cite{Brown:2011ik, Brown:2011mot} that the motivic multizeta algebra $\MZ$ is isomorphic to a certain Hopf-algebra comodule ${\cal F}$, known as the $f$-alphabet algebra, which has a very simple structure: it is a commutative algebra under the shuffle multiplication, multiplicatively generated by all monomials in an alphabet of letters $f_2$ and $f_3,f_5,f_7,\ldots$ which is free apart from the unique relation that $f_2$ commutes with all the other letters; thus we have
\beq
{\cal F}=\Q[f_2]\otimes_\Q \overline{\cal F} \, ,
\eeq
where $\overline{\cal F}$ is freely generated under the shuffle multiplication by all monomials in $f_3,f_5,\ldots$, and we sometimes write $f_{2n} = \frac{\zeta_{2n}}{(\zeta_2)^n} f_2^n$ for $n \in \mathbb N$. The space $\overline{\cal F}$ is a commutative Hopf algebra equipped with the shuffle multiplication and the deconcatenation coproduct, and ${\cal F}$ is a Hopf-algebra comodule equipped with the following extension of the deconcatenation coproduct to a coaction:
\begin{align}\label{deconc1}
    \Delta:{\cal F}&\rightarrow {\cal F}\otimes \overline{\cal F} \, ,\\
f_2^nf_{w_1}\cdots f_{w_r}&\mapsto
\sum_{i=0}^r f_2^nf_{w_1}\cdots f_{w_i}\otimes f_{w_{i+1}}\cdots f_{w_r}\,.
\nn
\end{align}
In \cite{Brown:2011ik,Brown:2011mot}, Brown identified the complete family of Hopf-algebra-comodule isomorphisms $\MZ\rightarrow {\cal F}$ normalized by $\zeta^{\mathfrak{m}}_w\mapsto f_w$, showing that it is parametrized by rational parameters indexed by any basis of non-single irreducible multizetas. In section~\ref{sec:falpha}, we display a canonical choice of one such isomorphism, uniquely determined as follows.

\begin{thm} There exists a canonical normalized Hopf algebra comodule isomorphism $\rho:\MZ\rightarrow {\cal F}$ whose definition depends only on the canonical decomposition \eqref{canondec0}; it is characterized by each of the two following properties, which are equivalent: 
\begin{itemize}
    \item $\rho$ satisfies
\beq
\rho(\xi)|_{f_w}=0\ \ \ \ \forall \ 
\xi\in I_w \, ,
\label{rhochar1}
\eeq
\item if $\Phi^{\mathfrak{m}}$ is written in a basis adapted to the canonical decomposition \eqref{canondec0}, then $\rho$ satisfies 
\beq
\rho(\Phi^{\mathfrak{m}})|_{f_w}=g_w\ \ \ \ \forall \ w  \ge 2 \, .
\label{rhochar2}
\eeq
\end{itemize}
This choice of isomorphism $\rho$ is canonical since the subspaces $I_w$ and the polynomials $g_w$ in
(\ref{rhochar1}) and (\ref{rhochar2}) are.
\end{thm}

\vspace{.3cm}
\noindent {\bf Example.}
The irreducible multizeta value of~\eqref{eq:Z35Intro} has the following $f$-alphabet image:
\begin{align}
    \rho(Z_{35})&=-\tfrac{20163}{2}f_3f_5+\tfrac{28743}{2}f_5f_3-3366f_2f_3f_3 \, .
\end{align}
The canonical choice of isomorphism is reflected in the absence of a term proportional to $f_8$.

\subsection{The canonical zeta generators in genus one}

Sections \ref{sec:sigman} to \ref{bigsecmd2} are dedicated to zeta generators $\sigma_w$ in genus one. These are derivations of the free graded Lie algebra 
${\rm Lie}[a,b]$ associated to the (pro-unipotent) fundamental group of the once-punctured torus. Based on earlier work in \cite{EnriquezEllAss,hain_matsumoto_2020, Ecalle, Schneps:2015mzv}, the action of the genus one generators~$\sigma_w$ on $a,b$ is determined in section \ref{sec:sigman.3} from the genus-zero polynomials $g_w$ via (with ${\rm B}_n$ the $n^{\rm th}$ Bernoulli number)
\begin{align}
\sigma_w(s_{12})&=0 \, ,
&\sigma_w(s_{01})&=  
\big[ s_{01} , g_w(s_{12},- s_{01})\big] \, , \label{sigintro} \\
s_{12} &= [b,a] \, ,  
&s_{01}&= -b-\sum_{n\ge 1} \frac{{\rm B}_n}{n!}\ad_a^n(b)\,,
\notag
\end{align}
together with the ``extension lemma'' 2.1.2 of~\cite{Schneps:2015mzv} reviewed in section \ref{sec:sigman.2}. In view of the canonical $g_w$ in the defining equation (\ref{sigintro}), we arrive at the first canonical choice of the zeta generators $\sigma_w$ in genus one at arbitrary odd $w\geq 3$.

Another important family of derivations on ${\rm Lie}[a,b]$ are Tsunogai's $\ep_{k}$ in even degree $k\geq 0$ (i.e.\ the combined homogeneity degrees in $a$ and $b$). The algebra $\mathfrak{u}$ generated by the $\ep_k$ is called the algebra of \textit{geometric derivations}. By work of Hain--Matsumoto \cite{hain_matsumoto_2020}, the $\sigma_w$ normalize the algebra $\mathfrak{u}$, i.e.\ commutators $[\sigma_w,\mathfrak{u}]$ are again contained in $\mathfrak{u}$. In fact, upon decomposing the zeta generators $\sigma_w$ into an infinite number of contributions at fixed even degree $\geq w+1$, all the terms lie in $\mathfrak{u}$ except for certain contributions at {\em key degree}~$2w$. The terms of $\sigma_w$ outside $\mathfrak{u}$ belong to yet another derivation known as the {\em arithmetic part} $z_w$ that furnishes a one-dimensional representation under the $\mathfrak{sl}_2$ spanned by the ${\rm Lie}[a,b]$-derivations $\ep_0, \ep_0^\vee$ and $\hhh \coloneqq [\ep_0, \ep_0^\vee]$ defined by
\beq
\ep_0(a) = b \, , \ \ \ \ \ep_0(b) = 0 \, , \ \ \ \ \ep_0^\vee(a) = 0 \, , \ \ \ \ \ep_0^\vee(b) = a\, .
\eeq
Even with the canonical definition of $\sigma_w$, the arithmetic derivations $z_w$ are not entirely characterized by requiring that they form an $\mathfrak{sl}_2$ singlet and that $\sigma_w-z_w \in \mathfrak{u}$. We arrive at canonical representatives of $z_w$ by additionally imposing that they exhaust the complete $\mathfrak{sl}_2$ singlet at key degree of $\sigma_w$. More specifically, the  $\ep_k^{(j)} \coloneqq {\rm ad}_{\ep_0}^j(\ep_k)$ with $j=0,1,\ldots,k{-}2$ composing $\sigma_w{-}z_w$ fall into $(k{-}1)$-dimensional representations of $\mathfrak{sl}_2$ because of $\ep_k^{(k-1)} \!=\! 0$. The canonical arithmetic derivations $z_w$ are then uniquely defined by imposing that any nested commutator $\ep_k^{(j)} $ at the key degree of $\sigma_w-z_w$ belongs to $\mathfrak{sl}_2$ representations of dimension~$\geq 3$.

Based on mould theory, we describe a first algorithm in section \ref{sec:prop2} to explicitly compute the action of $\sigma_w$ on $a$ and $b$ degree by degree and prove the following theorem:
\begin{thm}[see Theorem~\ref{thm:522} (iii)]\label{thmintro.3}
The genus-one zeta generators $\sigma_w$ are entirely determined by their parts of degree $<2w$.
\end{thm}
This remarkable property of $\sigma_w$ can be combined with the commutation relation~\cite{hain_matsumoto_2020}
\beq
[N,\sigma_w]= 0 \ \ {\rm with} \ \ N\coloneqq -\epsilon_0 + \sum_{k=2}^\infty (2k-1) \frac{{\rm B}_{2k}}{(2k)!} \epsilon_{2k}\,,
\label{nsigintro}
\eeq
to make $\sigma_w$ computationally accessible to all degrees.
By solving (\ref{nsigintro}) for $[\ep_0,\sigma_w]$, it relates contributions to $\sigma_w-z_w$ with different numbers of $\ep_{k_i}^{(j_i)}$ factors (with $0\leq j_i \leq k_i-2$) to be referred to as {\em modular depth}.\footnote{The Lie algebra $\mathfrak{u}$ is not free on the $\ep_k^{(j)}$ but satisfies relations~\cite{{IharaTakao,Schneps:2006,Pollack}} that are not homogeneous in modular depth and which, for this reason, only provides a filtration rather than a grading of $\mathfrak{u}$, see Remark~\ref{Pollackrel}.} On these grounds, we describe a second algorithm in section \ref{secmd2} based on (\ref{nsigintro}) to determine $\sigma_w-z_w$ recursively in modular depth, up to highest-weight vectors of $\mathfrak{sl}_2$ in each step which are defined to lie in the kernel of $\ad_{\ep_0}$. We will infer from the results of \cite{hain_matsumoto_2020} that there are no highest-weight vectors beyond key degree. From the viewpoint of (\ref{nsigintro}), it is thus sufficient to know the degree $\leq 2w$ parts (though Theorem \ref{thmintro.3} even guarantees that the complete information is available from degree $< 2w$) of $\sigma_w$. The infinity of terms at degree $\geq 2w+2$ follows from (\ref{nsigintro}) together with representation theory of $\mathfrak{sl}_2$.

This setup leads us to present a closed all-degree formula for $\sigma_w$ up to contributions in $\mathfrak{u}$ of modular depth $\geq 3$ (in the ellipsis),
\begin{align}
\sigma_w
 &= z_w - \frac{1}{(w-1)!} \ep_{w+1}^{(w-1)} \label{clmd2in} 
\\*
&\quad
-\frac{1}{2} \sum_{d=3}^{w-2} \frac{ \BF_{d-1} }{\BF_{w-d+2} }
\sum_{k=d+1}^{w-1} \BF_{k-d+1} \BF_{w-k+1} s^d(\ep_k,\ep_{w-k+d})
\notag\\
&\quad 
- \sum_{d=5}^w \BF_{d-1} s^d(\ep_{d-1},\ep_{w+1}) 
 - \frac{1}{2} \BF_{w+1} s^{w+2}(\ep_{w+1},\ep_{w+1})
\notag \\*
&\quad + \sum_{k=w+3}^\infty
\BF_k \sum_{j=0}^{w-2}  \frac{(-1)^j  \binom{k{-}2}{j}^{-1} }{j! (w{-}2{-}j)! } \, 
[  \ep_{w+1}^{(w-2-j)}  , \ep_k^{(j)} ] +\ldots \, ,
\notag
\end{align}
where we employ the shorthand $\BF_{k}\coloneqq \frac{{\rm B}_k}{k!}$ and we define
\beq
s^d(\ep_{k_1},\ep_{k_2}) \coloneqq 
\frac{(d{-}2)! }{(k_1{-}2)! (k_2{-}2)!} \sum_{i=0}^{d-2} (-1)^i   [ \ep_{k_1}^{(k_1-2-i)}, \ep_{k_2}^{(k_2-d+i)}] \, .
\label{anothersdkk}
\eeq
The highest-weight-vector contribution $\sim \ep_{w+1}^{(w-1)}$ in first line of (\ref{clmd2in}) is well-known and is used to determine the modular-depth two terms in the third and fourth line from (\ref{nsigintro}). The second line of (\ref{clmd2in}) is conjectural and features highest-weight vectors $s^d(\ep_k,\ep_{w-k+d})$ in each term -- they are not fixed by (\ref{nsigintro}) and confirmed by direct computation in a large number of examples. Moreover, the $d=3$ terms in the second line of (\ref{clmd2in}) reproduce the closed formula of Brown \cite{brown_2017} on depth-three terms in the terminology of the reference.

\vspace{.3cm}
\noindent {\bf Examples.}
As an illustration of Theorem~\ref{thmintro.3}, the first two zeta generators are determined fully by the following terms of their expansion in (\ref{clmd2in}) bounded by the respective key degrees 6 and 10,
\begin{align}
\sigma_3 = z_3- \frac{1}{2} \ep_4^{(2)}  +\ldots\,, \hspace{10mm}
\sigma_5 =  z_5 -\frac{1}{24} \epsilon_6^{(4)}
-\frac{5  }{48} [\epsilon_4^{(1)},\epsilon_4^{(2)}]+\ldots \, .
\end{align}
We refer to section~\ref{bigsecmd2} for more examples. 

\medskip

Finally, (\ref{nsigintro}) together with the terms of modular depth $d$ in $\sigma_w - z_w$ fix the explicit form of $[z_w, \ep_k] \in \mathfrak{u}$ up to and including modular depth $d+1$. Accordingly, the closed formula (\ref{clmd2in}) determines the terms of modular depth three beyond the well-known contributions \cite{hain_matsumoto_2020}
\beq
[z_w, \ep_k] = \frac{\BF_{w+k-1}}{\BF_k} 
\sum_{i=0}^{w-1} \frac{ (-1)^i(k+i-2)!}{i!(w+k-3)!}
[\ep_{w+1}^{(i)}, 
     \ep_{  w+k -1 }^{(w-i-1)} ] 
+\ldots
\eeq
and we give closed formulae for $[z_3,\ep_k]$ and $[z_5,\ep_k]$ at modular depth three in section \ref{secmd2.b}.

\subsection{Motivation and outlook}

A major motivation for our study of zeta generators stems from their relevance for periods of configuration spaces of Riemann surfaces with marked points. In genus zero, the canonical polynomials $g_w$ take center stage in the recent reformulation \cite{Frost:2023stm} of the motivic coaction \cite{goncharov2001multiple, Goncharov:2005sla, Brown:2011mot} and the single-valued map \cite{svpolylog, Broedel:2016kls, DelDuca:2016lad} of multiple polylogarithms on the sphere. The genus-one zeta generators $\sigma_w$ and their interplay with geometric derivations $\ep_k$ unlocked a fully explicit generating-series description of non-holomorphic modular forms in a companion paper \cite{Dorigoni:2024oft} to this work. 

As detailed in \cite{Dorigoni:2024oft}, the expansion of $\sigma_w$ in terms of the geometric derivations $\epsilon_k$ determines the appearance of (single-valued) multizeta values in so-called modular graph forms \cite{DHoker:2015wxz, DHoker:2016mwo}. 
The latter are non-holomorphic modular forms appearing in genus-one string scattering amplitudes.
At a computational level, the precise expressions for $\sigma_w$ in terms of $\epsilon_k$ presented in this work are crucial for an explicit realization of Brown's construction of non-holomorphic modular forms in \cite{Brown:2017qwo, Brown:2017qwo2} which was related to modular graph forms in \cite{Dorigoni:2022npe}. At a conceptual level, the intimate connection between zeta generators in genus zero and genus one presented in section \ref{sec:sigman} leads to a unified description of the single-valued map of multiple polylogarithms in one variable and iterated Eisenstein integrals \cite{Dorigoni:2024oft}.

These applications of zeta generators in genus zero and genus one lead us to expect that generalizations thereof to compact Riemann surfaces of arbitrary genus with any number of marked points may in fact exist. Our work sets the stage for two lines of follow-up research: 
\begin{itemize}
\item adapting zeta generators in genus one to systematic constructions of single-valued elliptic polylogarithms pioneered by Zagier \cite{Ramakrish} in any number of variables and which were more recently approached in the framework of ``elliptic modular graph forms'' in the string-theory literature \cite{DHoker:2018mys, Basu:2020pey, Dhoker:2020gdz, new:eMGF};
\item determining higher-genus incarnations of zeta generators from degenerations of the flat connections \cite{Enriquez:2011, Zerbini:2021, Zerbini:2022, DHoker:2023vax} used for constructions of polylogarithms on Riemann surfaces of arbitrary genus and applying them to non-holomorphic modular graph forms \cite{DHoker:2013fcx, Pioline:2015qha, DHoker:2017pvk, DHoker:2018mys, Basu:2018bde} and tensors \cite{Kawatalk1, Kawatalk2, DHoker:2020uid, kawazumi2022twisted}.
\end{itemize}

\section*{Acknowledgements}

We are grateful to David Broadhurst, Francis Brown, Emiel Claasen, Eric D'Hoker, Benjamin Enriquez, Hadleigh Frost, Deepak Kamlesh, Pierre Lochak, Franziska Porkert, Christophe Reutenauer, Carlos Rodriguez, Oliver Schnetz and Federico Zerbini for combinations of inspiring discussions and collaboration on related topics, and to Adam Keilthy for bringing results in his thesis \cite{Keilthy} to our attention. The authors would like to thank the organizers of the workshops ``Geometries and Special Functions for Physics and Mathematics'' at the BCTP Bonn and ``New connections between physics and number theory'' at the Pollica Physics Centre for creating a stimulating atmosphere. We are grateful to the Hausdorff Research Institute for Mathematics in Bonn for the hospitality and the vibrant atmosphere during the
follow-up workshop ``Periods in Physics, Number Theory and Algebraic Geometry''.
We thank the Galileo Galilei Institute (GGI) for Theoretical Physics in Florence for the hospitality and the INFN for partial support during the program ``Resurgence and Modularity in QFT and String Theory''.
OS is grateful to the Simons foundation for financial support during the GGI programme.
DD thanks Riken iTHEMS and the Yukawa Institute for
Theoretical Physics at Kyoto University for the hospitality and support during the iTHEMS-YITP Workshop ``Bootstrap, Localization and
Holography''.
The research of MD was supported by the IMPRS for Mathematical and Physical Aspects of Gravitation, Cosmology and Quantum Field Theory. JD is supported by the Royal Society (Spectral theory of automorphic forms: trace formulas and more). MH, OS and BV are supported by the European Research Council under
ERC-STG-804286 UNISCAMP. MH and BV are furthermore supported by the Knut and Alice Wallenberg Foundation under grants KAW 2018.0116 and KAW 2018.0162, respectively. The research of OS is partially supported by the strength area ``Universe and mathematical physics'' which is funded by the Faculty of Science and Technology at Uppsala University.

\section*{Data availability statement}

The ancillary files in the arXiv and journal submissions of this work provide supplemental data concerning the canonical polynomials $g_w$ and the arithmetic derivations $z_w$ mentioned in the introduction.

\section{Background on multizeta values}\label{sec:Background}

In this section, we review basic definitions on different types of multizeta values, their relations and their Hopf-algebraic properties. See \cite{BurgosFresan,ZigZag} for textbook introductions to the subject and \cite{IKZagier, Brown:2011mot}  for earlier references.

\subsection{Real and formal multizeta values}\label{sec:Background.1}
\subsubsection{Real multizeta values, shuffle and stuffle multiplication}

The {\em real multizeta values} are defined by the infinite sums 
\beq
\zeta_{k_1,k_2,\ldots,k_r} \coloneqq \sum_{1\leq n_1<n_2<\ldots<n_r}^{\infty} n_1^{-k_1} n_2^{-k_2}\ldots n_r^{-k_r} \, ,
\label{appMZV.01}
\eeq
where $k_1,\ldots,k_r \in \mathbb N$ and $k_r>1$ in order to ensure convergence
of the sum. The integers~$r$ and $\sum_{i=1}^r k_i$ in (\ref{appMZV.01}) are respectively referred to as the {\em depth} and {\em weight} of $\zeta_{k_1,k_2,\ldots,k_r}$.
Multizeta values (MZVs) satisfy a number of algebraic
relations over $\mathbb Q$ which we discuss further below. Let us first introduce the monomial notation
\beq
\zeta(x^{k_r-1}y\cdots x^{k_2-1}y x^{k_1-1}y)=\zeta_{k_1,k_2,\ldots,k_r}\, ,
\label{appMZV.02}
\eeq
where $x$ and $y$ are non-commutative indeterminates and the convergence property $k_r>1$ implies that the first letter on the left-hand side is $x$. We say that a non-trivial monomial in $x,y$ is {\em convergent} if 
it begins with $x$ and ends with $y$; all other monomials are {\em non-convergent}.
We extend the notation \eqref{appMZV.02} to the definition of the {\em regularized zeta 
values} $\zeta(w)$ for all non-convergent monomials $w=y^r v x^s$ with
$v$ convergent, by the explicit formula 
(established in Prop.~3.2.3 of \cite{Furusho2000TheMZ}, based on the regularization methods of \cite{LeMura})
\begin{align}
\zeta(w)=\sum_{a=0}^r \sum_{b=0}^s (-1)^{a+b}\zeta(y^a\shuffle y^{r-a}vx^{s-b}\shuffle x^b)\, ,
\label{appMZV.03}
\end{align}
an expression in which all the non-convergent $\zeta(w)$ cancel out so that
$\zeta(y^rvx^s)$ is expressed as a linear combination of convergent words 
only, and which ensures that for all pairs of (convergent or non-convergent) 
words $u,v$, the $\zeta$-values satisfy the shuffle relation
\beq
\zeta(u)\zeta(v)=\zeta(u\shuffle v) =\zeta(v\shuffle u)\, ,
\label{appMZV.04}
\eeq
where $\zeta$ is considered as a linear function on words, and we fix the values
$\zeta(x)=\zeta(y)=0$ and also $\zeta({\bf 1})=1$, where ${\bf 1}$ in the argument denotes the empty word.
We recall here that the shuffle product of monomials can be defined recursively
as follows: for any monomial $u$, we have ${\bf 1}\shuffle u=u\shuffle {\bf 1}=u$, 
and if $u,v\ne {\bf 1}$ we write $u=au'$ and $v=bv'$, where $a$ and $b$ are single letters (either $x$ or $y$), and we have
\beq
u\shuffle v=a(u'\shuffle v)+b(u\shuffle v')\, .
\label{appMZV.05}
\eeq
For example, writing $\zeta_2=\zeta(xy)$, we have \beq
\zeta_2^2=\zeta(xy)^2=\zeta(xy\shuffle xy)=4\zeta(xxyy)+2\zeta(xyxy)=
4\zeta_{1,3}+2\zeta_{2,2}\, .
\label{appMZV.06}
\eeq
This multiplication rule is called the {\em shuffle
multiplication} of real MZVs. 

There is a second multiplication, restricted to a subset of words $w$,
which arises when considering the MZVs written as infinite sums as in 
\eqref{appMZV.01}. Indeed, the result of multiplying two such series is 
itself a sum of such series, as can be seen on the first example: 
\begin{align}
\zeta_2^2&=\sum_{n_1\ge 1} n_1^{-2}\sum_{n_2\ge 1} n_2^{-2}\notag\\
&= \sum_{n_1>n_2\ge 1}n_1^{-2}n_2^{-2}+
\sum_{n_2>n_1\ge 1}n_1^{-2}n_2^{-2}+\sum_{n_1=n_2\ge 1}n_1^{-4}\notag\\
&= 2\zeta_{2,2}+\zeta_4\, .
\label{appMZV.06bis}
\end{align}
This product, called the {\em stuffle product}, can be defined for any pair of words $u,v$ ending in $y$ as follows: we first note that every monomial $u$ ending in $y$ can be rewritten in the free variables $y_i=x^{i-1}y$, with $i\ge 1$:  
\begin{align}
u=y_{i_1}\cdots y_{i_r}\, .
\end{align}
We stipulate that for all such monomials, we have $u\stuffle {\bf 1}={\bf 1}\stuffle u=u$. Then, in the case where $u,v\ne {\bf 1}$,
we peel off the first letter of each of the two words, writing $u=y_{i_1}u'$ and $v=y_{j_1}v'$ with $u'=y_{i_2}\cdots y_{i_r}$ and 
$v'=y_{j_2}\cdots y_{j_r}$,  and define the stuffle product by the recursive rule (first developed by Hoffman in \cite{Hoffman1})
\begin{align}\label{stuffledef}
u\stuffle v=y_{i_1}(u'\stuffle v)+y_{j_1}(u\stuffle v')+y_{i_1+j_1}(u'\stuffle v')\, .
\end{align}
The stuffle product is commutative and associative on words ending in $y$. 

Associated with the stuffle product, one can define a
{\em stuffle regularization} $\zeta_\stuffle(w)$ of MZVs for words ending in $y$. For convergent words $w$ (beginning with $x$ and ending in $y$) we set $\zeta_\stuffle(w) = \zeta(w)$.
The stuffle-regularized MZVs for 
non-convergent words ending in $y$ are defined as follows. First we deal with $\zeta_\stuffle(y^i)$ for $i\geq 0$  by writing the generating series
\begin{align}\label{appMZV.06ter}
 \sum_{n\ge 0}\zeta_\stuffle(y^n)  y^n \coloneqq{\rm exp}\bigg(\sum_{n\ge 2} \frac{(-1)^{n-1}}{n}\zeta(x^{n-1}y)y^n\bigg)\, ,
\end{align}
leading for instance to 
\begin{align}
\zeta_\stuffle({\bf 1})&=1 \, , \notag\\
\zeta_\stuffle(y)&=0 \, , \notag\\  \zeta_\stuffle(y^2)&= - \tfrac12 \zeta(xy)  = - \tfrac12 \zeta_2 \, , \label{stuffex}\\
\zeta_\stuffle(y^3)&=\tfrac13 \zeta(x^2y)= \tfrac13\zeta_3 \, , \notag \\
\zeta_\stuffle(y^4)&=-\tfrac{1}{4}\zeta(x^3y)+\tfrac{1}{8}\zeta(xy)^2=-\tfrac{1}{4}\zeta_4+\tfrac{1}{8}\zeta_2^2 \, .
\notag
\end{align}
Then for monomials $y^iv$ for a non-trivial convergent word $v$ we define the stuffle regularization~by
\beq\label{appMZV.06quater}
\zeta_\stuffle(y^iv)=\sum_{j=0}^i \zeta_\stuffle(y^j)\zeta(y^{i-j}v)\,,
\eeq
where the notation $\zeta(y^{i-j}v)$ refers to the shuffle regularization defined in \eqref{appMZV.03}.

The stuffle-regularized zeta values $\zeta_\stuffle(u)$ defined in this way satisfy the stuffle relations  
\beq
\zeta_\stuffle(u)\zeta_\stuffle(v)=\zeta_\stuffle(u\stuffle v)=\zeta_\stuffle(v\stuffle u)
\label{appMZV.07}
\eeq
for every pair of monomials $u,v$ both ending in $y$ as a direct consequence of their infinite sum expressions \eqref{appMZV.01} (see the original reference \cite{Hoffman1}). In 
particular the stuffle relations hold for ordinary MZVs $\zeta(u)$ and $\zeta(v)$ when $u$ and $v$ are convergent words; for example, we have
\beq
xy\stuffle xy=y_2\stuffle y_2=2y_2^2+y_4=2xyxy+xxxy \, ,
\label{appMZV.08}
\eeq
which corresponds to $\zeta_2^2=2\zeta_{2,2}+\zeta_4$ as in \eqref{appMZV.06bis} above. Note that if both $u$ and $v$ are convergent, then since $\zeta_\stuffle(u)=\zeta(u)$ and $\zeta_\stuffle(v)=\zeta(v)$, combining~\eqref{appMZV.04} and~\eqref{appMZV.07} implies that
\begin{align}\label{appMZ.09}
    \zeta(u) \zeta(v) = \zeta(u \shuffle v) = \zeta(u \stuffle v)
    \quad\quad \text{($u,v$ convergent)}\,.
\end{align}

The family of relations between MZVs consisting of the (``regularized'') shuffle relations
\eqref{appMZV.04} for all pairs of monomials $u,v$ and the (``regularized'') stuffle relations \eqref{appMZV.07} for all pairs of words $u,v$ both ending in $y$ is known as {\em the family of regularized double shuffle relations} on MZVs. These were studied fully in \cite{IKZagier}; for a standard reference text containing all the basic material on MZVs, see also \cite{BurgosFresan}.

\subsubsection{Formal MZVs} 
\label{sec:FMZV}

The formal MZVs, denoted by $\zeta^{\mathfrak{f}}(w)$, are symbols which by definition satisfy
only the (regularized) double shuffle relations explained above, as opposed to the real MZVs which may in theory satisfy any number of additional relations, even including the possibility of being rational numbers. General references for this material are~\cite{Ecalledim,IKZagier,Racinet2002,Burmesteretal2024, EspieNovelliRacinet+2003+1+16}. Let us introduce the notation for the ring of formal MZVs. 

For each $n\ge 0$, let $\Q_n[Z(w)]$ denote the vector space spanned
by formal symbols $Z(w)$ indexed by all degree $n$ monomials $w$ in two
non-commutative variables $x$ and $y$; in particular we have
 $\Q_0[Z(w)]=\Q$. We set 
\begin{align}
    \Q[Z(w)]\coloneqq \bigoplus_{n\ge 0} \Q_n[Z(w)]\, ,
\end{align}
and make this vector space into a commutative $\Q$-algebra by equipping
it with the (commutative) shuffle multiplication
\beq
Z(u)Z(v)=Z(u\shuffle v) \, .
\label{sec1eq.1}
\eeq

Let us introduce a second set of formal symbols $Z_\stuffle(w)$ for monomials $w$ ending in $y$, by 
\begin{itemize}[itemsep=0mm]
\item 
setting $Z_\stuffle(w)\coloneqq Z(w)$ for convergent $w$,
\item 
defining $Z_\stuffle(y^n)$ for $n\ge 1$ by the equation \eqref{appMZV.06ter} with $\zeta$ replaced by $Z$,
\item defining $Z_\stuffle(y^iv)$ for convergent words $v$ by equation \eqref{appMZV.06quater} with $\zeta$ replaced by $Z$.
\end{itemize}
Given that multiplying the symbols $Z(w)$ by the shuffle multiplication \eqref{sec1eq.1} reduces products to linear combinations, all of the new symbols $Z_\stuffle(w)$ can be expressed in terms of linear combinations of the symbols $Z(w)$.

\begin{dfn}
\label{dfn:IFZ} 
Let $\I_\FZ$ be the ideal of the ring $\Q[Z(w)]$ 
generated by the following two families of linear combinations: on the one hand the {\em regularizations} 
\begin{align}\label{sec1eq.2a}
Z(w)
-\sum_{a=0}^r \sum_{b=0}^s (-1)^{a+b}Z(y^a\shuffle y^{r-a}vx^{s-b}\shuffle x^b) \, ,
\end{align}
 for all words $w =y^r v x^s$  with $v$ convergent (adapted from \eqref{appMZV.03}),  and on the other hand the regularized stuffles given for all pairs of monomials $u$ and $v$ both ending in $y$ by
\begin{align}\label{sec1eq.2b}
Z_\stuffle(u)Z_\stuffle(v)&-Z_\stuffle(u\stuffle v)
\end{align}
(adapted from \eqref{appMZV.07}).
The expression \eqref{sec1eq.2b} is to be computed as a linear combination of symbols $Z(w')$ where the monomials $w'$ are all of homogeneous weight equal to the sum of the weights of $u$ and $v$ by (i) expanding out the right-hand term as a linear combination, (ii) replacing every occurrence of $Z_\stuffle$ by a polynomial expression in $Z$ using \eqref{appMZV.06ter} and \eqref{appMZV.06quater}, (iii) using the shuffle multiplication \eqref{sec1eq.1} to express all products $Z(w')Z(w'')$ as linear combinations $Z(w'\shuffle w'')$. Thus each of the expressions in \eqref{sec1eq.2a} and \eqref{sec1eq.2b} is a linear combination of fixed weight; we take them all together as the generators 
 of the ideal $\mathcal{I}_{\FZ}$. 
\end{dfn}

\vspace{.3cm}
\noindent {\bf Examples.} {\em Regularization:} the formula \eqref{sec1eq.2a} above for
the non-convergent word $w=yxy$ tells us to add the linear combination
\beq\label{exreg}
Z(yxy)-Z(yxy)+Z(y\shuffle xy)=Z(yxy)+2Z(xyy)
\eeq
to the ideal 
$\I_\FZ$.

\vspace{.2cm}
\noindent {\em Stuffle:} Let us compute the linear combination \beq \label{exstuf0}
Z_\stuffle(y^2)Z_\stuffle(xy)-Z_\stuffle(y^2\stuffle xy)
\eeq 
as a linear combination of $Z$-symbols using the three steps explained below \eqref{sec1eq.2b}. Using \eqref{stuffledef}, we have
\beq
yy\stuffle xy=y_1y_1\stuffle y_2=y_2y_1y_1+y_1y_2y_1+y_1y_2y_2+y_3y_1+y_1y_3=xyyy+yxyy+yyxy+xxyy+yxxy \, ,
\eeq
so by the first step, which consists of expanding out $Z_\stuffle(yy\stuffle xy)$, \eqref{exstuf0} can be rewritten as
\beq\label{exstuf1}
Z_\stuffle(yy)Z_\stuffle (xy)-Z_\stuffle(xyyy)-Z_\stuffle(yxyy)-Z_\stuffle(yyxy)-Z_\stuffle(xxyy)-Z_\stuffle(yxxy)\,.
\eeq
In the second step we replace each $Z_*$ by an expression in $Z$. For the three convergent words $xy$, $xyyy$ and $xxyy$ we have $Z_*=Z$; by \eqref{stuffex} we have $Z_\stuffle(y)=0$ and $Z_*(yy)=-\tfrac{1}{2}Z(xy)$, and finally by \eqref{appMZV.06quater} we have
\begin{align}
Z_\stuffle(yxyy)&=Z(yxyy)+Z_\stuffle(y)Z(xyy)=Z(yxyy) \notag \, , \\
Z_\stuffle(yyxy)&=Z(yyxy)+Z_\stuffle(y)Z(yxy)+Z_\stuffle(yy)Z(xy)=Z(yyxy)-\tfrac{1}{2}Z(xy)^2 \, ,\\
Z_\stuffle(yxxy)&=Z(yxxy)+Z_\stuffle(y)Z(xxy)=Z(yxxy)\, . \notag
\end{align}
Plugging these into \eqref{exstuf1} allows us to rewrite \eqref{exstuf0} as
\beq\label{exstuf2}
-\tfrac{1}{2}Z(xy)^2-Z(xyyy)-Z(yxyy)-Z(yyxy)+\tfrac{1}{2}Z(xy)^2-Z(xxyy)-Z(yxxy)\,.
\eeq
If necessary we could now expand out the products of $Z$-symbols using the shuffle, but since they cancel out we don't need to, so in the end we add the linear combination
\beq\label{exstuf3}
-Z(xyyy)-Z(yxyy)-Z(yyxy)-Z(xxyy)-Z(yxxy)
\eeq
to the ideal $\I_{\FZ}$.

\vspace{.3cm}

\begin{rmk}
Note that by \eqref{sec1eq.1}, for convergent words $u$ and $v$, the relations~\eqref{sec1eq.2b} of $\mathcal{I}_{\FZ}$ are of the ``shuffle=stuffle'' form $Z(u\shuffle v)=Z(u \stuffle v)$ since $Z_\stuffle(u)=Z(u)$ and $Z_\stuffle(v)=Z(v)$. 
A conjecture by Hoffman (cf.~\cite{Hoffman2} which is useful for practical computations in low weight and further discussed in \cite{IKZagier}) posits that the combinations
\begin{align}
    Z_\stuffle (u\stuffle v) -Z(u \shuffle v)
\end{align}
with both $u$ and $v$ convergent or $u=y$ and $v$ convergent suffice to generate the ideal $\mathcal{I}_{\FZ}$.

\end{rmk}

\begin{dfn}
    Let $\I_\Z$ be the ideal of $\Q[Z(w)]$ generated by all algebraic relations 
between real MZVs, as defined in
 (\ref{appMZV.01}). Since these are known to satisfy the regularized double 
shuffle relations, we have the inclusions
\beq
\I_\FZ\subset \I_\Z \subset \Q[Z(w)]\,.
\label{sec1eq.4}
\eeq
The space $\FZ$ of {\em formal MZVs} and the space $\Z$ of {\em real MZVs} are defined~by
\begin{align}
\label{sec1eq.8}
\FZ & \coloneqq \Q[Z(w)]/\I_\FZ\,  ,\nn\\
\Z &\coloneqq \Q[Z(w)]/\I_\Z\,,
\end{align}
so that there is a natural surjection
\beq\label{sec1eq.9}\FZ\rightarrow\!\!\!\!\!\rightarrow \Z \, .
\eeq
\end{dfn}

\vspace{0.4cm}
The space $\FZ$ is generated by the images
of the $Z(w)$ in the quotient modulo $\I_\FZ$, which we denote $\zeta^{\mathfrak{f}}(w)$; these
formal MZVs are subject by definition only to the regularized 
double shuffle relations coming from Definition \ref{dfn:IFZ}.
The elements of the $\Q$-algebra $\Z$ of real MZVs are denoted by $\zeta(w)$. 

The $\Q$-algebra $\FZ$ is weight-graded by definition since all of
its defining relations are weight-graded, while $\Z$ is conjectured but of course not known to 
be weight-graded; if it were, this would imply that all real MZVs are transcendental since any MZV which is algebraic would have to be the root of a $\mathbb Q$ polynomial in which each term would be of different weight. A standard conjecture asserts that the surjection \eqref{sec1eq.9} is an isomorphism.

\subsubsection{The Goncharov--Brown coaction}\label{Gonchcoprod}

Let $\overline{\FZ}$ denote the quotient of $\FZ$ modulo the ideal generated
by $\zeta^{\mathfrak{f}}_2$.
In \cite{goncharov2001multiple, Goncharov:2005sla}, Goncharov introduced a
coproduct on the space of motivic iterated integrals, making it into a Hopf algebra. 

In his thesis, published as \cite{Racinet2002}, Racinet showed that the space $\mathfrak{ds}$ of double-shuffle Lie polynomials is a Lie algebra under the Ihara bracket to be defined below in \eqref{eq:Ihara}, which implies that the universal enveloping algebra $\U\ds$ is a Hopf algebra. In \cite{Goncharov:2005sla}, Goncharov showed that his coproduct is dual to the (Poincar\'e--Birkhoff--Witt) multiplication on $\U\ds$. Therefore, $\overline{\FZ}$, which is the dual of $\U\ds$, is a Hopf algebra under the Goncharov coproduct. These objects and ideas are made fully explicit in the discussion around (\ref{Iharaderbrack}).
Brown's work~\cite{Brown:2011mot} allows us to define an extension of Goncharov's coproduct on the Hopf algebra $\overline\FZ$ to a coaction on the comodule
\beq
\label{def:ffmzv}
\FZ \cong \Q[\zeta^{\mathfrak{f}}_2]\otimes_\Q \overline{\FZ}
\eeq
as in (\ref{cozeta2}) below. See section 5 of \cite{Burmesteretal2024} for an introductory recapitulation of these facts.

There are in fact two different versions of the Goncharov--Brown coaction, which differ from each other only by the order of the tensor factors. We denote them by  
\beq
\begin{cases}
\label{eq:coact}
\Delta^{GB} : \FZ \to  \overline{\FZ}\otimes \FZ\, ,\\
    \Delta_{GB} : \FZ \to \FZ \otimes \overline{\FZ}\, .
    \end{cases}
\eeq
Both versions of the coaction are used regularly in the literature, so that it is important to keep track of which one is being used at all times. In the present paper, as we will specify,  
the coaction $\Delta^{GB}$ is implicitly used in numerous proofs in view of its compatibility with  double-shuffle theory and Hopf-algebra duals. The coaction $\Delta_{GB}$ entering explicit formulae (most notably in section~\ref{sec:falpha}) is used to
remain coherent with the recent literature\footnote{The coaction $\Delta^{GB}$ based on Goncharov's original coproduct was introduced in \cite{goncharov2001multiple, Goncharov:2005sla} and \cite{Brown:2011mot}. The coaction $\Delta_{GB}$ is used in the recent particle-physics, string-theory and mathematics literature such as~\cite{Duhr:2012fh, Drummond:2013vz, Brown:2018omk, Brown:2019jng, Abreu:2022mfk, Mafra:2022wml}.}.

Let us describe the construction of the Goncharov--Brown coaction $\Delta_{GB}$.

\begin{dfn}
\label{dfn:Gonch} 
Let $w$ be a convergent monomial in $x$ and
$y$, i.e.~starting with $x$ and ending with $y$. Write 
$w=x^{k_r-1}y\cdots
x^{k_1-1}y$ to match the monomial notation of
$\zeta^{\mathfrak{f}}_{k_1,\ldots,k_r}$ in (\ref{appMZV.02}), and associate to it the symbol 
\begin{align} 
\label{eq:Itozeta}
 I(0;1,0^{k_1-1},\ldots,1,
0^{k_r-1};1) = \zeta^{\mathfrak{f}}_{k_1,\ldots,k_r} \, .
\end{align} 
Let $n=k_1+\cdots+k_r$ denote the degree of $w$. Visualize the sequence 
$(0;1,0^{k_1-1},\ldots,1,0^{k_r-1};1)$
in order from left to right around a semi-circle as illustrated in Figure~\ref{fig:co}, with the terminal $0$ and $1$ at the outer edges and the middle $n$ points placed in clockwise order along the inner part of the semi-circle. 
To compute the coaction of the symbol $I(0;1,0^{k_1-1},\ldots,1,0^{k_r-1};1)$ associated with $\zeta^{\mathfrak{f}}_{k_1,\ldots,k_r}$,
draw every possible ``polygon'' inside the
half-circle starting with the outer 0 on the left and ending with the outer 1 on the right, with vertices at any subset of the inner letters (including the empty set). In the notation
\beq
(a_1,a_2,\ldots,a_n)=(1,0^{k_1-1},1,0^{k_2-1},\ldots,1,0^{k_r-1})
\label{acoprod}
\eeq
for the middle $n$ points (apart from the outer points 0 and 1), the contributing polygons are parametrized by
subsets 
$\{a_{i_1},a_{i_2},\ldots,a_{i_r}\}$ with $1\leq i_1<i_2< \cdots<i_r\leq n$ and all cardinalities in the range $0\leq r \leq n$; see Figure~\ref{fig:co} for the example of $r=2$. 

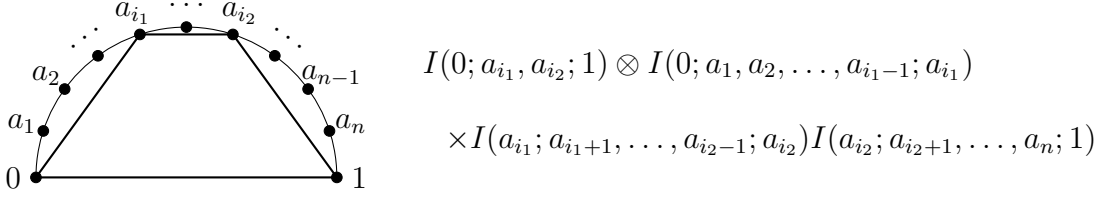
\vspace{.3cm}
\begin{figure}[ht!]
\centering
\begin{tikzpicture}
  \draw[-,draw=black] (2,0) arc (0:180:2);
  \filldraw[black] (-2,0) circle (2pt);
  \filldraw[black] (2,0) circle (2pt);
  \filldraw[black] ({2*cos(1*18)},{2*sin(1*18)}) circle (2pt);  
  \filldraw[black] ({2*cos(2*18)},{2*sin(2*18)}) circle (2pt);  
  \filldraw[black] ({2*cos(3*18)},{2*sin(3*18)}) circle (2pt);  
\filldraw[black] ({2*cos(4*18)},{2*sin(4*18)}) circle (2pt);
\filldraw[black] ({2*cos(5*18)},{2*sin(5*18)}) circle (2pt);  
\filldraw[black] ({2*cos(6*18)},{2*sin(6*18)}) circle (2pt);  
\filldraw[black] ({2*cos(7*18)},{2*sin(7*18)}) circle (2pt);  
\filldraw[black] ({2*cos(8*18)},{2*sin(8*18)}) circle (2pt);  
   \filldraw[black] ({2*cos(9*18)},{2*sin(9*18)}) circle (2pt);  
  \draw (-2.3,0)  node {0};
  \draw ({2.3*cos(9*18)},{2.3*sin(9*18)})  node {$a_1$};
  \draw ({2.3*cos(8*18)},{2.3*sin(8*18)})  node {$a_2$};
  \draw ({2.3*cos(7*18)},{2.3*sin(7*18)})  node[rotate=36] {$\ldots$};
  \draw ({2.3*cos(6*18)},{2.3*sin(6*18)})  node {$a_{i_1}$};
\draw ({2.3*cos(5*18)},{2.3*sin(5*18)})  node {$\ldots$};
  \draw ({2.3*cos(4*18)},{2.3*sin(4*18)})  node {$a_{i_2}$};
  \draw ({2.3*cos(3*18)},{2.3*sin(3*18)})  node[rotate=-36] {$\ldots$};
  \draw ({2.3*cos(2*18)},{2.3*sin(2*18)})  node {$\ a_{n-1}$};
  \draw ({2.3*cos(1*18)},{2.3*sin(1*18)})  node {$a_n$};
  \draw (2.3,0)  node {1};  
    \draw[black,thick] (-2,0) -- (2,0);      
  \draw[black,thick] (-2,0) --  ({2*cos(6*18)},{2*sin(6*18)});  
  \draw[black,thick] ({2*cos(6*18)},{2*sin(6*18)}) --  ({2*cos(4*18)},{2*sin(4*18)});    
  \draw[black,thick] ({2*cos(4*18)},{2*sin(4*18)}) --  (2,0);
  \draw (3,1.5) node[anchor=west] {$I(0;a_{i_1},a_{i_2};1) \otimes I(0;a_1,a_2,\ldots,a_{i_1-1};a_{i_1})$};
  \draw (3.3,0.5) node[anchor=west] {$\cdot\, 
  I(a_{i_1};a_{i_1+1},\ldots,a_{i_2-1};a_{i_2})
  I(a_{i_2};a_{i_2+1},\ldots,a_{n};1)$};
\end{tikzpicture}
\caption{\label{fig:co}\textit{Contributions to the coaction formula (\ref{Gonchco}) for $\Delta_{GB} I(0;a_{1},a_{2},\ldots,a_n;1)$ from polygons with inner vertices $a_{i_1},a_{i_2}$, i.e.\ quadrilaterals associated with subsets of $a_1,a_2,\ldots,a_n$ of cardinality $r=2$.}}
\end{figure}

The coaction is computed by adding up the contributions of all possible polygons: 
\begin{align}
&\Delta_{GB} I(0;a_{1},a_{2},\ldots,a_{n};1) = \sum_{r=0}^n \sum_{1\leq i_1 < i_2 < \ldots < i_r \leq n}
\! \! \! \! \! \! \! \! \! \! \! I(0;a_{i_1},a_{i_2},\ldots,a_{i_r};1)\otimes 
I(0;a_{1},a_2,\ldots, a_{{i_1}-1};a_{i_1}) \notag \\
&\quad \!\!\cdot
I(a_{i_1};a_{{i_1}+1},\ldots, a_{{i_2}-1};a_{i_2})
\cdots
I(a_{i_{r-1}};a_{{i_{r-1}}+1},\ldots, a_{{i_r}-1};a_{i_r})
 I(a_{i_r}; a_{{i_r}+1},\ldots,a_n;1)\, ,
\label{Gonchco}
\end{align}
where $I(0;a_{i_1},\ldots,a_{i_r};1)$ specializes to $I(0;1)=1$ in case of the empty subset at $r=0$ and $I(0; A;1) I(0; B;1)$ (or $I(0; A;1) \cdot I(0; B;1)$) is the shuffle product.
We simplify the expression (\ref{Gonchco}) according to the following rules:

\begin{itemize}[itemsep=0mm]

\item
$I(a;b)=1$ for all $a,b\in \{0,1\}$,

\item 
$I(a;b;c)=0$ for all $a,b,c\in \{0,1\}$,

\item 
$I(a;S;a)=0$ for $a\in \{0,1\}$ and any non-empty sequence $S$ of $0's$ and $1$'s,

\item
$I(1;S;0)=(-1)^nI(0;\overleftarrow{S};1)$ if $S$ is a sequence of $0$'s and
$1$'s of length $n$ and $\overleftarrow{S}$ denotes the sequence $S$
in the reversed order.
\end{itemize} 
We can also replace each term $I(0;S;1)$ by the formal (shuffle-regularized) MZV $\zeta^{\mathfrak{f}}(w_S)$, where if $S$ is
any sequence of $0$'s and $1$'s then $w_S$ is the monomial obtained by
reversing the order of $S$ and replacing every $0$ with an $x$ and every $1$
with a $y$. We finally  project the entries of the second factor of the tensor product modulo $\zeta^{\mathfrak{f}}(xy)=\zeta^{\mathfrak{f}}_2$ to $\overline{\FZ}$, so that the Goncharov--Brown coaction  takes values in $\FZ\otimes \overline{\FZ}$ as announced in \eqref{eq:coact}.
\end{dfn}

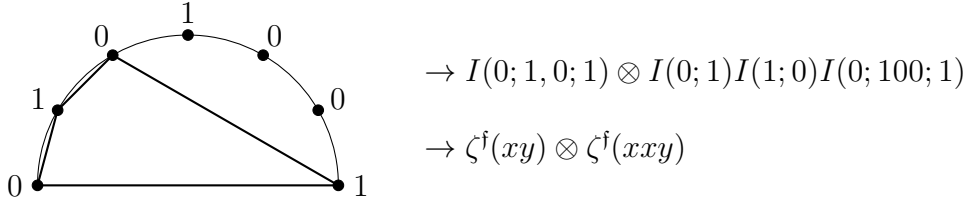
\begin{figure}[ht!]
\centering
\begin{tikzpicture}
  \draw[-,draw=black] (2,0) arc (0:180:2);
  \filldraw[black] (-2,0) circle (2pt);
  \filldraw[black] (2,0) circle (2pt);
  \filldraw[black] ({2*cos(5*30)},{2*sin(5*30)}) circle (2pt);  
  \filldraw[black] ({2*cos(4*30)},{2*sin(4*30)}) circle (2pt);  
  \filldraw[black] ({2*cos(3*30)},{2*sin(3*30)}) circle (2pt);    
  \filldraw[black] ({2*cos(2*30)},{2*sin(2*30)}) circle (2pt);    
  \filldraw[black] ({2*cos(1*30)},{2*sin(1*30)}) circle (2pt);
  \draw (-2.3,0)  node {0};
  \draw ({2.3*cos(5*30)},{2.3*sin(5*30)})  node {1};
  \draw ({2.3*cos(4*30)},{2.3*sin(4*30)})  node {0};
  \draw ({2.3*cos(3*30)},{2.3*sin(3*30)})  node {1};
  \draw ({2.3*cos(2*30)},{2.3*sin(2*30)})  node {0};
  \draw ({2.3*cos(1*30)},{2.3*sin(1*30)})  node {0};
  \draw (2.3,0)  node {1};          
  \draw[black,thick] (-2,0) -- (2,0);      
  \draw[black,thick] (-2,0) --  ({2*cos(5*30)},{2*sin(5*30)});  
  \draw[black,thick] ({2*cos(5*30)},{2*sin(5*30)}) --  ({2*cos(4*30)},{2*sin(4*30)});    
  \draw[black,thick] ({2*cos(4*30)},{2*sin(4*30)}) --  (2,0);
  \draw (3,1.5) node[anchor=west] {$\rightarrow  I(0;1,0;1) \otimes I(0;1)I(1;0)I(0;100;1)$};
  \draw (3,0.5) node[anchor=west] {$\rightarrow \zeta^{\mathfrak{f}}(xy)\otimes \zeta^{\mathfrak{f}}(xxy)$};
\end{tikzpicture}
\caption{
\label{fig:coact}
\textit{Example of a contribution to  
$\Delta_{GB} \zeta^{\mathfrak{f}}(xxyxy)$ as computed in~\eqref{eq:coex}. }}
\end{figure}

\vspace{0.2cm}\noindent
{\bf Example.}
The coaction on the convergent word $\zeta^{\mathfrak{f}}(xxyxy)$ is computed from the semi-circle drawn in Figure~\ref{fig:coact}, which shows one example of a contribution from a quadrilateral. The total result of the coaction is given by
\begin{align}
    \label{eq:coex}
    \Delta_{GB} \zeta^{\mathfrak{f}}(xxyxy) = 1 \otimes \zeta^{\mathfrak{f}}(xxyxy) + \zeta^{\mathfrak{f}}(xxyxy) \otimes 1 +3 \zeta^{\mathfrak{f}}(xy)  \otimes \zeta^{\mathfrak{f}}(xxy)\,.
\end{align}
The first term comes from the degenerate polygon consisting of the straight line from the outer 0 to the outer 1 with no inner vertices and the second to the  full polygon touching all the inner vertices. 
The term with factor $3$ arises from quadrilaterals involving the earliest $1$ (in clockwise direction) of the type shown in Figure~\ref{fig:coact}, and there are three such quadrilaterals which produce the same non-vanishing contribution.
All other polygons have  a vanishing contribution; in particular the polygons going from $0$ directly to the $1$ at the top  produce a $\zeta^{\mathfrak{f}}(xy)=\zeta^{\mathfrak{f}}_2$ to the right of the tensor product $\otimes$, which is projected to zero.

\vspace{0.4cm}

\begin{dfn} The coaction $\Delta^{GB}$ is obtained from $\Delta_{GB}$ by the identity
\beq\label{eq:coactb}
\Delta^{GB}=\iota\circ \Delta_{GB}\,,
\eeq
where $\iota$ exchanges the two tensor factors
\begin{align}
    \iota:\FZ\otimes \overline{\FZ}&\mapsto \overline{\FZ}\otimes \FZ \, ,\notag\\
    \alpha\otimes \beta&\mapsto \beta\otimes \alpha \, .
\end{align}
Reducing the $\FZ$ factor mod $\zeta^{\mathfrak{f}}_2$ in not just one but both factors of the image yields two coproducts 
\beq
\Delta_G,\Delta^G:\overline{\FZ}\rightarrow \overline{\FZ}\,,
\eeq
each of which confers a Hopf algebra structure on $\overline{\FZ}$. We will study the Hopf algebra $\overline{\FZ}$ equipped with $\Delta^G$ and its dual Hopf algebra $\overline{\FZ}^\vee$ further in section \ref{sec:Zmap.1}.
\end{dfn}

\subsection{Motivic MZVs}
\label{sec:MMZV}

We shall here recall the $\Q$-algebra of
{\em motivic MZVs}, which were constructed and studied in depth as a
subcategory of the category of mixed Tate motives ($MTM$) unramified over $\mathbb{Z}$ 
by Deligne, Goncharov, Manin and others, until Brown proved that 
the subcategory is equal to the full category (see \cite{Brown:2011mot}). Our definition follows from 
Brown's results. For further reading, see \cite{BurgosFresan,DeligneMultizetas,dupont2024motives}.
The motivic Lie algebra associated to motivic MZVs is known to be freely generated by one element in each odd degree; a fact that we will make crucial use of in this paper.

\subsubsection{Motivic versus formal multizetas, coproduct and coaction}

Let $\overline{\MZ}$ denote the space of motivic multizetas (modulo $\zeta^\mathfrak{m}_2$) as defined by Goncharov~\cite{Goncharov:2005sla}. It is known to be a Hopf algebra equipped with the coproduct that we introduced in section~\ref{Gonchcoprod}.

Goncharov showed that these motivic MZVs satisfy the double shuffle relations and surject via the period map to the quotient $\overline\Z\coloneq \Z/\langle \zeta_2\rangle $ of the $\mathbb Q$-algebra of real multizetas by the ideal generated by $\zeta_2$. We therefore have the surjections
\beq\label{secondcondition}
\overline{\FZ}\rightarrow\!\!\!\!\!\rightarrow\overline{\MZ} 
\rightarrow\!\!\!\!\!\rightarrow \overline{\Z} \, . 
\eeq
Let $\MZ$ be the Hopf-algebra comodule of motivic multizeta values defined by Brown~\cite{Brown:2011mot}; he showed that it has the structure
\beq
\label{def:mmzv}
\MZ \cong \Q[\zeta^{\mathfrak{m}}_2]\otimes_\Q \overline{\MZ}\,.
\eeq
This allows us to lift the first surjection in~\eqref{secondcondition} to a surjection from the space of formal multizetas $\mathcal{FZ}$ thanks to (\ref{def:ffmzv}) to the space of motivic multizetas $\mathcal{MZ}$ as follows:
for every word $w$ in $x,y$ of length $>2$, we map $\zeta^{\mathfrak{f}}(w)
\in \overline{\FZ}$ to $\zeta^{\mathfrak{m}}(w) \in \overline{\MZ}$ and similarly map $\zeta^{\mathfrak{f}}_2=\zeta^{\mathfrak{f}}(xy)$ to $\zeta_2^{\mathfrak{m}}=\zeta^{\mathfrak{m}}(xy)$. 

For any sequence $S$ of letters $0,1$, let $I^{\mathfrak{m}}(0;S;1)$ denote the image of the $I(0;S;1) \in \FZ$ in section \ref{Gonchcoprod}. The motivic MZVs surject down to the real MZVs by the period 
map 
\beq
I^{\mathfrak{m}}(0;S;1)\mapsto \zeta(w_S)
\eeq
(see \cite{Brown:2011mot}) with $w_S$ as in section \ref{Gonchcoprod}, so \eqref{secondcondition} lifts to the following sequence of $\Q$-algebra surjections
\beq
\FZ\rightarrow\!\!\!\!\!\rightarrow \MZ\rightarrow\!\!\!\!\!\rightarrow \Z \, ,
\eeq
with conjectured equality. 

Like $\FZ$, the Hopf algebra comodule $\MZ$ is 
graded by the weight of the MZVs, as is $\overline{\MZ}$.  We write
$\MZ_w$ (resp.~$\overline{\MZ}_w$\, $\FZ_w$, $\overline{\FZ}_w$) for the weight $w$ part of $\MZ$ 
(resp.~$\overline{\MZ}$\, $\FZ$, $\overline{\FZ}$). Note that we have
\begin{align}
\FZ_0&=\overline{\FZ}_0=\MZ_0=\overline{\MZ_0}=\Q\,,\\
\FZ_1&=\overline{\FZ}_1=\MZ_1=\overline{\MZ_1}=\{0\}\,.
\nn
\end{align}

The coactions $\Delta^{GB}$ and  $\Delta_{GB}$ reviewed in section \ref{Gonchcoprod} both descend directly to $\MZ$. Let us recall the notation for
$\Delta_{GB}$; it is identical to $\Delta^{GB}$ in (\ref{eq:coact}) up to exchanging the two factors of the tensor product. 
The descended coaction~\cite{Brown:2011mot} 
\begin{align}
    \Delta_{GB} : \MZ \to \MZ \otimes \overline{\MZ}\,,
\end{align}
makes $\MZ$ into a Hopf algebra comodule.
In particular we have 
\beq
\Delta_{GB}\bigl(\zeta^{\mathfrak{m}}_2\bigr)=\zeta^{\mathfrak{m}}_2\otimes 1\, .
\label{cozeta2}
\eeq
In analogy with~\eqref{eq:Itozeta} we write $\zeta^{\mathfrak{m}}_{k_1,\ldots,k_r}=I^{\mathfrak{m}}(0;1,0^{k_1-1},\ldots,1,
0^{k_r-1};1)\in \MZ$. We also use the notation  $\zeta^{\mathfrak{dr}}_{k_1,\ldots,k_r}=I^{\mathfrak{dr}}(0;1,0^{k_1-1},\ldots,1,
0^{k_r-1};1)\in \overline{\MZ}$ for the second tensor factor of $\Delta_{GB}$ whose reduction modulo $\zeta_2$ translates into $\zeta^{\mathfrak{dr}}_2=0$.\footnote{The superscript in $\zeta^{\mathfrak{dr}}$ refers to de Rham periods \cite{Goncharov:2005sla, Brown:2011mot, Brown2014MotivicPA, brown2015notes}. In the motivic coaction $\Delta_{GB}$, de Rham periods occur in the right entry of tensor products $A \otimes B$, i.e.\ $B$ is considered modulo $i\pi$. This is opposite to the coaction $\Delta^{GB}$ of \cite{Brown:2011mot} where de Rham periods are in the left entry, so that \eqref{cozeta2} would instead read $\Delta^{GB}(\zeta^{\mathfrak{m}}_2)=1\otimes \zeta^{\mathfrak{m}}_2$.}
The explicit form of the coaction for motivic MZVs $\zeta^{\mathfrak{m}}_{k_1,\ldots,k_r}$ is  encoded in symbols exactly as in (\ref{Gonchco}): we write
\begin{align}
&\Delta_{GB} I^{\mathfrak{m}}(0;a_{1},a_{2},\ldots,a_{n};1) = \sum_{r=0}^n \sum_{1\leq i_1 < i_2 < \ldots < i_r \leq n}
\! \! \! \! \! \! I^{\mathfrak{m}}(0;a_{i_1},\ldots,a_{i_r};1)\otimes 
I^{\mathfrak{dr}}(0;a_{1},\ldots, a_{{i_1}-1};a_{i_1}) \notag \\
& \quad\cdot
I^{\mathfrak{dr}}(a_{i_1};a_{{i_1}+1},\ldots, a_{{i_2}-1};a_{i_2})
\cdots
I^{\mathfrak{dr}}(a_{i_{r-1}};a_{{i_{r-1}}+1},\ldots, a_{{i_r}-1};a_{i_r})
I^{\mathfrak{dr}}(a_{i_r}; a_{{i_r}+1},\ldots,a_n;1)\, ,
\label{brco} 
\end{align}
where the rules detailed below (\ref{Gonchco}) apply in identical form to the terms $I^{\mathfrak{m}}$ and $I^{\mathfrak{dr}}$ on the right-hand side of \eqref{brco} and can be used to put all terms into the standard form $I^{\mathfrak{m}}(0;S;1)$ and $I^{\mathfrak{dr}}(0;S;1)$ for finite tuples $S$ of $0$'s and $1$'s.

\vspace{.3cm}\noindent
{\bf Examples.} 
When $w=x^{n-1}y$ for odd values of $n=2k+1$, the only polygons with a non-zero contribution are the degenerate one (going directly from $0$ to $1$) and the full polygon including every point on the semi-circle: thus we have
\begin{align}
    \Delta_{GB}\zeta^{\mathfrak{m}}_{2k+1} &= \zeta^{\mathfrak{m}}_{2k+1} \otimes 1 + 1\otimes \zeta^{\mathfrak{dr}}_{2k+1}\in {\MZ}\otimes \overline{\MZ}\, .
    \label{primcoact}
\end{align}
Such elements are said to be {\em primitive} for the coproduct. The counterparts of (\ref{primcoact}) for $w=x^{n-1}y$ at even $n=2k$ simplifies to $\Delta_{GB} \zeta^{\mathfrak{m}}_{2k} = \zeta^{\mathfrak{m}}_{2k}\otimes 1$ by the vanishing of $\zeta^{\mathfrak{dr}}_{2k}$.

We also give a few other illustrative instances:
\begin{align}
\label{eq:DG8}
    \Delta_{GB}(\zeta^{\mathfrak{m}}_3 \zeta^{\mathfrak{m}}_5) &=  \zeta^{\mathfrak{m}}_3 \zeta^{\mathfrak{m}}_5 \otimes 1 + 1\otimes \zeta^{\mathfrak{dr}}_3 \zeta^{\mathfrak{dr}}_5 + \zeta^{\mathfrak{m}}_3\otimes  \zeta^{\mathfrak{dr}}_5 + \zeta^{\mathfrak{m}}_5 \otimes \zeta^{\mathfrak{dr}}_3\,,\nn\\
    \Delta_{GB}(\zeta^{\mathfrak{m}}_{3,5}) &= \zeta^{\mathfrak{m}}_{3,5} \otimes 1 + 1\otimes \zeta^{\mathfrak{dr}}_{3,5} -5\, \zeta^{\mathfrak{m}}_3 \otimes \zeta^{\mathfrak{dr}}_5\,,\\
    \Delta_{GB}(\zeta^{\mathfrak{m}}_{2,6}) &= \zeta^{\mathfrak{m}}_{2,6} \otimes 1 + 1\otimes \zeta^{\mathfrak{dr}}_{2,6} +4\, \zeta^{\mathfrak{m}}_3 \otimes \zeta^{\mathfrak{dr}}_5+2\, \zeta^{\mathfrak{m}}_5 \otimes \zeta^{\mathfrak{dr}}_3\,.\nn
\end{align}
These relations are compatible with 
\begin{align}
    \zeta^{\mathfrak{m}}_{2,6} =  - \frac25 \zeta^{\mathfrak{m}}_{3,5} + 2\zeta^{\mathfrak{m}}_3 \zeta^{\mathfrak{m}}_5  - \frac{42}{125} (\zeta^{\mathfrak{m}}_2)^4\,,
\end{align}
where one has to use that the second entries of tensor products in $\MZ\otimes\overline{\MZ}$ are automatically projected modulo $\zeta^{\mathfrak{m}}_2$, so that $  \zeta^{\mathfrak{dr}}_{2,6} = -\frac25   \zeta^{\mathfrak{dr}}_{3,5} + 2   \zeta^{\mathfrak{dr}}_3\zeta^{\mathfrak{dr}}_5$. 

\subsubsection{Reducible motivic MZVs}\label{sec:Background.2}

General references for the following subsection are \cite{ZigZag,Furusho2000TheMZ}. Let $\fz$ denote the quotient of the $\Q$-algebra $\FZ$ given by
\beq\label{fz}
\fz\coloneqq \FZ/\bigl(\FZ_0\oplus \FZ_2\oplus (\FZ_{>0})^2\bigr)=
\overline{\FZ}/\bigl(\overline{\FZ}_0\oplus (\overline{\FZ}_{>0})^2\bigr)\, ,
\eeq
and analogously, let $\mz$ denote the quotient of the $\Q$-algebra $\MZ$ given by
\beq
\mz\coloneqq \MZ/\bigl(\MZ_0\oplus \MZ_2\oplus (\MZ_{>0})^2\bigr)=
\overline{\MZ}/\bigl(\overline{\MZ}_0\oplus (\overline{\MZ}_{>0})^2\bigr)\, .
\label{defmz}
\eeq
From the Hopf algebra structure on $\overline{\FZ}$ (resp.~$\overline{\MZ}$), the vector space $\fz$ (resp.~$\mz$) inherits the structure of a 
Lie coalgebra, dual to the Lie algebras that will be introduced in section~\ref{sec:Zmap.1}.
Note that by~\eqref{fz} and~\eqref{defmz}, the element $\zeta^{\mathfrak{f}}_2$ (resp.~$\zeta^{\mathfrak{m}}_2$) maps down to zero in $\fz$ (resp.~$\mz$).

\vspace{.2cm}
\begin{dfn}
\label{dfn:z2n} 
    For all even positive integers $w=2n$, let ${\rm B}_{2n}$ be the Bernoulli number, and set \cite{Brown:2011mot}
\beq\label{evens}
\zeta^{\mathfrak{m}}_{2n}  \coloneqq \frac{\zeta_{2n}}{\zeta_2^{n}}(\zeta^{\mathfrak{m}}_2)^{n}= (-1)^{n-1} \frac{(24)^n {\rm B}_{2n}}{2 (2n)!} \, (\zeta^{\mathfrak{m}}_2)^n\in \MZ_{2n}\, .
\eeq
\end{dfn}

\vspace{.3cm}
\begin{dfn}
\label{dfn:red}
For all $w\ge 3$, let $\hat R_w$ 
denote the canonical subspace of {\em reducible MZVs} 
in $\MZ_w$. The space $\hat R_w$ is the subspace generated by all total-weight $w$ products of lower-weight MZVs, or in other words by all weight 
$w$ elements of $({\MZ}_{>0})^2$. 
Note that $\hat R_3=\{0\}$, so there are actually non-trivial reducible 
subspaces only for $w\ge 4$, starting with $\hat R_4=\mathbb Q \zeta^{\mathfrak{m}}_4$ and $\hat R_5 = \mathbb Q \zeta^{\mathfrak{m}}_2 \zeta^{\mathfrak{m}}_3$.

\end{dfn}

\vspace{.2cm}
The Lie coalgebras $\fz$ and $\mz$ are weight-graded, and for each weight $w>1$ 
we have
\beq\label{fzandmz}
\fz_w=\FZ_w/\hat{R}_w\, , \qquad
\mz_w=\MZ_w/\hat{R}_w\, .
\eeq

\subsubsection{Irreducible MZVs}\label{sec:Background.3}

Let $\hat I_w$ be any supplementary subspace of $\hat R_w$ in $\MZ_w$ so that
\beq
\MZ_w=\hat R_w\oplus \hat I_w\, .
\eeq
Since the map $\MZ_w\rightarrow \mz_w$ is the quotient mod $\hat R_w$,
it induces an isomorphism $\hat I_w\rightarrow \mz_w$.
We will always choose $\hat I_w$ containing
$\zeta^{\mathfrak{m}}_w$ if $w$ is odd, see Lemma 3.2 of \cite{Brown:2011mot}. If $w$ is even, we set $I_w \coloneqq \hat I_w$, and if $w$ is odd we choose a supplementary subspace $I_w$ in $\hat I_w$ such that
$\hat I_w=\Q\zeta^{\mathfrak{m}}_w\oplus I_w$. Similarly, if $w$ is odd we set $R_w \coloneqq\hat R_w$ and if $w$ is even we choose a supplementary subspace $R_w\subset \hat R_w$ such that
$\hat R_w=\Q\zeta^{\mathfrak{m}}_w\oplus R_w.$
Then for all $w\ge 2$ we have the direct sum decomposition
\beq
\MZ_w=\Q\zeta^{\mathfrak{m}}_w\oplus I_w\oplus R_w\, .
\label{sec2eq.2}
\eeq

\section{The Z-map associating polynomials to MZVs}\label{sec:Zmap}

In this section we will introduce the {\em Z-map} (see \cite{Schneps:2013}), which provides a
family of canonical isomorphisms
between the MZV spaces studied in section \ref{sec:Background} (namely $\FZ$, $\overline{\FZ}$, $\MZ$, $\overline{\MZ}$, $\Z$, $\fz$ or $\mz$) 
and 
their dual spaces. As will be detailed in section \ref{sec:Zmap.1}, since all the MZV spaces are quotients of $\Q[Z(w)]$, all of their duals are subspaces of $\Q[Z(w)]^\vee$, which is nothing other than the polynomial algebra
$\Q\langle x,y\rangle$ in the non-commutative variables $x$ and $y$. 

Thanks to the fact that the double shuffle relations generate all relations satisfied by $\FZ$ (and in their linearized version, $\fz$), we can give an explicit description of the elements of the dual spaces $\FZ^\vee$ and $\fz^\vee$ in $\Q\langle x,y\rangle$. In the case of motivic and real MZVs we do not have an explicit description of this type since they may satisfy further, unknown relations. Still, thanks to Brown's theorem in \cite{Brown:2011mot}, we do know the structure and dimensions of the graded parts of the dual spaces $\MZ^\vee$ and $\mz^\vee$, which allows us to compute their elements explicitly in low weights (see section \ref{sec:Zmap.4a}).

\subsection{The double shuffle dual space of formal MZVs}\label{sec:Zmap.1}

Let $\Q\langle x,y\rangle$ denote the polynomial ring in two non-commutative variables $x$ and $y$,
equipped with its canonical basis of monomials $w$ in $x$ and $y$ (including the constant monomial ${\bf 1}$),
and let $\Q\langle\langle x,y\rangle\rangle$ denote its degree-completion, the power series ring in
$x$ and $y$. The space $\Q[Z(w)]$ introduced in section~\ref{sec:FMZV} can be identified with the graded dual of
$\Q\langle x,y\rangle$, equipped with the dual basis of symbols $Z(w)$ 
such that
\beq
\langle Z(u),v\rangle=\delta_{u,v}\, ,
\label{eq:Zmap.01}
\eeq
on monomials $u$ and $v$ and extended linearly to give a canonical pairing between $\Q\langle x,y\rangle$
and $\Q[Z(w)]$.

Recall from \eqref{sec1eq.8} that $\FZ$ is the quotient of $\Q[Z(w)]$ by the ideal $\I_{\FZ}$. The graded dual space $\FZ^\vee$ is thus the subspace of
$\Q\langle x,y\rangle$ that annihilates the elements of $\I_\FZ$; explicitly, $\FZ^\vee
\subset \Q\langle x,y\rangle$
is a weight-graded space in which $\FZ_0^\vee=\Q$,
$\FZ^\vee_1=0$ and for $w\ge 2$, $\FZ^\vee_w$ consists of all degree $w$ homogeneous 
polynomials $P\in \Q\langle x,y\rangle$ satisfying 
\beq
\langle L,P\rangle=0\ \ \hbox{for all }L\in\I_\FZ\, ,
\label{eq:Zmap.02}
\eeq
for the pairing in \eqref{eq:Zmap.01}
(see Definition~\ref{dfn:IFZ} for an explicit description of the elements $L$ of the ideal $\I_\FZ$). The subspace $\FZ^\vee$ is strictly smaller than $\Q\langle x,y\rangle$. In weight $w=2$, for
instance, since $Z(xy)+Z(yx) \in \I_{\FZ}$, we have $xy-yx \in \FZ^\vee_2$ 
whereas $xy$ and $yx$ are not individually contained in $\FZ^\vee_2$.

Similarly, the dual space of the quotient $\overline{\FZ}$ of $\FZ$ modulo $\zeta^{\mathfrak{f}}_2$
is a subspace 
$\overline{\FZ}^\vee\subset \FZ^\vee$.
We now consider $\overline{\FZ}$ with its Hopf algebra structure given by the coproduct $\Delta^G$; then the dual space $\overline{\FZ}^\vee$ is also a Hopf algebra. The coproduct on $\overline{\FZ}^\vee$ is inherited directly from the standard coproduct $\Delta_s$ on $\Q\langle x,y\rangle$, given by 
\beq\label{eq:Zmap.2a}
\Delta_s(x)=x\otimes {\bf 1}+{\bf 1}\otimes x\, , \ \ \ \ 
\Delta_s(y)=y\otimes {\bf 1}+{\bf 1}\otimes y\,;
\eeq
it satisfies 
\beq\label{eq:Zmap.2aa}
\langle \xi_1\otimes \xi_2, \Delta_s(g)\rangle=\langle \xi_1\shuffle \xi_2, g\rangle
\eeq
for $g\in \overline{\FZ}^\vee$, $\xi_1,\xi_2\in \overline{\FZ}$.
The multiplication on $\overline{\FZ}^\vee$, which we denote by $\circp$, is uniquely determined by the equality
\beq\label{eq:Zmap.2ab}
\langle \Delta^G(\xi),g\otimes h\rangle=\langle \xi,g\circp h\rangle
\eeq
for $\xi\in \overline{\FZ}$ and $g,h\in \overline{\FZ}^\vee$, see Prop.~3.18 of \cite{Burmesteretal2024} for a fully explicit proof of this duality formula. An explicit formula for $g\circp h$ in the restricted case of $g \in \fz^\vee$ can be found in (\ref{starmult}) below; the full formula (which we do not need here) can be
found in (23) of \cite{Burmesteretal2024}.

Let us now explain how to view $\overline{\FZ}^\vee$ as the universal enveloping algebra of the Lie algebra consisting of its primitive elements.
We begin by identifying
the subspace ${\rm Lie}[x,y]$ of Lie polynomials in $\Q\langle x,y\rangle$ as the subspace of primitive elements, which are those satisfying
\beq\label{eq:Zmap.2b}\Delta_s(g)=g\otimes {\bf 1}+{\bf 1}\otimes g\, .
\eeq
An equivalent formulation of this property is that $g$ is a Lie polynomial in $\Q\langle x,y\rangle$ if and only if 
\beq\label{eq:Zmap.2c}\langle Z(u\shuffle v)\,,\,g\rangle=0\,,
\eeq
for all pairs of non-empty words $u,v$. The Lie subalgebra of the Hopf algebra $\overline{\FZ}^\vee$ is likewise the space of elements $g\in \overline{\FZ}^\vee$ satisfying \eqref{eq:Zmap.2b}; the
Lie bracket is given by
\beq\label{eq:Zmap.03}
\{g,h\} \coloneqq g\circp h-h\circp g\,,
\eeq
for the multiplication $\circp$ of \eqref{eq:Zmap.2ab}. 

This Lie algebra is identified with the dual of the space
$\fz$ defined in \eqref{fz} above; indeed, the vector space $\fz$ inherits
the structure of a Lie coalgebra from the Hopf algebra structure on 
$\overline{\FZ}$, so its dual space $\fz^\vee\subset \overline{\FZ}^\vee$
thus forms a Lie algebra, which is precisely the Lie algebra of primitive elements of $\overline{\FZ}^\vee$.

Since $\fz$ is the quotient of $\overline{\FZ}$ modulo non-trivial products and the relations
\begin{align}\label{fulldoubleshuffle}
\zeta^{\mathfrak{f}}(u)\zeta^{\mathfrak{f}}(v)&=\zeta^{\mathfrak{f}}(u\shuffle v) \, ,\notag\\
\zeta^{\mathfrak{f}}_\stuffle(u)\zeta^{\mathfrak{f}}_\stuffle(v)&=\zeta^{\mathfrak{f}}_\stuffle(u\stuffle v)\,,
\end{align}
hold in $\FZ$ (the second equality being valid whenever $u,v$ both end in $y$), we see that the images of $\zeta^{\mathfrak{f}}(w)$ 
in the quotient $\fz$ satisfy
\beq\label{fzvee}
\zeta^{\mathfrak{f}}(u\shuffle v)=\zeta^{\mathfrak{f}}_\stuffle(u\stuffle v)=0 \,\ \ \ {\rm in}\ \fz\,.
\eeq
Thus the dual space $\fz^\vee$ is the subspace of polynomials $g\in \Q\langle x,y\rangle$
such that
\beq
\langle Z(u\shuffle v),g\rangle=
\langle Z_\stuffle (u\stuffle v),g_\stuffle\rangle=0\,,
\eeq
for all pairs of monomials $u$ and $v$ (ending in $y$ for the $\stuffle$ term), where
\beq
g_\stuffle=g+\sum_{n\ge 2} \frac{(-1)^{n-1}}{n}\zeta^{\mathfrak{f}}(x^{n-1}y) y^n\,, 
\eeq
(the term added to $g$ is the linearized version of \eqref{appMZV.06ter}).
We note in particular that by \eqref{eq:Zmap.2c}, the first equality $\langle Z(u\shuffle v),g\rangle=0$ shows that we have an inclusion of vector spaces (which is not a Lie algebra  morphism as the brackets are different)
\beq\label{eq:Zmap.2d}\fz^\vee\subset {\rm Lie}[x,y]\,.\eeq

The Lie algebra $\fz^\vee$ is known as the {\em double shuffle Lie algebra} and usually
denoted by $\ds$ for ``double shuffle'' (or $\dmr$ for ``double m\'elange r\'egularis\'e''
by French authors). The Lie bracket $\{\cdot,\cdot\}$ on $\ds$ corresponds to the Ihara bracket
\beq
\label{eq:Ihara}
\{g,h\}=[g,h]+D_g(h)-D_h(g)\,,
\eeq
where for each $g\in {\rm Lie}[x,y]$, the {\em Ihara derivation} $D_g$ of ${\rm Lie}[x,y]$
is defined by
\beq
\label{eq:IharaD}
D_g(x)=0 \, , \ \ \ \ D_g(y)=[y,g]\,,
\eeq
and the Lie bracket arises from the bracket of derivations
\beq\label{Iharaderbrack}
[D_g,D_h]=D_{\{g,h\}}\, .
\eeq
The Hopf algebra $\overline{\FZ}^\vee$ is identified with the universal enveloping algebra $\U\ds$. As such, the multiplication $\circp$ is identified with the Poincar\'e--Birkhoff--Witt multiplication (which exists for every universal enveloping algebra of a Lie algebra).

In the case where $g\in \ds$ and $h\in \U\ds$ the multiplication $\circp$ can be written succinctly as
\begin{equation}\label{starmult}
g\circp h=gh+D_{g}(h) \, ,
\end{equation}
which suffices for our purposes and implies that the two representations (\ref{eq:Zmap.03}) and (\ref{eq:Ihara}) of the Ihara bracket agree. 

In the rest of this article with the exception
of section \ref{sec:falpha}, we will consider the space $\FZ$ as a Hopf algebra comodule
equipped with the coaction $\Delta^{GB}$ over the Hopf algebra $\overline{\FZ}$
equipped with the coproduct $\Delta^G$;
the multiplication $\circp$ extends to $\FZ$ by the identity
\beq\label{starmultduality}\langle \Delta^{GB}(\xi),g\otimes h\rangle = \langle \xi,g\circp h\rangle
\eeq
for $\xi\in \FZ$ and $g,h\in \FZ^\vee$. The quotient space $\overline{\MZ}$ of
$\overline{\FZ}$ is then also a Hopf algebra 
equipped with the coproduct $\Delta^G$, and $\MZ$ equipped with $\Delta^{GB}$ is
a Hopf algebra comodule over it. The dual space 
\beq
\overline{\MZ}^\vee\subset \overline{\FZ}^\vee = {\cal U}\ds
\eeq
of $\overline{\MZ}$ is a Hopf algebra equipped with the standard 
coproduct $\Delta_s$ and the (restriction of the) multiplication $\circp$, and the Lie
algebra 
\beq
\mz^\vee\subset \fz^\vee =  \ds
\eeq
consists of the primitive elements for $\Delta_s$ in $\overline{\MZ}$, and
is equipped with the (restriction of the) Ihara bracket \eqref{eq:Zmap.03}.  

\subsection{The Z-map and dual spaces}\label{sec:Zmap.2}

\begin{dfn}
\label{dfn:Zmap}
We define the {\em Z-map} to be the canonical 
isomorphism
\beq
\xymatrix{\Q\langle x,y\rangle\ar[r]^Z&\Q[Z(w)]}
\label{eq:Zmap.04}
\eeq
mapping ${\bf 1}$ to $1$ and each non-trivial monomial $w$ to $Z(w)$, so that the notation $Z(w)$, previously just a symbol (see section \ref{sec:FMZV}), can now be interpreted as the image of
the monomial $w$ under the map $Z$.
The Z-map restricts to a canonical isomorphism on each (finite-dimensional) weight-graded part,
    and passes to corresponding isomorphisms (also called Z-maps) between
    any quotient of $\Q[Z(w)]$ (in particular the MZV spaces) and its dual viewed as a subspace
    of $\Q\langle x,y\rangle$.
\end{dfn}    
The situation is summarized in \eqref{diag0} below, in which all of the
    horizontal arrows are the canonical isomorphisms inherited from the top Z-map
    \beq
    Z:\Q\langle x,y\rangle\rightarrow \Q[Z(w)]\, ,
    \eeq
    all surjective maps are quotients, and all injective maps are inclusions of the dual spaces. The space $\overline{\Z}$ denotes the quotient of the $\Q$-algebra $\Z$ of real MZVs modulo the ideal generated by $\zeta_2$, and in analogy with $\fz$ and $\mz$, we denote the quotient of $\overline{\Z}$ mod constants and non-trivial products by $\z$. 
    For instance, the Z-map $Z(xy)$ is given by $\zeta_2^{\mathfrak{m}}$ in ${\cal MZ}$ and $0$ in $\overline{{\cal MZ}}$, respectively. More generally, we have
    \beq
Z(x^{k_r-1}y\cdots x^{k_2-1}y x^{k_1-1}y)
=\zeta^{\mathfrak{m}}_{k_1,k_2,\ldots,k_r}\ 
{\rm in} \ {\cal MZ}
    \eeq
for convergent words ($k_r\geq 2$), whereas the Z-map of divergent words follows from setting the combinations in (\ref{sec1eq.2a}) to zero.   
    
    Note that while both $\fz$ and $\mz$ are equipped with a Lie coalgebra structure inherited from the Hopf algebra structures on $\overline{\FZ}$ and $\overline{\MZ}$, we do not know that $\overline{\Z}$ is a Hopf algebra and therefore we do not know that $\z$ has a Lie coalgebra structure. However we still have vector space surjections $\fz\rightarrow\!\!\!\!\!\rightarrow \mz\rightarrow\!\!\!\!\!\rightarrow \z$ and the corresponding vector space inclusions of the dual spaces, all of which lie in the vector space ${\rm Lie}[x,y]$ by \eqref{eq:Zmap.2d}:
\beq
\label{eq:manyz}
\z^\vee\subset \mz^\vee\subset \fz^\vee\subset {\rm Lie}[x,y]\,. 
\eeq
We underline once more that all maps in the following diagram are to be viewed as vector space morphisms.

\vspace{.3cm}
\begin{equation}\label{diag0}
\xymatrixcolsep{3pc}
    \xymatrix{\phantom{\Q[Z(w)]}&\phantom{\Q\langle x,y\rangle}&\phantom{\Q[Z(w)]}&\phantom{\Q\langle x,y\rangle}&\Q[Z(w)]\ar@{->>}[d]&\Q\langle x,y\rangle\ar[l]_Z\\
              &&&&\FZ\ar@{-}[d(0.4)]\ar@{}[dd]^(.3){}="a"^(.9){}="b" \ar@{->>} "a";"b"\ar@{->>}[dll]&\FZ^\vee\ar[l]\ar@{^{(}->}[u]\\
              &&{\ \overline{\FZ}\ }\ar@{->>}[lld]\ar@{-}[d(0.4)]\ar@{}[dd]^(.3){}="a"^(.9){}="b" \ar@{->>} "a";"b"&{\ \overline{\FZ}^\vee\ }\ar[l]\ar@{^{(}->}[urr]&&\\
              {\ \fz\ }\ar@{->>}[dd]&{\ \fz^\vee\ }\ar[l]
              \ar@{^{(}->}[urr] &&&\MZ\ar@{-}[d(0.4)]\ar@{}[dd]^(.3){}="a"^(.9){}="b" \ar@{->>} "a";"b"\ar@{-}[dll(0.45)]\ar@{}[dll]^(.55){}="a"^(.9){}="b" \ar@{->>} "a";"b"&\MZ^\vee\ar[l]\ar@{^{(}->}[uu]\\
              &&{\ \overline{\MZ}\ }\ar@{-}[dll(0.45)]\ar@{}[dll]^(.55){}="a"^(.9){}="b" \ar@{->>} "a";"b"\ar@{-}[d(0.4)]\ar@{}[dd]^(.3){}="a"^(.9){}="b" \ar@{->>} "a";"b"&{\ \overline{\MZ}^\vee\ }\ar[l]\ar@{^{(}->}[urr]\ar@{^{(}->}[uu]&&\\
              {\ \mz\ }\ar@{->>}[dd]&{\ \mz^\vee\ }\ar[l]\ar@{^{(}->}[uu]\ar@{^{(}->}[urr]&&&\Z\ar@{-}[dll(0.45)]\ar@{}[dll]^(.55){}="a"^(.9){}="b" \ar@{->>} "a";"b"&\Z^\vee\ar[l]\ar@{^{(}->}[uu]\\
              &&{\ \overline{\Z}\ }\ar@{-}[dll(0.45)]\ar@{}[dll]^(.55){}="a"^(.9){}="b" \ar@{->>} "a";"b"&{\ \overline{\Z}^\vee\ }\ar[l]\ar@{^{(}->}[uu]\ar@{^{(}->}[urr]&&\\
              {\ \z\ }&{\ \z^\vee\ }\ar[l]\ar@{^{(}->}[uu]\ar@{^{(}->}[urr]&&&&\\
    }
\end{equation}

\vspace{.2cm}
We will make constant use of the Z-maps as well as the quotient maps and inclusions in this diagram
for our constructions below.

\subsection{The canonical decomposition of motivic MZV spaces and zeta generators in genus zero}
\label{sec:Zmap.4}

In this section we will define a specific canonical decomposition of $\MZ_w$ for each weight $w\ge 2$ into singles, irreducibles and reducibles of the type
\beq\label{canondecomp}\MZ_w=\Q\zeta^{\mathfrak{m}}_w\oplus I_w\oplus R_w
\eeq
introduced in \eqref{sec2eq.2}. 
The existence of this decomposition relies on working in the space of motivic multizeta values.
More generally, the main results of this work are stated only for
motivic multizetas as opposed to real multizetas in (\ref{appMZV.01}) since our arguments and proofs crucially rely on the freeness of the Lie algebra $\mz^\vee$ below which is tied to the {\it motivic} incarnation of multizetas.

\vspace{.3cm}
\begin{dfn}
\label{gwcan}
For each $w\ge 2$, let $\hat R_w\subset \MZ_w$ denote the subspace of reducible MZVs as in section \ref{sec:Background.2}, let
$\mz_w=\MZ_w/\hat R_w$ as in \eqref{fzandmz}, let $\mz_w^\vee\subset \MZ_w^\vee$ denote the dual space, and let $(\mz_w^\vee)^{\ge 2} 
\subset \mz_w^\vee$ denote the subspace of $\mz_w^\vee$ consisting of elements of depth $\ge 2$, where depth is the minimal $y$-degree of a polynomial. 
\begin{itemize}
    \item
Define the {\em canonical subspace of non-single irreducibles} $I_w$ of $\MZ_w$ by
\beq\label{defIw}
I_w=Z\bigl((\mz_w^\vee)^{\ge 2}\bigr)\subset \MZ_w\, .
\eeq
\item Define the {\em canonical subspace of non-single reducibles} $R_w$ as follows. For odd weights~$w$, set $R_w=\hat R_w$, and for even weights $w$, let $R_w\subset \hat R_w$ be the subspace 
spanned by all weight $w$ products of the elements: $\zeta^{\mathfrak{m}}_2$, the single zetas $\zeta^{\mathfrak{m}}_v$ for odd $v<w$, and all elements of $I_v$ with
$v<w$, excluding only the product $(\zeta^{\mathfrak{m}}_2)^{w/2}$.  Then since $\MZ =  \Q[\zeta^{\mathfrak{m}}_2]\otimes_\Q \overline{\MZ}$ (cf.~\eqref{def:mmzv}), using \eqref{evens}, 
we have $\hat R_w=\Q\zeta^{\mathfrak{m}}_{w}\oplus R_w$ when $w$ is even.
\item Define the {\em canonical decomposition} of $\MZ_w$ to be
\beq
\MZ_w=\Q\zeta^{\mathfrak{m}}_w\oplus I_w\oplus R_w
\eeq
for the canonical subspaces $R_w$ and $I_w$ defined above. 
\item Finally,
define the {\em canonical polynomial} $g_w\in \MZ_w^\vee$ for each $w\ge 2$ to be the unique polynomial in $x,y$ that 
\begin{itemize}
\item takes the value $1$ on $\zeta^{\mathfrak{m}}_w =\zeta^{\mathfrak{m}}(x^{w-1}y)$ in the sense that $\langle Z(x^{w-1}y),g_w\rangle=1$, and
\item annihilates $I_w$ and $R_w$ in the sense that $\langle \xi,g_w\rangle=0$ for any $\xi \in I_w$ and $\xi\in R_w$. 
\end{itemize}
That such polynomials exist follows from Lemma~\ref{lem:332}, but also from their alternative characterization in terms of the Drinfeld associator, see section~\ref{sec:Drin} below.
\end{itemize}
\end{dfn}

\noindent Examples of the polynomials $g_w$ will be given in section~\ref{sec:Zmap.4a} below.

\begin{lemma}\label{lem:332}
The canonical polynomials $g_w$ for $w\ge 2$ are uniquely characterized by the following properties:
\begin{itemize}
\item[(i)] The polynomial $g_w$ is normalized by $g_w|_{x^{w-1}y}=1$;
\item[(ii)] The polynomial $g_w$ lies in the subspace $(\MZ_w/R_w)^\vee\subset\MZ_w^\vee$; in particular for odd $w$ it lies in $\mz_w$ and is thus a Lie polynomial;
\item[(iii)] If we consider $g_w$ as lying in $(\MZ_w/R_w)^\vee$, the image $Z(g_w)$ of $g_w$ 
under the Z-map is a rational multiple 
of $\zeta^{\mathfrak{m}}_w\in \MZ_w/R_w$;
equivalently, if we consider $g_w$ as lying in $\MZ_w^\vee$, then $Z(g_w)$ does not contain any irreducible multizeta values in $I_w$
\beq
Z(g_w)\in \Q\zeta^{\mathfrak{m}}_w\oplus R_w\subset \MZ_w\,.
\label{char3rd}
\eeq
\end{itemize}
\end{lemma}

\vspace{.2cm}
\noindent {\bf Proof.} (i) is equivalent to $\langle Z(x^{w-1}y),g_w\rangle=1$. 

For (ii), saying that $g_w$ annihilates $R_w$ is equivalent to saying that $g_w$ lies in the dual space of $\MZ_w/R_w$,
namely $(\MZ_w/R_w)^\vee$; this space is equal to $\mz_w^\vee$ when $w$ is odd, so by \eqref{eq:manyz} $g_w$ is then in ${\rm Lie}[x,y]$.

For (iii),  we consider $g_w\in (\MZ_w/R_w)^\vee$ and for $\MZ_w/R_w =  \Q\zeta^{\mathfrak{m}}_w\oplus I_w$ we choose any basis consisting of
$\zeta^{\mathfrak{m}}_w$ and a basis for $I_w$. Then since $\langle g_w,I_w\rangle=0$ for all $\xi\in I_w$ we have $\langle Z(g_w),Z^{-1}(I_w)\rangle=0$, but $Z^{-1}(I_w)= (\mz_w^\vee)^{\ge 2}$, and the subspace of $\MZ_w/R_w$ annihilated by $(\mz_w^\vee)^{\ge 2}$ is
the 1-dimensional subspace generated by $\zeta^{\mathfrak{m}}_w$. Therefore if $g_w$ is considered as lying in $(\MZ_w/R_w)^\vee$ we have $Z(g_w)\in \Q\zeta^{\mathfrak{m}}_w\subset \MZ_w/R_w$, or equivalently, if $g_w$ is considered as lying in $\MZ_w^\vee$, we have
$Z(g_w)\in \Q\zeta^{\mathfrak{m}}_w\oplus R_w$. This construction proves the uniqueness of~$g_w$: the 1-dimensional subspace it generates annihilates the subspace $(\mz_w^\vee)^{\ge 2}$ of non-single zetas in the dual, and the specific choice of $g_w$ is given by the normalization in (i).
\qed

\begin{rmk}
The lemma shows that in order to compute the canonical polynomials $g_w$ for any $w\geq 2$, once conditions (i) and (ii) of Lemma \ref{lem:332} are fulfilled, the third defining condition of $g_w$, namely that it annihilates the subspace $I_w$, can be replaced by condition (iii) of the Lemma, which does not require computing the space $I_w$. Once $g_w$ is determined, it is then possible to recover the space $I_w$ as the image of a Lie subspace of $\MZ^\vee$ under $Z$ as in~\eqref{defIw} if needed. We will actually provide a very natural explicit basis for $I_w$, called {\em the semi-canonical basis}, in section~\ref{sec:Zmap.4b} below.
\end{rmk}

\begin{rmk}
    As mentioned in Remark~\ref{rmk:114} in the introduction, Keilthy has provided a construction of polynomials $g_w$ for odd $w$ in his dissertation~\cite{Keilthy}. The idea of his method is the following: assuming that the $g_w$ have been chosen for odd $w$ up to $2k{-}1$, one can then fix a choice of $g_{2k+1}$ by requiring it to be orthogonal to all weight $2k{+}1$ Ihara brackets of the previously chosen $g_w$; for example, $g_{11}$ is fixed by the unique condition that it must be orthogonal to $\{g_3,\{g_3,g_5\}\}$, where orthogonality is defined by the inner product on all pairs of monomials $u,v$ taking value $\delta_{u,v}$, analogously to our~\eqref{eq:Zmap.01}. 
    This is equivalent to our condition $\langle \xi, g_w\rangle =0$ for all irreducibles $\xi\in I_w$ since a basis of $I_w$ can be given using Ihara brackets of lower degree $g_{w'}$ (i.e.~with $w'<w$). 
    Note however that Lemma \ref{lem:332} provides canonical elements $g_w$ also for even $w$.
\end{rmk}

Given that we know from \cite{Brown:2011mot} that $\mz^\vee$ is free on one depth 1 generator in each odd weight $w\ge 3$ and the $g_w$ are such elements, the  set of $g_w$ for odd $w\ge 3$ form a canonical generating set for $\mz^\vee$.
By Lemma \ref{lem:332}, each $g_w$ is
characterized uniquely as the only depth 1 element of $\mz_w^\vee\subset \MZ^\vee$
normalized by $g_w|_{x^{w-1}y}=1$ such that $Z(g_w) \in \Q\zeta^{\mathfrak{m}}_w\oplus R_w\subset \MZ_w$.
\begin{dfn}
\label{dfn:mzvee}
 The Ihara derivations (\ref{eq:IharaD}) associated with the $g_w$ with $w\geq 3$ odd are referred to as {\em zeta generators in genus zero}. 
\end{dfn}

\vspace{.2cm}
The method of using the Z-map to produce canonical generators by taking the duals of the single zetas was initially developed in the framework of formal multizetas in \cite{Schneps:2013}.
 The family of polynomials $g_w$ will play a crucial role in the main results of this paper, namely
\begin{itemize}
\item the construction of a canonical isomorphism $\rho:\MZ\rightarrow {\cal F}$ from the motivic MZVs to the $f$-alphabet (section \ref{sec:falpha.2});
\item the construction of a canonical set of zeta generators in genus one (section \ref{sec:sigman.2}).
\end{itemize}
\noindent 
In the next subsection we give the explicit calculation of the canonical decomposition in weights up to $w=11$ and spell out the canonical polynomials $g_w$ up to $w=7$. We emphasize that our approach is unaffected by the irregular behaviour of the depth filtration of MZVs starting in weight 12: Our method to construct the canonical polynomials $g_w$ does not depend on the depth of the MZVs encountered in the canonical decompositions for $\MZ_w$, and the explicit form of $g_{12}$ can be found in the ancillary files of the arXiv submission. The canonical morphism to
the $f$-alphabet resulting from the discussion of section \ref{sec:falpha} below and relying on the canonical decomposition of $\MZ_w$ in intermediate steps is explicitly worked out up to and including weight 17 in the ancillary files of \cite{Dorigoni:2024oft}.

\subsection{\texorpdfstring{The canonical decomposition for $\MZ_w$ for $w\le 11$}{The canonical decomposition for MZ(w) for w<=11}}
\label{sec:Zmap.4a}

Since all MZVs in this subsection and the next one are motivic, we drop the superscript $\mathfrak{m}$ and simply write $\zeta_{k_1,\ldots,k_r}$ instead of $\zeta^{\mathfrak{m}}_{k_1,\ldots,k_r}$. We have 
\begin{align}\label{MZdecomp}
\MZ_2&=\langle\zeta_2\rangle \, , \notag\\
\MZ_3&=\langle\zeta_3\rangle \, , \notag\\
\MZ_4&=\langle\zeta_4\rangle\, ,\notag\\
\MZ_5&=\langle\zeta_5\rangle\oplus \langle \zeta_2\zeta_3\rangle=\Q\zeta_5\oplus R_5\, ,\notag\\
\MZ_6&=\langle\zeta_6\rangle\oplus \langle \zeta_3^2\rangle=\Q\zeta_6\oplus R_6\, ,\\
\MZ_7&=\langle\zeta_7\rangle\oplus \langle \zeta_2\zeta_5\,,\ \zeta_2^2\zeta_3\rangle=\Q\zeta_7\oplus R_7\, ,\notag\\
\MZ_8&=\langle\zeta_8\rangle\oplus \langle Z_{35}\rangle\oplus \langle \zeta_3\zeta_5\,,\ \zeta_2\zeta_3^2\rangle =\Q\zeta_8\oplus I_8\oplus R_8\, ,\notag\\
\MZ_9&=\langle \zeta_9\rangle\oplus \langle \zeta_3^3\,,\ \zeta_2\zeta_7\,,\ \zeta_4\zeta_5\,,\ \zeta_6\zeta_3\rangle=\Q\zeta_9\oplus R_9\, ,\nn\notag\\
\MZ_{10}&=\langle \zeta_{10}\rangle\oplus \langle Z_{37}\rangle\oplus
\langle \zeta_3\zeta_7\,,\ \zeta_5^2\,,\ \zeta_2\zeta_3\zeta_5\,,\ \zeta_2Z_{35} \,,\ \zeta_4\zeta_3^2
\rangle=\Q\zeta_{10}\oplus I_{10}\oplus R_{10}\, ,\notag\\
\MZ_{11}&=\langle \zeta_{11}\rangle\oplus \langle Z_{335}\rangle\oplus
\langle \zeta_3Z_{35}\,,\ \zeta_3^2\zeta_5\,,\ \zeta_2\zeta_9\,,\ \zeta_2\zeta_3^3\,,\ \zeta_4\zeta_7\,,\
\zeta_6\zeta_5\,,\ \zeta_8\zeta_3\nn\rangle=\Q\zeta_{11}\oplus I_{11}\oplus R_{11}\,,\notag
\end{align}
where the irreducibles $Z_{35}$, $Z_{37}$ and $Z_{335}$ are the Z-map images of the generators $\{g_3,g_5\}$, $\{g_3,g_7\}$
and $\{g_3,\{g_3,g_5\}\}$ of $(\mz_w^\vee)^{\ge 2}$ for $w=8,10$ and $11$, respectively (see (\ref{eq:Ihara}) for the definition of the Ihara bracket): they are explicitly given in terms of a common (arbitrary) choice of MZVs $\zeta_{3,5}$, $\zeta_{3,7}$ and $\zeta_{3,3,5}$ by
\begin{align}
\label{Z35Z37Z335}
Z_{35} &\coloneqq Z(\{g_3,g_5\})=
-\tfrac{1105181}{80}\zeta_8
+\tfrac{24453}{5}\zeta_{3,5}
+\tfrac{28743}{2}\zeta_3\zeta_5-1683\,\zeta_2\zeta_3^2\,,\notag\\
Z_{37} &\coloneqq Z(\{g_3,g_7\})=
\tfrac{6614309}{112} \zeta_{3,7} + \tfrac{7796217}{16}\zeta_3\zeta_7 + \tfrac{26525967}{112} \zeta_5^2-\tfrac{2159}{627}\zeta_2 Z_{35} \nn\\
&\hspace{40mm}
-\tfrac{3203187}{76}\zeta_2\zeta_3\zeta_5 
-  \tfrac{60072829}{608} \zeta_4 \zeta_3^2 
-\tfrac{408872741707}{680960} \zeta_{10}
\,,\\
Z_{335}&\coloneqq Z(\{g_3,\{g_3,g_5\}\})
= -\tfrac{3683808}{5}\zeta_{3,3,5} +\tfrac{1119631493}{20}\zeta_{11}-\tfrac{28597725}{38}\zeta_3^2\zeta_5\nn\\
&\hspace{40mm}
+\tfrac{296304}{2717} \zeta_3 Z_{35}-\tfrac{198893689}{6}\zeta_2\zeta_9 +\tfrac{25828428}{247}\zeta_2\zeta_3^3
- \tfrac{90515817}{40} \zeta_4\zeta_7
\nn\\
&\hspace{40mm}
+ \tfrac{6826931}{4} \zeta_6\zeta_5
+ \tfrac{1953356831 }{23712}
\zeta_8\zeta_3\, . \notag
\end{align}
We observe here that the products listed above spanning the spaces of reducibles $R_w$ actually form bases for these spaces. 
This is a general result valid for all $w$, which will be proven in the following section \ref{sec:Zmap.4b}, in which we actually determine an explicit basis for $\MZ$ adapted to the canonical decomposition of Definition \ref{gwcan}. 

Up to $w=7$, the canonical polynomials $g_w$ are given by
\begin{align}
g_2&=[xy] \,,\notag\\
g_3&=[x[xy]]+[[xy]y] \,, \notag \\
g_4&=[x[x[xy]]]+\tfrac{1}{4}[x[[xy]y]]+[[[xy]y]y]+\tfrac{5}{4}(xyxy-xyyx-yxxy+yxyx) \,,\notag\\
%
g_5&=[x[x[x[xy]]]]{+}2[x[x[[xy]y]]]{-}\tfrac{3}{2}[[x[xy]]\, [xy]]{+}2[x[[[xy]y]y]]{+}\tfrac{1}{2}[[xy]\, [[xy]y]]{+}[[[[xy]y]y]y]  
\,,\notag \\
g_6&=[x[x[x[x[xy]]]]]{+}\tfrac{3}{4}[x[x[x[[xy]y]]]]{+}\tfrac{1}{6}[x[[x[xy]]\,[xy]]]{+}\tfrac{23}{16}[x[x[[[xy]y]y]]]{+}\tfrac{1}{12}[x[[xy]\,[[xy]y]]]\notag\\*
&\ \ \ \ \ -\tfrac{89}{48}[x[[[xy]y]\,[xy]]]+\tfrac{3}{4}[x[[[[xy]y]y]y]]
+\tfrac{5}{3}[[xy]\,[[[xy]y]y]]+[[[[[xy]y]y]y]y] \notag\\*
&\quad+\tfrac{7}{4}(xyxxxy-xyyxxx+xyyyxy-xyyyyx-yxxxxy+yxyxxx-yyyxxy+yyyxyx)\notag\\*
&\quad+\tfrac{21}{4}
(xyxyxx-xyxxyx+yxxxyx-yxxyxx-yxyyxy+yxyyyx+yyxyxy-yyxyyx)\notag\\*
&\quad+\tfrac{7}{16}(xyxxyy-xyyyxx-yxxxyy+yxyyxx)+\tfrac{7}{48}(yxxyxy-xyxyxy)\notag\\*
&\quad+\tfrac{35}{48}(yxxyyx+yxyxxy-xyxyyx-xyyxxy)+\tfrac{77}{48}(xyyxyx-yxyxyx) \,,\notag\\
g_7&=[x[x[x[x[x[xy]]]]]]+3[x[x[x[x[[xy]y]]]]]
-5 [x [x[[x[x, y]]\,[x, y]]]] 
+2[[x[x[xy]]\,[x[xy]]]
\notag\\
&\quad
+5 [x[x[x[[[x y] y] y]]]] 
+\tfrac{19}{16}[x[x[[xy]\,[[xy]y]]]]
-\tfrac{173}{16} [x[[x[[x y] y]]\, [x y]]] 
-2[[x[xy]]\,[x[[xy]y]]]\notag\\
&\quad+\tfrac{17}{16}[[[x[xy]]\,[xy]]\,[xy]]+5[x[x[[[[xy]y]y]y]]]
+\tfrac{99}{16}[x[[xy]\,[[[xy]y]y]]]
-\tfrac{61}{16}[[x[[xy]y]]\,[[xy]y]]\notag\\
&\quad-\tfrac{109}{16}[[x[[[xy]y]y]]\,[xy]]
+\tfrac{65}{16}[[xy]\,[[xy]\,[[xy]y]]]
+3[x[[[[[xy]y]y]y]y]]+4[[xy]\,[[[[xy]y]y]y]]\notag\\
&\quad+3[[[xy]y]\,[[[xy]y]y]]
+[[[[[[xy]y]y]y]y]y]\, .
\label{explgws}
\end{align}
In these expressions, we have omitted the separating comma between the two arguments of the Lie bracket in ${\rm Lie}[x,y]$ to condense the formulas. 
The odd degree (Lie) polynomials satisfy the symmetry property $g_{2k+1}(x,y) = g_{2k+1}(y,x)$ that follows from the arguments in footnote \ref{gwsymm}. This is easy to see for $g_3$, but requires also the use of the Jacobi identity to make it manifest for $g_5$ and $g_7$. Our expressions are  chosen to be adapted to the Lyndon basis of ${\rm Lie}[x,y]$ that we introduce in the next section.

For $w\ge 8$ the polynomials $g_w$ become too unwieldy to write down, although they can be calculated on a computer easily (either by the methods presented here, or from the Drinfeld associator as in (\ref{Phiexp}) below). The explicit form of all $g_w$ at $w \leq 12$ can be found in machine-readable form in an ancillary file of the arXiv submission of this work.
However, since the Z-map is an isomorphism, no information is lost in giving their Z-map images, which determine them completely and are much shorter to write down: 
\begin{align}
Z(g_2)&=2 \zeta_2\,,\notag\\
Z(g_3)&=12\zeta_3\,,\notag\\
Z(g_4)&=\tfrac{375}{8}\zeta_4\,,\notag\\
Z(g_5)&=385\zeta_5-105\zeta_2\zeta_3\,,\notag\\
Z(g_6)&=\tfrac{251797}{288}\zeta_6-\tfrac{679}{4}\zeta_3^2\,,\notag\\
Z(g_7)&=\tfrac{49203}{4}\zeta_7
    -\tfrac{14091}{4}\zeta_2\zeta_5
    -\tfrac{11865}{4}\zeta_4\zeta_3\,,\\
Z(g_8)&=\tfrac{769152355481}{40974336}\zeta_8-\tfrac{18246083}{1824}\zeta_3\zeta_5+\tfrac{74974943}{71136}\zeta_2\zeta_3^2\,,\notag\\
Z(g_9)&=\tfrac{373659143}{864}\zeta_9-\tfrac{264398849}{3456}\zeta_6\zeta_3-\tfrac{3702413}{36}\zeta_4\zeta_5-\tfrac{70513729}{576}\zeta_2\zeta_7+\tfrac{133133}{16}\zeta_3^3\,,\notag\\
Z(g_{10})&=\tfrac{22565838727030761032761}{48180785666457600} \zeta_{10} +  
+ \tfrac{23603271373}{184515876480} \zeta_2 Z_{35} -\tfrac{70504768535925229}{227096463360} \zeta_3\zeta_7 -\tfrac{66965094752611}{436723968} \zeta_5^2  \notag\\
&\quad +\tfrac{21865877274704331}{321719989760} \zeta_2\zeta_3\zeta_5 + \tfrac{3916397111572098571}{100376636805120} \zeta_4 \zeta_3^2\,,\notag\\
Z(g_{11})&=\tfrac{1316030287522093}{78587904}\zeta_{11}
+\tfrac{67235}{1227936}\zeta_3Z_{35}
+\tfrac{4632642114815}{4911744}\zeta_3^2\zeta_5
-\tfrac{824237896586533}{176822784}\zeta_2\zeta_9
\notag\\
&\quad 
-\tfrac{470709526441}{4911744}\zeta_2\zeta_3^3
-\tfrac{3026492983085}{818624}\zeta_4\zeta_7
-\tfrac{218501860145855}{78587904}\zeta_6\zeta_5
-\tfrac{3190686062952839}{1414582272}\zeta_8\zeta_3\,.\nn
\end{align}
Note that, in agreement with the third characterizing property (\ref{char3rd}) of $g_w$, the non-single irreducibles $Z_{35}\in I_8$, $Z_{37}\in I_{10}$ and $Z_{335}\in I_{11}$ are absent in $Z(g_{8})$, $Z(g_{10})$ and $Z(g_{11})$, respectively. The contributions $\zeta_2Z_{35}$ and $\zeta_3Z_{35}$ to $Z(g_{10})$ and $Z(g_{11})$ lie in $R_{10}$ and $R_{11}$, respectively, and are therefore compatible with (\ref{char3rd}).

\subsection{\texorpdfstring{The semi-canonical basis for $\MZ_w$}{The semi-canonical basis for MZ(w)}}
\label{sec:Zmap.4b}

In this section we determine an explicit basis for $\MZ$ which is adapted to the canonical decomposition. The basis of the irreducible parts $I_w$ is given by the Z-map images of the Lyndon brackets of the canonical free generators $g_w$ of $\mz_w^\vee$. The basis of the reducible part $R_w$ in turn consists of all weight $w$ products of elements of the set 
given by $\zeta_2$, $\zeta_v$ for
all odd $v<w$, and the chosen basis elements for $I_v$
for $v<w$, which form a linearly independent
set as proven in Corollary \ref{cor:342} at the end of this subsection. Because the Lyndon basis for
a free Lie algebra,
although very natural and practical, 
cannot justifiably be called canonical, we refer to our basis
as the {\em semi-canonical} basis for the canonical decomposition of $\MZ_w$.

Let us recall the definition and the basic result we need concerning Lyndon bases.

\begin{dfn}
    Let $B=\{b_1,b_2,\ldots\}$ be an ordered set of letters. A {\em Lyndon word} in the alphabet $B$ is a word $W_1=b_{i_1}b_{i_2}\cdots b_{i_r}$ that has the property that 
    every right subword $W_j=b_{i_j} b_{i_{j+1}}\cdots b_{i_r}$ with $j>1$ is lexicographically larger than $W_1$.
\end{dfn}

The following classic theorem was discovered simultaneously in 1958 by Chen--Fox--Lyndon and Shirshov (cf.~\cite{Lyndon}, \cite{Shirshov}, or \cite{Reutenauer} for a comprehensive introduction). 

\begin{thm}
\label{lem:lynd}
Let $B=\{b_1,b_2,\ldots\}$ be an ordered set of letters and let ${\rm Lie}[B]$ be the free Lie algebra generated by $B$ (over a field which we take to be $\Q$). Then a basis of ${\rm Lie}[B]$ is given by the individual letters $b_i$ and the set of {\em Lyndon brackets}
\beq
\label{eq:lynd}
[b_{i_1}b_{i_2}\ldots b_{i_r}]\,,
\eeq
where the word $b_{i_1}b_{i_2}\ldots b_{i_r}$ is a Lyndon word, and the rule for making it into a Lie bracket is to place the comma at the leftmost position such that it divides the Lyndon word into two shorter Lyndon words: 
\begin{align}
    [b_{i_1}b_{i_2}\ldots b_{i_r}] = \left[ [b_{i_1}\ldots b_{i_{k-1}}], [b_{i_k}\ldots b_{i_r}]\right]
    \label{reclyn}
\end{align}
and to proceed recursively until it is a multiple bracket of single letters for which we set $[b_i] \coloneqq b_i$.
\end{thm}

\vspace{.3cm}
\noindent {\bf Examples.} The first few Lyndon brackets in the free Lie algebra ${\rm Lie}[x,y]$ are given by
\begin{align}
    [xy] = [x,y]\,,\quad
    [xxy] = [x,[x,y]]\,,\quad
    [xyy] = [[x,y],y]\,,\quad
    [xxyy] = [x,[[x,y],y]]]\,.
\end{align}

The first few Lyndon brackets in the free Lie algebra $\mz^\vee$  on one generator $g_w$ for each odd $w\ge 3$ (see Definition~\ref{dfn:mzvee}) equipped with its Ihara Lie bracket $\{\cdot,\cdot\}$ from~\eqref{eq:Ihara} are given by
\begin{align}
    \{g_3g_5\} = \{g_3,g_5\}\,,\ \ 
    \{g_3g_7\} = \{g_3,g_7\}\,,\ \ 
    \{g_3g_3g_5\} = \{g_3,\{g_3,g_5\}\}\,.
\end{align}

\begin{dfn}\label{defbasIw}
Since $\mz^\vee$ is freely generated by the canonical Lie polynomials $g_3,g_5,\ldots$, the Lyndon brackets in these generators form a basis. Every such Lyndon bracket corresponds as above to a Lyndon word $g_{v_1}\cdots g_{v_r}$ with $r>1$.
We write the corresponding Lyndon bracket as 
\beq
L_{v_1 v_2 \cdots v_r} \coloneqq \{ g_{v_1}g_{v_2} \cdots g_{v_r} \}\in \mz^\vee 
\label{Lsubdef}
\eeq
with odd $v_1,\ldots,v_r \geq 3$. For example, $L_{335}$ denotes the Lyndon bracket $\{g_3,\{g_3,g_5\}\}$. We denote the Z-map images of the Lyndon bracket by 
\beq
Z_{v_1\cdots v_r} \coloneqq Z(L_{v_1\cdots v_r})\, ,
\label{Zsubdef}
\eeq
consistently with (\ref{Z35Z37Z335}).
These elements with $v_1+\cdots+v_r=w$ form the {\em semi-canonical basis}
for the canonical subspace of weight 
$w$ non-single irreducibles $I_w\subset \MZ_w$.
\end{dfn}

\vspace{.3cm}
Our next task is to establish a basis for the spaces $R_w$. 

\begin{prop}\label{propbasRw} Let $C_w\subset \MZ$ be the set consisting of $\zeta_2$, the $\zeta_v$ for odd $3\leq v<w$, and the Z-map images $Z_{v_1\cdots v_r}$ of Lyndon brackets $L_{v_1\cdots v_r}\in \mz^\vee$ with $r>1$,  {$v_1{+}\cdots{+}v_r<w$}.
Then, the set of weight $w$ products of elements of $C_w$ forms a linearly independent set. If $w$ is odd (resp.~even) all of these products (resp.~all of these products except for $(\zeta_2)^{w/2}$) form a basis for $R_w$. 
\end{prop}

This proposition follows from the general result on Hopf algebras given in the following theorem (see Corollary \ref{cor:342}). It seems like this result should be well-known, however it appears to have only been written down in an unpublished note by
Perrin and Viennot \cite{PerrinViennot}.

\vspace{.3cm}
\begin{thm}
\label{thm:PV}
Let $X$ denote an alphabet of 
weighted letters having the 
property that the number of letters in each weight is finite. Let $A^\vee$ 
denote the graded associative $\Q$-algebra on~$X$, considered as a 
Hopf algebra equipped with a multiplication denoted $\circp$ and the standard 
coproduct $\Delta_s$ for which the letters of $X$ are primitive. 
Let $A$ denote the graded dual space of~$A^\vee$, let $L^\vee\subset A^\vee$ 
denote the subspace of primitive elements for $\Delta_s$, and let 
$B=\{b_1,b_2,\ldots\}$ be a vector space basis for $L^\vee$. Then,
\begin{enumerate}
\item[(i)]
$L^\vee$ forms a Lie algebra whose bracket is given by $[g,h]=g\circp h-h\circp g$.

\item[(ii)] Both $A$ and $A^\vee$ have bases given by the
monomials $w$ in the letters of $X$, which we denote by $w\in A$ and 
$w^\vee\in A^\vee$. The map $w^\vee\mapsto w$ provides an isomorphism
of graded vector spaces from $A^\vee$ to $A$. As a $\Q$-algebra, however,
$A$ is commutative, equipped with the shuffle multiplication.

\item[(iii)]
Let $\xi_i$ denote the images of the elements $b_i\in A^\vee$ under the
isomorphism in (ii). The $\xi_i$ then form a multiplicative set of generators 
for $A$ under the shuffle multiplication.

\item[(iv)]
The ordered monomials $\xi_{i_1}\shuffle \xi_{i_2}\shuffle \cdots\shuffle 
\xi_{i_m}$ with $i_1\leq i_2\leq \ldots \leq i_m$ form a linear basis for $A$; those with $m>1$ form a basis for 
the subspace $S\subset A$ annihilating $L^\vee$.
\end{enumerate}
\end{thm}

\vspace{.3cm}
\noindent {\bf Proof.} (i) follows directly from the Milnor--Moore theorem~\cite{MilnorMoore}. The 
vector space part of (ii) follows from the fact that each graded part is 
finite-dimensional, so has a dual that is isomorphic to it and equipped with a 
dual basis; the notation $w^\vee$ for the basis of $A^\vee$ simply defines a 
dual basis to the basis of monomials $w\in A$.  The fact that the 
multiplication on $A$ is the shuffle is standard, corresponding to the fact that
an element of~$A^\vee$ is a Lie element if and only if it satisfies the 
shuffle relations (see \eqref{eq:Zmap.2c}), completing the proof of (ii). This is the same as saying
that the subspace $S\subset A$ spanned by all 
shuffles of monomials is the subspace that annihilates
the Lie algebra $L^\vee$. For this reason, the quotient space 
$L=A/S$ is the Lie coalgebra dual to $L^\vee$, and the linear isomorphism in (ii) induces a linear isomorphism between $L$ and $L^\vee$.  Hence, the $\xi_i\in A$ form a basis for
a subspace $\tilde{L}\subset A$ isomorphic to $L$, restricted to which the
quotient map $A\rightarrow A/S = L$ is an isomorphism. Thus we have
$A = S\oplus \tilde{L}$, completing the proof of (iii).

The final point (iv) follows from the 
Poincar\'e--Birkhoff--Witt theorem~\cite{Bourbaki:Lie}, which states that the universal enveloping 
algebra of a Lie algebra is generated by the ordered monomials in elements 
of a basis, and the only relations come from relations in the Lie algebra. 
We consider $L=A/S$ as a Lie algebra with the trivial bracket, so that the
only multiplicative relations between the generators $\xi_i$ of $L$ are
given by the fact that they commute. By the Poincar\'e--Birkhoff--Witt theorem,
the ordered monomials $\xi_{i_1}\shuffle \cdots \shuffle\xi_{i_m}$
with $m\ge 1$ then form a basis for the universal enveloping algebra $A$ of $L$,
and the monomials with $m>1$ form a basis for the kernel of the map 
$A\rightarrow L$, so in fact they form a basis for~$S$, proving (iv).\qed

\begin{rmk}
Essentially what this proof expresses is that the usual basis of the free associative algebra $A^\vee$ on the alphabet $X$, given by the monomials in the letters of $X$, can be replaced by a different basis consisting of the basis $b_i$ of Lie elements on the one hand, spanning the Lie algebra $L^\vee\subset A^\vee$, completed by the space $S^\vee$ spanned by shuffles of monomials on the other, so that $A^\vee=L^\vee\oplus S^\vee$. In the dual space $A$, this corresponds to an equivalent decomposition $A=L\oplus S$ where $L$ is the subspace whose basis is the $\xi_i$ and $S$ is 
the subspace spanned by all  non-trivial shuffles of the $\xi_i$, which are in fact linearly independent by (iv).
\end{rmk}

\begin{cor}\label{cor:1} 
Let $A^\vee=\overline{\MZ}^\vee$, which by Brown's theorem \cite{Brown:2011mot} is freely generated by $g_3,g_5,\ldots$ under the $\circp$ multiplication. Then the elements $Z(g_w)$ for odd $w\ge 3$ together with the shuffles
\beq\label{shufwords}
Z(g_{w_1})\shuffle Z(g_{w_2})\shuffle\cdots\shuffle Z(g_{w_r})\ \ \ {\rm with}\ \ \ 
w_1\le w_2\le \cdots \le w_r
\eeq
(called {\em ordered} shuffle products) form a basis for $\overline{\MZ}=A$; in particular the ordered shuffles are linearly independent.
\end{cor}
 
We now pass from $\overline{\MZ}$ to $\MZ$ by using the isomorphism (\ref{def:mmzv}).

\begin{cor}
\label{cor:342}
Let $g_3,g_5,\ldots$ denote the canonical generators of $\mz^\vee$. Then a basis for $\MZ$ is given by
the following elements:
\begin{enumerate}
\item[(i)]
the single motivic zeta values $\zeta_w$ for $w\ge 2$;

\item[(ii)]
the Z-map images $Z_{w_1\cdots w_r}$ of the basis of $\mz^\vee$ given by the Lyndon brackets $L_{w_1\cdots w_r}$ with $r>1$ of the canonical generators $g_3,g_5,\ldots$; the weight $w=w_1{+}\ldots{+}w_r$ elements of this type give a basis of $I_w$;

\item[(iii)]
the ordered shuffle products of all the basis elements in (i) and (ii) above, excluding the products of even single zetas (since these products are equal to rational multiples of powers of $\zeta_2$); the weight $w$ elements of this type form a basis for  $R_w$.
\end{enumerate}
\end{cor}

\vspace{.2cm}
\noindent {\bf Proof.} 
A basis of $\Q[\zeta_2]$ is given by the powers of $\zeta_2$, so by \eqref{evens} the single zeta values $\zeta_w$ for all even $w\ge 2$ also give a basis. A basis for $\overline{\MZ}$ is given in Corollary~\ref{cor:1}. Thanks to \eqref{def:mmzv}, a basis for 
the tensor product is given by the products of the basis elements of each 
of the two vector spaces, which is precisely as described by (i), (ii) and (iii) of the statement.\qed

\subsection{Canonical polynomials from the Drinfeld associator}
\label{sec:Drin}

In this section we introduce the
 Drinfeld associator \cite{Drinfeld:1989st, Drinfeld2} which offers an alternative method of computing the canonical polynomials $g_w$. The Drinfeld associator is given by the power series \cite{LeMura}
\beq
\Phi_{{\rm KZ}}(x,y)
\coloneqq
{\bf 1}+\sum_w (-1)^{d(w)}\zeta(w)w\in  {\cal Z}  \otimes_\Q \Q\langle\langle x,y
\rangle\rangle\, ,
\label{defdrin}
\eeq
where $\Q\langle\langle x,y\rangle\rangle$ denotes the degree completion of the
polynomial ring $\Q\langle x,y\rangle$, the sum runs over non-trivial monomials $w$ in $x$ and $y$, and
for each such $w$, $d(w)$ denotes the depth of the monomial, i.e.\ the
number of $y$'s contained in it.\footnote{The subscript ``KZ'' in $\Phi_{{\rm KZ}}(x,y)$ stems from the fact that the Drinfeld associator can be constructed by solving the Knizhnik--Zamolodchikov equation, see appendix~\ref{app:deg}.}
Removing the signs in front of each term produces a power series that we call the {\it modified Drinfeld associator}, given by\footnote{This definition gives an a posteriori explanation of the stuffle-regularized MZVs defined in~\eqref{appMZV.06quater}: for words $w$ ending with $y$ the value $\zeta_\stuffle(w)$ is nothing other than the coefficient of $w$ in the product $C  \Phi$ of formal power series, where $C$ is the power series defined in~\eqref{appMZV.06ter}.} 
\beq
\label{eq:Phi}
\Phi(x,y)  \coloneqq \Phi_{{\rm KZ}}(x,-y)={\bf 1}+\sum_w \zeta(w)w\in \Z \, \hat\otimes \, \Z^\vee\, ,
\eeq
where $\hat\otimes$ denotes the completed tensor product (allowing infinite sums).
We also have formal and motivic versions \beq\label{canonicalelements}
\Phi^{\mathfrak{f}}\in \FZ \, \hat\otimes \, \FZ^\vee\ \ \ {\rm and}\ \ \ \Phi^{\mathfrak{m}}\in \MZ \, \hat\otimes\, \MZ^\vee\, ,
\eeq
obtained by replacing
$\zeta(w)$ by $\zeta^{\mathfrak{f}}(w)$ and $\zeta^{\mathfrak{m}}(w)$, respectively. The coefficients of all three power
series $\Phi$, $\Phi^{\mathfrak{f}}$ and $\Phi^{\mathfrak{m}}$ satisfy the regularized double shuffle relations.

\begin{dfn} \label{basind} 
Let $V=\bigoplus_w V_w$ be a graded vector space for which each graded part is finite-dimensional, and let $V^\vee$ denote the graded dual (the direct sum of the duals of the graded parts of $V$). Choose any basis $e_1,e_2,\ldots$ for $V$ respecting the grading decomposition, and let $e_1^\vee,e_2^\vee,\ldots$ denote the dual basis of $V^\vee$, with 
$\langle e_i^\vee,e_j\rangle=\delta_{ij}$.
Let 
\beq
\Psi=\sum_{i=1}^\infty e_i\otimes e_i^\vee\in V\, \hat\otimes\, V^\vee\, .
\label{psiinf}
\eeq
 We call $\Psi$ the {\em canonical element} of $V\, \hat\otimes \, V^\vee$.
\end{dfn}

Note that the element $\Psi$  is independent of the choice
of basis of $V$ due to the use of dual bases.

\begin{prop} \label{canonel}
Let $V$ be as in Definition~\ref{basind} and let $\phi:V\rightarrow W$ denote any surjective linear morphism and
$\phi^\vee:W^\vee\rightarrow V^\vee$ denote
the dual morphism.
Let $\Psi$ be the canonical element
of $V\,\hat\otimes\, V^\vee$. Then $\bigl(\phi\otimes(\phi^\vee)^{-1}\bigr)(\Psi)$ (in the sense specified in the proof) is the canonical element
of $W\,\hat\otimes \,W^\vee$.
\end{prop}

\noindent {\bf Proof.} We may assume that $V$ is finite-dimensional by working with a fixed graded piece. 
Since $\phi$ is surjective, we have that $V/{\rm Ker}\,\phi\cong W$. Choose a basis of $V$ adapted to this quotient, i.e.\ linearly independent elements $\tilde{w}_1,\ldots, \tilde{w}_m\in V$  that get mapped to a basis $\{w_i=\phi(\tilde{w}_i)\}$ of $W$ under $\phi$ and a basis $k_1,\ldots,k_n$ of ${\rm Ker}\,\phi$. Write the canonical element $\Psi$ in this basis:
\beq
\Psi=\sum_{i=1}^m \tilde w_i\otimes \tilde w_i^\vee+\sum_{j=1}^n k_j\otimes k_j^\vee\,.
\eeq
We now apply the map $\phi\otimes (\phi^\vee)^{-1}$ to $\Psi$, with the understanding that this map is interpreted as the composition
\beq\label{compos}
\bigl({\rm id}\otimes (\phi^\vee)^{-1}\bigr)\circ (\phi\otimes {\rm id}\bigr)\, ,
\eeq
which avoids appearing to apply $(\phi^\vee)^{-1}$ to elements not in $\phi^\vee(W^\vee)$.
We thus obtain
\beq
\label{eq:phiphivee}
\bigl(\phi\otimes (\phi^\vee)^{-1}\bigr)(\Psi)=\sum_{i=1}^m w_i\otimes (\phi^\vee)^{-1}(\tilde w_i^\vee)=\sum_{i=1}^m w_i\otimes w_i^\vee\,,
\eeq
which is the canonical element of $W\otimes W^\vee$.
\qed

\vspace{.3cm}
Recall from diagram \eqref{diag0} that $\Q[Z(w)]$ is the graded dual of the power series ring $\Q\langle \langle x,y\rangle\rangle$. Then, the element
\beq
\Phi^Z =
{\bf 1}+\sum_w Z(w) \otimes w\in \Q[Z(w)]\,  \hat\otimes_\Q \, \Q\langle\langle x,y\rangle\rangle
\eeq
is the canonical element of the tensor product $\Q[Z(w)] \, \hat\otimes_\Q \,\Q\langle\langle x,y\rangle\rangle$.
Since $\Z,\FZ$ and $\MZ,$ are all quotients of $\Q[Z(w)]$ (see diagram~\eqref{diag0}), Proposition \ref{canonel} then implies that $\Phi$, $\Phi^{\mathfrak{f}}$ and $\Phi^{\mathfrak{m}}$ are the canonical elements for the respective rings $\Z\,\hat\otimes\, \Z^\vee$, $\FZ\,\hat\otimes\, \FZ^\vee$ and $\MZ\,\hat\otimes\, \MZ^\vee$. In particular, the choice of basis in which to
express $\Phi^{\mathfrak{m}}$ is of little significance in general. However, writing
$\Phi^{\mathfrak{m}}$ in the semi-canonical basis does have one convenient advantage: it provides another method to compute the canonical polynomials $g_w$.

In our semi-canonical basis of $\MZ$ (see (\ref{Z35Z37Z335}) for $Z_{35}, Z_{37}$ and $Z_{335}$), the expansion of the modified Drinfeld associator $\Phi$ to weight 11 reads as follows, see \cite{Drummond:2013vz} for the analogous expansion of the Drinfeld associator and its significance for the motivic coaction:\footnote{The product $\circ$ among $h\in \ds$ and $g\in \U\ds$ in \cite{Drummond:2013vz} is related to the Poincar\'e--Birkhoff--Witt multiplication~$ \circp$ in (\ref{starmult}) via $\overleftarrow{g \circ h} = \overleftarrow{h} \circp \overleftarrow{g}$, where $\overleftarrow{w}$ is obtained by reversing the letters $x,y$ of $w \in \mathbb Q\langle x,y \rangle$. This is a consequence of $D_{\overleftarrow{h}}(\overleftarrow{g}) = -\overleftarrow{D_h (g)}$ which can be proven by induction. The Drinfeld associator in the conventions of \cite{Drummond:2013vz} is obtained from the series $\Phi(x,y)$ in the present work by reversing the words $w \mapsto \overleftarrow{w}$ in~(\ref{eq:Phi}).}
\begin{align}
\Phi&={\bf 1}+\zeta_2 g_2+\zeta_3 g_3+\zeta_4 g_4+\zeta_5 g_5+\zeta_2\zeta_3 g_3\circp g_2
+\zeta_6 g_6+\tfrac{1}{2}\zeta_3^2g_3\circp g_3+\zeta_7g_7 \label{Phiexp}\\
&\quad +\zeta_3\zeta_4g_3\circp g_4+\zeta_2\zeta_5g_5\circp g_2
+\zeta_8g_8
+\zeta_2\zeta_3^2\bigl(\tfrac{1}{2}g_3\circp g_3\circp g_2+\tfrac{17}{247}\{g_3,g_5\}\bigr)\notag\\
&\quad +\tfrac{1}{24453}Z_{35}\{g_3,g_5\} 
+\zeta_3\zeta_5\Bigl(\tfrac{47}{114}g_3\circp g_5+\tfrac{67}{114}g_5\circp g_3\Bigr)\nn\\
&\quad +
\zeta_9g_9+\tfrac{1}{6}\zeta_3^3g_3\circp g_3\circp g_3 +\zeta_2\zeta_7g_7\circp g_2+\zeta_4\zeta_5g_5\circp g_4+\zeta_6\zeta_3g_3\circp g_6\notag\\
&\ \ \ +\zeta_{10}g_{10}
+\tfrac{8}{6614309}Z_{37}\{g_3,g_7\}+\zeta_3\zeta_7\bigl(\tfrac{24581}{59858}g_3\circp g_7+\tfrac{35277}{59858}g_7\circp g_3\bigr)\notag\\
&\quad +\zeta_5^2\bigl(\tfrac{1}{2}g_5\circp g_5-\tfrac{2160}{29929}\{g_3,g_7\}\bigr)+\zeta_2Z_{35}\bigl(\tfrac{1016}{243951279}\{g_3,g_7\}+\tfrac{1}{24453}\{g_3,g_5\}\circp g_2\bigr)\notag\\
&\quad+\zeta_2\zeta_3\zeta_5\bigl(\tfrac{47}{114}g_3\circp g_5\circp g_2+\tfrac{67}{114}g_5\circp g_3\circp g_2+\tfrac{492798}{9667067}\{g_3,g_7\}\bigr)\notag\\
&\quad +\zeta_4\zeta_3^2\bigl(\tfrac{85}{494}\{g_3,g_5\}\circp g_2
+\tfrac{60072829}{502687484}
\{g_3,g_7\}+\tfrac{1}{2}g_3\circp g_3\circp g_4\bigr)\notag\\
&\quad +\zeta_{11}g_{11}
+\tfrac{1}{3683808}Z_{335}\{g_3,\{g_3,g_5\}\}\notag\\
&\quad+\zeta_3Z_{35}\bigl(\tfrac{7063}{625556646}g_3\circp g_3\circp g_5+\tfrac{5728}{312778323}g_3\circp g_5\circp g_3-\tfrac{6173}{208518882}g_5\circp g_3\circp g_3\bigr)\notag\\
&\quad +\zeta_3^2\zeta_5\bigl(\tfrac{5439455}{46661568}g_3\circp g_3\circp g_5+\tfrac{4179377}{23330784}g_3\circp g_5\circp g_3+\tfrac{3177525}{15553856}g_5\circp g_3\circp g_3\bigr)\notag\\
&\quad +\zeta_2\zeta_9\bigl(-\tfrac{31943}{22102848}\{g_3,
\{g_3,g_5\}\}+g_9\circp g_2\bigr)
+\zeta_4\zeta_7
\bigl(\tfrac{46765}{3274496}
\{g_3,\{g_3,g_5\}\}+g_7\circp g_4\bigr)\notag\\
&\quad +\zeta_2\zeta_3^3\bigl(
\tfrac{3066359}{75825048}  g_3\circp g_3\circp g_5-\tfrac{456995}{37912524}  g_3\circp g_5\circp g_3-\tfrac{ 2152369}{75825048}  g_5\circp g_3\circp g_3+\tfrac{1}{6}g_3\circp g_3\circp g_3\circp g_2\bigr)\notag\\
&\quad +\zeta_8\zeta_3\bigl(
-\tfrac{1953356831}{87350455296}
\{g_3, \{ g_3, g_5 \}\} 
+g_3\circp g_8\bigr)
+\zeta_6\zeta_5\bigl( \tfrac{540685}{14735232}
\{g_3,\{g_3,g_5\}\}+  g_5\circp g_6\bigr)
+\ldots \notag
\end{align}

\subsubsection*{Computational remarks}

In order to write motivic MZVs in a given basis in weight $w$ we need to know the linear relations between motivic MZVs in that weight. While these are not known in general, we have several possible approaches: (i) in weights up to $w=22$ (and also at weight $w=23$ modulo a 31-bit prime), it is known by dimension arguments that $\MZ_w=\FZ_w$ \cite{Blumlein:2009cf} so we can use the double shuffle relations, (ii) since Brown gave the dimension of $\MZ_w$ in all weights, if we reached any weight where $\MZ_w$ is not equal to $\FZ_w$ (in spite of the conjecture that they are equal) we could write the real MZVs as real numbers, seek for enough linear relations between them with rational coefficients to reach the correct dimension and then prove that these relations are motivic \cite{Blumlein:2009cf}. In practice, the latter method has been used to create the available datamines, making the decomposition particularly easy by computer as it is enough to enter an MZV into the datamine to automatically obtain its decomposition. Note that the $\mathbb Q$-bases of \cite{Blumlein:2009cf} were extended from weight 22 to weight 34 in the HyperlogProcedures of Schnetz~\cite{Schnetz:www}.

\vspace{.2cm}\noindent
(1) In computing the expression \eqref{Phiexp}, we have written multiple $\circp$-products without parentheses with the understanding that we can evaluate them as $g_{w_1}\circp(g_{w_2}\circp\cdots(g_{w_{r-1}}\circp(g_{w_r}\circp g_k))\cdots )$ with $w_i$ odd and $k$ odd or even. In this way, the 
left factor of each $\circp$ multiplication is a Lie polynomial, i.e.\ a $g_w$ with $w$ odd, which allows us to use the simplified expression \eqref{starmult} for
the multiplication $\circp$ in $\MZ^\vee$.

\vspace{.3cm}\noindent (2)
This gives us three ways to recursively compute the $g_w$, of which we saw the first two earlier:
\begin{itemize}
\item[(i)] from the properties in Lemma \ref{lem:332} that uniquely characterize the $g_w$,
\item[(ii)] get the semi-canonical basis for $I_w$ using the Lyndon words and then compute the unique normalized polynomial $g_w\in \MZ_w^\vee$ annihilating the basis elements of $R_w$ and $I_w$, or 
\item[(iii)] decompose $\Phi$ into the semi-canonical basis of $\MZ$; then
\beq
g_w = \Phi |_{\zeta_w}
\label{niceeq}
\eeq
\end{itemize}
The equivalence of the third approach with the others is a direct consequence of Proposition \ref{canonel}, which implies that the polynomial appearing in $\Phi$ with coefficient $\zeta_w$ must be the element of the dual basis of the semi-canonical basis taking the value 1 on $\zeta_w$ and annihilating $I_w$ and $R_w$. 

\vspace{.3cm}\noindent (3) 
As an advantage of the first method (i) over methods (ii) and (iii), the conditions of Lemma \ref{lem:332} make it clear that the canonical $g_w$ do not depend on any basis choice for $\MZ_w$. For those weights $w$ where the expansion of the Drinfeld associator is available (e.g.\ from \cite{Blumlein:2009cf, Schnetz:www}), the third approach (iii) enjoys the computational advantage that ans\"atze and solutions of linear equation systems can be bypassed.


\section{\texorpdfstring{The canonical morphism from motivic MZVs to the $f$-alphabet}{The canonical morphism from motivic MZVs to the f-alphabet}}
\label{sec:falpha}

In~\cite{Brown:2011ik, Brown:2011mot}, Brown proved a remarkable theorem showing that the
motivic MZV Hopf algebra comodule $\MZ$ is isomorphic to a certain Hopf algebra comodule $\F$ with a particularly simple structure that we recall below. However, Brown did not display a canonical isomorphism, but rather showed the existence and described the construction of a family of isomorphisms $\rho_{\vec{c}\,} : \MZ \to \F$ parametrized by free rational parameters $\vec{c}$ associated to a chosen basis of non-single irreducible motivic MZVs. The goal of this section is to use the canonical  polynomials $g_w$ of Definition~\ref{gwcan} to fix a canonical choice of isomorphism
\beq
\rho:\MZ\rightarrow {\cal F}\, .
\eeq
As in section \ref{sec:Zmap.4a}, we will allow ourselves to simplify the notation
by writing $\zeta$ instead of $\zeta^{\mathfrak{m}}$ throughout the present section, which will deal uniquely with motivic MZVs. Furthermore, in order for this section to remain coherent with the literature (see footnote 3 above) we will consider $\MZ$ as a Hopf algebra comodule with the structure conferred on it by the choice of coaction $\Delta_{GB}$ and not $\Delta^{GB}$ (see \eqref{eq:coact} and \eqref{eq:coactb}). This change also modifies the structure of the dual Hopf algebra $\MZ^\vee$, which instead of being equipped with the multiplication $\circp$ satisfying \eqref{starmultduality}, becomes equipped with the multiplication $\bullet$ defined by
\beq\label{bulletmult1}
h\bullet g \coloneqq g\circp h\,,
\eeq
satisfying
\beq\label{bulletmult2}
\langle \Delta_{GB}(\xi),g\otimes h\rangle=\langle \xi,g\bullet h\rangle
\eeq
for all $\xi\in \MZ$, $g,h\in \MZ^\vee$. Moreover, the simple expression
(\ref{starmult}) for $g\circp h$ in case of $g\in \ds$ translates into
\beq
h\bullet g = gh + D_g(h) 
\eeq
with the Ihara derivation $D_g$ defined by (\ref{eq:Ihara}).
The Lie subspace of $\MZ^\vee$ is then equipped with the Lie bracket associated to $\bullet$, defined by
\beq
[\![g,h]\!] \coloneqq g\bullet h-h\bullet g\, .
\label{defbull}
\eeq
(Note that this Lie bracket satisfies $[\![g,h]\!] =-\{g,h\}$ in relation to the Ihara bracket~\eqref{eq:Ihara}.) 

\subsection{\texorpdfstring{Definition of the $f$-alphabet}{Definition of the f-alphabet}}
\label{sec:falpha.1}

We begin by defining the Hopf algebra comodule
${\cal F}$, familiarly called the $f$-alphabet \cite{Brown:2011ik, Brown:2011mot}. To start with,
let $\overline{\cal F}^\vee \coloneqq \Q\langle f^\vee_3,f^\vee_5,\ldots\rangle$ be the free associative
Hopf algebra on one non-commutative indeterminate $f^\vee_w$ in each odd weight $w\ge 3$, with the usual
(concatenation) multiplication and the standard coproduct defined by
\beq
\Delta_s(f_w^\vee)=f^\vee_w\otimes 1+1\otimes f^\vee_w
\eeq
for all odd $w\ge 3$. The subspace of Lie polynomials ${\cal L}^\vee \coloneqq {\rm Lie}[f^\vee_3,f^\vee_5,\ldots]
\subset \overline{\cal F}^\vee$ is the space of primitive elements $f^\vee\in 
\overline{\cal F}^\vee$, i.e.~elements satisfying
\beq
\Delta_s(f^\vee)=f^\vee\otimes 1+1\otimes f^\vee \, .
\eeq

Now let $\overline{\cal F}$ denote the Hopf algebra dual to
$\overline{\cal F}^\vee$.  The underlying vector space of $\overline{\cal F}$ 
is isomorphic to that of $\Q\langle f_3,f_5, \ldots\rangle$, the free
associative algebra spanned by all monomials 
$f_{i_1}\cdots f_{i_r}$ in the free non-commutative indeterminates $f_i$ for odd $i\ge 3$; these monomials
form a dual basis to the basis of monomials $f_{i_1}^\vee\cdots f_{i_r}^\vee$
of $\overline{\cal F}^\vee$ in the sense that $\langle f_{i_1}^\vee\cdots f_{i_r}^\vee ,  f_{j_1}\cdots f_{j_r}\rangle \!= \delta_{i_1,j_1}\cdots \delta_{i_r,j_r}$.
The Hopf algebra structure
of $\overline{\cal F}$ is given by equipping $\overline{\cal F}$ with the
(commutative) shuffle multiplication on the monomials $f_{i_1}\cdots f_{i_r}$
and the {\em deconcatenation coproduct} $\Delta$
defined by
\beq
\Delta(f_{i_1}\cdots f_{i_r})=\sum_{j=0}^r f_{i_1}\cdots f_{i_j}\otimes
f_{i_{j+1}}\cdots f_{i_r}\, .
\label{sec3eq.1}
\eeq

Following Brown, let us now define the comodule ${\cal F}$ to be the tensor product
\beq
{\cal F} \coloneqq \Q[f_2]\otimes_\Q \overline{\cal F} \, ,
\label{sec3eq.2}
\eeq
where $f_2$ is a new commutative indeterminate of weight 2 and the factor $\Q[f_2]$ denotes
the polynomial ring over $\Q$ in the single indeterminate $f_2$. The algebra 
structure of $\overline{\cal F}$ extends to ${\cal F}$ by letting $f_2$ 
commute with $\overline{\cal F}$; the general rule is
\beq\label{shuffleonF}
(f_2^mf_{i_1}\cdots f_{i_r})\shuffle (f_2^n f_{j_1}\cdots f_{j_s})=
f_2^{m+n} \bigl(f_{i_1}\cdots f_{i_r}\shuffle f_{j_1}\cdots f_{j_s}\bigr) 
\eeq
for odd $i_1,\ldots,i_r,j_1,\ldots,j_s\geq 3$.
By a slight abuse of terminology, we continue to call this product on all of 
${\cal F}$ the {\em shuffle product} on ${\cal F}$.

The $\Q$-algebra ${\cal F}$ is made into a
$\overline{\cal F}$-comodule by defining a coaction
\begin{equation}
\Delta:{\cal F}\rightarrow {\cal F}\otimes \overline{\cal F}
\label{sec3eq.3}
\end{equation}
on ${\cal F}$ by \eqref{sec3eq.1} above together with
\beq
\Delta(f_2)=f_2\otimes 1 \, .
\label{sec3eq.4}
\eeq
Thus the general formula for this coaction is given by
\begin{equation}
\Delta(f_2^n f_{i_1} f_{i_2} \ldots f_{i_r}) = 
\sum_{j=0}^r f_2^n f_{i_1} \ldots f_{i_j}
\otimes  f_{i_{j+1}} \ldots f_{i_r}
\label{sec3eq.5}
\end{equation}
with integer $n,r \geq 0$ and odd $i_1,\ldots,i_r \geq 3$.

\vspace{.2cm}
Now let ${\cal F}^\vee$ denote the dual of ${\cal F}$. The underlying vector space
of ${\cal F}^\vee$ is a tensor product of two vector spaces
\beq
\langle f^\vee_2,f^\vee_4,\ldots\rangle\otimes_\Q \overline{\cal F}^\vee \, ,
\eeq
where $\overline{\cal F}^\vee$ is as defined at the beginning of this section, and
the left-hand factor denotes the vector space (not ring) dual of 
$\Q[f_2]$, with basis $f^\vee_{2n}\in {\cal F}^\vee$ satisfying
\beq
\langle f^\vee_{2n},f_2^m\rangle=\delta_{m,n}
\frac{\zeta_2^n}{\zeta_{2n}} \, .
\eeq
By analogy with  Definition~\ref{dfn:z2n} we set
\beq
f_{2m} \coloneqq \frac{ \zeta_{2m}}{\zeta_2^m} f_2^m\in {\cal F} \, ,
\label{sec3eq.6}
\eeq
so that
\beq\label{needalabel}
\langle f^\vee_{2m},f_{2n}\rangle=\delta_{m,n} \, .
\eeq
The fact that ${\cal F}$ is a Hopf algebra comodule and not a Hopf algebra
is reflected in the dual space by the fact that ${\cal F}^\vee$ is not a Hopf algebra
but a Hopf algebra module over the Hopf algebra $\overline{\cal F}^\vee$. 
Thus, the concatenation multiplication does not extend from the subspace
$\overline{\cal F}^\vee$ to all of ${\cal F}^\vee$; instead we only
have an action of $\overline{\cal F}^\vee$ on ${\cal F}^\vee$, which we 
write as
\beq
a(f^\vee_{2n}b)=f^\vee_{2n}ab\in {\cal F}^\vee
\eeq
for $n\ge 1$ and $a,b\in \overline{\cal F}^\vee$. This action can be considered 
as a multiplication of an element of the space $\Q[f^\vee_2,f^\vee_4,\ldots]$ 
with an element of $\overline{\cal F}^\vee$, but the $f^\vee_{2n}$ 
cannot be multiplied together.  Thus every element of ${\cal F}^\vee$ is 
a sum of monomials which can be written uniquely in the form
$f^\vee_{2n}b$ for some $n\ge 0$ (with the convention $f_0^\vee=1$) and 
some $b\in \overline{\cal F}^\vee$.

\subsection{\texorpdfstring{A canonical choice of normalized isomorphism from $\MZ$ to ${\cal F}$}{A canonical choice of normalized isomorphism from MZ to F}}
\label{sec:falpha.2}

\begin{dfn}
\label{dfn:normmor}
A morphism $\phi: \MZ\to\F$ is a {\em normalized morphism} if the following conditions hold \cite{Brown:2011ik, Brown:2011mot}:
\begin{itemize}
\item[(i)] normalization: $\phi\big( \zeta_n\big) = f_n$ for all $n \geq 2$, where $f_{n}$ for even values $n=2m$ was defined in~\eqref{sec3eq.6}. 
\vspace{-.2cm}
\item[(ii)] compatibility with the shuffle multiplication (\ref{shuffleonF}) on ${\cal F}$,
\beq
 \phi\big( \zeta(w_1) \zeta(w_2)\big) = \phi\big( \zeta(w_1)\big) \shuffle \phi\big( \zeta(w_2)\big)\, . \label{falpha.01}
\eeq

\vspace{-.2cm}
\item[(iii)] compatibility with coactions $\Delta$ in (\ref{sec3eq.5}) and $\Delta_{GB}$ in (\ref{brco}), 
given by the following formula for all monomials $w$ in $x$ and $y$: 
\beq
\Delta \phi\big( \zeta(w)\big)
= \phi\big( \Delta_{GB} \zeta(w)\big) \, .
\label{falpha.02}
\eeq
It is understood that $\phi$ acts on each factor of the tensor product, with
an additional projection from ${\cal F}$ to $\overline{\cal F}$ in the
second factor, meaning that each term involving a power of $f_2$ in the 
second factor will be projected to zero. 
\end{itemize}
\end{dfn}

\begin{rmk}
The third property (\ref{falpha.02}) 
translates the Goncharov--Brown coaction $\Delta_{GB}$, which is expressed
by the complicated procedure given in Definition~\ref{dfn:Gonch},
into the considerably simpler deconcatenation coaction \eqref{sec3eq.5} in the $f$-alphabet.
\end{rmk}

The results summarized in the next theorem follow directly from the results of
Brown in~\cite{Brown:2011ik, Brown:2011mot} that we state here in a version adapted to the semi-canonical basis of Definition~\ref{defbasIw}.

\begin{thm}[Brown] 
\label{thm:421}
Let $w\ge 2$, let $\MZ_w=\Q\zeta_w\oplus I_w\oplus R_w$ 
be the canonical decomposition of Definition~\ref{gwcan} and choose the semi-canonical  basis of $I_w$ expressed via Lyndon words $Z_{v_1\ldots v_r}$ with odd $v_1,\ldots,v_r \geq 3$ introduced in Definition~\ref{defbasIw}. 
Let $\vec{c}=\{ c_{v_1\ldots v_r}\}$ denote an infinite family of rational parameters
indexed by the same Lyndon words. Then for any choice of rational values for the parameters $\vec{c}$,
there exists a normalized Hopf algebra comodule isomorphism 
\beq
\rho_{\vec{c}\,}:\MZ\to \F\, .
\eeq
Furthermore, any normalized Hopf algebra comodule isomorphism in the sense of Definition \ref{dfn:normmor} corresponds to a specific choice of rational values of the parameters in $\vec{c}$. 
\end{thm}

\begin{rmk}
    We have used our choice of semi-canonical basis to state Brown's theorem, but the result is in fact independent of the choice of basis and even of the choice of subspace $I_w$ of non-single irreducibles. For any such choice of $I_w$ equipped with any basis, we can use that basis to index a set of rational numbers $\vec{c}$ parametrizing the inequivalent normalized isomorphisms from $\MZ$ to the $f$-alphabet, with the same constructive proof as the one indicated below for our particular choice.
\end{rmk}

Essentially, the proof of this result comes down to actually constructing
the isomorphisms $\MZ\rightarrow {\cal F}$ inductively weight by weight
\cite{Brown:2011ik, Brown:2011mot}. 
We sketch the procedure here and work
it out explicitly for small weights. 

We saw in section \ref{sec:Zmap.4b} that for weights $w\le 7$ we have $I_w=\{0\}$. Thus for these weights the theorem 
says that the normalized isomorphism is uniquely fixed up to $w\leq 7$; it is in fact determined solely by properties (i) and (ii) of Definition~\ref{dfn:normmor}.
For $w=2,3,4$, we must have
\begin{align}
\rho_{\vec{c}\,}:\MZ_w&\rightarrow {\cal F}_w \, , \notag\\
\zeta_w&\mapsto f_w\, ,
\end{align}
since the weight spaces $\MZ_w$ are $1$-dimensional for these values. For weight 5, $\MZ_5$ is $2$-dimensional spanned by
$\zeta_5$ and $\zeta_2\zeta_3$, so by (i) and (ii) we have
\begin{align}
\rho_{\vec{c}\,}: \MZ_5&\rightarrow {\cal F}_5 \, , \notag\\
\zeta_5&\mapsto f_5 \, ,\notag\\
\zeta_2\zeta_3&\mapsto f_2f_3\, .
\end{align}
For weight 6, $\MZ_6$ is $2$-dimensional spanned by $\zeta_6=\frac{35}{8}\zeta_2^3$ and $\zeta_3^2$, so all $\rho_{\vec{c}\,}$ are given by
\begin{align}
\rho_{\vec{c}\,}:\MZ_6&\rightarrow {\cal F}_6 \, , \notag\\
\zeta_6&\mapsto f_6 \, , \notag\\
\zeta_3^2&\mapsto f_3\shuffle f_3=2f_3f_3 \, .
\end{align}
Finally, in weight 7, $\MZ_7$ is $3$-dimensional, spanned by $\zeta_7$,
$\zeta_2\zeta_5$ and $\zeta_4\zeta_3$, so we have
\begin{align}
\rho_{\vec{c}\,}:\MZ_7&\rightarrow {\cal F}_7 \, , \notag\\
\zeta_7&\mapsto f_7 \, , \notag\\
\zeta_2\zeta_5&\mapsto f_2f_5 \, , \notag\\
\zeta_4\zeta_3&\mapsto f_4f_3 \, .
\end{align}

\vspace{.3cm}
Starting from weight $w=8$, the presence of non-trivial spaces of non-single irreducibles $I_w\subset \MZ_w$ requires additional input from the coaction property (\ref{falpha.02}) in (iii). 

\vspace{.3cm}
\noindent {\bf Example.} Let us illustrate this for the case of weight $w=8$, where we use the element $Z_{35}$ defined in (\ref{Z35Z37Z335}) appearing in our semi-canonical basis constructed in section \ref{sec:Zmap.4b}.
The image under $\rho_{\vec{c}\,}$ of this element is not fixed by (i) and (ii) alone, so we make the most general ansatz 
\begin{align}
\label{eq:Z35ans}
\rho_{\vec{c}\,}(Z_{35}) =  a_1 f_3 f_5 + a_2 f_5f_3 + a_3 f_2 f_3 f_3  + c_{35} f_8
\end{align}
with rational parameters $a_i,c_{35}$
and then impose (iii). 
By combining~\eqref{Z35Z37Z335} and~\eqref{eq:DG8} we find that
\begin{align}
\Delta_{GB}(Z_{35}) =  Z_{35}\otimes 1 + 1\otimes Z_{35}  - \tfrac{20163}{2} \zeta_3\otimes \zeta_5 + \tfrac{28743}{2} \zeta_5\otimes \zeta_3 -3366 \zeta_2\zeta_3\otimes \zeta_3\,,
\end{align}
whose $\rho_{\vec{c}\,}$-image is
\begin{align}
\rho_{\vec{c}\,}\big(\Delta_{GB}(Z_{35})\big) = 
\rho_{\vec{c}\,}(Z_{35}) \otimes 1 + 1\otimes \rho_{\vec{c}\,}(Z_{35}) 
- \tfrac{20163}{2} f_3\otimes f_5 + \tfrac{28743}{2} f_5\otimes f_3 -3366 f_2f_3\otimes f_3\,.
\end{align}
To impose (iii) we have to compare this with the deconcatenation coaction~\eqref{sec3eq.5} applied to the ansatz~\eqref{eq:Z35ans}, which is
\begin{align}
\Delta\big(\rho_{\vec{c}\,}(Z_{35})\big) &=  \rho_{\vec{c}\,}(Z_{35}) \otimes 1 + 1\otimes \rho_{\vec{c}\,}(Z_{35}) + a_1 f_3\otimes f_5 + a_2 f_5 \otimes f_3 + a_3 f_2f_3\otimes f_3 \,.
\end{align}
Comparing coefficients fixes the parameters $a_i$ but leaves $c_{35}$ undetermined, so for \eqref{eq:Z35ans} we obtain
\beq\label{c35example}\rho_{\vec{c}\,}(Z_{35})=-\tfrac{20163}{2}f_3f_5+\tfrac{28743}{2}f_5f_3-3366f_2f_3^2+c_{35}f_8\, .\eeq
This is the first appearance of a rational parameter of $\vec{c}$ from Theorem~\ref{thm:421}.
Analogous free parameters appear as the coefficient of $f_w$ in the image under $\rho_{\vec{c}\,}$ of  each basis element of $I_w$. 
In the semi-canonical basis the parameter $c_{v_1\ldots v_r}$ corresponds to the coefficient of $f_{v_1+\ldots+v_r}$ in $\rho_{\vec{c}\,}(Z_{v_1\ldots v_r})$.

\vspace{.3cm}
\begin{dfn}\label{def1rho} 
For $w\ge 8$, let 
$\MZ_w=\Q\zeta_w\oplus I_w\oplus R_w$ denote the canonical decomposition constructed in section \ref{sec:Zmap.4}. 
Let $\rho_{\vec{c}}$ be the family of normalized Hopf algebra comodule isomorphisms established in the semi-canonical basis as in Theorem~\ref{thm:421} such that its rational parameters $\vec{c}=\{ c_{v_1\ldots v_r}\}$ are indexed by Lyndon words with odd $v_1,\ldots,v_r\geq 3$. Then we define the {\em canonical $f$-alphabet isomorphism} 
\beq\label{canonrho}
\rho:\MZ\rightarrow {\cal F} \quad \text{by} \quad \rho \coloneqq \rho_{\vec{0}}\,.
\eeq
\end{dfn}

The definition of the canonical isomorphism implies immediately 
\beq
\rho(Z_{v_1\ldots v_r})|_{f_{w}}=0 
\label{olisfavorite}
\eeq
for all $v_1{+}\ldots{+}v_r=w\geq 8$ (with $r>1$), which is an alternative unique characterization of~$\rho$. 
This leads for instance to
\begin{align}
\label{morerhos}
\rho(Z_{35})&=-\tfrac{20163}{2}f_3f_5+\tfrac{28743}{2}f_5f_3-3366f_2f_3f_3 \, ,\notag\\
\rho(Z_{37}) &=-\tfrac{5432401}{16}f_3f_7+\tfrac{7796217}{16}f_7f_3
+119340 f_5 f_5
 -\tfrac{2698111}{16}f_4f_3 f_3
-\tfrac{29731}{4}f_2f_3f_5\notag\\
&\quad  -\tfrac{366535}{4}f_2f_5f_3
\,,\nn\\
\rho(Z_{335}) &= 1629441 f_5f_3f_3 -1037295 f_3f_5f_3 -20223 f_3f_3f_5 +\tfrac{31943}{6} f_2 f_9-473832 f_2 f_3 f_3 f_3\nn\\
&\quad 
- \tfrac{420885}{8} f_4 f_7
- \tfrac{540685}{4} f_6 f_5
+\tfrac{1953356831}{23712}  f_8 f_3
\,.
\end{align}

\begin{prop}
The isomorphism $\rho$ is uniquely characterized by the property:
\beq
\rho(\xi)|_{f_w}=0
\ \textrm{for all} \ \xi \in I_w \, .
\label{canrhodef}
\eeq
Equivalently, one can characterize $\rho$ as the unique isomorphism $\MZ \rightarrow {\cal F}$ that preserves the
relation (\ref{niceeq}) between the canonical polynomial $g_w$ and the modified Drinfeld associator $\Phi$, i.e.\ 
\beq
 \rho(\Phi) |_{f_w} = g_w\,.
\eeq
\end{prop}

\noindent 
{\bf Proof.}
Since the $Z_{v_1\ldots v_r}$ at $v_1{+}\ldots{+}v_r=w$ with $r>1$ form a basis of $I_w$, we also have from~\eqref{olisfavorite} for all $w\geq 2$ that $\rho(\xi)|_{f_w}=0$ for any $\xi \in I_w$.
Therefore, writing $\Phi$ in the semi-canonical basis,
no irreducible MZV can contribute to the coefficient of $f_w$ in $\rho(\Phi)$ and the property~\eqref{niceeq} is preserved.

Note that even though the semi-canonical basis appears when defining $\rho=\rho_{\vec{0}}$ in (\ref{canonrho}), $\rho$ is characterized by the property~\eqref{canrhodef} which refers only to the {\em canonical} subspace $I_w$ and therefore $\rho$ can be defined canonically in this way. 
\qed

\begin{rmk}
We end this section with a  brief observation about the specific MZVs $\zeta_{3,5}$, $\zeta_{3,7}$ and $\zeta_{3,3,5}$, that 
are widely used in the physics literature as a basis for a non-canonical choice of (1-dimensional) subspace of non-single irreducibles in $I_w\subset \MZ_w$ for $w=8,10,11$. 
Using~\eqref{Z35Z37Z335} and~\eqref{morerhos}, the
canonical parameter choice $c_{35}=c_{37}=c_{335}=0$ translates into the $f$-alphabet images
\begin{align}
\rho( \zeta_{3,5}) &= - 5 f_3 f_5 + \tfrac{100471}{35568} f_8 \, , 
\label{fimag} \\
\rho( \zeta_{3,7}) &= - 14 f_3 f_7 - 6 f_5 f_5 + \tfrac{408872741707}{40214998720} f_{10} \, , \notag \\
\rho( \zeta_{3,3,5}) &= - 5 f_3 f_3 f_5 - 45 f_2 f_9 
- 3f_4 f_7
+ \tfrac{5}{2}f_6 f_5
+ \tfrac{1119631493}{14735232} f_{11}
\notag
\end{align}
for these elements.
The analogous $\rho$-images of all irreducible higher-depth motivic MZVs of weights $\leq 17$ in the basis choice of \cite{Blumlein:2009cf} can be found in the ancillary files of \cite{Dorigoni:2024oft}. 
\end{rmk}

\section{\texorpdfstring{Canonical zeta generators $\sigma_w$ in genus one}{Canonical zeta generators sigma(w) in genus one}}
\label{sec:sigman}

In this section we show how the canonical polynomials $g_w$ associated with zeta generators in genus zero
as defined in section~\ref{sec:Zmap.4} induce canonical zeta generators $\sigma_w$ in genus one. The construction  also includes a canonical split of $\sigma_w$ into an arithmetic and a geometric part. 

\subsection{\texorpdfstring{The Tsunogai derivations $\epsilon_k$}{The Tsunogai derivations epsilon(k)}}
\label{sec:sigman.0}

In this section we write ${\rm Lie}[a,b]$ for the fundamental Lie algebra
associated to a once-punctured torus. This is a free Lia algebra on two generators
and thus isomorphic to ${\rm Lie}[x,y]$, but we prefer to distinguish the
letters used because the topological fundamental group
of a thrice-punctured sphere maps non-trivially to that of a once-punctured torus
when two of the holes are joined together. We also have a natural map
between the pro-unipotent fundamental groups, which gives a natural but
highly non-trivial Lie algebra morphism
\beq
{\rm Lie}[x,y]\rightarrow {\rm Lie}[a,b]
\eeq
between the associated graded Lie algebras (see \eqref{sec4eq.2} below and appendix \ref{app:deg}).

 We write
${\rm Der}^0{\rm Lie}[a,b]$ for the subspace of Lie algebra derivations of
${\rm Lie}[a,b]$ which annihilate the bracket $[a,b]=ab-ba$, where the last expression is valued in  $\Q\langle a,b\rangle$.
A derivation
in ${\rm Der}^0{\rm Lie}[a,b]$ is entirely determined by its value on $a$ (see for example Thm.~2.1 of \cite{Schneps:2012} giving an explicit formula for the value of such a derivation on $b$).

\begin{dfn} 
\label{abdegrs} Let $\delta\in {\rm Der}^0{\rm Lie}[a,b]$.
 We say that
$\delta$ is  of {\em homogeneous degree $n$} if $\delta(a)$ (and thus also $\delta(b)$)
is a Lie polynomial of homogeneous degree $n+1$, i.e.\ if $\delta$ adds $n$ to the
degree of any polynomial it acts on. We furthermore assign {\em $a$-degree $k$ and $b$-degree $\ell$} to $\delta$ if $\delta(a)$ is a Lie polynomial of homogeneous degree $k+1$ in $a$ and $\ell$ in $b$, in which case $\delta(b)$ is necessarily of $a$-degree $k$ and $b$-degree $\ell+1$ (unless it vanishes). The $b$-degree of a derivation and the homogeneous $b$-degree of a polynomial in $a,b$ is also referred to as the {\em depth}. The (homogeneous) degree of $\delta$ is equal to the sum of its $a$- and its $b$-degree.
\end{dfn}

\vspace{.3cm}
We now need to introduce the {\em Tsunogai derivations} which were introduced by Tsunogai in 1995 \cite{Tsunogai}, also see \cite{Tsunogai1}.

\begin{dfn}
\label{defmfu}
For all $i\ge 0$, let $\epsilon_{2i}$ denote the derivation of ${\rm Lie}[a,b]$
defined by 
\beq
\epsilon_{2i}(a)=\ad_a^{2i}(b)
\, ,\ \ \ \
\epsilon_{2i}([a,b])=0
\, ,\ \ \ \ i\ge 0\,.
\label{epsaction}
\eeq
These two conditions determine $\epsilon_{2i}$ completely: its action on $b$ is given explicitly by
\begin{align}
\epsilon_0(b) =0
\quad\text{and}\quad
    \epsilon_{2i}(b) = \sum_{j=0}^{i-1} (-1)^j  \left[ \ad_a^j(b), \ad_a^{2i-1-j}(b)
    \right] \, ,\ \ \ \ i\geq 1\,.
    \label{ep2ionb}
\end{align}
We write $\mathfrak{u}$ for the Lie algebra of derivations of ${\rm Lie}[a,b]$ generated by the $\epsilon_{2i}$ for $i\ge0$; the $\epsilon_{2i}$ are also called {\em geometric derivations}. 
\end{dfn}

\vspace{.3cm}
The Lie algebra $\mathfrak{u}$ of geometric derivations $\epsilon_{2i}$ has a rich history dating back to pioneering work of Ihara \cite{Ihara:1990}, with detailed studies in the work of Tsunogai \cite{Tsunogai, Tsunogai1}. They have become ubiquitous in the theory of elliptic MZVs as reviewed in appendix \ref{app:eMZV}, see for example \cite{KZB, EnriquezEllAss, Enriquez:Emzv,  Broedel:2015hia, hain_matsumoto_2020, LMSonEMZV} and \cite{Pollack}, with numerous references in the recent mathematics and string-theory literature.
The derivations $\epsilon_0$ and $\epsilon_2$ defined in~\eqref{epsaction} play a special role. The derivation $\epsilon_0$ is nilpotent on the $\epsilon_k$ (with even $k\geq 2$) in the sense that $\ad_{\epsilon_0}^{k-1} (\epsilon_k) =0$, see part (i) of Lemma \ref{lem:513} below.
The derivation $\epsilon_2$ is central in ${\rm Der}^0{\rm Lie}[a,b]$ and will play no role in our construction.

We will also make essential use of the following $\mathfrak{sl}_2$-subalgebra of ${\rm Der}^0{\rm Lie}[a,b]$:
\begin{dfn}
    \label{dfn:sl2}
    Define derivations $\epsilon_0^\vee, \hhh \in {\rm Der}^0{\rm Lie}[a,b]$ by
\begin{align}
    \epsilon_0^\vee(a) = 0 \,,\ \ \ \  \epsilon_0^\vee(b) = a  \,,\ \ \ \ 
    \hhh= [\epsilon_0,\epsilon_0^\vee]
    \,.
    \label{ep0vee}
\end{align}
The derivations $\epsilon_0$, $\epsilon_0^\vee$ and $\hhh$ generate the Lie subalgebra of ${\rm Der}^0{\rm Lie}[a,b]$ denoted $\mathfrak{sl}_2$. The generator $\hhh$ satisfies $\hhh(a) = -a$ and $\hhh(b)=b$. We refer to vectors that are annihilated by $\ep_0$ as {\em highest-weight vectors} and vectors that are annihilated by $\ep_0^\vee$ as {\em lowest-weight vectors}, respectively.
\end{dfn}

\begin{dfn}
    \label{dfn:switch}
We will also need to introduce the {\em switch} operator $\theta$, which can be considered as the automorphism of $\Q\langle\langle a,b\rangle\rangle$ that exchanges $a$ and $b$, mapping a polynomial $f=f(a,b)$ to $\theta(f)$ with $[\theta(f)](a,b)=f(b,a)$, but also acts on derivations $\delta$  of $\Q\langle a,b\rangle$ by conjugation via the formula
\beq
\theta(\delta)\coloneqq \theta\circ\delta\circ\theta^{-1}\, ,
\eeq
i.e.
\beq
[\theta(\delta)](a)=\theta\bigl(\delta(b)\bigr) \, ,\ \ \ \  
[\theta(\delta)](b)=\theta\bigl(\delta(a)\bigr) \, .
\eeq
Notice that  $\theta(\epsilon_0)=\epsilon_0^\vee$ and therefore $\theta(\hhh)=-\hhh$.
\end{dfn}

\vspace{.3cm}
The interplay of the derivations $\ep_k$ with the $\mathfrak{sl}_2$-algebra and the switch operation $\theta$ in the previous definitions is reviewed in the following lemma (see for instance \cite{Tsunogai, hain_matsumoto_2020, Pollack}).

 \begin{lemma}
 \label{lem:513} 
 For even values $k\ge 2$ and even or odd $j\ge 0$, set 
\beq\label{defepsilonkj}
\epsilon_k^{(j)} \coloneqq  \ad_{\epsilon_0}^j(\epsilon_k)
\eeq
including $\ep_k^{(0)}=\ep_k$. Then the $\epsilon_k^{(j)}$ for $k\geq 2$ together with the generators $\ep_0,\ep_0^\vee,\hhh$ of the $\mathfrak{sl}_2$ in Definition \ref{dfn:sl2} satisfy the following properties:
\begin{itemize}
\item[(i)] The derivation $\epsilon_k^{(j)}$ is of $a$-degree $k-j-1$ and $b$-degree $j+1$ for $0\le j\le k-2$ (in other words, $\epsilon_k^{(j)}(a)$ is a polynomial of homogeneous $a$-degree $k-j$ and $b$-degree $j+1$) and thus of homogeneous degree $k$. We have the nilpotency property
\beq
\epsilon_k^{(j)}=0 \quad \forall \ j>k-2
\,.
\label{nilpote0}
\eeq
The $\ep_k^{(k-2)}$ at maximum value of $j$ are highest-weight vectors of the $\mathfrak{sl}_2$.
\item[(ii)] The derivations $\epsilon_k$ with $k\geq 2$ commute with
$\ep_0^\vee$:
\beq
[\ep_0^\vee , \ep_k ] = 0 \quad \forall \ k\geq 2\, ,
\label{lwproperty}
\eeq
i.e.\ they furnish lowest-weight vectors of $\mathfrak{sl}_2$.
\item[(iii)] The generator $\hhh$ of
$\mathfrak{sl}_2$ satisfies the following commutation relations:
\beq
[\hhh , \ep_k ] = (2-k) \ep_k \quad \forall \ k\geq 0 \, , \ \ \ \
[\hhh , \ep_0^\vee ] = -2 \ep_0^\vee
\,.
\label{hbrak.01}
\eeq
In particular this implies that the $\epsilon_k^{(j)}$ are all eigenvectors
for $\hhh$, with eigenvalues given~by
\beq
[\hhh , \ep_k^{(j)} ]
= (2+2j-k)  \ep_k^{(j)} 
\quad \forall \ k\geq 2 \, , \ 0\leq j \leq k-2\, .
\label{hbrak.02}
\eeq
\item[(iv)] The commutation relations of the $\mathfrak{sl}_2$ generators with $ \ep_k^{(j)}$ at $k\geq 2$ and
$0\leq j \leq k-2$ are
$[\ep_0,  \ep_k^{(j)}]=  \ep_k^{(j+1)}$
by definition, $[\hhh , \ep_k^{(j)} ] = (2j+2-k)\ep_k^{(j)}$ by the previous point and
\begin{align}
[\ep_0^\vee , \ep_k^{(j)} ] &= j (k-1-j)\ep_k^{(j-1)}\, . \label{epjkep} 
\end{align}
\item[(v)] The switch operator in Definition \ref{dfn:switch} acts on the $\ep_k^{(j)}$ with $k\geq 2$ and $0\leq j \leq k-2$~via
\begin{align}
    \theta\bigl(\epsilon_k^{(j)}\bigr) &=-\frac{j!}{(k-2-j)!}\epsilon_k^{(k-2-j)}\, .
   \label{lemepjk.b}
\end{align}
\end{itemize}

\end{lemma}

\vspace{.2cm}
\noindent {\bf Proof.} (i) The derivation $\epsilon_k$ is of $a$-degree $k-1$ and $b$-degree $1$ by definition, and each application of $\ad_{\epsilon_0}$ increases the $b$-degree by 1 without changing the total degree, so it decreases the $a$-degree by $1$, proving the first statement. For the second statement, it is enough to show that $\epsilon_k^{(k-1)}=0$ even though since $\ep^{(j)}_k$ shifts the $(a,b)$ degrees of any polynomial in $a,b$ by $(k-1-j,1+j)$, the case $\ep^{(k-1)}_k$ of interest has $(a,b)$ degrees $(0,k)$ as a derivation, meaning that a priori $\epsilon_k^{(k-1)}(a)$ could be a polynomial of $a$-degree $1$ and $b$-degree $k$. Since the only Lie polynomial with these degrees is $\ad_b^k(a)$ up to scalar multiple, we must have 
\beq
\epsilon_k^{(k-1)}(a)=c\cdot \ad_b^k(a)
\eeq
for some constant $c$, and $\epsilon_k^{(k-1)}(b)=0$. However, the derivation $\ep^{(k-1)}_k$ must annihilate the commutator $[a,b]$ since both $\ep_k$ and $\ep_0$ do, so by the above, we have $\ep^{(k-1)}_k([a,b])= c\cdot[ \ad_b^k(a),b]$ which only vanishes for $c=0$. Thus $c=0$, so the derivation $\epsilon_k^{(k-1)}=0$.

(ii) is readily established by evaluating $[\ep_0^\vee, \ep_{2i}]=\ep_0^\vee \ep_{2i} - \ep_{2i} \ep_0^\vee$ on $a$ and $b$. The least straightforward part of the computation is to note that $\ep_0^\vee \sum_{j=0}^{i-1}(-1)^j [\ad_a^j(b) , \ad_a^{2i-1-j}(b)]$ receives a single contribution from the $j=0$ term, resulting in $[\ep_0^\vee(b), \ad_a^{2i-1}(b)]= \ep_{2i}(a)$.

(iii) Any monomial in $a,b$ is an eigenvector for $\hhh$, with the difference of the $b$-degree minus the $a$-degree as its eigenvalue. Since $\ep_k$ at $k\geq 0$ and $\ep_0^\vee$ shift the $(a,b)$-degrees by $(k-1,1)$ and $(1,-1)$, respectively, the associated differences ``$b$-degree minus $a$-degree'' are shifted by $2-k$ in case of $\ep_k$ and $-2$ in case of $\ep_0^\vee$. This implies both identities in (\ref{hbrak.01}) as eigenvalue equations. The second claim (\ref{hbrak.02}) is a corollary which can for instance be inferred from $\ep_k^{(j)}$ shifting the $(a,b)$-degrees by $(k-1-j,j+1)$.

(iv) One can conveniently prove (\ref{epjkep}) by induction in $j$, starting with $[\ep_0^\vee,\ep_k^{(0)}]=0$ as a base case which follows from (ii). The inductive step relies on the Jacobi identity $[\ep_0^\vee,\ep_k^{(j)}]= [\ep_0^\vee,[\ep_0,\ep_k^{(j-1)}]]
=[ [\ep_0^\vee,\ep_0],\ep_k^{(j-1)}]
+[\ep_0,[\ep_0^\vee,\ep_k^{(j-1)}]]$ as well as (\ref{hbrak.02}) to evaluate the first term $[ [\ep_0^\vee,\ep_0],\ep_k^{(j-1)}]=- [ \hhh,\ep_k^{(j-1)}]$.

(v) We proceed by induction in $j$, first proving $\theta(\ep_k)= -\frac{1}{(k-2)!} \ep_k^{(k-2)}$ as a base case of (\ref{lemepjk.b}) at $j=0$.

{\em Base case:} If a derivation of degree $>0$ annihilates the bracket $[a,b]$, then knowing its value on one of the variables $a$ or $b$ determines it completely. Hence, it suffices to show that $\theta(\ep_k)$ and $ -\frac{1}{(k-2)!} \ep_k^{(k-2)}$ have the same action on $b$ to establish their equality as derivations in ${\rm Der}^0{\rm Lie}[a,b]$. For this purpose, we successively simplify
\begin{align}
\ep_{k}^{(k-2)}(b) &= (\ep_0)^{k-2}\ep_{k}(b) =\sum_{j=0}^{\frac{k}{2}-1} (-1)^j     (\ep_0)^{k-2}  \left[ \ad_a^j(b), \ad_a^{k-1-j}(b)
    \right] \notag \\
    &= (\ep_0)^{k-2}  \left[ b, \ad_a^{k-1}(b)
    \right] 
    = - \left[ b,  (\ep_0)^{k-2} \ad_a^{k-2}([b,a])
    \right] \notag \\
    &=-(k-2)! \left[ b,   \ad_b^{k-2}([b,a])
    \right] =  -(k-2)!\,\ad_b^{k}(a)\, .
    \label{ep2idual}
\end{align}
In the first step, we have used $\ep_0(b)=0$ to remove all contributions to $\ep_{k}^{(k-2)}(b)$ with an $\ep_0$ on the right of $\ep_{k}$. The second step makes use of the expression (\ref{ep2ionb}) for $\ep_{k}(b)$ and $k$ even. The third step relies on the fact that for $m\ge 1$,
$\ad_a^m(b)$ is annihilated by $(\ep_0)^m$ such that $[ \ad_a^j(b), \ad_a^{k-1-j}(b)]$ is annihilated by $(\ep_0)^{k-2} $ unless $j=0$. After redistributing the $(k-1)$-fold action of $\ad_a$ in the fourth step, we note in the fifth step that the $k-2$ factors of $\ep_0$ can act on the $k-2$ exposed powers of $\ad_a$ (besides $[b,a]$ which is annihilated by $\ep_0$) in $(k-2)!$ different permutations, converting $\ad_a^{k-2}$ to $\ad_b^{k-2}$ in all cases. The end result of (\ref{ep2idual}) after repackaging the powers of $\ad_b$ is equivalent to 
\beq
\ep_{k}^{(k-2)}(b)
= -(k-2)! \,\ad_b^{k}(a) 
= -(k-2)!\, \theta\big( \ep_k(a) \big)
\eeq
by virtue of (\ref{epsaction}). As a consequence, $\theta(\ep_k)$ and $ -\frac{1}{(k-2)!} \ep_k^{(k-2)}$ have the same action on $b$ and must agree as derivations since they both annihilate $[a,b]$ and have degree $>0$.

{\em Inductive step:}
Now we can take care of (\ref{lemepjk.b}) at values $j>0$ by induction as follows:
\begin{align}
\theta(\ep_k^{(j)}) &=\theta\big( [\ep_0,\ep_k^{(j-1)}] \big) = \big[\theta(\ep_0), \theta(\ep_k^{(j-1)} )\big] 
=- \frac{(j-1)!}{(k-1-j)!} [ \ep_0^\vee,  \ep_k^{(k-1-j)} \big]\notag \\
& =
-\frac{(j-1)!}{(k-1-j)!} j (k-1-j)\ep_k^{(k-2-j)} =  -\frac{j!}{(k-2-j)!} \ep_k^{(k-2-j)} \, ,
\label{indeps}
\end{align}
where we used $\theta(\ep_0) = \ep_0^\vee$ and the induction hypothesis
$\theta(\ep_k^{(j-1)} ) = -\frac{(j-1)!}{(k-1-j)!} \ep_k^{(k-1-j)}$ in the third step and (\ref{epjkep}) proven as (iv) in passing to the second line.  
\qed

\begin{rmk}
\label{Pollackrel}
Note that the $\epsilon_k^{(j)}$ are by no means free generators of $\mathfrak{u}$; commutators of two or more of them obey a number of relations related to period polynomials of holomorphic cusp forms on ${\rm SL}_2(\mathbb{Z})$, the first of which were noticed by Ihara and Takao (cf.~\cite{IharaTakao}). The relations between brackets of two $\epsilon_k$'s were classified in \cite{Schneps:2006} where the connection with cusp forms was made explicit; subsequently Pollack in \cite{Pollack} unearthed many more relations, and made a general conjecture about the full set of relations between the $\epsilon_k^{(j)}$. These relations, which we call {\it Pollack's relations}, were proved to be motivic in \cite{hain_matsumoto_2020}. They appear in many works related to elliptic MZVs, such as for example
\cite{LNT} and \cite{Broedel:2015hia}, see appendix \ref{app:eMZV} for a brief recap. The lowest-degree Pollack relations  arise in degrees 14 and 16, and are given by
\begin{align}
0 &= [\ep_4,\ep_{10}] - 3 [\ep_6,\ep_8] \, ,
\label{tsurels} \\
0 &= 80 [\ep_4^{(1)}, \ep_{12} ] + 16 [\ep_{12}^{(1)},\ep_4]
- 250 [\ep_6^{(1)},\ep_{10}] -125 [\ep_{10}^{(1)},\ep_6] + 280 [\ep_8^{(1)},\ep_8] \notag \\*
&\quad
- 462[\ep_4, [\ep_4,\ep_8]] - 1725 [\ep_6,[\ep_6,\ep_4]] \,.\label{tsurels2}
\end{align}
\end{rmk}

\subsection{The genus one motivic Lie algebra}\label{sec:sigman.1}

In \cite{hain_matsumoto_2020}, Hain and Matsumoto define a Tannakian category
$MEM$ of {\em mixed elliptic motives} and study its fundamental Lie algebra.
We do not recall their construction here, but restrict ourselves to giving
the main result of their article that we will use here. Let ${\rm Lie}\,\pi_1(MEM)$ denote the graded Lie algebra associated to the unipotent radical of the fundamental group of the category $MEM$. Let $\mathfrak{sl}_2$ denotes the Lie subalgebra of ${\rm Der}^0{\rm Lie}[a,b]$ from Definition~\ref{dfn:sl2}.
\begin{thm}[Hain--Matsumoto]
\label{thm:HM}
There is a Lie algebra
morphism (the ``monodromy representation'', see section 22 of \cite{hain_matsumoto_2020})
\beq
{\rm Lie}\,\pi_1(MEM)\rightarrow {\rm Der}^0{\rm Lie}[a,b]\, 
\eeq
whose image ${\cal L}$ is generated by the derivations $\epsilon_k^{(j)}$ for even $k>0$ and $0\le j\le k-2$ together with derivations $\sigma_w$ for each odd $w\ge 3$, and has the following properties:
\begin{itemize}
    \item[(i)] The Lie subalgebra ${\cal S} \coloneqq {\rm Lie}[\sigma_3,\sigma_5,\ldots]\subset {\cal L}$ is free,
    \item[(ii)] The Lie subalgebra $\mathfrak{u}$ generated by the $\epsilon_k^{(j)}$ is normal in ${\cal L}$, i.e.\ ${\cal L}= \mathfrak{u} \rtimes \sigalg$,
    \item[(iii)] ${\cal L}$ is an 
$\mathfrak{sl}_2$-module,  and $\mathfrak{u}$ is also an $\mathfrak{sl}_2$-module,
\item[(iv)] the Lie subalgebra $\mathfrak{u} \rtimes \mathfrak{sl}_2$ is normal inside ${\cal L} \rtimes \mathfrak{sl}_2$.
\end{itemize}
\end{thm}

\vspace{.2cm}
\begin{rmk}\label{framework}  Although entirely phrased in terms of the monodromy representation of the fundamental Lie algebra of the category $MEM$, this theorem reflects essential geometric/arithmetic content. The quotient of ${\cal L}$ by the normal Lie subalgebra $\mathfrak{u}$ is isomorphic to~${\cal S}$, which is itself free on one generator in each odd rank $\geq 3$, i.e.~isomorphic to ${\rm Lie}\,\pi_1(MTM)$ the fundamental Lie algebra of the category of mixed Tate motives unramified over $\mathbb{Z}$. This in turn reflects the geometric situation where an elliptic curve parametrized by $\tau$ degenerates to the nodal elliptic curve when $\tau$ tends to $i\infty$ (see appendix~\ref{app:deg}). 

To be more precise, if one considers the universal elliptic curve ${\cal E}$ as a fibration over the Deligne--Mumford compactification $\overline{\cal M}_{1,1}$ of the moduli space of elliptic curves ${\cal M}_{1,1}$ (viewed as the usual fundamental domain for the action of ${\rm SL}_2(\mathbb{Z})$ on the Poincar\'e upper half-plane, parametrized by the variable $\tau$), then the fiber over $\tau=i\infty$ is the so-called nodal (or degenerate) elliptic curve $E_\infty$.  Let $\pi_1$ denote the fundamental group of the punctured torus, freely generated by loops $\alpha$ and $\beta$ through and around the genus hole, and let $\hat\pi_1$ be its profinite completion. Then there is a canonical {\em arithmetic} outer Galois action of the absolute Galois group ${\rm Gal}(\overline{\Q}/\Q)$ on  $\hat\pi_1(E_\infty)$. Furthermore, since ${\cal E}$ is a fibration over the base ${\cal M}_{1,1}$ with an elliptic curve as a fiber, $\pi_1({\cal E})$ fits into a short exact sequence whose kernel is free on two generators (the $\pi_1$ of the fiber) and whose quotient is ${\rm SL}_2(\mathbb{Z})$ (the $\pi_1$ of the base). Hence, there is a second, {\em geometric} outer group action on  $\pi_1(E_\infty)$ by the group ${\rm SL}_2(\mathbb{Z})$, which extends to an action of the profinite completion $\widehat{\rm SL}_2(\mathbb{Z})$ on $\hat\pi_1(E_\infty)$.  Thus we have two disjoint profinite groups, $\widehat{\rm SL}_2(\mathbb{Z})$ and the absolute Galois group ${\rm Gal}(\overline{\Q}/\Q)$~\cite{SGA1}, acting as automorphism groups of  $\hat\pi_1(E_\infty)$. 

The pro-unipotent version of this situation, or rather the associated Lie algebra version, has $\sigalg = {\rm Lie}\,\pi_1(MTM)$ playing the role of ${\rm Gal}(\overline{\Q}/\Q)$ and $\mathfrak{u}\rtimes {\rm sl}_2$ playing the role of ${\rm SL}_2(\mathbb{Z})$, both acting as derivation Lie algebras (the Lie algebra version of automorphism groups) of  
${\rm Lie}[a,b]$, the free Lie algebra on two generators which plays the role of $\hat\pi_1(E_\infty)$. The fact that $\sigalg$ acts on $\mathfrak{u}\rtimes\mathfrak{sl}_2$  
reflects the fact that ${\rm Gal}(\overline{\Q}/\Q)$ acts not only on $\hat\pi_1(E_\infty)$ but also on $\widehat{\rm SL}_2(\mathbb{Z})$, since the latter group is also a fundamental group, namely of ${\cal M}_{1,1}$.
\end{rmk}

\vspace{.2cm}
Hain and Matsumoto conjecture that the surjective morphism
from ${\rm Lie}\,\pi_1(MEM)$ to ${\cal L}$ is actually an isomorphism, but
this is still an open question.
They further explain that there is a natural surjection from ${\rm Lie}\,\pi_1(MEM)$ to ${\rm Lie}\,\pi_1(MTM)$, the fundamental Lie algebra of the category of mixed Tate motives unramified over $\mathbb{Z}$. Since this category was shown by Brown to be generated by the motivic MZVs, we have the isomorphism
\beq
{\rm Lie}\,\pi_1(MTM) = \mz^\vee\, ,
\eeq
where $\mz^\vee$ is the Lie algebra associated to the motivic MZVs. Hain and Matsumoto further proved the existence of a section map
\beq\label{sectionmap}
{\rm Lie}\,\pi_1(MTM)\hookrightarrow {\rm Lie}\,\pi_1(MEM)\, ,
\eeq
which explains the semi-direct product structure in Theorem \ref{thm:HM} (ii), with the image of ${\rm Lie}\, \pi_1(MTM)$ identified with $\sigalg \subset {\rm Lie} \,\pi_1(MEM)$.
The section map was defined explicitly in independent parallel work by Enriquez in \cite{EnriquezEllAss}, working with the {\em Grothendieck--Teichm\"uller Lie algebra} $\grt$. Thanks to this work, $\sigalg$ is  identified as a canonical Lie subalgebra of ${\cal L}$. However, neither Hain--Matsumoto nor Enriquez gave a canonical choice of the actual generators $\sigma_w$ for odd $w\ge 3$; a priori, the choice of generator $\sigma_w$ is
only defined up to adding on brackets of $\sigma_u$ with smaller $u<w$. This exactly
parallels the fact that no special set of free generators of the motivic Lie algebra ${\rm Lie}\,\pi_1(MTM) = \mz^\vee$ was defined prior to the canonical family of $g_w$ in genus zero defined in section~\ref{sec:Zmap.4}.

Our main purpose in this section is to point out that, thanks to the canonical genus zero generators $g_w$ and the existence of the section map \eqref{sectionmap}, we can now define a canonical choice of genus one generators $\sigma_w$ simply as the images of the $g_w$ under the section map. More precisely, we will construct an explicit Lie algebra morphism
\beq
\tilde\gamma:\mz^\vee \rightarrow {\rm Lie}[\sigma_3,\sigma_5,\ldots,]\subset {\rm Der}^0{\rm Lie}[a,b]
\label{sec4eq.2old}
\eeq
and use it to define the $\sigma_w$ (as images of the $g_w$), to compute them and to determine many of their properties. In the same way as the Ihara derivations of $g_w$ are called zeta generators in genus zero, we will refer to the $\sigma_w$ as {\it zeta generators in genus one}. The tight interplay of zeta generators in genus zero and one can for instance be seen from (\ref{g01.11}) below where the action of $\sigma_w$ is computed from $g_w$. Additional facets of the relation between zeta generators in genus zero and genus one can be found in appendix~\ref{app:deg}.

Let us show how the map $\tilde\gamma$ in (\ref{sec4eq.2old}) relates to the Grothendieck--Teichm\"uller section map defined by Enriquez. We do not need to give the definition of $\grt$ here, but only to mention two essential properties that we need: firstly, there is an injective morphism 
\begin{align}\label{minusy}
\mz^\vee&\hookrightarrow\grt \, ,\notag\\
h(x,y)&\mapsto h(x,-y) \, ,
\end{align}
(this is a direct consequence of the fact that Goncharov's motivic MZVs satisfy the associator relations, see for example \cite{Andre}) and secondly, Enriquez \cite{EnriquezEllAss} defined an injective map
\beq\label{enriquezsection}
\grt\hookrightarrow  {\rm Der}^0{\rm Lie}[a,b]
\eeq
which was shown in \cite{Schneps:2015mzv} to be equivalent to the Hain--Matsumoto section, using methods from \'Ecalle's mould theory that will be explained in section \ref{sec:prop} below. Let
\beq\label{gammamap}
\gamma:\mz^\vee\hookrightarrow {\rm Der}^0{\rm Lie}[a,b]
\eeq
denote the composition of \eqref{minusy} with \eqref{enriquezsection}. 
The explicit isomorphism $\tilde\gamma$ announced in \eqref{sec4eq.2old} is given by
\beq\label{gammatilde}\tilde\gamma=\theta\circ\gamma\, ,
\eeq
where $\theta$ is the switch automorphism of $\Q\langle\langle a,b\rangle\rangle$ exchanging $a$ and $b$, see Definition \ref{dfn:switch}. 

\begin{dfn}\label{thesigmaw}
Let $g_w$ for odd $w\ge 3$ denote the family of canonical free generators of $\mz^\vee$ given in Definition~\ref{dfn:mzvee}. Set
\beq
\tau_w \coloneqq \gamma(g_w)\, , \ \ \ \
\sigma_w \coloneqq \tilde\gamma(g_w)\, ,
\label{deftausig}
\eeq
where $\gamma$ is as in \eqref{gammamap} and $\tilde\gamma$ as in \eqref{gammatilde}. This definition accomplishes the second goal of this article of giving a canonical choice for the zeta generators $\sigma_w$ in genus one for odd $w\ge 3$.
\end{dfn}

\vspace{.3cm}
The remainder of section \ref{sec:sigman} and all of sections \ref{sec:prop} and \ref{bigsecmd2} are devoted to the study of the canonical zeta generators $\sigma_w$ in genus one. Section \ref{sec:sigman.2} gives an explicit step-by-step construction of the Enriquez map \eqref{enriquezsection}, and in Theorem~\ref{thm:522} of section \ref{sec:sigman.3} we list several properties of the zeta generators $\sigma_w$ and their switch images $\tau_w$. Section \ref{sec:sigman.4} contains the low-degree parts of $\sigma_w$ for $w=3,5,7,9$. The proofs of some of the properties in Theorem \ref{thm:522} rely on a second, mould theoretic construction of the map $\gamma$, which is given in section \ref{sec:prop1} along with a necessary introduction to mould theory; the full proof of the theorem is contained in section \ref{sec:prop2} (using mould theory), section \ref{sec:prop3} (using the $\mathfrak{sl}_2$ subalgebra of Definition \ref{dfn:sl2}) and section \ref{newsec:71} (summarizing the essential argument of \cite{Hain, hain_matsumoto_2020}). Section \ref{secmd2} introduces a recursive procedure to compute high-degree contributions to $\sigma_w$ in terms of $\ep_k$ which leads to a variety of explicit results beyond the state-of-the-art in section \ref{newsec:74}.

\subsection{Genus one derivations from genus zero polynomials}\label{sec:sigman.2}

Since in this section we will work only in odd weights $w$, we can work entirely
mod $\zeta_2$, in the $\Q$-algebras $\overline{\FZ}$ and $\overline{\MZ}$.

The surjection $\FZ\rightarrow\!\!\!\!\!\rightarrow \MZ$ from section~\ref{sec:MMZV}
induces a surjection $\overline{\FZ}\rightarrow\!\!\!\!\!\rightarrow\overline{\MZ}$ and a surjection $\fz\rightarrow\!\!\!\!\!\rightarrow\mz$. As we saw in the previous sections, we can pass to the dual spaces using the Z-map and these surjections induce
injections $\mz^\vee\hookrightarrow\fz^\vee=\ds$ and $\overline{\MZ}^\vee\hookrightarrow\overline{\FZ}^\vee=\U\ds$
in the dual spaces. The complete situation combining all the surjections, dual inclusions and Z-maps is summarized in the diagram \eqref{diag0}.

\vspace{.2cm}
The map from $g_w$ to $\sigma_w$ is to be viewed as a map from genus zero
to genus one, see appendix \ref{app:deg}. The genus zero situation here is represented by the Lie algebra ${\rm Lie}[x,y]$,
which is identified with the graded Lie algebra associated to the pro-unipotent completion of the fundamental group~$\pi_1$ of
the sphere with three punctures (which is free on two generators). The genus one situation
is represented by the completion $\widehat{\rm Lie}[a,b]\subset \Q\langle\langle a,b\rangle\rangle$ of the free Lie algebra on two generators ${\rm Lie}[a,b]$, the graded Lie algebra of the pro-unipotent
fundamental group of the once-punctured torus. The topological map from the sphere to the torus obtained by
joining two of the punctures passes to the topological fundamental groups, their unipotent completions and then via formality isomorphisms to the 
corresponding graded Lie algebras, yielding the following Lie algebra morphism: 
\begin{align}
\psi:{\rm Lie}[x,y]&\rightarrow \widehat{\rm Lie}[a,b]\, ,\notag\\
x&\mapsto t_{12}\, ,\notag\\
y&\mapsto t_{01} \, ,
\label{sec4eq.2}
\end{align}
where letting ${\rm B}_n$ denote the standard Bernoulli numbers,
\begin{align}
t_{01}&\coloneqq \frac{\ad_b}{e^{\ad_b}-1}(-a)=-a-\sum_{n\ge 1} \frac{{\rm B}_n}{n!}\ad_b^n(a)
=-a + \tfrac{1}{2} \ad_b(a)
- \tfrac{1}{12} \ad_b^2(a)
+ \tfrac{1}{720} \ad_b^4(a) +\ldots \, ,
\notag \\
t_{12}&\coloneqq [a,b] \, .
\label{g01.01}
\end{align}
The map $\psi$ in \eqref{sec4eq.2} also arises when computing the Knizhnik--Zamolodchikov--Bernard connection on a degeneration limit of the torus (corresponding topologically to the degenerate torus obtained by joining two punctures of the thrice-punctured sphere), and matching the result with the Knizhnik--Zamolodchikov connection on the sphere. This calculation is spelled out in detail in appendix \ref{app:deg}.

In order to explicitly define the map $\gamma$ in \eqref{gammamap}, we will make use of the notion of a {\it partner}~\cite{Schneps:2015mzv}: for any $g(a,b)\in {\rm Lie}[a,b]$, we write $g=g_aa+g_bb$ and
define the partner of $g$ by the formula
\beq
g' \coloneqq \sum_{i\ge 0} \frac{(-1)^{i-1}}{i!} a^ib\,\partial^i_a\bigl(g_a\bigr)\in \Q\langle a,b\rangle\,,
\label{g01.03}
\eeq
where $\partial_a$ is the derivation of $\Q\langle\langle a,b\rangle\rangle$
defined by $\partial_a(a)=1$ and $\partial_a(b)=0$. It is shown in Lemma 2.1.1 of \cite{Schneps:2015mzv} that the derivation $a\mapsto g$, $b\mapsto g'$ lies in ${\rm Der}^0{\rm Lie}[a,b]$ if and only if $g$ has a certain property called {\it push-invariance} to which we will return in section \ref{sec:prop} (see \eqref{pushpoly}).

\vspace{.2cm}
We can now proceed to the explicit definition of the map $\gamma$ of \eqref{gammamap}. Define $\tau_h \coloneqq\gamma(\candh)\in {\rm Der}^0\widehat{\rm Lie}[a,b]$
to be the derivation obtained from $\candh \in \mz^\vee$ by the following procedure:
\begin{itemize}
    \item Let $\candh=\candh(x,y)$ be in $\mz^\vee$ and define a derivation $\kappa_h$ of the Lie subalgebra
 ${\rm Lie}[t_{12},t_{01}]\subset \widehat{\rm Lie}[a,b]$ 
 by\footnote{The minus sign in front of $t_{01}$ in \eqref{g01.04} is present because if $h(x,y)\in \mz^\vee\subset\ds$, then as in
 \eqref{minusy}, the polynomial $h(x,-y)$ lies in $\grt$. Since the process described in the present section is an explicit version of Enriquez's map \eqref{enriquezsection} from $\grt$ to ${\rm Der}^0{\rm Lie}[a,b]$, the starting point of the map is the $\grt$ polynomial $h(x,-y)$, or more precisely, the associated Ihara derivation which maps
 $x\mapsto 0$ and $y\mapsto [y,h(x,-y)]$. The first step in the explicit construction of the Enriquez map is transporting this Ihara derivation to a derivation on ${\rm Lie}[t_{01},t_{12}]$ via the map \eqref{sec4eq.2}, which is what is expressed in \eqref{g01.04}.}
\beq
\kappa_{\candh}(t_{12})=0\, , \ \ \ \ \kappa_{\candh}(t_{01})= [t_{01},\candh(t_{12},-t_{01})]\, .
\label{g01.04}
\eeq

\item 
By the ``extension lemma'' 2.1.2 of~\cite{Schneps:2015mzv}, there exists a unique derivation $\tau_{\candh}$ of $\Q\langle\langle a,b\rangle\rangle$ having the following two properties:
firstly
\beq
\tau_{\candh}(t_{01})=\kappa_{\candh}(t_{01}) \,, 
\label{g01.05b}
\eeq
and secondly $\tau_{\candh}(b)$ is (in each degree) the {\it partner} of $\tau_{\candh}(a)$ as defined in \eqref{g01.03}.

Specifically, the action of the derivation $\tau_h$ on $a$ can be inferred from (\ref{g01.05b}) degree by degree as follows. Suppose $h(x,y)$ is homogeneous of degree $w$ in $x,y$. We have from~\eqref{g01.01}
\beq
\tau_h(t_{01})=\tau_h\bigl(-a+\tfrac{1}{2}[b,a]-\tfrac{1}{12}[b,[b,a]]+\cdots\bigr)
\label{tauhabearlier}
\eeq
so 
\beq
\label{tauhab}
\tau_h(a)=-\kappa_h(t_{01})+\tfrac{1}{2}\tau_h\big([b,a]\big)-\tfrac{1}{12}\tau_h\big([b,[b,a]]\big)+\cdots
\eeq
since $\tau_\candh(t_{01})=\kappa_\candh(t_{01})$.
In particular, the lowest degree part of $\tau_h(a)$ is equal to the lowest degree part of $-\kappa_h(t_{01})$, which is equal to $[a,h^d([a,b],a)]$ from~\eqref{g01.04} and where~$d$ denotes the minimal $x$-degree of $h$ and $h^d(x,y)$ are the contributions to $h(x,y)$ of $x$-degree $d$; the term $[a,h^d([a,b],a)]$ is of degree $w+d+1$ in $a,b$. So we have
\beq\label{lowestpart}\tau_h(a)_{w+d+1}=-\kappa_h(t_{01})_{w+d+1}=
[a,h^d([a,b],a)]
\eeq
in lowest degree, where $g_d$ denotes the degree-$d$ contributions to polynomials $g$ in $a$ and $b$.
We set $\tau_h(b)_{w+d+1}$ to be the partner of $\tau_h(a)_{w+d+1}$ using the formula \eqref{g01.03}.

 We then use \eqref{tauhab} to recursively compute $\tau_h(a)$ in successive degrees $w+d+i$ ($i>1$):
\begin{align}
\tau_h(a)_{w+d+2}&=-\kappa_h(t_{01})_{w+d+2}+\tfrac{1}{2}[\tau_h(b)_{w+d+1},a]+\tfrac{1}{2}[b,\tau_h(a)_{w+d+1}] \, , \notag\\
\tau_h(a)_{w+d+3}&=-\kappa_h(t_{01})_{w+d+3}+\tfrac{1}{2}[\tau_h(b)_{w+d+2},a]+\tfrac{1}{2}[b,\tau_h(a)_{w+d+2}]\notag\\
&\ \ \ \ \ \ \ \
-\tfrac{1}{12}[\tau_h(b)_{w+d+1},[b,a]] 
-\tfrac{1}{12}[b,[\tau_h(b)_{w+d+2},a]]-\tfrac{1}{12}[b,[b,\tau_h(a)_{w+d+1}]]
\, ,
\notag\\
\hbox{etc.}\, ,\label{recurs}
\end{align}
defining $\tau_h(b)_{w+d+i}$ to be the partner of $\tau_h(a)_{w+d+i}$ at each successive degree via \eqref{g01.03}.
This process yields a unique Lie series $\tau_h(a)$. As observed just after \eqref{g01.03}, if $\tau_h(a)$ has the property of {\it push-invariance} then $\tau_\candh\in{\rm Der}^0\widehat{\rm Lie}[a,b]$, so in particular $\tau_\candh$ annihilates $[a,b]=t_{12}$, and thus $\tau_h$ is an extension of $\kappa_h$ to all of ${\rm Der}^0\widehat{\rm Lie}[a,b]$. The fact that $\tau_h(a)$ does indeed possess the necessary property of push-invariance  is proved in Theorem \ref{thm:613} (iii) below.

\item For each $h\in \mz^\vee$, we define $\sigma_h\in {\rm Der}^0\widehat{\rm Lie}[a,b]$ to be the derivation obtained from $\tau_h$ by the switch operator in Definition \ref{dfn:switch}:
we set 
\beq\label{g01.06bis}\sigma_{\candh}=\theta(\tau_{\candh})\, ,
\eeq
or equivalently, $\sigma_{\candh}$ acts on $a$ and $b$ via
\beq
\sigma_{\candh}(a)=\theta\bigl(\tau_{\candh}(b)\bigr)\, , \ \ \ \
\sigma_{\candh}(b)=\theta\bigl(\tau_{\candh}(a)\bigr)\, .
\label{g01.07}
\eeq
\end{itemize}

\vspace{.3cm}
\noindent Combining all the steps of the  process above then yields explicit versions  
\begin{align}
\gamma:\mz^\vee&\hookrightarrow {\rm Der}^0{\rm Lie}[a,b]\, ,
&\tilde \gamma:\mz^\vee&\hookrightarrow {\rm Der}^0{\rm Lie}[a,b] \, ,
\notag \\
 h &\mapsto \tau_{\candh} \, ,
 &h &\mapsto \sigma_{\candh} \, ,
 \label{g01.10}
\end{align}
of the maps $\gamma$ from \eqref{gammamap} and $\tilde\gamma$ from \eqref{sec4eq.2old}. 

\subsection{\texorpdfstring{The canonical genus one derivations $\sigma_w$}{The canonical genus one derivations sigma(w)}}
\label{sec:sigman.3}

We shall now specialize the above construction of $\gamma(h)$ and $\tilde \gamma(h)$ for general $h \in \mz^\vee$ to the canonical polynomials $h \rightarrow g_w$ of Definition~\ref{dfn:mzvee} for odd $w\geq 3$. The concrete realization of the maps $\gamma,\tilde \gamma$ in (\ref{g01.10}) provided by the previous section allows for an explicit computation of the zeta generators $\sigma_w,\tau_w$ in (\ref{deftausig}).
By (\ref{g01.04}) and (\ref{g01.07}), the action of the genus one zeta generators $\sigma_w =\tilde \gamma(g_w)$ on the smaller Lie subalgebra ${\rm Lie}[t_{01},t_{12}]\subset \widehat{{\rm Lie}}[a,b]$ is given~by
\beq
\sigma_w(t_{12})=0  \, , \ \ \ \
\sigma_w\big( \theta (t_{01}) \big) =  
\theta\big(\big[ t_{01} , g_w(t_{12},- t_{01})\big] \big)
\, ,
\label{g01.11}
\eeq
obtained from applying the switch $\theta$ to
\beq
\tau_w(t_{12})=0  \, , \ \ \ \
\tau_w(t_{01})=  
\big[  t_{01} , g_w( t_{12},-  t_{01})\big] \, .
\label{g01.11tau}
\eeq
With the notation $s_{01} = \theta(t_{01})$ and $s_{12} = \theta(t_{12})$ for the images under the
switch, (\ref{g01.11}) takes the more compact form $\sigma_w(s_{12})=0 $ and
$\sigma_w(  s_{01}) =  \big[ s_{01} , g_w(s_{12},- s_{01})\big] $ as advertised in
(\ref{sigintro}) in the introduction.

By the discussion in section \ref{sec:Drin}, the canonical polynomials $g_w$ are determined by the (modified) Drinfeld associator and the $\mathbb Q$ relations among MZVs. Hence, the information from iterated integrals in genus zero already fixes the defining relations (\ref{g01.11tau}) of zeta generators in genus one. Further discussions of the tight interplay between genus zero and genus one can be found in  appendix \ref{app:deg}.

In the previous section, we explained how to infer $\tau_h(a)$ and $\tau_h(b)$ from $\tau_h(t_{01})$ and $\tau_h(t_{12})$ for general $h \in \mz^\vee$ from (\ref{g01.04}) by the extension lemma 2.1.2 of \cite{Schneps:2015mzv}. To compute $\sigma_w(a)$ and $\sigma_w(b)$, we can either apply that method with $h=g_w$ and use the switch $\theta$ or use the same method directly from (\ref{g01.11}).

\vspace{0.3cm}
The derivations $\tau_w$ and $\sigma_w$ associated to $g_w$ for odd $w$ have many remarkable properties, of which a number are listed in the following theorem. Several of these are statements for the different degree parts of $\tau_w$ and $\sigma_w$ (where degree refers to the degree as a derivation). The degree $2w$ parts of $\tau_w$ and $\sigma_w$ turn out to play a special role and are called the {\em key degree} parts  $\tau^{\rm key}_w$ and $\sigma^{\rm key}_w$.  In section~\ref{sec:prop1}
we will present a brief introduction to mould theory which will enable us to
prove the first three of these in section~\ref{sec:prop2}; the others are proved in section~\ref{sec:prop3}. 
Part (i) and (ii) of the theorem below are already known from \cite{Schneps:2015mzv} and implicitly from \cite{EnriquezEllAss, hain_matsumoto_2020}.
Part (iv) follows straightforwardly from Theorem \ref{thm:HM} \cite{hain_matsumoto_2020}. Part (v) is essentially in  \cite{hain_matsumoto_2020}, see for instance Remark 20.4.
The last two sentences of part (vi) readily follow from Theorem \ref{thm:HM} as can be seen from their proof in section \ref{prfrest} below.
Part (vii) was proven in section 27 of \cite{hain_matsumoto_2020} as will be reviewed in section \ref{newsec:71} below.

\begin{thm}
\label{thm:522}
For odd $w\geq 3$, the zeta generators $\tau_w$ and $\sigma_w$ in Definition~\ref{thesigmaw} satisfy:
\begin{enumerate}
\item[(i)]
Both $\tau_w$ and $\sigma_w$ lie in ${\rm Der}^0\widehat{\rm Lie}[a,b]$.

\item[(ii)]
The minimal degree of $\tau_w$ and $\sigma_w$ 
is $w+1$, and all odd-degree terms are equal to zero. All terms of the power series $\tau_w(a)$ are of constant $a$-degree $w+1$, or equivalently (thanks to the switch), all terms of the power series $\sigma_w(a)$ have constant $b$-degree $w$.

\item[(iii)]
Both $\tau_w$ and $\sigma_w$ are entirely determined by their parts of
degree $<2w$.
\item[(iv)] There are no highest-weight vectors of $\mathfrak{sl}_2$ in $\sigma_w$ beyond key degree.
\item[(v)]
All contributions to $\tau_w$ and $\sigma_w$ of degree different from $2w$ lie in $\mathfrak{u}$. The key-degree parts $\tau^{\rm key}_w$ and $\sigma^{\rm key}_w$ do not lie in $\mathfrak{u}$.
\item[(vi)] Define the \emph{arithmetic part} $z_w \in {\rm Der}^0\widehat{\rm Lie}[a,b]$ of the derivation $\sigma_w$ to be the one-dimensional component of $\sigma^{\rm key}_w$ as an $\mathfrak{sl}_2$ representation, i.e.\ which commutes
with~the generators $\ep_0,\ep_0^\vee$ of $\mathfrak{sl}_2 \subset {\rm Der}^0 \widehat{\rm Lie}[a,b]$ in Definition~\ref{dfn:sl2}. Then, the difference  $\sigma^{\rm key}_w - z_w$ and by (v) in fact all of $\sigma_w {-} z_w$ lies in $\mathfrak{u}$. Moreover, while the $z_w$ themselves do not lie in $\mathfrak{u}$, the brackets $[z_w,\ep_k]$ for any even $k\geq 0$ lie in~$\mathfrak{u}$.
\item[(vii)] 
$\sigma_w$ commutes with the infinite series $N$ in geometric derivations defined by 
\beq
N\coloneqq -\epsilon_0 + \sum_{k=2}^\infty (2k-1) \frac{{\rm B}_{2k}}{(2k)!} \epsilon_{2k}\, .
\label{g01.12}
\eeq
\end{enumerate}
\end{thm}

\begin{rmk}
As pointed out in \cite{hain_matsumoto_2020, Brown:2017qwo2}, the characterization of the arithmetic parts $z_w$ in the earlier literature as commuting with $\mathfrak{sl}_2$ and not lying in $\mathfrak{u}$ does not identify the $z_w$ uniquely; ambiguities remain for $w\geq 7$, since one can modify $z_w$ by adding on $\mathfrak{sl}_2$-invariant combinations of $\ep_k^{(j)}$ in $\sigkey_w-z_w$ while keeping the overall $\sigma_w$ unchanged (see for instance Remark 20.3 (ii) of \cite{hain_matsumoto_2020}). In order to eliminate this ambiguity, we added the defining property in Theorem \ref{thm:522} (vi) that $z_w$ exhausts the one-dimensional irreducible $\mathfrak{sl}_2$ representations of $\sigkey_w$ (or equivalently, $\sigkey_w-z_w$ contains no one-dimensional irreducible representations of $\mathfrak{sl}_2$). Moreover, the canonical zeta generators $\sigma_w$ established with the help of the polynomials $g_w(x,y)$ resolve an independent class of earlier ambiguities in $z_w$, namely it is no longer possible to add on nested brackets of lower-weight $z_v$ with $v<w$ (e.g.\ for example, we cannot add a multiple of $[z_3,[z_3,z_5]]$ to $z_{11}$). Hence, the properties in part (vi) of Theorem \ref{thm:522} single out unique canonical arithmetic derivations $z_w$ at each odd $w\geq 3$.
\end{rmk}

\subsection{\texorpdfstring{Expansions of $\sigma_w$ in low degree}{Expansions of sigma(w) in low degree}}
\label{sec:sigman.4}

In this section we spell out the explicit low-degree parts of the $\sigma_w$ up
to $w=9$, in order to give a feel for their appearance. 
For this purpose, we rewrite the expansion of $\sigma_w(a)$ and $\sigma_w(b)$ resulting from (\ref{g01.11}) and the extension lemma in terms of the geometric derivations $\epsilon_{k}^{(j)}$ in (\ref{defepsilonkj}) acting on $a$ and $b$, up to the arithmetic parts $z_w$ at key degree described in Theorem \ref{thm:522} (vi). Note that according to Theorem \ref{thm:522} (v), the ``key
degree'' part $\sigkey_w$ of $\sigma_w$, which is the part in degree $2w$ (as a derivation) is the only part not consisting of
brackets of $\epsilon_{k}^{(j)}$. 

In section~\ref{sec:compasp} below we give a more detailed description of a first computation algorithm, but begin by presenting a few examples to convey an impression of the structure of the $\sigma_w$.
In the following
examples for $w=3,5,7$,
we decompose $\sigkey_w$ into the unique choice of its $\mathfrak{sl}_2$ invariant part $z_w$ in Theorem \ref{thm:522} (vi) and nested brackets of $\ep_k^{(j)}$ in $(\geq 3)$-dimensional irreducible representations of $\mathfrak{sl}_2$. An alternative algorithm for the computation of $\sigma_w$ --- in particular how to determine the infinity of terms beyond key degree $2w$ --- is described in section \ref{secmd2}.

\subsubsection{The case $w=3$} 

In this situation, we first give the complete calculation of the derivations $\tau_3$ and $\sigma_3$ related by the switch, and specify the arithmetic derivation $z_3$ by directly giving its values on $a$ and $b$. Recall that the switch maps $\tau_w$ to $\sigma_w$ via (\ref{g01.07}) and
 acts on the derivations $\epsilon_k^{(j)}$ according to (\ref{lemepjk.b}).
Direct computation based on (\ref{g01.11tau}) shows that 
\begin{align}
    \tau_3&=\epsilon_4+\tau_3^ {\rm key}-\tfrac{1}{960}[\epsilon_4^{(1)},\epsilon_4^{(2)}]+\tfrac{1}{725760}[\epsilon_4^{(1)},\epsilon_6^{(4)}]-\tfrac{1}{1451520}[\epsilon_4^{(2)},\epsilon_6^{(3)}]
 \notag\\
    &\quad   +\tfrac{1}{1741824000}[\epsilon_4^{(2)},\epsilon_8^{(5)}]
  -\tfrac{1}{870912000}[\epsilon_4^{(1)},\epsilon_8^{(6)}]
    +\tfrac{1}{2786918400}
    [\epsilon_4^{(2)},[\epsilon_4^{(2)},\epsilon_6^{(4)}]] \notag \\
    &\quad
    +\tfrac{1}{1931334451200}[\epsilon_4^{(1)},\epsilon_{10}^{(8)}]
    -\tfrac{1}{3862668902400}[\epsilon_4^{(2)},\epsilon_{10}^{(7)}]+\ldots \, ,
       \label{tau3stuff} 
\end{align}
with an infinite series in nested brackets of $\epsilon_{k_i}^{(j_i)}$
of total degree $\sum_i k_i \geq 16$ in the ellipsis.
Here and in section \ref{sec:ldgr} below, we have made a choice on how the Pollack relations of Remark \ref{Pollackrel} are used to represent the degree $\geq 14$ terms of $\sigma_w$ and $\tau_w$.

The key-degree part $\tau_3^{\rm key}$ concentrated in degree 6 is given explicitly by
\begin{align}
 \tau_3^{\rm key}(a)&=
    -\tfrac{1}{4}[aaababb]
-\tfrac{1}{4}[aaabbab]
-\tfrac{1}{12}[aababab] \, ,\notag\\
\tau_3^{\rm key}(b)&=\tfrac{1}{4}[aababbb]+\tfrac{1}{4}[aabbabb]+\tfrac{1}{4}[aabbbab]+\tfrac{1}{12}[abababb]\, ,
\end{align}
where we employ the Lyndon-bracket notation introduced in Theorem \ref{lem:lynd}.

Applying the switch 
(\ref{g01.07}) and (\ref{lemepjk.b})
to $\tau_3$ and $\tau_3^{\rm key}$, we obtain the following explicit formula for $\sigma_3$ (again skipping an infinity of contributions at degree $\sum_i k_i \geq 16$):
\begin{align}
\sigma_3 &= -\tfrac{1}{2} \epsilon_4^{(2)}+z_3+\tfrac{1}{480} [\epsilon_4,\epsilon_4^{(1)}]+\tfrac{1}{30240} [\epsilon_4^{(1)},\epsilon_6] -\tfrac{1}{120960} [\epsilon_4,\epsilon_6^{(1)}]+\tfrac{1}{7257600} [\epsilon_4,\epsilon_8^{(1)}]
\label{firsts3} \\
&\quad -\tfrac{1}{1209600} [\epsilon_4^{(1)},\epsilon_8]
-\tfrac{1}{58060800} [\epsilon_4,[\epsilon_4,\epsilon_6]]+\tfrac{1}{47900160} [\epsilon_4^{(1)},\epsilon_{10}]-\tfrac{1}{383201280} [\epsilon_4,\epsilon_{10}^{(1)}] + \ldots  \, .
\notag
\end{align}
For $w=3$ it turns out that the key-degree part  $\sigma_3^{\rm key}$ is already $\mathfrak{sl}_2$ invariant and therefore coincides with the arithmetic derivation $z_3$ whose action on $a$ and $b$ is given by
\begin{align}
z_3(a)&=\tfrac{1}{4}[aaababb]
+\tfrac{1}{4}[aaabbab]
+\tfrac{1}{12}[aababab]\, , \label{z3act}\\*
z_3(b)&=-\tfrac{1}{4}[aababbb]
-\tfrac{1}{4}[aabbabb]
-\tfrac{1}{4}[aabbbab]
-\tfrac{1}{12}[abababb] \, .\notag
\end{align}
An exact expression for the whole of the power series $\sigma_3$ will be given as a closed formula in section 
\ref{sec:sigman.8} below.

\subsubsection{The case $w=5,7,9$}
\label{sec:ldgr}

Now we give the lowest-degree contributions to the expansions of $\sigma_5$, $\sigma_7$ and $\sigma_9$: 
\begin{align}
\sigma_5 &= -\tfrac{1}{24} \epsilon_6^{(4)}
-\tfrac{5  }{48} [\epsilon_4^{(1)},\epsilon_4^{(2)}]
+z_5
+\tfrac{1}{5760} [\epsilon_4,\epsilon_6^{(3)}]
-\tfrac{1}{5760} [\epsilon_4^{(1)},\epsilon_6^{(2)}] +\tfrac{1}{5760} [\epsilon_4^{(2)},\epsilon_6^{(1)}]
\notag \\
&\quad
+\tfrac{1}{3456} [\epsilon_4,[\epsilon_4,\epsilon_4^{(2)}]]
+\tfrac{1}{6912} [\epsilon_4^{(1)},[\epsilon_4^{(1)},\epsilon_4]]
+\tfrac{1}{145152} [\epsilon_6^{(1)},\epsilon_6^{(2)}] -\tfrac{1}{145152} [\epsilon_6,\epsilon_6^{(3)}]
\notag \\
&\quad
-\tfrac{1}{2073600}
[\epsilon_4,[\epsilon_4,\epsilon_6^{(2)}]]
+\tfrac{139 }{72576000} [\epsilon_4^{(1)},[\epsilon_4,\epsilon_6^{(1)}]]
-\tfrac{23 }{24192000} [\epsilon_4,[\epsilon_4^{(1)},\epsilon_6^{(1)}]]
\notag \\
&\quad
-\tfrac{1007 }{145152000} [\epsilon_4^{(2)},[\epsilon_4,\epsilon_6]]
-\tfrac{1}{4147200}
[\epsilon_4^{(1)},[\epsilon_4^{(1)},\epsilon_6]]
+\tfrac{289 }{48384000}
[\epsilon_4,[\epsilon_4^{(2)},\epsilon_6]]
\notag \\
&\quad
+\tfrac{1}{145152000}
[\epsilon_6,\epsilon_8^{(3)}]
-\tfrac{1}{36288000} [\epsilon_6^{(1)},\epsilon_8^{(2)}]
+\tfrac{1}{14515200}
[\epsilon_6^{(2)},\epsilon_8^{(1)}]
-\tfrac{1}{7257600}
[\epsilon_6^{(3)},\epsilon_8]
+ \ldots \label{sig5exp} \\
\sigma_7 &=
-\tfrac{1}{720} \epsilon_8^{(6)}
+\tfrac{7 }{1152} [\epsilon_4^{(2)},\epsilon_6^{(3)}]
-\tfrac{7 }{1152}
[\epsilon_4^{(1)},\epsilon_6^{(4)}]
-\tfrac{661 }{57600}
[\epsilon_4^{(1)},[\epsilon_4^{(1)},\epsilon_4^{(2)}]]
-\tfrac{661 }{57600}
[\epsilon_4^{(2)},[\epsilon_4^{(2)},\epsilon_4]]
\notag \\
&\quad
+\tfrac{1}{172800}
[\epsilon_4,\epsilon_8^{(5)}]
-\tfrac{1}{172800}
[\epsilon_4^{(1)},\epsilon_8^{(4)}]
+\tfrac{1}{172800}
[\epsilon_4^{(2)},\epsilon_8^{(3)}]
+\tfrac{1}{13824}
[\epsilon_6^{(1)},\epsilon_6^{(4)}]
-\tfrac{1}{13824}
[\epsilon_6^{(2)},\epsilon_6^{(3)}]
\notag \\
&\quad
+z_7
-\tfrac{1}{4354560}
[\epsilon_6,\epsilon_8^{(5)}]
+\tfrac{1}{4354560}
[\epsilon_6^{(1)},\epsilon_8^{(4)}]
-\tfrac{1}{4354560}
[\epsilon_6^{(2)},\epsilon_8^{(3)}]
+\tfrac{1}{4354560}
[\epsilon_6^{(3)},\epsilon_8^{(2)}]
\notag \\
&\quad
-\tfrac{1}{4354560}
[\epsilon_6^{(4)},\epsilon_8^{(1)}]
+\tfrac{7 }{552960}
[\epsilon_4,[\epsilon_4,\epsilon_6^{(4)}]]
+\tfrac{7 }{552960}
[\epsilon_4,[\epsilon_4^{(1)},\epsilon_6^{(3)}]]
+\tfrac{7 }{184320} [\epsilon_4^{(1)},[\epsilon_4^{(2)},\epsilon_6^{(1)}]]
\notag \\
&\quad
+\tfrac{7 }{552960}
[\epsilon_4^{(2)},[\epsilon_4,\epsilon_6^{(2)}]]
-\tfrac{7 }{184320}
[\epsilon_4,[\epsilon_4^{(2)},\epsilon_6^{(2)}]]
-\tfrac{7 }{276480}
[\epsilon_4^{(2)},[\epsilon_4^{(2)},\epsilon_6]]
\notag \\
&\quad
-\tfrac{7 }{552960}
[\epsilon_4^{(1)},[\epsilon_4,\epsilon_6^{(3)}]]
-\tfrac{7 }{552960}
[\epsilon_4^{(2)},[\epsilon_4^{(1)},\epsilon_6^{(1)}]]
+ \ldots \label{sig7exp} \\
\sigma_9 &=
-\tfrac{1}{40320} \epsilon_{10}^{(8)}
-\tfrac{1}{5184} [\epsilon_4^{(1)},\epsilon_8^{(6)}] +\tfrac{1}{5184}
[\epsilon_4^{(2)},\epsilon_8^{(5)}]
-\tfrac{7 }{20736} [\epsilon_6^{(3)},\epsilon_6^{(4)}]
+\tfrac{1}{9676800}
[\epsilon_4,\epsilon_{10}^{(7)}]
\notag \\
&\quad
-\tfrac{1}{9676800}
[\epsilon_4^{(1)},\epsilon_{10}^{(6)}]
+\tfrac{1}{9676800}
[\epsilon_4^{(2)},\epsilon_{10}^{(5)}]
+\tfrac{7 }
{4147200}
[\epsilon_6^{(1)},\epsilon_8^{(6)}]
-\tfrac{7 }{4147200}
[\epsilon_6^{(2)},\epsilon_8^{(5)}]
\notag \\
&\quad
+\tfrac{7 }{4147200}
[\epsilon_6^{(3)},\epsilon_8^{(4)}]
-\tfrac{7 }{4147200}
[\epsilon_6^{(4)},\epsilon_8^{(3)}]
-\tfrac{529 }{691200}
[\epsilon_4,[\epsilon_4^{(2)},\epsilon_6^{(4)}]]
+\tfrac{2959 }{2419200}
[\epsilon_4^{(1)},[\epsilon_4^{(2)},\epsilon_6^{(3)}]]
\notag \\
&\quad
+\tfrac{5891 }{6220800}
[\epsilon_4^{(2)},[\epsilon_4,\epsilon_6^{(4)}]]
-\tfrac{443 }{967680}
[\epsilon_4^{(1)},[\epsilon_4^{(1)},\epsilon_6^{(4)}]]
-\tfrac{799 }{1088640}
[\epsilon_4^{(2)},[\epsilon_4^{(2)},\epsilon_6^{(2)}]]
\notag\\
&\quad
-\tfrac{10651 }{21772800}
[\epsilon_4^{(2)},[\epsilon_4^{(1)},\epsilon_6^{(3)}]] + \ldots
 \label{sig9exp}
\end{align}
In all cases, the ellipsis refers to an infinite series in nested brackets of $\epsilon_{k_i}^{(j_i)}$
of total degree $\sum_i k_i \geq 16$, and the expansion of $\sigma_9$ additionally involves
an arithmetic contribution $z_9$ at key degree 18. The action of the arithmetic derivation $z_5$ on the generators $a$ is given by
\begin{align}
    z_5(a)&=-\tfrac{[aaaaababbbb]}{240}-\tfrac{[aaaaabbbbab]}
   {240}+\tfrac{[aaaabaabbbb]}{120}+\tfrac{[aaaabababbb]}{80}-\tfrac{[aaaababbabb]}{30}\nn\\
&\quad +\tfrac{[aaaababbbab]}{60}   +\tfrac{[aaaabbaabbb]}{80}
   -\tfrac{7[aaaabbababb]}{120}-\tfrac{[aaaabbabbab]}{30}
   +\tfrac{[aaaabbbaabb]}{80}\nn\\
&\quad   +\tfrac{[aaaabbbabab]}{240}+\tfrac{[aaaabbbbaab]}{240} -\tfrac{[aaabaababbb]}{24}-\tfrac{3 [aaabaabbabb]}{80}-\tfrac{7[aaabaabbbab]}{240}\nn\\
&\quad -\tfrac{[aaababaabbb]}{240}+\tfrac{73 [aaababababb]}{240}+\tfrac{49[aaabababbab]}{80}+\tfrac{3[aaababbaabb]}{80}+\tfrac{149[aaababbabab]}{240}\nn\\
&\quad +\tfrac{[aaababbbaab]}{240}-\tfrac{[aaabbaababb]}{240}-\tfrac{[aaabbaabbab]}{60}+\tfrac{[aaabbabaabb]}{240}+\tfrac{5
   [aaabbababab]}{16}\nn\\
&\quad-\tfrac{[aaabbabbaab]}{240}
   +\tfrac{[aaabbbaabab]}{240}+\tfrac{[aaabbbabaab]}{120}+\tfrac{[aabaabaabbb]}{240}+\tfrac{[aabaabababb]}{240}\nn\\
&\quad   -\tfrac{[aabaababbab]}{30}+\tfrac{[aabaabbaabb]}{120}-\tfrac{[aabaabbabab]}{30}-\tfrac{3 [aababaababb]}{80}-\tfrac{3[aababaabbab]}{80}\nn\\
&\quad-\tfrac{[aabababaabb]}{240}    +\tfrac{[aababababab]}{16}\,,
\label{z5ona}
\end{align}
again using the Lyndon bracket notation of Theorem \ref{lem:lynd}. A similar expression for $z_5(b)$ can be reconstructed from (\ref{z5ona}) by virtue of the following observation:

\begin{rmk}
The Lie polynomials $z_w(a)$ and $z_w(b)$ at $w=3$, $w=5$ and $w=7$ are related by the switch $\theta$ via
\beq
z_w(b) = - \theta \big( z_w(a) \big) \, , \ \ \ \  w\leq 7 \, .
\label{funswitch}
\eeq
\end{rmk}

Note that an alternative method for the computation of $z_3(a),z_3(b),z_5(a),z_5(b)$
was given by Pollack in \cite{Pollack}, though the approach in that reference has not yet led to explicit results for $z_{w\geq 7}$. Machine-readable expressions for $z_w(a)$ and $z_w(b)$ at $w=3,5,7$ can be found in an ancillary file of the arXiv submission of this work.

\subsubsection{Computational aspects}
\label{sec:compasp}

We close this section by giving more details on the practical implementation of Definition~\ref{thesigmaw}  to determine the canonical zeta generators $\sigma_w$ and their arithmetic parts~$z_w$. 

The starting point of the construction is to solve the conditions~\eqref{g01.11} degree by degree following~\eqref{lowestpart} and~\eqref{recurs} and the partner condition. We recall from Theorem~\ref{thm:522} that at degree $d$ the derivation $(\sigma_w)_d$ has $a$-degree $d{-}w$ and $b$-degree $w$. 

For the example of $\sigma_3$ the extension lemma leads at lowest degree to\footnote{Note that, as a derivation, the lowest degree of $\sigma_3$ is $4$, but here we are writing the degree of the image of $a$ and $b$ as the subscript.}
\begin{align}
\label{eq:s35}
    \big(\sigma_3(a) \big)_{5} & = -[aabbb] +[ababb]\,,\quad
    \big(\sigma_3(b) \big)_{5}  = -[abbbb] \,,
\end{align}
by using $g_3$ presented in~\eqref{explgws} as well as~\eqref{lowestpart}. We here employ Lyndon bracket notation in the Lie algebra ${\rm Lie}[a,b]$. From~\eqref{recurs} we then obtain at the next degree (which is here already key degree):
\begin{align}
\label{eq:s37}
   \big(\sigma_3(a) \big)_{7} &= \tfrac14 [aaababb] + \tfrac14 [aaabbab] +\tfrac1{12} [aababab]\,,\nn\\
   \big(\sigma_3(b) \big)_{7} &=  -\tfrac14 [aababbb] -\tfrac14[aabbabb] -\tfrac14 [aabbbab] -\tfrac1{12} [abababb]\,.
\end{align}

Since~\eqref{eq:s35} is not at key degree, we know from Theorem~\ref{thm:522} that it must be possible to rewrite it completely as the action of a geometric derivation, i.e.\ an element of $\mathfrak{u}$. We know moreover from part (ii) of that theorem that the {\em total depth}, meaning the total number of $\epsilon_i$ (for $i\geq 0$) of any term is equal to $w=3$. Together with the information on the degree, computable from Lemma~\ref{lem:513}, this leaves very few possible terms. 
For any nested backet of  the form 
$\ep_{k_1}^{(j_1)} \cdots \ep_{k_r}^{(j_r)}$ (with $k_i\geq 4$ and any allowed placement of brackets) the conditions to be allowed at degree $d$ in $\sigma_w$ are 
\begin{align}
\label{eq:ansconst}
   r+\sum_{i=1}^r j_i &= w  &&\text{for the total depth and}\nn\\[-1mm]
   \sum_{i=1}^r k_i&= d &&\text{for the degree.}
\end{align}
For example, for the lowest degree $d=4$ in~\eqref{eq:s35}, the only possible term in $\sigma_3$ is proportional to $\ep_4^{(2)}$ and the constant of proportionality $c_1$ is fixed by
\begin{align}
    \big(\sigma_3(a) \big)_5 = \left[\big(c_1 \ep_4^{(2)} \big)(a)\right]_5 = c_1\left( 2[aabbb]-2[ababb]\right)
\end{align}
to the value $c_1=-\tfrac12$ when comparing to~\eqref{eq:s35}, in agreement with~\eqref{firsts3} and a general formula to be derived in Corollary~\ref{cor:lowdeg}.

The next-to-lowest degree in $\sigma_3$, given by~\eqref{eq:s37}, is the key degree $d=2w=6$ and therefore contains both the arithmetic $z_3$ part, transforming in an $\mathfrak{sl}_2$ singlet,  as well as possible geometric contributions. The most general ansatz compatible with~\eqref{eq:ansconst} is
\begin{align}
    \sigkey_3 = z_3 + c_2 \ep_6^{(2)}\,.
\end{align}
In order to separate out the geometric from the arithmetic term, we use that $z_3$ is a singlet under $\mathfrak{sl}_2$ and thus commutes with $\ep_0$. 
The general relations
\begin{align}
\label{eq:gencond}
\epsilon_0\big(\sigma_w(a) \big) -\sigma_w(b) = [\epsilon_0,\sigma_w](a)\,,\quad\quad
\epsilon_0\big(\sigma_w(b) \big)=[\epsilon_0,\sigma_w](b)
\end{align}
at key degree depend only on the geometric part due to $[\ep_0,\sigkey_w] = [\ep_0,\sigkey_w-z_w]$. Moreover, the commutator $[\ep_0,\sigkey_w-z_w]$ of the geometric term can be evaluated easily according to general representation theory as in Lemma~\ref{lem:513}.
The left-hand sides of the general conditions~\eqref{eq:gencond} only depend on $\sigma_w(a)$ and $\sigma_w(b)$ that are furnished by 
\eqref{g01.11} whereas the geometric contribution on the right-hand sides can be computed using the ansatz. 

In the case of~\eqref{eq:s37} we can use the second equation of~\eqref{eq:gencond} and find for the left-hand side
\begin{align}
    \ep_0\big(\sigkey_3(b) \big) = 0
\end{align}
as well as
\begin{align}
    c_2 \ep_6^{(3)}(b) =  12c_2 \big( 2[aabbbb]+5 [ababbbb] + 2[abbabbb]\big)
    \label{ep62ex}
\end{align}
for the right-hand side, implying $c_2=0$ and that the action of $z_3$ is given by~\eqref{eq:s37}, which agrees with the expression already presented in~\eqref{z3act}. 

The ans\"atze for the degree $d$ parts of $\sigma_w$ rapidly grow with $d$ and $w$. For instance, the candidate terms for $(\sigma_7)_{12}$ compatible with (\ref{eq:ansconst}) are given by
\begin{align}
(\sigma_7)_{12} &=
c_{1} \ep_{12}^{(6)}
+c_{2}
[\epsilon_4,\epsilon_8^{(5)}]
+c_{3}
[\epsilon_4^{(1)},\epsilon_8^{(4)}]
+c_4
[\epsilon_4^{(2)},\epsilon_8^{(3)}] \label{si7d12} \\
&\quad+c_{5}
[\epsilon_6^{(1)},\epsilon_6^{(4)}]
+c_{6}
[\epsilon_6^{(2)},\epsilon_6^{(3)}]
+c_{7}
[\epsilon_4^{(1)},[\epsilon_4^{(1)},\epsilon_4^{(2)}]]
+c_{8}
[\epsilon_4^{(2)},[\epsilon_4^{(2)},\epsilon_4]]
\, .\notag 
\end{align}
By matching the action of this ansatz on $a$ with $\big(\sigma_7(a)\big)_{13}$ computed from (\ref{g01.11}), we find the values of the above $c_i$ noted in the degree 12 parts of (\ref{sig7exp}) including a vanishing coefficient $c_1$ of $\ep_{12}^{(6)}$.
The absence of terms in $\sigma_w$ with a single $\ep_k^{(j)}$ at any degree besides the minimal degree $w+1$ will follow from Proposition \ref{cor:isomd} (i) below.

\medskip
In summary, the strategy for converting the result of the extension lemma construction of $\sigma_w$ into expressions in terms of geometric and arithmetic derivations is to make an ansatz for the geometric terms at a given degree subject to the constraints~\eqref{eq:ansconst}.\footnote{It can be useful, although not necessary, to group these terms according to $\mathfrak{sl}_2$ representations.} Away from key degree, evaluating this ansatz on $a$ and $b$ and equating it with the explicit form of $\sigma_w$ then fixes the ansatz (modulo free parameters that are in one-to-one correspondence with the Pollack relations defining $\mathfrak{u}$). At key degree one can separate the geometric from the arithmetic part of $\sigma_w$ using~\eqref{eq:gencond} by first computing the geometric part;  then the arithmetic $z_w$ is simply the difference $z_w=\sigkey_w-\left(\sigkey_w |_{\mathfrak{u}}\right)$.

In section~\ref{secmd2}, we will provide additional calculational tools that recursively determine~$\sigma_w$ up to highest-weight vectors of $\mathfrak{sl}_2$ (see Definition \ref{dfn:sl2}). In case of (\ref{si7d12}), the ansatz contains two highest-weight vectors $[\epsilon_6^{(1)},\epsilon_6^{(4)}]-[\epsilon_6^{(2)},\epsilon_6^{(3)}]$ and $[\epsilon_4^{(1)},[\epsilon_4^{(1)},\epsilon_4^{(2)}]]
+[\epsilon_4^{(2)},[\epsilon_4^{(2)},\epsilon_4]]$, and the method of section~\ref{secmd2} can efficiently determine 6 out of the 8 parameters $c_i$. By Theorem \ref{thm:522} (iv), there are no highest-weight vectors in $\sigma_w$ beyond key degree. Hence, a major virtue of the method in section \ref{secmd2} is that the evaluation of {\em infinitely} many contributions $\sigma_w(a)_{d>2w+1}$ via (\ref{g01.11}) can be bypassed, i.e.\ that the extension lemma construction of section \ref{sec:sigman.2} only needs to be applied to a {\em finite} range of degrees where it fixes {\em all} terms.


\section{\texorpdfstring{Properties of $\tau_w$ and $\sigma_w$}{Properties of tau(w) and sigma(w)}}
\label{sec:prop}

In this section, we prove the properties of the derivation $\tau_w,\sigma_w$ or zeta generators in genus one listed in Theorem~\ref{thm:522} (i) to (vi). One of the key tools for parts (i)-(iii) will be \'Ecalle's theory of moulds developed in \cite{Ecalle} (see also \cite{ZigZag} for an exposition of the basic theory), and the proof of parts (iv)-(vi) will make use of the $\mathfrak{sl}_2$ algebra in Definition \ref{dfn:sl2}.

\subsection{Introduction to moulds}\label{sec:prop1}

For the reader's convenience, we first review a few basic definitions and facts about moulds,
and one fundamental theorem due to \'Ecalle (cf.~\cite{Ecalle},\cite{ZigZag}).

\subsubsection{Moulds and power series}

\begin{dfn}
A {\em rational mould} over a ring $R$ is a family of rational functions
$F=\bigl(F_r\bigr)_{r\ge 0}=(F_0,F_1,F_2,\ldots)$ such that
\beq
F_r(u_1,\ldots,u_r)\in R(u_1,\ldots,u_r)\, ,
\label{mld.01}
\eeq
i.e.\ $F_r$ is a function of $r$ commutative variables $u_i$. The {\em constant term} of the mould
$F_0$ lies in the ring $R$. We will generally refer to a rational mould simply as a ``mould'',
and most of the time we will work over the base field $\Q$.
Also, when there is no possibility of confusion, we often write $F(u_1,\ldots,u_r)$ instead of
$F_r(u_1,\ldots,u_r)$. The function $F_r$ or $F(u_1,\ldots,u_r)$ is called the {\em depth $r$ part
of the mould $F$}.  When the rational functions $F_r$ are polynomials for all $r>0$, we say that $F$
is a {\em polynomial mould}. Moulds can be added componentwise and multiplied by a constant in $R$ componentwise. The moulds with constant term $0$ thus form a vector space, denoted $ARI$; its
vector subspace of polynomial moulds is denoted $ARI^{pol}$. The names of the various objects, morphisms and properties are due to \'Ecalle \cite{Ecalle}.
\end{dfn}

Let $c_i=\ad_x^{i-1}y$ for $i\ge 1$. From now on unless otherwise stated we will work with $R=\Q$. The power series in $\Q\langle\langle x,y\rangle\rangle$
that can be written as power series in the $c_i$ are exactly the ring of power series
$p$ satisfying $\partial_x(p)=0$, where $\partial_x$ is the derivation defined by $\partial_x(x)=1$,
$\partial_x(y)=0$. These power series are in bijection with the free ring $\Q\langle\langle c_1,c_2,\ldots
\rangle \rangle$ of power series on the non-commutative variables $c_i$. All Lie-like and
group-like power series in $\Q\langle \langle x,y\rangle\rangle$ belong to $\Q\langle\langle c_1,c_2,\ldots
\rangle \rangle$ and indeed, with the exception of the element $x$, all Lie polynomials in $x,y$ are in bijection with the Lie polynomials in the $c_i$.  There is a
simple bijection between power series $p\in \Q\langle\langle c_1,c_2,\ldots\rangle\rangle$ and polynomial moulds,
given by letting $p^r$ denote the part of $p$ of homogeneous degree $r$ in the $c_i$ (i.e.\ homogeneous degree $r$ in $y$) and mapping
$p^r$ to the space of polynomial moulds of depth $r$ by the map on monomials
\beq
ma:c_{i_1}\dots c_{i_r}\mapsto (-1)^{r+i_1+\cdots+i_r}u_1^{i_1-1}\cdots u_r^{i_r-1}\, ,
\label{mld.08}
\eeq
extended by linearity. We often use the notation $P=ma(p)$ for the polynomial mould associated to a power
series $p\in\Q\langle\langle c_1,c_2,\ldots\rangle\rangle$ under the map $ma$. The vector space of power series without constant term maps isomorphically under $ma$ to the vector space $ARI^{pol}$.

\subsubsection{Basic operators on moulds}

The space of moulds $ARI$ is equipped with many operations. All those given in the following list are natural extensions to moulds of familiar operations on power series in $x$ and $y$ (see \cite{Ecalle} or \cite{ZigZag} for
complete definitions and details).

\begin{itemize}
\item  Mould multiplication is defined by: 
\beq
mu(G,H)(u_1,\ldots,u_r)=\sum_{i=0}^r G(u_1,\ldots,u_i)H(u_{i+1},\ldots,u_r)\, .
\label{mld.06}
\eeq
This multiplication is valid for moulds with non-zero constant term as well, and is compatible with power series multiplication in the sense that if $G=ma(g)$ and $H=ma(h)$ for $g,h\in \Q\langle\langle c_1,c_2,\ldots,\rangle\rangle$, then 
\beq
ma\bigl(gh\bigr)=mu(G,H)\, .
\eeq

\item
The Lie bracket $lu$ on $ARI$ is defined by 
\beq
\label{eq:lu}
lu(G,H)=mu(G,H)-mu(H,G)\, ,
\eeq
and when $ARI$ is considered as a Lie algebra under this bracket, it is denoted
$ARI_{lu}$. Again, for $G=ma(g)$ and $H=ma(h)$ as above, we have
\beq
ma\bigl([g,h]\bigr)=lu(G,H)\, .
\eeq
\item 
For each mould $G\in ARI$, there is a
derivation $arit(G)$ of the Lie algebra $ARI_{lu}$ which generalizes the
Ihara derivation $D_g$ for $g\in {\rm Lie}[x,y]$ defined by (\ref{eq:IharaD}) in the sense that
if $G=ma(g)$ and $H=ma(h)$ for $g,h\in {\rm Lie}[x,y]$ then
\beq
arit(G)\cdot H=-ma\bigl(D_g(h)\bigr)\, .
\label{mld.17}
\eeq
(The minus sign is due to the original definition of $arit$ by \'Ecalle).

\item 
The {\em ari-bracket} is another Lie bracket on the space $ARI$ (besides $lu$ introduced in~\eqref{eq:lu}), defined by
\beq
ari(G,H)=arit(H)\cdot G-arit(G)\cdot H+lu(G,H)\, .
\label{mld.18}
\eeq
The ari-bracket generalizes the Ihara bracket~\eqref{eq:Ihara} on the underlying vector space ${\rm Lie}[x,y]$ in the sense that if
$G=ma(g)$ and $H=ma(h)$ for $g,h\in {\rm Lie}[x,y]$ then
\beq
ari(G,H)=ma\bigl(\{g,h\}\bigr)\, .
\label{mld.19}
\eeq
We denote the Lie algebra formed by the vector space $ARI$ equipped with the $ari$-bracket by $ARI_{ari}$.

\item 
The universal enveloping algebra $\mathcal{U} ARI_{ari}$ of the Lie algebra $ARI_{ari}$ is nothing other than the space
of all (rational in the context of this article) moulds; these are essentially
the same moulds as in $ARI$ except that arbitrary constant terms are allowed.
By the Poincar\'e--Birkhoff--Witt theorem, this universal enveloping algebra
is equipped with an associative multiplication law which we denote by $\circp$.
The expression for this multiplication
$G\circp H$ simplifies in the case where $G\in ARI$, in which situation it is given for $G$ in $ARI_{ari}$ and $H$ in $\mathcal{U}ARI_{ari}$ by 
\beq
G\circp H=mu(G,H)-arit(G)\cdot H\, ,
\label{mld.20}
\eeq
which thanks to \eqref{mld.17} generalizes the $\circp$ multiplication introduced in \eqref{starmult}:
\beq\label{mld.20a}
G\circp H=ma\bigl(g\circp h\bigr).
\eeq

\item
The {\em ari-exponential} map from $ARI_{ari}$ to the group-like elements
in the universal enveloping algebra is defined for $F\in ARI$ by
\beq
\exp_{ari}(F)=Id+\sum_{n\ge 1} \frac{1}{n!} \bigl(\underbrace{F\circp F\circp\ldots\circp F}_n\bigr)\, ,
\label{mld.21}
\eeq
where the $\circp$ multiplication must be applied from right to left so that
the leftmost element being multiplied is always $F$, and $Id$ denotes the $mu$- and $\circp$-identity mould $(1,0,0,\ldots)$. The image of the space
$ARI$ under the map $\exp_{ari}$ is called $GARI$, and it consists precisely
of the set of all (here rational) moulds with constant term 1. The set $GARI$ forms a group with respect to the multiplication obtained from lifting the $ari$ Lie bracket to $GARI$ using the Baker--Campbell--Hausdorff formula. The
ari-exponential has an inverse map, the {\em ari-logarithm} 
\beq
\log_{ari}:GARI\rightarrow ARI\, .
\label{mld.22}
\eeq

\item 
The group $GARI$ acts
on the Lie algebra $ARI_{ari}$ via the {\em adjoint action}, under which each mould $P\in GARI$ gives an isomorphism
of the Lie algebra $ARI_{ari}$ via the {\em adjoint operator} $\Ad_{ari}(P)$. 
Let $L\coloneqq\log_{ari}(P)$, so $L\in ARI$. Then the adjoint action of $P$ on a mould
$A\in ARI$ can be expressed and computed explicitly by the standard formula
\begin{align}
\Ad_{ari}(P)(A)
&=A+ari(L,A)+\tfrac{1}{2}ari(L,ari(L,A))+\tfrac{1}{6}ari(L,ari(L,ari(L,A)))+
\cdots 
\label{adjoint}
\end{align}
by exponentiating the $ari$ bracket $ari(L,\cdot)$.
\item We define an operator $dur$ acting on all moulds by $dur(F)(\emptyset)=F(\emptyset)$ and the following
formula for $r\ge 1$:
\beq
dur(F)(u_1,\ldots,u_r)=(u_1+\cdots+u_r)F(u_1,\ldots,u_r)\, .
\label{mld.02}
\eeq
If $F=ma(f)$ for a power series $f\in \Q\langle\langle
c_1,c_2,\ldots\rangle \rangle$ (considered as a function $f(x,y)$), then
\beq
dur(F)=ma([x,f]) \, .
\eeq
\item 
We will also need the mould operator $\Delta$ defined by $\Delta(F)(\emptyset)=F(\emptyset)$ and 
\beq
\Delta(F)(u_1,\ldots,u_r)=u_1\cdots u_r(u_1+\cdots+u_r) F(u_1,\ldots,u_r)\, .
\label{mld.03}
\eeq
If $F=ma(f)$ as above, we have
\beq\label{Deltaonpoly}\Delta(F)=ma\bigl([x,f(x,[x,y])]\bigr)\, .
\eeq
The inverse operator of $\Delta$ is given by
\beq
\Delta^{-1}(F)(u_1,\ldots,u_r) = \frac{1}{u_1\cdots u_r(u_1+\cdots+u_r)} F(u_1,\ldots,u_r)\, .
\label{mld.04}
\eeq
Of course, the operator $\Delta$ on power series given in \eqref{Deltaonpoly} cannot always be inverted in the world of non-commutative power series.
\item
The {\it push-operator} acts on moulds $F$ by the formula
$push(F)(\emptyset)=F(\emptyset)$ and for $r\ge 1$,
\beq\label{mld.05}
push(F)(u_1,\ldots,u_r)=F(-u_1-\cdots-u_r,u_1,u_2,\ldots,u_{r-1})\, .
\eeq
The push-operator corresponds to an operation on power series (also called push) monomial by monomial defined as follows:
\beq\label{pushpoly}
push(x^{a_1}yx^{a_2}y\cdots yx^{a_r-1}yx^{a_r})=x^{a_r}yx^{a_1}y\cdots yx^{a_{r-2}}yx^{a_{r-1}}
\eeq
in the sense that if $h\in \Q\langle\langle c_1,c_2,\ldots\rangle\rangle$ then
\beq\label{pushequal}
ma\bigl(push(h)\bigr)=push\bigl(ma(h)\bigr)\, ,
\eeq
where the left-hand push is as in \eqref{pushpoly} and the right-hand one is as in \eqref{mld.05} (for this equivalence, see \cite{SchnepsRaphael}, section 3.3). In particular, $h$ is push-invariant if and only if $ma(h)$ is.
\item 
The {\em swap} operator on moulds is defined by the formula
$swap(F)(\emptyset)=F(\emptyset)$ and
\beq
swap(F)(v_1,v_2,\ldots,v_r)=F(v_r,v_{r-1}-v_r,\ldots,v_1-v_2)\, .
\label{mld.06bis}
\eeq
We could write the mould $swap(F)$ in the variables $u_i$ instead of $v_i$, of course, but to keep apart a mould
and its swap it is convenient to consider the swapped mould parts $swap(F)_r$ as lying in
$\Q(v_1,\ldots,v_r)$.
\item 
Finally, we need to define the {\em alternality} property on moulds. A mould $P\in ARI$ is said to be {\em alternal} if for all $r\ge 2$ we have
\beq
\sum_{w\in u\shuffle v} P(w)=0
\label{mld.10}
\eeq
for all pairs of non-empty words $u=(u_1,\ldots,u_i)$,
$v=(u_{i+1},\ldots,u_r)$. (There is no condition at $r=1$.)
When $P=ma(p)$ for a power series $p\in \Q\langle\langle c_1,c_2,\ldots\rangle\rangle$ without constant term, then $P$ is alternal if and only if $p$ is a Lie element in the $c_i$, or equivalently, if and only if $p(x,y)\in {\rm Lie}[x,y]$.
\end{itemize}

\vspace{.3cm}
\noindent {\bf Example.}
Recall that the first non-trivial element of $\mz^\vee$ is given by
\beq
g_3=[x,[x,y]]+[[x,y],y]=c_3+[c_2,c_1]=c_3+c_2c_1-c_1c_2\, .
\eeq
By \eqref{mld.08}, the associated mould $G_3=ma(g_3)\in ARI$ is given by 
\begin{align}
0&\mapsto 0=G_3(\emptyset)\ &\hbox{in depth 0} \, , \notag\\
c_3&\mapsto u_1^2=G_3(u_1)\ 
&\hbox{in depth 1}\, ,\\
c_2c_1-c_1c_2&\mapsto
-u_1+u_2=G_3(u_1,u_2)\ 
&\hbox{in depth 2} \, . \notag
\end{align}
The fact that $g_3$ is a Lie polynomial is reflected in the alternality condition satisfied by $G_3$:
\beq
\sum_{w\in (u_1\shuffle u_2)} G_3(w)=G_3(u_1,u_2)+G_3(u_2,u_1)=0\, .
\eeq

\subsubsection{The fundamental operator $\Ad_{ari}(pal)$ and \'Ecalle's theorem}

\'Ecalle defined a remarkable pair of inverse moulds in the group $GARI$, called $pal$ and $invpal$, which have the following property: when acting on $ARI$ via the adjoint action, $invpal$ transforms the double shuffle property into a much simpler property known as {\em bialternality}, where a bialternal mould is an alternal mould with alternal swap, and $pal$ does the opposite (this is a major result due to \'Ecalle, see \cite{Ecalle,EcalleVietnam} and an expository version in section 4.6 of \cite{ZigZag}). The isomorphisms $\Ad_{ari}(invpal)$ and $\Ad_{ari}(pal)^{-1}$ are mutually inverse. The action of $\Ad_{ari}(invpal)$ on a double shuffle Lie polynomial mould introduces certain denominators, but these are eliminated by the operator $\Delta$ in (\ref{mld.03}), yielding a polynomial mould once again (cf.~\cite{BaumardSchneps:2015}); in other words, restricted to $ma(\ds)$, the composition 
$\Delta\circ \Ad_{ari}(invpal)$ takes polynomial moulds to polynomial moulds.
The key result for our purposes here is that when restricted to the subspace
$ma(\mz^\vee)\subset ma(\ds)$, the
map $\Delta\circ \Ad_{ari}(invpal)$
is directly related to the morphism
\beq
\gamma:\mz^\vee\rightarrow
{\rm Der}^0{\rm Lie}[a,b]
\eeq
of \eqref{gammamap} by the following formula: if
$h\in \mz^\vee$, then 
\beq\label{mainconnection}\Delta\circ \Ad_{ari}(invpal)\bigl(ma(h)\bigr)=ma\bigl(\gamma(h)(a)\bigr)\, ,
\eeq
where 
\beq
\gamma(h)\in {\rm Der}^0{\rm Lie}[a,b]
\eeq
and $\gamma(h)(a)$ denotes the Lie series obtained by applying that derivation to $a$ (cf.~\cite{Schneps:2015mzv}, Thm. 1.3.1).  The connection \eqref{mainconnection} 
enables us to apply the known properties of the operator $\Ad_{ari}(invpal)$ to prove properties of the derivations $\tau_w$ and $\sigma_w$ in view of their relation to $\gamma$ in (\ref{deftausig}).

We now proceed to the definition of the moulds $pal$ and $invpal$. 
\begin{dfn}
Let $dupal$ be the mould defined explicitly by
$dupal(\emptyset)=0$ and for $r>0$ by 
\beq
dupal(u_1,\ldots,u_r)=\frac{{\rm B}_r}{r!} \frac{1}{u_1\cdots u_r}
\Biggl(\sum_{j=0}^{r-1}(-1)^j \binom{r-1}{j}
u_{j+1}\Biggr)\, .
\label{mld.11}
\eeq
\end{dfn}

\begin{lemma} 
The mould $dupal$ is related to $t_{01}$ in (\ref{g01.01}) by the equation
\beq\label{dupalt01} 
dupal(u_1,\ldots,u_r)=\frac{1}{u_1\cdots u_r}ma(t_{01}^r)
\eeq
for all $r\ge 1$, where $t_{01}^r$ is the part of $t_{01}$ of $b$-degree $r$.
\end{lemma}

\noindent {\bf Proof.} The map $ma$ maps power series in $a,b$ to moulds exactly like those in $x,y$, namely via \eqref{mld.08} with $c_i=\ad_a^{i-1}(b)$.
 To prove \eqref{dupalt01}, notice that since we have 
\beq
\ad_b^{r-1}(a)=-\ad_b^{r-2}([a,b])=-\ad_{c_1}^{r-1}(c_2)=-\sum_{j=0}^{r-1}(-1)^j \binom{r-1}{j}
c_1^jc_2c_1^{r-1-j}\, ,
\eeq
the associated mould is
\beq\label{mabra}
ma\bigl(\ad_b^{r-1}(a)\bigr)=-\sum_{j=0}^{r-1}(-1)^j \binom{r-1}{j}
u_{j+1}\, .
\eeq
Hence, since the part $t_{01}^r$ of $b$-degree $r$ of $t_{01}$ is just given by $-\tfrac{{\rm B}_r}{r!}\ad_b^r(a)$, \eqref{dupalt01} follows from comparing \eqref{mld.11} and \eqref{mabra}.\qed

\begin{dfn}
Let $pal$ be the mould defined recursively
by $pal(\emptyset)=1$ and the formula
\beq
dur(pal)=mu(pal,dupal)
\label{mld.12}
\eeq
with $dur$ defined in (\ref{mld.02}) and $dupal$ in~\eqref{mld.11}.
\end{dfn}

This formula might look circular but in fact it defines each depth of
$pal$ successively thanks to the fact that $dupal(\emptyset)=0$. For example,
in depth 1, we have
\begin{align}
dur(pal)(u_1)&=u_1pal(u_1)\nn\\
&=mu(pal,dupal)(u_1)\nn\\
&=pal(\emptyset)dupal(u_1)+pal(u_1)dupal(\emptyset)\nn\\
&=dupal(u_1)\nn\\
&=-\frac{1}{2} \, ,
\label{mld.13}
\end{align}
so
\beq
pal(u_1)=-\frac{1}{2u_1}\, .
\label{mld.14}
\eeq
Then in depth 2, we have
\begin{align}
dur(pal)(u_1,u_2)&=(u_1+u_2)pal(u_1,u_2)\nn\\
&=mu(pal,dupal)(u_1,u_2)\nn\\
&=pal(\emptyset)dupal(u_1,u_2)+pal(u_1)dupal(u_2)\nn\\
&=\frac{u_1-u_2}{12u_1u_2} +\frac{1}{4u_1}\nn\\
&= \frac{u_1+2u_2}{12u_1u_2} \, ,
\label{mld.15}
\end{align}
so
\beq
pal(u_1,u_2)=\frac{u_1+2u_2}{12u_1u_2(u_1+u_2)} \, .
\label{mld.16}
\eeq

\begin{dfn}
Let $lopal=\log_{ari}(pal)$ using the ari-logarithm map defined in \eqref{mld.22}, and recall that $invpal$ is the inverse of $pal$ in the group $GARI=\exp_{ari}(ARI)$, equipped with the Baker--Campbell--Hausdorff multiplication law, so that we have
\beq
\log_{ari}(pal)=-\log_{ari}(invpal)\,.
\eeq
\end{dfn}

In lowest depths we have
\begin{align}\label{invpallopal}
    lopal &=\bigg(0,\, -\frac1{2u_1}, \, \frac{u_1-u_2}{12u_1u_2(u_1+u_2)},\, \ldots\biggr)\,,\nn\\
    invpal&=\biggl(1,\, \frac1{2u_1},\,  \frac{-u_1+4u_2}{12u_1 u_2(u_1+u_2)},\, \ldots\biggr) \, .
\end{align} 
Both of these moulds will be used below in our computations of $\sigma_w$. 

\vspace{.2cm}
The following theorem summarizes the key results from mould theory needed for the proof of Theorem \ref{thm:522} (i) to (iii). 

\begin{thm}
\label{thm:613}
Let $h\in \ds$ and let $H=ma(h)$ denote the associated mould. Let $\tau_h$ be the derivation of $\widehat{\rm Lie}[a,b]$
constructed from $h$ as in section \ref{sec:sigman.2} and write
$T_h=ma\bigl(\tau_h(a)\bigr)$. Then, 

\begin{enumerate}
\item[(i)] The mould $\Ad_{ari}(invpal)(H)$ is bialternal, i.e.~it is alternal and its swap is alternal (cf.~\cite{Ecalle,EcalleVietnam} and \cite{ZigZag}, Thm. 4.6.1);

\item[(ii)]
We have the following equality of moulds in $ARI$ (cf.~\cite{Schneps:2015mzv}, Thm.\ 1.3.1):
\beq\label{keyequality}
T_h=\Delta\circ\Ad_{ari}(invpal)(H) \, ; 
\eeq

\item[(iii)] All bialternal moulds are push-invariant (cf.~\cite{Ecalle}, \cite{ZigZag} Lemma 2.5.5); in particular $\Ad_{ari}(invpal)(H)$ is push-invariant, and so is $T_h$ since $\Delta$ does not modify push-invariance;

\item[(iv)] A bialternal rational mould $A$ satisfies
\beq
A(-u_1,\ldots,-u_r)=A(u_1,\ldots,u_r)
\eeq
for all $r\ge 1$. In particular if $A(u_1,\ldots,u_r)$ is of odd total degree then it is equal to zero (cf.~\cite{ZigZag}, Lemma 2.5.5).
\end{enumerate}
\end{thm}
Note that the push invariance of $T_h$ and therefore $\tau_h(a)$ established in part (iii) is crucial to obtain extensions of derivations of the Lie subalgebra ${\rm Lie}[t_{12},t_{01}]\subset \widehat{\rm Lie}[a,b]$ to all of ${\rm Der}^0\widehat{\rm Lie}[a,b]$, see the discussion around (\ref{recurs}).

\subsection{Proof of Theorem~\ref{thm:522} (i)-(iii)}\label{sec:prop2}

For all $h\in \mz^\vee$, let $\tau_h$ denote the associated derivation in ${\rm Der}^0{\rm Lie}[a,b]$ constructed in section~\ref{sec:sigman.2}. Let $g_w$ for odd $w\ge 3$ be the canonical free generators of $\mz^\vee$; recall that we write $\tau_w$ and $\sigma_w$ for the zeta generators in genus one rather than $\tau_{g_w}$ and $\sigma_{g_w}$. The results of Theorem \ref{thm:613} are valid for all elements $h\in\ds$, in particular for elements of the subspace $\mz^\vee\subset \ds$, but in this section we will apply them specifically to the elements $g_w$.

\begin{cor}[Theorem~\ref{thm:522} (i)]
\label{cor:621}
The derivations $\tau_w$ and $\sigma_w$ satisfy
\beq\label{mld.32}
\tau_w([a,b])=
\sigma_w([a,b])=0\, ,
\eeq
i.e.\ $\tau_w$ and $\sigma_w$ lie in ${\rm Der}^0 {\rm Lie}[a,b]$.
\end{cor}

\noindent {\bf Proof.} The mould $T_w=ma\bigl(\tau_w(a)\bigr)$ is push-invariant by Theorem \ref{thm:613} (iii), and we saw in \eqref{pushequal} that push-invariance for moulds is equivalent to push-invariance of power series. Thus $\tau_w(a)$ is push-invariant. It is shown in Lemma 2.1.1 of \cite{Schneps:2015mzv} that for any derivation $\delta$ of ${\rm Lie}[a,b]$ such that $\delta(b)$ is the partner of $\delta(a)$ as defined in \eqref{g01.03}, then $\delta([a,b])=0$ if and only if $\delta(a)$ is push-invariant. Since $\tau_w(b)$ is the partner of $\tau_w(a)$ by construction (i.e.~\eqref{g01.03} with $g=\tau_w(a)$ and $g' = \tau_w(b)$) and $\tau_w(a)$ is push-invariant, we thus have $\tau_w([a,b])=0$ as desired. Then 
\beq
\sigma_w([a,b])=\theta\circ\tau_w\circ\theta([a,b])=\theta\circ\tau_w([b,a])=0
\eeq
as well.\qed

\begin{prop} The mould $T_w$ is zero in all even depths, and in odd depths $r\ge 1$, $T_w(u_1,\ldots,u_r)$ is a  polynomial of homogeneous degree $w+1$ in the variables $u_i$. In particular
\beq
T_w(u_1)=u_1^{w+1}\, .
\eeq
\end{prop}

\vspace{.2cm}
\noindent {\bf Proof.} We first show that the mould $T_w$ is of constant degree $w+1$ in $u_1,\ldots,u_r$ in every depth. For this, we begin by noting that the Lie series 
\beq
\tau_w(t_{01})=[t_{01},g_w(t_{12},-t_{01})]
\eeq
has constant $a$-degree equal to $w+1$ since $g_w$ is a polynomial of homogeneous degree $w$ and both $t_{01}$ and $t_{12}$ have $a$-degree $1$. Then, using the degree-by-degree computation of $\tau_w(a)$ given in \eqref{tauhabearlier} to \eqref{recurs} (with $h=g_w$), we see that $\tau_w(a)_n$ is a Lie polynomial of constant $a$-degree $w+1$ in every degree $n$ 
since the $a$-degree of the partner $\tau_w(b)$ is one less than that of $\tau_w(a)$ at every degree. 
By the defining property $g_w(x,y)|_{x^{w-1}y}=1$ of the canonical polynomials in genus zero and their symmetry property $g_w(x,y) = g_w(y,x)$,\footnote{This symmetry property follows from the odd degree $w$ of $g_w$ together with the facts that $g_w(x,-y)\in \grt$ by \eqref{minusy} and that one of the defining properties of elements $h\in \grt$ is $h(x,y)+h(y,x)=0$.\label{gwsymm}} the monomial $y^{w-1}x$ also appears in $g_w(x,y)$ with coefficient~1. Since $g_w(x,y)$ for odd $w$ is a Lie polynomial this implies that the Lie word $\ad_y^{w-1}(x)$ appears in $g_w$ with coefficient 1. Thus the minimal $x$-degree in $g_w$ is 1 and by \eqref{lowestpart} we have 
\beq\label{lowestpart2}
\tau_w(a)_{w+2}=[a,\ad_a^{w-1}([a,b])]=\ad_a^{w+1}(b)\, ,
\eeq
where the sign in $t_{01} = -a + \ldots$ disappears since $w$ is odd.

Under the map $ma$ from power series to commutative variables $u_1,\ldots,u_r$ defined in \eqref{mld.08} (with $c_i=\ad_a^{i-1}b$ for $i\ge 1$), we see that the $a$-degree corresponds to the degree in $u_1,\ldots,u_r$ while the $b$-degree corresponds to the mould depth $r$; thus for all $r\ge 1$, the depth $r$ part of the mould $T_w=ma\bigl(\tau_w(a)\bigr)$ is a polynomial in $u_1,\ldots,u_r$ of degree $w+1$. Furthermore, the lowest depth part of $T_w$ appears in depth 1 and is given by
\beq
T_w(u_1)=ma\bigl(\ad_a^{w+1}(b)\bigr)=u_1^{w+1}\, .
\eeq
It remains only to prove that $T_w(u_1,\ldots,u_r)=0$ for all even $r$. For this, we apply Theorem \ref{thm:613} to the case $h=g_w$ and $H=G_w=ma(g_w)$. By (ii) of that theorem, we have
\beq
T_w=\Delta\circ \Ad_{ari}(invpal)(G_w)\, .
\eeq
Therefore for each $r\ge 1$ we have
\beq
\Delta^{-1}(T_w)(u_1,\ldots,u_r)=\frac{T_w(u_1,\ldots,u_r)}{u_1\cdots u_r(u_1+\cdots+u_r)}=\Ad_{ari}(invpal)(G_w)(u_1,\ldots,u_r)\, .
\eeq
By (i) of Theorem \ref{thm:613}, the mould $\Ad_{ari}(invpal)(G_w)$ is bialternal, so the rational mould in the middle term is bialternal. The total degree of this rational function is $w-r$, which is odd whenever $r$ is even. Thus, by Theorem \ref{thm:613} (iv), the mould $T_w$ is zero in all even depths $r$.
This
concludes the proof of the Proposition.
\qed

\begin{cor} [Theorem~\ref{thm:522} (ii)] 
\label{cor:lowdeg}

\noindent
\begin{itemize}
\item[(i)] The minimal degree part of the Lie
series $\tau_w(a)$ is equal to $\ad_a^{w+1}(b)$, so the minimal degree part of $\tau_w$ is $\epsilon_{w+1}$. The minimal degree part of $\sigma_w$ is given by $-\tfrac{1}{(w-1)!}\epsilon_{w+1}^{(w-1)}$. 
\item[(ii)] There are no terms of degree $<w+2$ and no terms of even degree in the Lie series $\tau_w(a)$, $\sigma_w(a)$ and their partners. 
For all odd $n\ge w+2$, the degree-$n$ terms of $\tau_w(a)$ (resp.~$\sigma_w(b)$) all have $b$-degree (resp.~$a$-degree) equal to $n-w-1$ and constant $a$-degree (resp.~constant $b$-degree) equal to $w+1$.
\end{itemize}
\end{cor} 

\vspace{.1cm}
\noindent {\bf Proof.} (i) We saw in \eqref{lowestpart2} that the lowest degree of $\tau_w(a)$ is $w+2$ and $(\tau_w(a))_{w+2}=\ad_a^{w+1}(b)$, which is also equal to
$\epsilon_w(a)$ by (\ref{epsaction}). The switch formula is given in \eqref{lemepjk.b}.

(ii) The statement is a direct translation of the corresponding statement of the previous proposition into terms of the non-commutative variables $a,b$. The minimal degree of $\tau_w$ and $\sigma_w$ as a derivations is $w+1$ by part (i), so the minimal degree of the Lie series $\tau_w(a)$ and $\sigma_w(a)$ is $w+2$. For the other terms, the map $ma$ sends a polynomial $h\in \Q\langle c_1,c_2,\ldots\rangle$ (with $c_i=\ad_a^{i-1}(b)$) of homogeneous  degree $n$ in $a,b$ and homogeneous depth $r$ to a mould $ma(h)$ concentrated in depth $r$ of homogeneous degree $n-r$ in the variables $u_1,\ldots,u_r$. Since the degree of $T_w(u_1,\ldots,u_r)$ is always $w+1$ by the previous Proposition, the $a$-degree of every term of $\tau_w(a)$ is $w+1$. The depth $r$ part of the mould $T_w$ corresponds to the $b$-degree $r$ part of the power series $\tau_w(a)$. We first observe that if $r$ is even then $T_w(u_1,\ldots,u_r)=0$ by the previous proposition, so all terms of $\tau_w(a)$ of even $b$-degree $r$ are zero, but these are precisely all the terms of total degree $w+1+r$, which are all of the even-degree terms. If we have a term $\tau_w(a)$ of odd total degree $n$, then since it has $a$-degree $w+1$ its  $b$-degree is equal to $n-w-1$. This concludes the proof for $\tau_w(a)$ and the switch gives the analogous result for $\sigma_w(b)$ with $b$-degree $w+1$ and $a$-degree $n-w-1$. \qed

\begin{prop}
For each odd $w\ge 3$, the mould
$T_w=ma\bigl(\tau_w(a)\bigr)$ is entirely determined by its parts of depth $r\le w-1$.
\end{prop}

\vspace{.2cm}
\noindent {\bf Proof.} 
By Theorem~\ref{thm:613} (ii), the mould $\Delta^{-1}T_w$ is equal to
$\Ad_{ari}(invpal)(G_w)$ where $G_w=ma(g_w)$ and $g_w$ is the canonical polynomial in genus zero.  For any moulds $P\in GARI$ and $A\in ARI$, set $L=\log_{ari}(P)$ and recall the adjoint operator formula \eqref{adjoint}. Since $L$ has no constant term, taking the $ari$-bracket with $L$ increases the depth, so the adjoint operator formula shows that for any given depth $r$, only the terms of $A$ of depth $\le r$ contribute to the depth $r$ part of $\Ad_{ari}(P)(A)$.  
Now let $A=\Ad_{ari}(invpal)(G_w)$ and $P=pal$, so that 
\beq\label{getGw}
\Ad_{ari}(P)(A)=\Ad_{ari}(pal)\bigl(\Ad_{ari}(invpal)(G_w)\bigr)=G_w\, .\eeq 
Since $g_w$ is a Lie polynomial of degree $w$ it has no terms of depth $\ge w$, so the same is true for the associated mould $G_w=ma(g_w)$. Thus, $G_w$ is determined entirely by its parts of depth $\le w-1$, which in turn by the adjoint action formula are determined entirely by the parts of $A=\Ad_{ari}(invpal)(G_w)$ in depths $\le w-1$. The parts of $T_w$ of depth $\le w-1$ determine those of $A=\Ad_{ari}(invpal)(G_w)$ by applying $\Delta^{-1}$, and the parts of $A$ of depths $\le w-1$ then determine $G_w$ up to depth $w-1$ by the adjoint action formula \eqref{getGw} -- but this is all of $G_w$, which then in turn determines all of $T_w$ by the formula
\beq
T_w=\Delta\circ \Ad_{ari}(invpal)(G_w)\, ,
\eeq 
concluding the proof of the proposition.\qed

\begin{cor}[Theorem~\ref{thm:522} (iii)] 
Both of the derivations $\tau_w$ and $\sigma_w$ are entirely determined by their parts of degree $\le 2w-1$ (as derivations).
\end{cor}

\vspace{.2cm}
\noindent {\bf Proof.} By the above Proposition, $T_w$ is entirely determined
by its parts of depth $\le w-1$, so the same holds for the Lie series
$\tau_w(a)$. But
we saw above that for all $r\ge 1$ the $b$-degree $r$ part of the Lie series $\tau_w(a)$ is of polynomial
degree $w+r+1$ in $a$ and $b$, so in particular the $b$-degree $w-1$ part of $\tau_w(a)$
is of degree $2w$. Saying that $\tau_w(a)$ is determined by its parts of $b$-degree $\le w-1$ is equivalent to saying that it is determined by its parts of total degree $\le 2w$. Since $\tau_w([a,b])=0$ by Corollary \ref{cor:621}, knowing $\tau_w(a)$ determines $\tau_w$ completely. The part of $\tau_w(a)$ of given polynomial degree $n$ corresponds to the part of $\tau_w$ of degree $n-1$ as a derivation; thus the derivation $\tau_w$ is entirely determined by its parts of degree $\le 2w-1$, and the same holds for $\sigma_w$.\qed

\subsection{Proof of Theorem~\ref{thm:522} (iv)-(vi)}\label{sec:prop3}

In this section, we use properties of the $\mathfrak{sl}_2$ algebra in Definition \ref{dfn:sl2} with generators $\ep_0,\ep_0^\vee,\hhh$ to prove parts (iv)-(vi) of Theorem~\ref{thm:522}.

Since the element $\hhh=[\ep_0,\ep_0^\vee]\in\mathfrak{sl}_2\subset {\rm Der}^0{\rm Lie}[a,b]$ acts by $\hhh(a)=-a$ and $\hhh(b)=b$, any derivation $\delta$ of $\widehat{\rm Lie}[a,b]$ of homogeneous $a$-degree $\alpha$ and $b$-degree $\beta$ is an eigenvector for
$\hhh$, with eigenvalue given by
\beq\label{heigenvalue}
[\hhh,\delta]=(\beta-\alpha)\delta \, .
\eeq
In particular, for the action of $\hhh$ on $\mathfrak{u}$, we have $[ \hhh , \ep_k^{(j)}] = (2j{+}2{-}k) \ep_k^{(j)}$ from \eqref{hbrak.02}, so  $\hhh$ has eigenvalues covering the spectrum of values ${-}k{+}2,\, {-}k{+}4,\ldots,-2,0,2,\ldots,k{-}4,\, k{-}2$ within the $(k-1)$-dimensional irreducible representations $\{ \ep_k^{(j)}, \ j=0,1,\ldots,k-2\}$ of $\mathfrak{sl}_2$ at fixed $k$. Similarly, $(r-1)$-dimensional irreducible subrepresentations in $\mathfrak{u}$ built from brackets of $ \ep_{k_1}^{(j_1)}
\ep_{k_2}^{(j_2)}\ldots \ep_{k_m}^{(j_m)}$
will have the spectrum of $\hhh$-eigenvalues ${-}r{+}2,\, {-}r{+}4,\ldots,-2,0,2,$ $\ldots,r{-}4,\, r{-}2,$ always including the eigenvalue zero since $r$ is even as will become clear from the discussion around (\ref{sl2tens}). 

By $[\hhh,\ep_0]= 2 \ep_0$ and $[\hhh,\ep_0^\vee]= - 2 \ep_0^\vee$, adjoint action of $\ep_0$ and $\ep_0^\vee$ shifts the $\hhh$ eigenvalue of any derivation $\delta \in {\rm Der}^0\widehat{\rm Lie}[a,b]$ (not necessarily $\delta \in\mathfrak{u}$) by $2$ and $-2$, respectively (except for highest- and lowest-weight vectors annihilated by $\ad_{\ep_0}$ and $\ad_{\ep^\vee_0}$, respectively). 

\begin{lemma}
\label{cor:rec1}
By the above spectra of $\hhh$ eigenvalues in irreducible representations of $\mathfrak{sl}_2$ and the action (\ref{epjkep}) as well as the fact that $\ad_{\ep_0}\ep_k^{(j)}=\ep_k^{(j+1)}$ and $\ep_k^{(k-1)}=0$,
we have:
\begin{itemize}
\item[(i)] for any $Y \in \ad_{\ep_0}\mathfrak{u}$, the equation $\ad_{\ep_0}X=Y$ has a unique solution $X\in\ad_{\ep_0^\vee}\mathfrak{u}$. In particular, ${\rm ad}_{\ep_0}$ has no kernel within eigenspaces at negative eigenvalues of $\hhh$.
\item[(ii)] for any $Y \in \ad_{\ep_0^\vee}\mathfrak{u}$, the equation $\ad_{\ep_0^\vee}X=Y$ has a unique solution $X\in \ad_{\ep_0}\mathfrak{u}$. In particular, $\ad_{\ep_0^\vee}$ has no kernel at positive eigenvalues of $\hhh$.
\end{itemize}
\end{lemma}

\subsubsection{Proof of Theorem \ref{thm:522} (iv)}
\label{sec:prfv}

For any term of  $\sigma_w$ of total degree $n$, since by Theorem \ref{thm:522} (ii) the $b$-degree is $w$, the $a$-degree must be $n-w$, and thus by \eqref{heigenvalue} this term is an $\hhh$-eigenvector with $\hhh$-eigenvalue equal to $2w - n$. Thus any term of $\sigma_w$ of bihomogeneous degree in $a$ and $b$ and total degree $n$ is an eigenvector for $\hhh$, and we have:  
\begin{align}
&\hbox{if $n<2w$, the eigenvalue of $\hhh$ is strictly positive,} \notag \\
&\hbox{if $n=2w$, the eigenvalue of $\hhh$ is zero,}\\
&\hbox{if $n>2w$, the eigenvalue of $\hhh$ is negative}. \notag
\end{align}

\begin{lemma}[Theorem \ref{thm:522} (iv)] 
\label{deghev} 
The derivation $\sigma_w$ has no highest-weight vectors in degrees $n>2w$.
\end{lemma}

\noindent {\bf Proof.} Since ${\rm ad}_{\ep_0}$  has no kernel at negative $\hhh$-eigenvalues by Lemma \ref{cor:rec1} (i), the infinite Lie series of geometric contributions to $\sigma_w$ above key degree $2w$ does not involve any highest-weight vectors. \qed

\subsubsection{Proof of Theorem \ref{thm:522} (v) and (vi)}
\label{prfrest}

We shall next prove parts (v) and (vi) of Theorem \ref{thm:522} based on
Theorem \ref{thm:HM}. In a notation where 
\begin{align}
\mathfrak{g}&\coloneqq \mathfrak{u} \rtimes \mathfrak{sl}_2  \, ,
\label{semidirect1} 
\end{align}
and $\sigalg$ denotes the free Lie algebra of zeta generators $\sigma_w$, Theorem \ref{thm:HM} implies that
\beq
[\mathfrak{g},\sigalg]\subset\mathfrak{g} \, .
\label{preHMmagic}
\eeq
Following the notation $p_d$ for degree-$d$ parts of polynomials $p$ in $a,b$, we shall write $(\sigma_w)_d$ for the degree-$d$ part of genus one zeta generators, so that in particular $\sigma_w^{\rm key} = (\sigma_w)_{2w}$.

\begin{prop}[Theorem \ref{thm:522} (v) and (vi)]
\label{deghev.bis}

\noindent
\begin{itemize}
\item[(i)] All terms of $\sigma_w$ in degrees $\neq 2w$ lie in $\mathfrak{u}$, but $\sigma_w^{\rm key}\notin\mathfrak{u}$.
\item[(ii)] The terms of $\sigma_w$ in key degree $2w$ that lie in irreducible $\mathfrak{sl}_2$ representations of dimension $\geq 3$ lie in $\mathfrak{u}$. 
\item[(iii)] The brackets $[z_w,\ep_k]$ of the $\mathfrak{sl}_2$-invariant part $z_w$ of $\sigma_w$ lie in $\mathfrak{u}$.
\end{itemize}
\end{prop}

\vspace{.2cm}
\noindent {\bf Proof.} (i) Recall from Theorem \ref{thm:522} (ii) that every term of $\sigma_w$ is of $b$-degree $w$ and that the minimum total degree of any term is given by $n=w+1$. Let 
\beq
\sigma_w = \sum_{n=w+1}^\infty 
(\sigma_w)_n
\eeq
denote the expansion of $\sigma_w$ according to total degree. Then
by \eqref{heigenvalue}, for each $n\ge w+1$, we have
\beq
[\hhh , (\sigma_w)_n] = (2w-n) (\sigma_w)_n \, .
\label{exp:sideg}
\eeq
Note that,
instead of (\ref{preHMmagic}),
we actually have the stronger statement
\beq
[\mathfrak{g},\sigalg]\subset\mathfrak{u}
\label{HMmagic}
\eeq
since the brackets on the left-hand cannot have any terms of degree zero and $\mathfrak{u}$ is the part of $\mathfrak{g}$ of degree $>0$.
Thus, the bracket $[\hhh,\sigma_w]$ must lie in $\mathfrak{u}$ and indeed each separate term $[\hhh,(\sigma_w)_n ]$ must already lie in $\mathfrak{u}$ since there are no linear relations between terms of different degree. Hence, by (\ref{exp:sideg}), we must have 
\beq
(2w-n) (\sigma_w)_n  \in \mathfrak{u}
\label{hhequ}
\eeq
for all $n\ge w+1$, i.e.~for all terms of $\sigma_w$. In particular, whenever $2w-n\ne 0$, (\ref{hhequ}) implies that $(\sigma_w)_n \in \mathfrak{u}$. Terms of $\sigma_w$ not in $\mathfrak{u}$ can thus only occur when $n=2w$, i.e.~in key degree. The fact that $\sigma_w^{\rm key}\notin \mathfrak{u}$ follows directly from Theorem \ref{thm:HM}, since if $\sigma_w^{\rm key}$ lied in  $\mathfrak{u}$ then we would have $\sigma_w \in\mathfrak{u}$, so~$\mathfrak{u}$ together with the $\sigma_w$ could not generate a semi-direct product as in Theorem \ref{thm:HM} (ii).

\vspace{.2cm}
\noindent (ii) Once again, by \eqref{HMmagic}, any bracket of $\mathfrak{sl}_2$ elements and $\sigma_w$, and therefore in particular $[\ep_0,(\sigma_w)_{2w}]$ must lie in $\mathfrak{u}$. 
If we decompose
\beq
(\sigma_w)_{2w}= \sum_{{\rm odd}\ d\ge 1}
(\sigma_w)_{2w}^{(d)}
\, ,
\eeq
where $(\sigma_w)_{2w}^{(d)}$ collects the key-degree terms in $\sigma_w$ that lie in $d$-dimensional irreducible representations of $\mathfrak{sl}_2$, we must then have 
\beq
[\ep_0, (\sigma_w)_{2w}^{(d)}] \in \mathfrak{u}
\label{neednumber}
\eeq
separately for each (odd) $d\geq 1$.
When $d\ge 3$, the terms  $[\ep_0, (\sigma_w)_{2w}^{(d)} ]\in \mathfrak{u}$  are non-zero since highest-weight vectors of $(d\geq 3)$-dimensional $\mathfrak{sl}_2$ representations have $\hhh$-eigenvalue $\geq 2$. 
Then, thanks to the equality\footnote{The prefactor follows from the fact that the $\mathfrak{sl}_2$ properties of $[\ep_0, (\sigma_w)_{2w}^{(d)}]$ are identical to $\ep_{d+1}^{(\frac{d+1}{2})}$, where the action of the lowering operator $\ad_{\ep_0^\vee}$ yields $\frac{1}{4}(d-1)(d+1)  \ep_{d+1}^{(\frac{d-1}{2})}$ by (\ref{epjkep}).}
\beq\label{useful}
(\sigma_w)_{2w}^{(d)} = \frac{4}{(d-1)(d+1)} [\ep_0^\vee,[\ep_0, (\sigma_w)_{2w}^{(d)}]]\, ,
\eeq
we see that for $d\ge 3$,
the term $(\sigma_w)_{2w}^{(d)}$ itself lies in $\mathfrak{u}$ since $\mathfrak{u}$ is an $\mathfrak{sl}_2$-module by Theorem \ref{thm:HM}.

When $d=1$, the term $[\ep_0,(\sigma_w)_{2w}^{(1)}]=0$ and therefore we cannot use \eqref{neednumber} to conclude that $(\sigma_w)_{2w}^{(1)}$ lies in $\mathfrak{u}$; indeed we know that it cannot lie in $\mathfrak{u}$ since otherwise all of $\sigma_w$ would, contradicting (i). This proves that the arithmetic terms $z_w$ of $\sigma_w$ form a one-dimensional $\mathfrak{sl}_2$ representation in key degree. 

Finally, (iii) follows directly from \eqref{HMmagic}, since this shows that $[\ep_k,\sigma_w]\in \mathfrak{u}$ and $z_w$ is the only term of $\sigma_w$ not already in $\mathfrak{u}$.
\qed

\section{\texorpdfstring{Recursive high-order computations of $\sigma_w$ and $[z_w,\ep_k]$}{Recursive high-order computations of sigma(w) and [z(w),epsilon(k)]}}
\label{bigsecmd2}

In this section, we combine representation theory of $\mathfrak{sl}_2$ with Theorem \ref{thm:522},  particularly part (vii) recalled below, to perform explicit high-order computations of $\sigma_w$ and $[z_w,\ep_k]$ in terms of nested brackets of $\ep_k^{(j)}$.

\subsection{Proof and first consequences of Theorem~\ref{thm:522} (vii)}
\label{newsec:71}

\begin{prop}[Theorem \ref{thm:522} (vii)] Let
$\BF_k \coloneqq  \frac{ {\rm B}_k}{k!}$ for $k\ge 2$, and set
\beq\label{N}
N \coloneqq  -\ep_0 + \sum_{k=4}^{\infty} (k-1)\BF_k \ep_k\, .
\eeq 
Then for all odd $w\ge 3$ we have
\beq
[N,\sigma_w]=0\in {\rm Der}^0\widehat{\rm Lie}[a,b] \, .
\label{cons.01}
\eeq
\end{prop}

\noindent {\bf Proof.} The proof of this result is given in section 27 of \cite{hain_matsumoto_2020} based on sections 12 and 13 of \cite{Hain}, so we simply indicate the essential argument here. In the framework set forth in Remark
\ref{framework}, we noted that the two profinite groups $\widehat{{\rm SL}}_2(\mathbb{Z})$ and the absolute Galois group ${\rm Gal}(\overline{\Q}/\Q)$ both act naturally as automorphisms on the profinite fundamental group $\hat\pi_1(E_\infty)$ of the nodal elliptic 
curve, where ${\rm SL}_2(\mathbb{Z})$ is identified with the fundamental group of the moduli space ${\cal M}_{1,1}$. There is a distinguished element in $\widehat{\rm SL}_2(\mathbb{Z})$ on which ${\rm Gal}(\overline{\Q}/\Q)$ acts via
its abelian quotient ${\rm Gal}(\overline{\Q}^{\rm ab}/\Q)$: this is the element corresponding to a small loop around the degenerate point $\tau=i\infty$ in the moduli space (or as
Hain--Matsumoto describe it, a small loop around $q=0$ in the $q$-disk where $q=e^{2\pi i\tau}$). Thus in the pro-unipotent version, or rather the associated Lie algebra version, the arithmetic part $\sigalg$ corresponding to the Galois action commutes with the image of this element in the Lie algebra $\mathfrak{u}\rtimes \mathfrak{sl}_2$. There are various ways of showing that this image is equal to the element $N$ defined in \eqref{N}; the method used in section 12 of \cite{Hain} is to identify it as the residue at $q=0$ of the restriction of the KZB connection (see appendix \ref{app:deg}) to a first order neighborhood of the degenerate nodal curve. \qed

In the remainder of this section, the commutation relation \eqref{cons.01} will be applied to recursively determine the {\em infinite} series expansions of $\sigma_w$ as in (\ref{firsts3}) 
to (\ref{sig9exp}) from the {\em finitely} many terms in degree $\le 2w$.
The finitely many contributions to $\sigma_w$ not yet determined by \eqref{cons.01} are precisely the highest-weight vectors of $\mathfrak{sl}_2$, i.e.~the elements in the kernel of $\ad_{\epsilon_0}$. By Theorem \ref{thm:522} (iv), these highest-weight-vector contributions to $\sigma_w$ occur only up to and including key degree $2w$ which explains the finite number of them for each $w$.

For example, when $w=3$, the key degree is $6$ and feeding the highest-weight vector contributions $-\frac{1}{2} \epsilon_4^{(2)}$ and $z_3$ into \eqref{cons.01} determines all of $\sigma_3$, see (\ref{exsig3}) below for the exact result. When $w=5$, the highest-weight vectors
$-\frac{1}{24} \epsilon_6^{(4)}, \, -\frac{5  }{48} [\epsilon_4^{(1)},\epsilon_4^{(2)}]$ and $z_5$ occurring in the low-degree part of  $\sigma_5$ feed into \eqref{cons.01} and determine all of $\sigma_5$.\footnote{The analogous highest-weight vectors in the expansion
(\ref{sig7exp}) that completely determine
$\sigma_7$ are given by $-\frac{1}{720} \epsilon_8^{(6)}$ at degree 8, by $\frac{7 }{1152} ( [\epsilon_4^{(2)},\epsilon_6^{(3)}] -[\epsilon_4^{(1)},\epsilon_6^{(4)}])$ at degree 10, by $\frac{1}{13824}([\epsilon_6^{(1)},\epsilon_6^{(4)}] -
[\epsilon_6^{(2)},\epsilon_6^{(3)}])$ and $-\frac{661}{57600} ( [\epsilon_4^{(1)},[\epsilon_4^{(1)},\epsilon_4^{(2)}]] +
 [\epsilon_4^{(2)},[\epsilon_4^{(2)},\epsilon_4]])$ at degree 12 and finally $z_7$ at key degree 14.}

Our construction of $\sigma_w$ from finitely many highest-weight vectors will be recursive in the {\em modular depth} of its geometric contributions which we define as follows: 
\begin{dfn}
    \label{dfn:moddepth}
Nested brackets $[[\ldots [[\epsilon_{k_1}^{(j_1)} ,\epsilon_{k_2}^{(j_2)} ], \epsilon_{k_3}^{(j_3)} ] ,\ldots ], \epsilon_{k_r}^{(j_r)}]$ of $r$ derivations $\epsilon_k^{(j)}$ in $\mathfrak{u}$ are said to have {\em modular depth} $r$. The modular depth forms a natural increasing filtration on $\mathfrak{u}$, but not a grading, as shown for example by the Pollack relation \eqref{tsurels2} which can be viewed as an equality between linear combinations of
terms of modular depth 2 with two terms of modular depth 3.
\end{dfn}

\vspace{.3cm}
In addition to the infinitely many terms in the series expansion of $\sigma_w$ above key degree, 
the recursive method of section \ref{secmd2} will completely determine the
explicit form of the brackets $[z_w,\epsilon_k]$ of the arithmetic contributions $z_w$ to the zeta generators. We reiterate that, by Theorem \ref{thm:522} (v) and (vi), the non-geometric part $z_w$ of $\sigma_w$ is concentrated in a one-dimensional $\mathfrak{sl}_2$ representation at key degree $2w$ and gives rise to brackets $[z_w,\epsilon_k] \in \mathfrak{u}$.

\subsection{\texorpdfstring{$\mathfrak{sl}_2$ prerequisites}{sl2 prerequisites}}
\label{sec:sl2p}

We start by organizing $\mathfrak{u}$ into representations of the subalgebra $\mathfrak{sl}_2$ of ${\rm Der}^0{\rm Lie}[a,b]$, and describing its irreducible pieces; in particular we determine the highest- and lowest-weight vectors of each one. 
 
In view of the nilpotency $\ad_{\epsilon_0}^{k-1} \epsilon_k =0$ (see \eqref{nilpote0}), the non-zero $\ep_k^{(j)}=\ad_{\epsilon_0}^{j} \epsilon_k $ for fixed even $k$ and $j=0,1,\ldots,k-2$ form a $(k-1)$-dimensional irreducible representation of $\mathfrak{sl}_2$, which we denote by $V(\ep_k)$. The generators $\epsilon_0,\epsilon_0^\vee,\hhh$ of $\mathfrak{sl}_2$ permute the elements of $V(\ep_k)$ simply by $\ad_{\epsilon_0}\ep_k^{(j)}=\ep_k^{(j+1)}$, (\ref{hbrak.02}) and (\ref{epjkep}), identifying $\ad_{\epsilon_0}$ and $\ad_{\epsilon_0^\vee}$ as the raising and lowering operators for the eigenvalues of $\hhh$, respectively. All irreducible representations of $\mathfrak{sl}_2$ inside $\mathfrak{u}$ are formed from nested commutators of the $\ep_k^{(j)}$, and they are all isomorphic (as $\mathfrak{sl}_2$-representations) to some $V(\ep_k)$ for even $k\geq 2$. Note that each odd-dimensional $\mathfrak{sl}_2$-representation occurs infinitely many times in $\mathfrak{u}$, and they can be arranged by modular depth.

The collections of commutators $[\ep_{k_1}^{(j_1)},\ep_{k_2}^{(j_2)}]$
for fixed $k_1,k_2$ and $j_i=0,1,\ldots,k_i-2$ sit inside the reducible tensor-product representations $V(\ep_{k_1})\otimes V(\ep_{k_2})$ of $\mathfrak{sl}_2$ which can be decomposed into
the following $(r-1)$-dimensional irreducible representations $V(\ep_r)$ of $\mathfrak{sl}_2$:
\beq
V(\ep_{k_1})\otimes V(\ep_{k_2}) = \bigoplus_{\substack{r= |k_1-k_2|+2\\ r \in 2\mathbb Z}}^{k_1+k_2-2}  V(\ep_r)\, .
\label{sl2tens}
\eeq
Since $r$ is restricted to even values, the dimensions of the irreducible representations of $\mathfrak{sl}_2$ in iterated tensor products of $V(\ep_{k_i})$ are always odd.

\subsubsection{Projectors to lowest-weight vectors}

The projection of the commutators $[\ep_{k_1}^{(j_1)},\ep_{k_2}^{(j_2)}]$ at modular depth two 
into the irreducible representations $V(\ep_r)$ on the right-hand side of (\ref{sl2tens}) is implemented by
\begin{align}
\label{eq:deft}
 t^d(\ep_{k_1},\ep_{k_2}) \coloneqq  \frac{(d{-}2)!}{(k_1{-}2)!(k_2{-}2)!}\sum_{i=0}^{d-2} (-1)^i \frac{(k_1{-}2{-}i)! (k_2{-}d{+}i)!}{i! (d{-}2{-}i)!} [ \ep_{k_1}^{(i)} , \ep_{k_2}^{(d-2-i)} ] 
\end{align}
with $d = \frac{1}{2}(k_1+k_2-r+2)$ and therefore $2\leq d \leq {\rm min}(k_1,k_2)$. In case of $k_1=k_2$, the $ t^d(\ep_{k},\ep_{k})$ at even values of $d$ vanish.

The outcomes $t^d(\ep_{k_1},\ep_{k_2})$ of the projectors in (\ref{eq:deft}) are lowest-weight vectors (see Definition \ref{dfn:sl2}) of the $V(\ep_r)$ in the tensor product (\ref{sl2tens}). The rest of the $(r-1)$-dimensional irreducible representations in $\mathfrak{u}$ at modular depth two is obtained from ${\rm ad}_{\ep_0}^j t^d(\ep_{k_1},\ep_{k_2})$ with $j=0,1,\ldots,r-2$ and terminates due to ${\rm ad}_{\ep_0}^{r-1} t^d(\ep_{k_1},\ep_{k_2})=0$. 

Since  $t^{d_1}(\ep_{k_1},\ep_{k_2})$ is a lowest-weight vector it can be inserted on the same footing as $\epsilon_r$ with $r=k_1+k_2-2d_1+2$ into another operation (\ref{eq:deft}). For instance,
\begin{align}
\label{itdeft}
 t^{d_2}(\ep_{k_3},t^{d_1}(\ep_{k_1},\ep_{k_2})) &= \frac{(d_2{-}2)!}{(k_3{-}2)!(r{-}2)!}\sum_{i=0}^{d_2-2} (-1)^i \frac{(k_3{-}2{-}i)! (r{-}d_2{+}i)!}{i! (d_2{-}2{-}i)!}\notag \\
 &\quad \times [ \ep_{k_3}^{(i)} ,  {\rm ad}_{\ep_0}^{d_2-2-i} t^{d_1}(\ep_{k_1},\ep_{k_2})] 
\end{align}
is the lowest-weight vector of a $(k_1+k_2+k_3-2d_1-2d_2+3)$-dimensional irreducible $\mathfrak{sl}_2$ representation in the triple tensor product $V(\ep_{k_1})\otimes V(\ep_{k_2})\otimes V(\ep_{k_3})$ which
may be decomposed into irreducibles by iterating (\ref{sl2tens}). Iterations of the $t^d$ projectors (\ref{eq:deft}) as exemplified in (\ref{itdeft}) are instrumental for compactly representing the contributions to $[z_w,\ep_k]$ at modular depth three in section \ref{secmd2.b} below.

\subsubsection{Projectors to highest-weight vectors}

One can similarly generate highest-weight vectors of the the irreducible representations $V(\ep_r)$ in $V(\ep_{k_1})\otimes V(\ep_{k_2})$
and tensor products at higher modular depth~via
\beq
s^d(\ep_{k_1},\ep_{k_2}) \coloneqq 
\frac{(d{-}2)! }{(k_1{-}2)! (k_2{-}2)!} \sum_{i=0}^{d-2} (-1)^i   [ \ep_{k_1}^{(k_1-2-i)}, \ep_{k_2}^{(k_2-d+i)}]
\label{defspq}
\eeq
 where again $d = \frac{1}{2}(k_1+k_2-r+2)$, as long as $2\leq d \leq {\rm min}(k_1,k_2)$. Nevertheless, we will see that an extension of 
(\ref{defspq}) to $d>{\rm min}(k_1,k_2)$ will be useful to bring certain contributions to $\sigma_w$
into a convenient form, though the highest-weight vector property $[\ep_0, s^d(\ep_{k_1},\ep_{k_2})]=0$ 
only holds for $d \leq {\rm min}(k_1,k_2)$. Since the entries $\ep_{k_1},\ep_{k_2}$ of the $s^d$-operation in (\ref{defspq}) are lowest-weight 
vectors, the nested brackets relevant to modular depth $m\geq 3$ are generated by $m$ iterations of
$t^{d_i}$ and a single $s^d$ operation for the outermost bracket. For instance,
\begin{align}
\label{itdefs}
 s^{d_2}(\ep_{k_3},t^{d_1}(\ep_{k_1},\ep_{k_2})) &= \frac{(d_2{-}2)! }{(k_3{-}2)! (r{-}2)!} \sum_{i=0}^{d_2-2} (-1)^i    [ \ep_{k_3}^{(k_3-2-i)},{\rm ad}_{\ep_0}^{k_2-d_2+i}  t^{d_1}(\ep_{k_1},\ep_{k_2}) ]
\end{align}
at suitable values for $d_1,d_2$ (with $r=k_1+k_2-2d_1+2$) generate all highest-weight vectors 
of the irreducible $\mathfrak{sl}_2$ representations in $V(\ep_{k_1})\otimes V(\ep_{k_2})
\otimes V(\ep_{k_3})$. In general, iterations of $s^{d_{m-1}} t^{d_{m-2}}\ldots  t^{d_1}$
conveniently capture the highest-weight-vector 
contributions to $\sigma_w$ at each modular depth that are not yet determined
by the recursion below based on $[N,\sigma_w]=0$ (see Theorem \ref{thm:522} (vii)).

\subsubsection{$\mathfrak{sl}_2$ representations of Pollack relations} 

The Pollack relations among $\ep_k^{(j)}$ with $k\geq 4$ and $ 0\leq j \leq k-2$ in Remark \ref{Pollackrel} fall into irreducible $\mathfrak{sl}_2$ representations of dimension $\geq 11$.\footnote{More specifically, Pollack relations whose relative factors in the modular-depth-two contributions $[ \ep_{k_1}^{(j_1)} , \ep_{k_2}^{(j_2)}]$ are governed by holomorphic cusp forms of modular weight $w$ \cite{Pollack} fall into irreducible $\mathfrak{sl}_2$ representations of
dimension $w-1$.} As exemplified by the second relation in (\ref{tsurels}), Pollack relations generically mix contributions of different modular depths $\geq 2$.

\subsection{\texorpdfstring{Recursive higher-order computations of $\sigma_w$ and $[z_w,\ep_k]$}{Recursive higher-order computations of sigma(w) and [z(w),epsilon(k)]}}
\label{secmd2}

Based on the vanishing of $[N,\sigma_w]$ in section \ref{newsec:71} and the $\mathfrak{sl}_2$ prerequisites of section \ref{sec:sl2p}, we shall now set up the recursive high-order computations of $\sigma_w$ and $[z_w,\ep_k]$ in terms of nested brackets of $\ep_k^{(j)}$. For this purpose, we parametrize the desired expressions according to modular depth.

\begin{dfn}\label{split} 
Given that $\sigma_w - z_w $
and $[z_w,\epsilon_k]$ both lie in $\mathfrak{u}$ for any odd $w\geq 3$ and even $k \geq 4$ by Theorem \ref{thm:522} (v) and (vi), we expand
\begin{align}
\sigma_w &= z_w + \sigma_w^{\{ 1 \}}+ \sigma_w^{\{ 2 \}}+ \sigma_w^{\{ 3 \}}+\ldots+ \sigma_w^{\{ w \}} \, ,\label{mdexp.01} \\
[z_w,\epsilon_k] &= [z_w,\epsilon_k]^{\{ 1 \}}+[z_w,\epsilon_k]^{\{ 2 \}}+ [z_w,\epsilon_k]^{\{ 3 \}} + \ldots + [z_w,\epsilon_k]^{\{ w+1 \}} \, ,
\notag
\end{align}
where $\sigma_w^{\{ m\}}$ and $[z_w,\epsilon_k]^{\{ m \}}$ refer to combinations of $[[\ldots [[\ep_{k_1}^{(j_1)},\ep_{k_2}^{(j_2)}],\ep_{k_3}^{(j_3)}],\ldots], \ep_{k_m}^{(j_m)} ]\in \mathfrak{u}$ at modular depth $m=1,2,\ldots,w+1$. The properties of the arithmetic derivations $z_w \in {\rm Der}^0\widehat{\rm Lie}[a,b]$ outside $\mathfrak{u}$ can be found in Theorem \ref{thm:522} (vi) --- $a$- and $b$-degree $w$ and vanishing commutators $[z_w,\ep_0]=[z_w,\ep_0^\vee]=0$. 
\end{dfn}

\vspace{.3cm}
\begin{rmk}
\label{nozep1}
The maximum modular depth $w$ of $\sigma_w$ and $w+1$ of $[z_w,\ep_k]$ in (\ref{mdexp.01}) both follow from the fact that each $\ep_m$ with $m\geq 0$ has $b$-degree 1: the $b$-degrees $w$ of $\sigma_w$ (see Theorem \ref{thm:522} (ii)) and
$w+1$ of $[z_w,\epsilon_k]$ are incompatible with modular depths 
$\sigma_w^{\{ m\geq w+1\}}$ and $[z_w,\epsilon_k]^{\{ m\geq w+2 \}}$.
The well-known vanishing of $[z_w,\epsilon_k]^{\{ 1 \}}$ \cite{Pollack, brown_2017, hain_matsumoto_2020} follows from the fact that only expression in $\mathfrak{u}$ compatible with its $a$- and $b$-degrees is $\ep_{2w+k}^{(w)}$ which violates the lowest-weight-vector property of $z_w$ and $\epsilon_k$.
\end{rmk}

\vspace{.3cm}
\begin{rmk}
\label{polrmk}
We recall that generic Pollack relations among $\ep_k^{(j)}$ with $k\geq 4$ and $0\leq j \leq k-2$ in Remark \ref{Pollackrel} relate nested brackets of different modular depth $\geq 2$. 
Accordingly, the individual contributions $\sigma_w^{\{ m\geq 2 \}}$ and $[z_w,\epsilon_k]^{\{ m\geq 2 \}}$ to the right-hand side of (\ref{mdexp.01}) are usually not well-defined before specifying a scheme of applying those Pollack relations that mix modular depths.\footnote{For instance, the image of the second relation in (\ref{tsurels}) under $\ad_{\ep_0}^{10}$ can be used to convert contributions $\sim s^3(\ep_4,\ep_{12}),s^3(\ep_6,\ep_{10})$ and $s^3(\ep_8,\ep_{8})$ to $\sigma_{13}^{ \{2\} }$ into contributions $\sim [\ep_4^{(2)},[\ep_4^{(2)},\ep_8^{(6)}]]$ and $[\ep_6^{(4)},[\ep_6^{(4)},\ep_4^{(2)}]]$ to $\sigma_{13}^{ \{3\} }$. Similarly, the coefficient of $t^4(\ep_4, \ep_{14})$ in $[z_3, \ep_{12}]^{\{2\}}$ can be modified through Pollack relations of degree 18 at the cost of extra terms in all of $[z_3, \ep_{12}]^{\{2\}}$, $[z_3, \ep_{12}]^{\{3\}}$ and $[z_3, \ep_{12}]^{\{4\}}$.} We will specify a choice of $\sigma_{w}^{ \{2\} }$ and $[z_{w},\ep_k]^{ \{2\} }$ for all odd $w\geq 3$ in (\ref{cl.md2}) and (\ref{clmd2}) below which eliminates some of the ambiguities in $\sigma_{w}^{ \{3\} }$ and $[z_{w},\ep_k]^{ \{3\} }$ (those that descend from Pollack relations involving terms of modular depth two). Nevertheless, the recursive relations among $\sigma_{w}^{ \{m\} }$ to be derived below are valid for any scheme of applying Pollack relations that mix different modular depths as long as the same choice is consistently applied to all modular depths $m\geq 2$.

In the companion paper \cite{Dorigoni:2024oft}, we study uplifts of zeta generators $\sigma_w \rightarrow \hat\sigma_w$ which no longer act on $\widehat{\rm Lie}[a,b]$ and where the $\ep_k^{(j)}$ in their series expansion in $\mathfrak{u}$ are promoted to free-algebra generators $ \eee_k^{(j)}$ with $k\geq 4$ and $0\leq j \leq k-2$. The expansion of the uplifted $\hat\sigma_w$ in terms of $\eee_k^{(j)}$ is determined from considerations of non-holomorphic modular forms and does not share the ambiguities from Pollack relations. Accordingly, the uplifted $\hat\sigma_w$ induce preferred representations of the 
$\sigma_{w}^{ \{m\} }$ and $[z_w,\ep_k]^{ \{m\} }$ at $m= 2$ and partially at $m=3$ which will be followed in section \ref{newsec:74}.
\end{rmk}

\vspace{.3cm}
With the notation of Definition \ref{split} for the contributions of fixed modular depth $m$, we organize the property $[N,\sigma_w]=0$ as written in (\ref{cons.01}) according to modular depth
\begin{align}
0 = [N,\sigma_w] &= {-} [\ep_0, \sigma_w^{\{ 1 \}}+\sigma_w^{\{ 2 \}}+\ldots +\sigma_w^{\{ w \}}]
\label{nsigcond} \\
&\quad + \sum_{k=4}^{\infty} (k-1) \BF_k
\Big([\ep_k,\sigma_w^{\{ 1 \}}]
+[\ep_k,\sigma_w^{\{ 2 \}}]+\ldots+[\ep_k,\sigma_w^{\{ w \}}] \notag \\
&\quad\quad\quad\quad\quad
-[z_w,\epsilon_k]^{\{ 1 \}}-[z_w,\epsilon_k]^{\{ 2 \}}-\ldots-[z_w,\epsilon_k]^{\{ w+1 \}}
\Big) \, , \notag
\end{align}
where $\BF_k \coloneqq  \frac{ {\rm B}_k}{k!}$, and we have used $\mathfrak{sl}_2$ invariance $[\ep_0,z_w]=0$.

\begin{prop}
\label{cor:isomd}
Upon isolating the contributions to (\ref{nsigcond}) at fixed modular depth $m=1,2,\ldots,{w+1}$, we deduce
\beq
[\ep_0,\sigma_w^{\{ m \}}] +  \sum_{k=4}^{\infty} (k-1) \BF_k 
[z_w,\epsilon_k]^{\{ m\}} = 
\sum_{k=4}^{\infty} (k-1) \BF_k [\ep_k,\sigma_w^{\{ m-1 \}}]\, .
\label{mdexp.02}
\eeq
In particular: 
\begin{itemize}
\item[(i)] By $\sigma_w^{\{ 0 \}}=0$ and $[z_w,\epsilon_k]^{\{ 1\}}=0$ (see Remark \ref{nozep1}), the $m=1$ instance of (\ref{mdexp.02}) enforces $[\ep_0,\sigma_w^{\{ 1 \}}]=0$. 
Hence, the only term in $\sigma_w^{\{ 1 \}}$ of modular depth one compatible with the $b$-degree $w$ of $\sigma_w$ and (\ref{mdexp.02}) is the
highest-weight vector  $\sigma_w^{\{ 1 \}}= - \frac{1}{(w-1)!} \ep_{w+1}^{(w-1)}$ identified in Corollary \ref{cor:lowdeg} (i).
\item[(ii)] Applying $\ad_{\ep_0^\vee}$ to both sides of (\ref{mdexp.02}) implies ($m=2,3,\ldots,w+1$)
\beq
[\ep_0^\vee,[\ep_0,\sigma_w^{\{ m \}}] ] = 
\sum_{k=4}^{\infty} (k-1) \BF_k [\ep_k,[\ep_0^\vee,\sigma_w^{\{ m-1 \}}]]
\label{mdexp.03}
\eeq
since both $z_w$ and $\ep_k$ are annihilated by $\ad_{\ep_0^\vee}$. This is the recursive approach announced earlier on to determine $\sigma_w^{\{ m \}}$ from its precursor at lower modular depth $\sigma_w^{\{ m-1 \}}$ up to the kernel of $\ad_{\ep_0^\vee} \ad_{\ep_0}$. Since $\ad_{\ep_0^\vee}$ is invertible on the image of $\ad_{\ep_0}$, see (ii) of Corollary \ref{cor:rec1} with $Y \in \ad_{\ep_0^\vee}\mathfrak{u}$
on the right-hand side composed of $[\ep_k,[\ep_0^\vee,\sigma_w^{\{ m-1 \}}]]=[\ep_0^\vee,[\ep_k,\sigma_w^{\{ m-1 \}}]]$, the only part of $\sigma_w^{\{ m \}}$ which is not yet determined by (\ref{mdexp.03}) is in the kernel of $\ad_{\ep_0}$, i.e.\ a combination of highest-weight vectors of $\mathfrak{sl}_2$. By Theorem \ref{thm:522} (iv) proven in section \ref{sec:prfv}, the highest-weight vectors in $\sigma_w$ all occur below or at key degree. In fact, $z_w$ gathers all highest-weight vectors in $\sigma^{\rm key}_w$ by definition, so $\sigma_w^{\{ m \}}$ at degree $2w$ is free of highest-weight vectors. Hence, the missing information on $\sigma_w^{\{ m \}}$ inaccessible from (\ref{mdexp.03}) amounts to a finite number of terms at degree $\leq 2w-2$.
\item[(iii)] By inserting the expression for $\sigma_w^{\{ m \}}$ modulo highest-weight vectors found in (ii) into (\ref{mdexp.02}) and isolating terms of degree $2w+k$, one can solve for
$[z_w,\epsilon_k]^{\{ m\}}$. Note that contributions to $[z_w,\epsilon_k]$ of modular depth $m$ determined from $[N,\sigma_w]=0$ only depend on the highest-weight vectors in $\sigma_w$ up to and including modular depth $m-1$.
\item[(iv)] Given that $\sigma_w^{\{ 1 \}}= - \frac{1}{(w-1)!} \ep_{w+1}^{(w-1)}$, the $m=2$ instances of (\ref{mdexp.02}) and (\ref{mdexp.03}) can be written more explicitly as
\beq
[\ep_0,\sigma_w^{\{ 2 \}}] +  \sum_{k=4}^{\infty} (k-1) \BF_k 
[z_w,\epsilon_k]^{\{ 2\}} = -\frac{1}{(w-1)!}
\sum_{k=4}^{\infty} (k-1) \BF_k [\ep_k,\ep_{w+1}^{(w-1)}]
\label{mdexp.04}
\eeq
and
\beq
[\ep_0^\vee,[\ep_0,\sigma_w^{\{ 2 \}}]]  = -\frac{1}{(w-2)!}
\sum_{k=4}^{\infty} (k-1) \BF_k [\ep_k,\ep_{w+1}^{(w-2)}]\, .
\label{mdexp.05}
\eeq
Inverting the operation $\ad_{\ep_0^\vee} \ad_{\ep_0}$ determines
\begin{align}
\sigma_w^{\{2\}}
&= - \sum_{d=5}^w \BF_{d-1} s^d(\ep_{d-1},\ep_{w+1}) 
 - \frac{1}{2} \BF_{w+1} s^{w+2}(\ep_{w+1},\ep_{w+1})
\notag \\
&\quad + \sum_{k=w+3}^\infty
\BF_k \sum_{j=0}^{w-2}  \frac{(-1)^j  \binom{k{-}2}{j}^{-1} }{j! (w{-}2{-}j)! } \, 
[  \ep_{w+1}^{(w-2-j)}  , \ep_k^{(j)} ]
\ {\rm mod} \ {\rm Ker}(\ad_{\ep_0})
\, ,
\label{s2modhw}
\end{align}
where ${\rm mod} \ {\rm Ker}(\ad_{\ep_0})$ refers to highest-weight vectors to be proposed in (\ref{clmd2}) below. All instances of the brackets $s^d(\ep_{k_1},\ep_{k_2})$ defined by (\ref{defspq}) that occur in (\ref{s2modhw}) have $d> {\rm min}(k_1,k_2)$ and are therefore not highest-weight vectors. Upon insertion of (\ref{s2modhw}) into (\ref{mdexp.04}) and isolating terms of degree $2w+k$, we reproduce the closed-form expression at modular depth two known form \cite{hain_matsumoto_2020}
\beq
[z_w,\epsilon_k]^{\{ 2\}} = \frac{{\rm BF}_{w+k-1}}{ {\rm BF}_{k} } t^{w+1}(\ep_{w+1},\ep_{w+k-1}) \, .
\label{cl.md2}
\eeq
\item[(v)] The instance of (\ref{mdexp.02}) at the maximum value $m=w+1$ simplifies to
\beq
\sum_{k=4}^{\infty} (k-1) \BF_k 
[z_w,\epsilon_k]^{\{ w+1\}} = 
\sum_{k=4}^{\infty} (k-1) \BF_k [\ep_k,\sigma_w^{\{ w \}}]
\label{mdexp.09}
\eeq
by $\sigma_w^{\{ w+1\}}=0$. Hence, the contribution to $[z_w,\epsilon_k]$ of 
highest modular depth $w+1$ can simply be determined from the highest-modular depth terms in $\sigma_w$ by isolating the parts of degree $2w+k$ in (\ref{mdexp.09}). Validity of (\ref{mdexp.02}) at $m=1,2,\ldots,w+1$ --- finitely many steps in the recursion in the modular depth --- is sufficient for $[N,\sigma_w]=0$, see (\ref{nsigcond}).
\end{itemize}
\end{prop}
Note that parts (ii) and (iii) of Proposition \ref{cor:isomd} can also be unified by the decomposition of $[\ep_k,\sigma_w^{\{ m-1 \}}]$ on the right-hand side of (\ref{mdexp.02}) into the  image of $\ad_{\ep_0}$ and the kernel of $\ad_{\ep_0^\vee}$,
\begin{align}
[\ep_0,\sigma_w^{\{ m \}}] &= 
\sum_{k=4}^{\infty} (k-1) \BF_k [\ep_k,\sigma_w^{\{ m-1 \}}] \, \big|_{{\rm Im}(\ad_{\ep_0})}\, ,
\notag \\
 \sum_{k=4}^{\infty} (k-1) \BF_k 
[z_w,\epsilon_k]^{\{ m\}} &= 
\sum_{k=4}^{\infty} (k-1) \BF_k [\ep_k,\sigma_w^{\{ m-1 \}}] \, \big|_{{\rm Ker}(\ad_{\ep_0^\vee})} \, .
\end{align}
This decomposition is unique since ${\rm Ker}(\ad_{\ep_0^\vee})$ projects the individual terms of $[\ep_k,\sigma_w^{\{ m-1 \}}] $ to lowest-weight vectors which do not occur in the image of $\ad_{\ep_0}$.

\subsection{\texorpdfstring{Applying the recursion for $\sigma_w^{\{ m\}}$ and $[z_w,\ep_k]^{\{ m\}}$}{Applying the recursion for sigma(w) and [z(w),epsilon(k)]}}
\label{newsec:74}

In this section, we gather explicit results for zeta generators and commutators $[z_w,\ep_k]$ at modular depth $2\leq m \leq 4$  that go considerably beyond the state of the art and found fruitful applications in the construction of non-holomorphic modular forms \cite{Dorigoni:2024oft}. 

\subsubsection{Zeta generators at modular depth two}
\label{secmd2.1}

The relation (\ref{mdexp.05}) for the modular-depth-two contributions $\sigma_w^{\{2\}}$ to the zeta generators determines the infinite series of terms in (\ref{s2modhw}) that are not highest-weight vectors. We shall now augment these terms by a conjectural closed formula for the highest-weight vectors in $\sigma_w^{\{2\}}$ given by the first line of 
\begin{align}
\sigma_w^{\{2\}}
&=
-\frac{1}{2} \sum_{d=3}^{w-2} \frac{ \BF_{d-1} }{\BF_{w-d+2} }
\sum_{k=d+1}^{w-1} \BF_{k-d+1} \BF_{w-k+1} s^d(\ep_k,\ep_{w-k+d})
\notag\\
&\quad 
- \sum_{d=5}^w \BF_{d-1} s^d(\ep_{d-1},\ep_{w+1}) 
 - \frac{1}{2} \BF_{w+1} s^{w+2}(\ep_{w+1},\ep_{w+1})
\notag \\
&\quad + \sum_{k=w+3}^\infty
\BF_k \sum_{j=0}^{w-2}  \frac{(-1)^j  \binom{k{-}2}{j}^{-1} }{j! (w{-}2{-}j)! } \, 
[  \ep_{w+1}^{(w-2-j)}  , \ep_k^{(j)} ]\, .
\label{clmd2}
\end{align}
This conjecture for the complete parts $\sigma_w^{\{2\}}$ of modular depth two is readily checked to
reproduce the terms $[\epsilon_{k_1}^{(j_1)},\epsilon_{k_2}^{(j_2)}]$ in the examples (\ref{firsts3}) to (\ref{sig9exp}) at $w\leq 9$.
The first line of (\ref{clmd2}) gathers highest-weight vectors such as $-\frac{5  }{48} [\epsilon_4^{(1)},\epsilon_4^{(2)}]$ in $\sigma_5^{\{2\}}$ and $\frac{7 }{1152} ( [\epsilon_4^{(2)},\epsilon_6^{(3)}] -[\epsilon_4^{(1)},\epsilon_6^{(4)}])+\frac{1}{13824}([\epsilon_6^{(1)},\epsilon_6^{(4)}] -
[\epsilon_6^{(2)},\epsilon_6^{(3)}])$ in $\sigma_7^{\{2\}}$\footnote{The analogous highest-weight vectors in $\sigma_9^{\{2\}}$ resulting from the first line of (\ref{clmd2}) are given by $\frac{1}{5184}([\epsilon_4^{(2)},\epsilon_8^{(5)}] - [\epsilon_4^{(1)},\epsilon_8^{(6)}])$ and $-\frac{7 }{20736}[\epsilon_6^{(3)},\epsilon_6^{(4)}]$ at degree 12, $\frac{7 }
{4147200} ( [\epsilon_6^{(1)},\epsilon_8^{(6)}]
- [\epsilon_6^{(2)},\epsilon_8^{(5)}] +[\epsilon_6^{(3)},\epsilon_8^{(4)}]-[\epsilon_6^{(4)},\epsilon_8^{(3)}])$ at degree 14 and $ -\frac{1}{26127360}([\epsilon_8^{(1)},\epsilon_8^{(6)}]-[\epsilon_8^{(2)},\epsilon_8^{(5)}]+[\epsilon_8^{(3)},\epsilon_8^{(4)}])$ at degree 16.} which have been tested for all cases of degree $\leq 22$ and are in general conjectural. Note that the highest-weight-vector contributions to $\sigma_w^{\{2\}}$ in the first line of (\ref{clmd2}) are in one-to-one correspondence with the $\tau \rightarrow i\infty$ asymptotics of the generalized Eisenstein series ${\rm F}^{+(s)}_{m,k}$ in \cite{Dorigoni:2021jfr,Dorigoni:2021ngn} at $m+k+s=w+1$ upon assembling their iterated-integral representations from the generating series of \cite{Dorigoni:2024oft}.

The images of the terms $s^d(\ep_{k_1},\ep_{k_2})$ under the switch operation in Definition \ref{dfn:switch} have $b$-degree or depth $d$, and their $d=3$ instances line up with Brown's general formula for the depth-three contributions to $\tau_w$ \cite{brown_2017}. However, the choice of $\tau_{w\geq 11}$ in the reference does not match the {\it canonical} zeta generators in this work since redefinitions via nested brackets of $\tau_v$ at $v<w$ have been used in \cite{brown_2017} to remove contributions of modular depth and $b$-degree three. The second and third line of (\ref{clmd2}) are rigorously derived by solving (\ref{mdexp.05}) and, together with the conjectural highest-weight vectors at depth $d\geq 5$ in the first line, furnish a partial generalization of Brown's result beyond depth three: On the one hand,  (\ref{clmd2}) is claimed to capture all contributions $[\ep_{k_1}^{(j_1)},\ep_{k_2}^{(j_2)}]$ to~$\sigma_w$, regardless of their values of $j_1,j_2,k_1,k_2$ or depth in the sense of \cite{brown_2017}. On the other hand, terms in $\sigma_w$ at depth or $b$-degree $d$ involve contributions of modular depth up to and including $d$, and closed formulae for $\sigma_w^{\{ m\geq 3\}}$ akin to (\ref{clmd2}) are currently out of reach.

Note that, following the comments below (\ref{defspq}), the $s^{d}(\ep_{k_1},\ep_{k_2})$ in the second line of (\ref{clmd2}) have $d>{\rm min}(k_1,k_2)$ and are therefore not highest-weight vectors. Moreover, the expression 
(\ref{clmd2}) for contributions to $\sigma_w$ of modular depth two can be rewritten in a variety of ways via Pollack relations among $\ep_k^{(j)}$, see Remark \ref{polrmk}. Hence, the closed formula (\ref{clmd2}) for $\sigma_w^{\{ 2 \}}$ realizes a specific choice of distributing terms between different modular depths.

\subsubsection{Commutators of arithmetic derivations at modular depth three}
\label{secmd2.b}

By Proposition \ref{cor:isomd} (iii), the highest-weight vectors in $\sigma_w$ at modular depth $m$ determine the contributions to the brackets $[z_w,\ep_k]$ at modular depth $m+1$ via (\ref{mdexp.02}). The conjectural expressions (\ref{clmd2}) for $\sigma_w^{\{2\}}$ therefore translate into expressions for $[z_w,\ep_k]^{\{3\}}$ that generalize the simple closed formula
(\ref{cl.md2}) for terms of modular depth two.

Contributions to $[z_3,\epsilon_k]$ and $[z_5,\epsilon_k]$ at modular depth $\geq 3$ and low values of $k$ have been firstly reported in \cite{Pollack} and the ancillary files of \cite{Dorigoni:2022npe}, respectively. Moreover, the combinatorial tools developed in \cite{Pollack} can be used to determine more general expressions for $[z_w,\epsilon_k]$. Our conjecture (\ref{clmd2}) for $\sigma_w^{\{2\}}$ gives access to arbitrary $[z_w,\ep_k]^{\{3\}}$, but the expressions resulting from the representation-theoretic manipulations become increasingly unwieldy with growing~$w$. Hence, we content ourselves to giving the following two infinite families of commutation relations beyond the state of the art with arbitrary even $k\geq 4$ (see (\ref{itdeft}) for the iteration of the projector $t^d$ to lowest-weight vectors),
\begin{align}
[z_3,\ep_k]^{\{3\}} &= 
 \frac{3 \BF_4 \BF_{k-2} }{\BF_k} \bigg\{
 {-} \frac{(k{-}3) }{(k{-}1)} \, t^2 (\ep_4, t^3(\ep_4,\ep_{k-2}))
 + \frac{(k{-}2) }{k} \, t^3 (\ep_4, t^2(\ep_4,\ep_{k-2}))
 \bigg\} \label{zwep3.56} \\*
 &\quad +  \frac{1}{(k{-}1) \BF_k} \sum_{\ell=6}^{k-4} (\ell{-}1) \BF_\ell \BF_{k+2-\ell}
\notag \\*
&\quad\quad\quad
\times \bigg\{
{-} \frac{2(k{-}\ell{+}1) }{(k{-}\ell{+}2) } \, t^2 (\ep_\ell, t^3(\ep_4,\ep_{k+2-\ell}))
+\frac{\ell{-}2 }{k} \,  t^3(\ep_\ell, t^2(\ep_4,\ep_{k+2-\ell}))
\bigg\}
 \notag
\end{align}
and
\begin{align}
[z_5,\ep_k]^{\{3\}}&= \frac{ \BF_{k+2} \BF_2^3 }{2 \BF_4 \BF_k } \, t^4(\ep_{k+2},t^3(\ep_4,\ep_4)) \label{zwep3.57} \\
&\quad + \frac{5 \BF_6 \BF_{k-2} }{\BF_k} \bigg\{
 {-}\frac{ (k{-}5)}{(k{-}1)} \, t^2 (\ep_6, t^5(\ep_6,\ep_{k-2})) 
 +\frac{2(k{-}3)(k{-}4) }{k(k{-}1)} \, t^3 (\ep_6, t^4(\ep_6,\ep_{k-2})) \notag \\
 &\quad\quad\quad\quad 
 - \frac{2 (k{-}2)(k{-}3) }{k (k{+}1)} \, t^4 (\ep_6, t^3(\ep_6,\ep_{k-2}))
 +\frac{ (k{-}2) }{(k{+}2)} \, t^5 (\ep_6, t^2(\ep_6,\ep_{k-2}))
\bigg\} \notag \\
&\quad + \BF_4 \bigg\{
{-}\frac{12(k{-}3) }{k(k{-}1)}\, t^2(\ep_4, t^5(\ep_6,\ep_k))
+\frac{ 36 (k{-}2) }{k^2(k{+}1)} \,  t^3(\ep_4, t^4(\ep_6,\ep_k)) \notag \\
&\quad\quad\quad\quad 
- \frac{24}{k(k{+}1)(k{+}2)} \,  t^4(\ep_4, t^3(\ep_6,\ep_k))
-\frac{9(k{-}2) }{5k} \,  t^3(\ep_k, t^4(\ep_4,\ep_6)) \notag \\
&\quad\quad\quad\quad 
- \frac{2(k{-}2)(k{-}3) }{k(k{+}1)} \, t^4(\ep_k, t^3(\ep_4,\ep_6))
- \frac{(k{-}2)(k{-}3)(k{-}4) }{k(k{+}1)(k{+}2) } t^5(\ep_k, t^2(\ep_4,\ep_6))
\bigg\}
\notag \\
&\quad+ \frac{1}{(k{-}1) \BF_k} \sum_{\ell=8}^{k-4} (\ell{-}1) \BF_\ell \BF_{k+4-\ell}
\bigg\{
{-}\frac{4(k{-}\ell{+}1) }{(k{-}\ell{+}4)} \, t^2(\ep_\ell, t^5(\ep_6,\ep_{k+4-\ell})) \notag \\
&\quad\quad\quad\quad\quad\quad\quad\quad\quad
+ \frac{6(\ell{-}2) (k{-}\ell{+}2)(k{-}\ell{+}3) }{k(k{-}\ell{+}4)(k{-}\ell{+}5) } \, t^3(\ep_\ell, t^4(\ep_6,\ep_{k+4-\ell})) \notag \\
&\quad\quad\quad\quad\quad\quad\quad\quad\quad
-\frac{4 (\ell{-}3)(\ell{-}2) (k{-}\ell{+}3) }{k(k{+}1)(k{-}\ell{+}6) } \, t^4(\ep_\ell, t^3(\ep_6,\ep_{k+4-\ell})) \notag \\
&\quad\quad\quad\quad\quad\quad\quad\quad\quad
+\frac{(\ell{-}2)(\ell{-}3)(\ell{-}4) }{k(k{+}1)(k{+}2) } \, t^5(\ep_\ell, t^2(\ep_6,\ep_{k+4-\ell}))
\bigg\}\, .
\notag
\end{align}
The remaining brackets $[z_w,\ep_k]^{\{3\}}$ at degree $\leq 20$ are given by
\begin{align}
[z_7,\ep_4]^{\{3\}} &= 
\frac{ \BF_8 \BF_2^2 }{\BF_6} \, t^6(\ep_8,t^3(\ep_4,\ep_6))
+ \frac{ \BF_6 \BF_2^2 }{2 \BF_4} \, t^4(\ep_6,t^5(\ep_6,\ep_6))
  \label{zwep3.53}\\
&\quad - \BF_6 \bigg\{
\frac{15}{14} \, t^3(\ep_4, t^6(\ep_6,\ep_8)) 
+ \frac{5}{14} \, t^4(\ep_4, t^5(\ep_6,\ep_8)) \notag \\
&\quad\quad\quad\quad
+\frac{5}{7} \, t^5(\ep_4, t^4(\ep_6,\ep_8))
+\frac{3}{28} \,t^6(\ep_4, t^3(\ep_6,\ep_8))
\bigg\}\, ,
\notag \\
[z_7,\ep_6]^{\{3\}}&= 
 \frac{ \BF_{10} \BF_4 \BF_2^2 }{\BF_6^2} \, t^6(\ep_{10},t^3(\ep_4,\ep_6))
 + \frac{ \BF_8 \BF_2^2 }{2 \BF_6} \, t^4(\ep_8,t^5(\ep_6,\ep_6))   \notag\\
 &\quad - \frac{ \BF_4 \BF_8 }{\BF_6} \bigg\{
 \frac{5}{2}\, t^5(\ep_8,t^4(\ep_4,\ep_8))
 + \frac{7}{2}\,  t^6(\ep_8,t^3(\ep_4,\ep_8))
 +\frac{14}{5}\,  t^7(\ep_8,t^2(\ep_4,\ep_8))
 \bigg\} \notag \\
 &\quad - \BF_6 \bigg\{
 \frac{10}{7} \,  t^3(\ep_6,t^6(\ep_6,\ep_8)) 
 + \frac{50}{49} \,  t^4(\ep_6,t^5(\ep_6,\ep_8)) \notag \\
 &\quad\quad\quad\quad
 +\frac{25}{84} \,  t^5(\ep_6,t^4(\ep_6,\ep_8))
 +\frac{1}{42} \,  t^6(\ep_6,t^3(\ep_6,\ep_8))
 \bigg\} \, .\notag 
\end{align}

\subsubsection{Exact results for $\sigma_3$ and $z_3$}
\label{sec:sigman.8}

Once the complete set of highest-weight vectors for a given $\sigma_w$ is available, then the recursion (\ref{mdexp.02}) determines all-degree expressions for both $ \sigma_w^{\{ 2 \}},\sigma_w^{\{ 3 \}},\ldots,\sigma_w^{\{ w \}}$ and $[z_w,\ep_k]^{\{2\}},[z_w,\ep_k]^{\{3\}},$ $\ldots, [z_w,\ep_k]^{\{w+1\}}$. With the highest-weight vectors for $\sigma_3,\sigma_5,\sigma_7$ noted in section \ref{newsec:71}, there is no obstruction to algorithmically assembling the exact results for the expansions of $\sigma_w$ and $[z_w,\ep_k]$ at $w\leq 7$.

We shall here display the exact results for $\sigma_3$ and $[z_3,\ep_k]$ which terminate with modular depth three and four, respectively. The all-order expansion of $\sigma_3$ is given by,
\begin{align}
 \sigma_3 &= - \frac{1}{2} \ep_4^{(2)} + z_3 + \frac{1}{480} [\ep_4,\ep_4^{(1)}]
+ \sum_{k=6}^\infty \BF_k \bigg( [\ep_4^{(1)} , \ep_k ] -\frac{ [ \ep_4 , \ep_k^{(1)} ] }{k{-}2} \bigg) \notag \\
&\quad + \sum_{m=4}^{\infty} \sum_{r=6}^{\infty}
\frac{(m{-}1)\BF_m \BF_r}{m{+}r{-}2} \big[ \ep_m ,[\ep_4, \ep_r]\big]\, ,
\label{exsig3}
\end{align}
where the second line is obtained by solving (\ref{mdexp.03}) at $m=w=3$ for $\sigma_3^{\{3\}}$ with the expression for $\sigma_3^{\{2\}}$ determined by the first line. The action of the arithmetic part $z_3$ on $a,b$ can be found in (\ref{z3act}). The expression for $[z_3,\ep_{k}]$ resulting from $[N,\sigma_3]=0$ can be assembled by combining $[z_3,\ep_k]^{\{2\}} = \frac{{\rm BF}_{k+2}}{ {\rm BF}_{k} } t^{4}(\ep_{4},\ep_{k+2})$ from (\ref{cl.md2}) with the expression (\ref{zwep3.56}) for $[z_3,\ep_k]^{\{3\}}$ and the degree-$(2w+k)$ parts of
\begin{align}
&\sum_{k=4}^{\infty}(k-1)\BF_k
[z_3,\ep_k]^{\{4\}} = \sum_{k=4}^{\infty}(k-1)\BF_k  \sum_{m=4}^{\infty}
\sum_{r=6}^{\infty} \frac{ (m-1) \BF_m \BF_r }{(m + r - 2) } [\ep_k,[\ep_m,[\ep_4, \ep_r]]]
\end{align}
which follows from (\ref{mdexp.09}) at $w=3$. The lowest-degree examples of $[z_3,\ep_k]^{\{4\}}$ occur in
\begin{align}
[z_3,\ep_{12}]   &=  \frac{\BF_{14}}{\BF_{12}} t^4(\ep_{4}, \ep_{14})
+\frac{ \BF_4 \BF_{10} }{\BF_{12}}
\bigg\{
{-}\frac{27}{11} t^2(\ep_{4}, t^3(\ep_4,\ep_{10}))  + \frac{ 5}{2}  t^3(\ep_{4}, t^2(\ep_4,\ep_{10})) \bigg\} \notag \\*
&\quad+  \frac{ \BF_6 \BF_8 }{\BF_{12}}
\bigg\{
{-}\frac{35}{44}   t^2(\ep_{6}, t^3(\ep_4,\ep_{8}))
+\frac{5}{33}  t^3(\ep_{6}, t^2(\ep_4,\ep_{8})) \notag \\*
&\quad\quad\quad\quad
- \frac{35}{33}  t^2(\ep_{8}, t^3(\ep_4,\ep_{6}))
+\frac{7}{22}  t^3(\ep_{8}, t^2(\ep_4,\ep_{6}))
\bigg\}  \notag \\*
&\quad+\frac{ 9 \BF_{4}^2 \BF_{6} }{88 \BF_{12} } [\ep_{4}, [\ep_{4}, [\ep_{4}, \ep_{6} ] ] ]
\notag 
\end{align}
as well as
\begin{align}
[z_3,\ep_{14}]   &=  \frac{ \BF_{16}}{\BF_{14}} t^4(\ep_4,\ep_{16})
+ \frac{ \BF_{4} \BF_{12}}{\BF_{14}} \bigg\{  
 \frac{18}{7} t^3( \ep_4, t^2(\ep_4,\ep_{12}))  - \frac{33}{13}t^2( \ep_4, t^3(\ep_4,\ep_{12}))  \bigg\} \notag \\
&\quad + \frac{ \BF_{6}\BF_{10}}{\BF_{14}} \bigg\{   
 \frac{10}{91}  t^3( \ep_6, t^2(\ep_4,\ep_{10}))  
- \frac{9}{13}   t^2( \ep_6, t^3(\ep_4,\ep_{10}))  
\notag \\
&\quad\quad\quad\quad  + \frac{36}{91}  t^3( \ep_{10}, t^2(\ep_4,\ep_{6}))  
  - \frac{15}{13}  t^2( \ep_{10}, t^3(\ep_4,\ep_{6}))   
  \bigg\} \notag \\
&\quad + \frac{ \BF_{8}^2}{ \BF_{14} } \bigg\{    
 \frac{3}{13}  t^3( \ep_8, t^2(\ep_4,\ep_{8}))  
 - \frac{49}{52}  t^2( \ep_8, t^3(\ep_4,\ep_{8}))  
  \bigg\} \notag \\
&\quad+  \frac{ 9 \BF_{4}^2 \BF_{8} }{ 130 \BF_{14}}  [\ep_{4}, [\ep_{4}, [\ep_{4}, \ep_{8} ] ] ]
+ \frac{ 27 \BF_{4} \BF_{6}^2 }{ 104 \BF_{14}} [\ep_{4}, [\ep_{6}, [\ep_{4}, \ep_{6} ] ] ]\, ,
\notag 
\end{align}
also see appendix E.1 of \cite{Dorigoni:2024oft} for
$[z_3,\ep_{k}]$ at $k=4,6,8,10$.

\subsubsection{Highest-weight vectors at modular depth three}
\label{sec:sigman.9}

While a comprehensive study of highest-weight vector contributions to $\sigma_w^{\{m\geq 3\}}$ is left for the future, their instances at $w\leq 11$ are accessible from the ancillary files of \cite{Dorigoni:2024oft}. The simplest highest-weight vector at modular depth three occurs in the expansion (\ref{sig7exp}) of $\sigma_7$ at degree 12 and can be compactly written as $- \frac{661 }{14400} s^{3}(\ep_ 4, t^{3}(\ep_4, \ep_4)) $ through the combination (\ref{itdefs}) of $s^d$ and $t^d$ operations. This shorthand also streamlines the expansions of $\sigma_9,\sigma_{11}$ to 
\begin{align}
\sigma_9  &=  -\frac{\ep_{10}^{(8)}}{ 8! }  + \frac{5 s^3(\ep_{4},\ep_{8})}{18}  + \frac{7 s^3(\ep_{6},\ep_{6}) }{72} 
+ \frac{s^5(\ep_{4},\ep_{10})}{720}   - \frac{7 s^5(\ep_{6},\ep_{8})}{1440} 
\label{sigexp.0} \\*
&\quad 
+ \frac{34921 s^{2}(\ep_ 4, t^{4}(\ep_4, \ep_6))}{1134000}
+ \frac{ 2587 s^{3}(\ep_ 4, t^{3}(\ep_4, \ep_6)) }{37800} - 
\frac{ 529 s^{4}(\ep_ 4, t^{2}(\ep_4, \ep_6)) }{14400} 
\notag\\*
&\quad   -  \frac{s^7(\ep_{6},\ep_{10})}{30240}  + \frac{s^7(\ep_{8},\ep_{8})}{12096}
+\frac{ s^{5}(\ep_ 4, t^{3}(\ep_4, \ep_8)) }{2592}
+\frac{ 7s^{5}(\ep_ 4, t^{3}(\ep_6, \ep_6))  }{51840 } \notag \\*
&\quad - 
\frac{ 34921 s^{4}(\ep_ 6, t^{4}(\ep_6, \ep_4)) }{47628000} - 
\frac{ 2587 s^{5}(\ep_ 6, t^{3}(\ep_6, \ep_4)) }{1587600} + 
\frac{ 529 s^{6}(\ep_ 6, t^{2}(\ep_6, \ep_4)) }{604800 }  
\notag \\
  %
&\quad
\frac{149 s^{3}(\ep_4,t^{3}(\ep_4,t^{3}(\ep_4,\ep_4)))}{13824}-\frac{149s^{4}(\ep_4,t^{2}(\ep_4,t^{3}(\ep_4,\ep_4)))}{69120}
+\ldots
  \notag \\
\sigma_{11} &= -\frac{\ep_{12}^{(10)} }{ 10! } 
+  \frac{11 s^3(\ep_{4},\ep_{10})}{40} + \frac{11 s^3(\ep_{6},\ep_{8}) }{60} 
 + \frac{ 242407 s^{2}(\ep_ 4, t^{2}(\ep_4, \ep_6))}{14735232} 
  +   \frac{s^5(\ep_{4},\ep_{12} )}{720}   \notag \\
&\quad -\frac{s^5(\ep_{6},\ep_{10})}{216}   - \frac{7 s^5(\ep_{8},\ep_{8})}{4320}
  +\frac{11090423 s^{2}(\ep_ 4, t^{4}(\ep_4, \ep_8)) }{309439872} + 
\frac{ 3197 s^{3}(\ep_ 4, t^{3}(\ep_4, \ep_8))  }{57600} \notag \\
&\quad - 
\frac{ 2983 s^{4}(\ep_ 4, t^{2}(\ep_4, \ep_8))  }{86400}  + 
\frac{ 148753 s^{3}(\ep_ 4, t^{3}(\ep_6, \ep_6))  }{7367616} + 
\frac{ 490853 s^{3}(\ep_ 6, t^{3}(\ep_6, \ep_4))  }{17191104} 
 \notag \\
&\quad
+
\frac{ 156805 s^{4}(\ep_ 6, t^{2}(\ep_6, \ep_4))  }{14735232}
+ c \, s^{2}(\ep_ 4, t^{2}(\ep_4, t^{3} (\ep_4,\ep_4))) + \ldots \,,\notag 
\end{align}
where the ellipsis refers to all
contributions of degree $\geq 18$, and the coefficient $c\in \mathbb Q$ of the first modular-depth-four contribution to $\sigma_{11}$ in the last line has not yet been computed. It is, however, a highest-weight vector and entirely fixed by our construction. Note that the $ s^{d_2}(\ep_ {k_3}, t^{d_1}(\ep_{k_1}, \ep_{k_2}))$ only furnish highest-weight vectors if $d_2\leq {\rm min}(k_3,r)$, where $r=k_1+k_2-2d_1+2$. Accordingly, all the terms $s^{d_2}(\ep_ {k_3}, t^{d_1}(\ep_{k_1}, \ep_{k_2}))$ of modular depth three in (\ref{sigexp.0}) are highest-weight vectors with the exception of the contributions
$s^{5}(\ep_ 4, t^{3}(\ep_4, \ep_8)) $ and $s^{5}(\ep_ 4, t^{3}(\ep_6, \ep_6))$ to $\sigma_9$. The ancillary files of \cite{Dorigoni:2024oft} provide all contributions to $\sigma_w^{\{m\leq 3\}}$ at degree $\leq 20$ in machine-readable form which determines all the highest-weight vectors of $\sigma_9^{\{3 \}}$ and~$\sigma_{11}^{\{3 \}}$.

\appendix

\section{Deriving the topological map from the sphere to the torus}
\label{app:deg}

The goal of this appendix is to derive the explicit form of the map (\ref{sec4eq.2}) and (\ref{g01.01}) between the generators $x,y$ and $a,b$ of the fundamental groups in genus zero and genus one, respectively. 
Our derivation will be based on a formulation of the zeta generators in terms of Knizhnik--Zamolodchikov (KZ) connections in genus zero and Knizhnik--Zamolodchikov--Bernard (KZB) connections in genus one. The form of the KZ connection obtained from the degeneration limit of the KZB connection then relates the generators $x,y$ of the fundamental group of the thrice punctured sphere to the generators $a,b$ of the fundamental group of the once-punctured torus.

\subsection{Zeta generators in terms of the KZ connection}
\label{app:deg.1}

In this appendix we assume the conjecture that the surjection from motivic to real MZVs is an isomorphism, and thus identify the motivic version $\Phi^{\mathfrak{m}}(x,y)$ of the modified Drinfeld associator with $\Phi(x,y)$ as defined in \eqref{eq:Phi}. We will systematically assume that $\Phi(x,y)$ is written in the semi-canonical basis defined in section \ref{sec:Zmap.4b}, and use the notation
\beq
g_w=\Phi(x,y)|_{\zeta_w}
\eeq
for the canonical polynomial $g_w$ that then appears in $\Phi$ with coefficient $\zeta_w$ for odd $w\ge 3$  (see Definition \ref{dfn:mzvee}). The power series $\Phi(x,y)$ in \eqref{eq:Phi} can be obtained as the path-ordered exponential of the modified KZ connection $J$ defined by\footnote{The connection $J(x,y;z)$ differs from the classical KZ connection $J_{\rm KZ}(x,y;z)= \bigl(\tfrac{x}{z}+\tfrac{y}{z-1}\bigr) dz$ by changing $y$ to $-y$, corresponding to the relation $\Phi(x,y)=\Phi_{\rm KZ}(x,-y)$ between the power series $\Phi$ and the classical Drinfeld associator (\ref{defdrin}) obtained by path-ordered integration of $J_{\rm KZ}$.} 
\begin{align}
J(x,y;z) &\coloneqq  \, \bigg( \frac{x}{z} + \frac{y}{1-z} \bigg)dz \, , \ \ \ \ z \in \mathbb C \setminus \{0,1 \} \, ,
\notag \\
\Phi(x,y) &= {\rm Pexp} \bigg( \int^1_0 J(x,y;z) \bigg) \, , \label{kzapp.01} \\
g_w(x,y) &= \Phi(x,y) \, \big|_{\zeta_w}
= {\rm Pexp} \bigg( \int^1_0  \bigg[ \frac{x}{z} + \frac{y}{1-z} \bigg]dz \bigg) \, \big|_{\zeta_w}\, ,
\notag
\end{align}
where the iterated integration is taken over the simplex $0< z_1< \cdots< z_r< 1$, and the convention for expanding path-ordered exponentials is
\beq\label{kzapp.02}
{\rm Pexp}\Biggl(\int_0^1 J(z)\Biggr)=1+\sum_{r= 1}^{\infty} \int_0^1 J(z_r) \int_0^{z_r} J(z_{r-1})
\cdots \int_0^{z_3} J(z_2) \int_0^{z_2} J(z_1) \, .
\eeq
The endpoint divergences in (\ref{kzapp.02}) are understood to be regularized by passing to shuffle-regularized versions (\ref{appMZV.03}) of the MZVs in the expansion of $\Phi(x,y)$. 

\begin{figure}[ht]
	\begin{center}
\begin{tikzpicture}[line width=0.30mm]
\def\thyck{0.5mm}
\begin{scope}[xshift=-8.5cm]
\begin{scope}[yshift=4cm]
  \draw  (0,1.2) arc (90:270:1.2cm);
  \draw [arrows={-Stealth[width=1.6mm, length=1.8mm]}, color=gray] (-2.5,0) -- (3.5,0) node[above]{$\quad{\rm Re}(z)$};
    \draw [arrows={-Stealth[width=1.6mm, length=1.8mm]}, color=gray] (0,-1.5) -- (0,1.8) node[right]{${\rm Im}(z)$};
\draw[arrows={-Stealth[width=1.6mm, length=1.8mm]}](0,0)--(0.3,0);
\draw(0,0)--(0.2,0);
\draw (2.5,0)node{$\bullet$} ;
  \draw(2.5,-0.4)node{$z=1$};
\draw (0.3,0) .. controls (1.3,0.1) and (0.8,1.2) .. (0,1.2);
\draw (0.3,0) .. controls (1.3,-0.1) and (0.8,-1.2) .. (0,-1.2);
%
%
%
\draw(-1.25,1.25)node{${\cal C}_x$};
\draw[arrows={-Stealth[width=1.6mm, length=1.8mm]}](-0.58,-1.055)--(-0.57,-1.062);
\draw[arrows={-Stealth[width=1.6mm, length=1.8mm]}](0.803,0.697)--(0.802,0.7);
\draw[arrows={-Stealth[width=1.6mm, length=1.8mm]}](-0.8485,0.8485)--(-0.8485-0.001,0.8485-0.001);
\end{scope}
 \draw [arrows={-Stealth[width=1.6mm, length=1.8mm]}, color=gray] (-2.5,0) -- (3.5,0) node[above]{$\quad{\rm Re}(z)$};
  \draw [arrows={-Stealth[width=1.6mm, length=1.8mm]}, color=gray] (0,-1.5) -- (0,1.8) node[right]{${\rm Im}(z)$};
%
\draw (0.433,0.25) arc (30:330:0.5); 
\draw (0.3,0) .. controls (0.55,0.0) and (0.5,0.15) .. (0.433,0.25);
\draw (0.3,0) .. controls (0.55,0.0) and (0.5,-0.15) .. (0.433,-0.25);
\draw[arrows={-Stealth[width=1.6mm, length=1.8mm]}]
  (-0.25,-0.433)--(-0.24,-0.44);
  \draw[arrows={-Stealth[width=1.6mm, length=1.8mm]}]
  (0.25,0.433)--(0.24,0.44);
\draw[arrows={-Stealth[width=1.6mm, length=1.8mm]}](0,0)--(0.3,0);
\draw(0,0)--(0.2,0);
  \draw (2.5,0)node{$\bullet$} ;
  \draw(2.5,-0.4)node{$z=1$};
\end{scope}
\begin{scope}[yshift=4cm]
  \draw(2.5,-1) arc (-90:90:1cm);
  \draw [arrows={-Stealth[width=1.6mm, length=1.8mm]}, color=gray] (-2.5,0) -- (4.0,0) node[below]{$\quad{\rm Re}(z)$};
    \draw [arrows={-Stealth[width=1.6mm, length=1.8mm]}, color=gray] (0,-1.5) -- (0,1.8) node[right]{${\rm Im}(z)$};
\draw[arrows={-Stealth[width=1.6mm, length=1.8mm]}](0,0)--(0.3,0);
\draw(0,0)--(0.2,0);
  \draw (2.5,0)node{$\bullet$} ;
  \draw(2.5,-0.4)node{$z=1$};
\draw  (0.3,0) .. controls (1.3,-0.0) and (1.5,-1) .. (2.5,-1);
\draw  (0.3,0) .. controls (1.3,0.0) and (1.5,1) .. (2.5,1);
\draw(2.2,1.4)node{${\cal C}_y$};
 \draw[arrows={-Stealth[width=1.6mm, length=1.8mm]}](1.4,-0.5)--(1.41,-0.51);
 \draw[arrows={-Stealth[width=1.6mm, length=1.8mm]}](3.28,-0.61)--(3.295,-0.59);
 \draw[arrows={-Stealth[width=1.6mm, length=1.8mm]}](2.51,1)--(2.5,1);
\end{scope}
  \draw[fill=white] (2.5,0) circle (0.5cm);
  \draw[white, fill=white] (2,0) circle (0.1cm);
 \draw [arrows={-Stealth[width=1.6mm, length=1.8mm]}, color=gray] (-2.5,0) -- (3.5,0) node[below]{$\quad{\rm Re}(z)$};
  \draw [arrows={-Stealth[width=1.6mm, length=1.8mm]}, color=gray] (0,-1.5) -- (0,1.8) node[right]{${\rm Im}(z)$};
%
\draw[arrows={-Stealth[width=1.6mm, length=1.8mm]}](0,0)--(0.3,0);
\draw(0,0)--(0.2,0);
\draw(0.2,0) .. controls (0.35,0) and (0.35,0.1) .. (0.5,0.1);
\draw(0.2,0) .. controls (0.35,0) and (0.35,-0.1) .. (0.5,-0.1);
 \draw(2.5-0.48,0.1) -- (0.5,0.1);
  \draw(2.5-0.48,-0.1) -- (0.5,-0.1);
%
  \draw[arrows={-Stealth[width=1.6mm, length=1.8mm]}](1.35,0.1)--(1.34,0.1);
  \draw[arrows={-Stealth[width=1.6mm, length=1.8mm]}](1.38,-0.1)--(1.39,-0.1);
  \draw[arrows={-Stealth[width=1.6mm, length=1.8mm]}]
  (2.25,-0.433)--(2.26,-0.44);
  \draw[arrows={-Stealth[width=1.6mm, length=1.8mm]}]
  (2.75,0.433)--(2.74,0.44);
  \draw (2.5,0)node{$\bullet$} ;
  \draw(2.5,-0.8)node{$z=1$};
\end{tikzpicture}
\end{center}
\caption{\textit{The loops ${\cal C}_x$ and ${\cal C}_y$ around $z=0$ and $z=1$ anchored at the origin (upper half) and their homotopy deformation to infinitesimal circles along with straight paths between zero and one in case of ${\cal C}_y$ (lower half) \cite{Ihara:1990}. Strictly speaking, all the contours start and end at the tangential base point from 0 to 1 as indicated by the arrows at the origin pointing along the positive real axis. The straight line portions of the path in the lower-right panel should be viewed as running along the real axis between 0 and 1; they have been slightly separated for visual convenience.}
\label{figsphere}}
\end{figure}
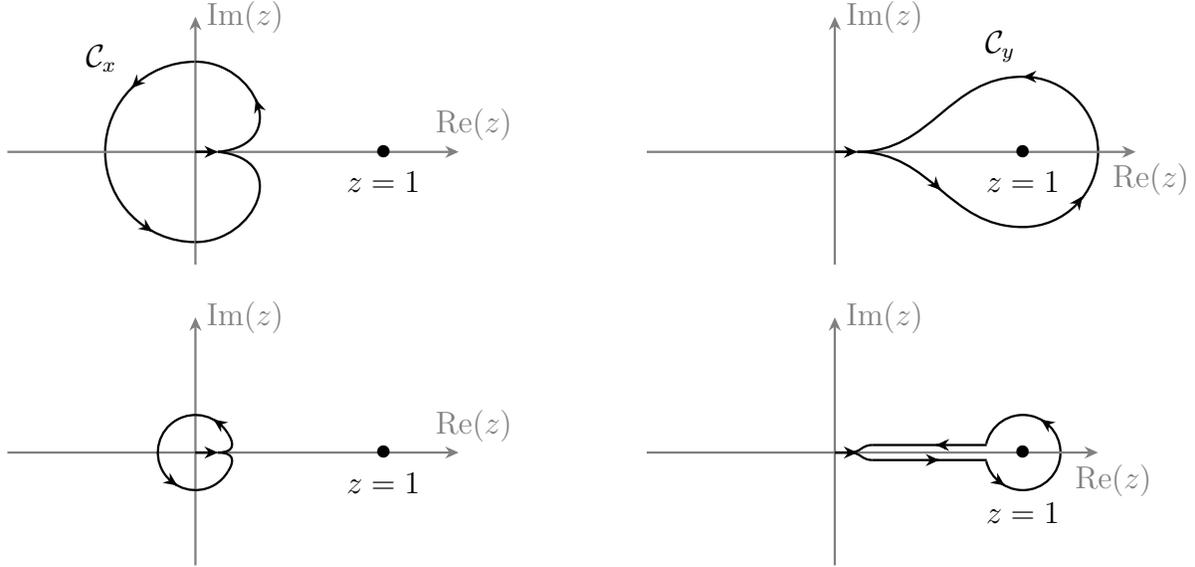

The zeta generators in genus zero are given by the Ihara derivations $D_{g_w}$ associated to the polynomials $g_w$, which act on the free Lie algebra ${\rm Lie}[x,y]$ via 
\beq
D_{g_w}(x) = 0 \, , \ \ \ \
D_{g_w}(y) = [y,g_w(x,y)]\, ; 
\label{kzapp.03}
\eeq
they can be interpreted as the coefficient of $\zeta_w$ in the holonomies of $J(x,y;z)$ w.r.t.\ the loops around $z=0$ and $z=1$, respectively. More specifically, (\ref{kzapp.03}) extracts the {\em linearized} monodromy of the loops ${\cal C}_x$ and ${\cal C}_y$ around $z=0$ and $z=1$ anchored at the 
tangential base point from 0 to 1 as drawn in Figure \ref{figsphere}, where only the first power of $2\pi i$ is retained:
\begin{align}
D_{g_w}(x) &= - {\rm Pexp} \bigg( \int_{{\cal C}_x} J(x,y;z) \bigg) \, \big|_{2\pi i \zeta_w}\, , \notag \\
D_{g_w}(y) &= - {\rm Pexp} \bigg( \int_{{\cal C}_y} J(x,y;z) \bigg) \, \big|_{2\pi i \zeta_w}
\label{kzapp.04}
\end{align}
Equivalence to (\ref{kzapp.03}) can be seen as follows:
\begin{itemize}
\item The path-ordered exponentials of $J(x,y;z)$ associated with the infinitesimal circles around $0$ and $1$ in counter-clockwise orientation are given by $e^{2 \pi i x}$ and $e^{-2 \pi i y}$, respectively.
\item Since ${\cal C}_x$ is homotopic to an infinitesimal circle around $z=0$, we have
\beq
{\rm Pexp} \bigg( \int_{{\cal C}_x} J(x,y;z) \bigg) = e^{2 \pi i x}
\label{kzapp.05}
\eeq
which does not contain any odd Riemann zeta values, thereby reproducing $D_{g_w}(x) =0$.
\item The path ${\cal C}_y$ is homotopic to the composition of the path $(0,1)$ followed by an infinitesimal circle around $z=1$ and the inverse path $(1,0)$ as seen in the lower-right panel of Figure \ref{figsphere}. Hence, the path-ordered exponential can be decomposed into
\beq
{\rm Pexp} \bigg( \int_{{\cal C}_y} J(x,y;z) \bigg) = \Phi(x,y)^{-1} e^{-2 \pi i y} \Phi(x,y) \,.
\label{kzapp.06}
\eeq
By the conventions (\ref{kzapp.02}) for path-ordered exponentials, the last segment $(1,0)$ of the deformed path ${\cal C}_y$ translates into the leftmost factor $\Phi(x,y)^{-1}$.
\item Extracting the coefficient of $\zeta_w$ from (\ref{kzapp.06}) leads to
\beq
{\rm Pexp} \bigg( \int_{{\cal C}_y} J(x,y;z) \bigg) \, \big|_{\zeta_w }= e^{-2 \pi i y}  g_w(x,y) - g_w(x,y) e^{-2 \pi i y}  
\label{kzapp.07}
\eeq
which upon linearization in $2\pi i$ reduces to $-2\pi i [y, g_w(x,y)]$ and
reproduces the action of $D_{g_w}$ on $y$ in (\ref{kzapp.03}).
\end{itemize}
We emphasize that it will be the formulation (\ref{kzapp.04}) of zeta generators in terms of linearized monodromies which generalizes from genus zero to genus one.

\subsection{Degenerating the KZB connection}
\label{app:deg.2}

In the same way as the (modified) KZ connection (\ref{kzapp.01}) can be used to generate multiple polylogarithms in genus zero, the Brown--Levin formulation of elliptic polylogarithms in genus one \cite{BrownLev} is based on the KZB connection
\begin{align}
J_{\rm KZB}(A,B;z|\tau) &\coloneqq   \, {\rm ad}_B F(z,{\rm ad}_B| \tau)A \,dz\, , \ \ \ \
F(z,\alpha | \tau) \coloneqq \frac{\theta'_1(0|\tau) \theta_1(z+\alpha|\tau)}{\theta_1(z|\tau) \theta_1(\alpha|\tau)}\, ,
\label{kzapp.08} \\
\theta_1(z|\tau) &\coloneqq 2 q^{1/8} \sin(\pi z) \prod_{n=1}^{\infty} (1-q^n) (1-e^{2\pi i z} q^n) (1-e^{-2\pi i z} q^n) \, , \ \ \ \  q\coloneqq e^{2\pi i \tau}\, ,
\notag
\end{align}
where $F(z,\alpha | \tau)$ is known as the Kronecker--Eisenstein series.
The modular parameter $\tau \in \mathbb H$ of the torus takes values in the upper half plane $\mathbb H\coloneqq \{ \tau \in \mathbb C\, , \ \Im \tau>0\}$, and $z,\alpha \in \mathbb C$ live on the universal cover of the torus $\mathbb C/(\mathbb Z+\tau \mathbb Z)$. The KZB connection $J_{\rm KZB}$ depends on non-commutative indeterminates $A,B$, and the adjoint actions of $B$ in $ {\rm ad}_B F(z,{\rm ad}_B| \tau)$ are performed after series expansion in the second argument of $F$. Note that the elliptic associators of \cite{KZB, EnriquezEllAss, Hain} are obtained from (regularized) path-ordered exponentials of (\ref{kzapp.08}), integrated over the homology cycles of the torus.

The degeneration $\tau \rightarrow i\infty$ of the Kronecker--Eisenstein series and its expansion coefficients w.r.t.\ the second argument $\alpha= {\rm ad}_B$ is well-known to yield \cite{ZagierF}
\beq
\lim_{\tau \rightarrow i\infty}F(z,\alpha | \tau) = \frac{1}{\alpha}+ \pi \cot(\pi z) 
-2 \sum_{n=1}^{\infty} \alpha^{2n-1} \zeta_{2n} \, .
\label{kzapp.09}
\eeq
The limit $\tau \rightarrow i\infty$ degenerates the torus to a nodal sphere. In the coordinate $\sigma \coloneqq e^{2\pi i z}$ of the nodal sphere, the pinched homology cycle of the degenerate torus translates into the identification of the points $\sigma=0$ with $\sigma=\infty$. Based on $dz = \frac{d\sigma}{2\pi i \sigma}$ and (\ref{kzapp.09}), the degeneration of the KZB connection (\ref{kzapp.08}) is readily found to be
\beq
\lim_{\tau \rightarrow i\infty} J_{\rm KZB}(A,B;z|\tau) = \bigg\{
A+2\pi i\bigg({-}\frac{1}{2}+ \frac{\sigma}{\sigma-1} \bigg) [B,A]
+ \sum_{n=1}^{\infty} (2\pi i)^{2n} \frac{{\rm B}_{2n}}{(2n)!} {\rm ad}_B^{2n}A
\bigg\} \frac{d\sigma}{2\pi i \sigma} \, .
\label{kzapp.10}
\eeq
In order to make contact with the images $t_{01}$ and $t_{12}$ in (\ref{g01.01}) of the genus-zero generators $x,y$, we redefine the non-commutative $A,B$ in (\ref{kzapp.08}) in terms of the generators $a,b$ introduced in section \ref{sec:sigman.0}
\beq
A = - 2\pi i a \, , \ \ \ \ B= \frac{b}{2\pi i}
\label{kzapp.11}
\eeq
and obtain the (modified) KZ connection (\ref{kzapp.01}) at $x=t_{01}$ and $y= -t_{12}$ from the degeneration (\ref{kzapp.10}),
\beq
\lim_{\tau \rightarrow i\infty} J_{\rm KZB}(A,B;z|\tau) =
\bigg( \frac{t_{01}}{\sigma}+ \frac{-t_{12}}{1-\sigma} \bigg) d\sigma= 
J(t_{01},-t_{12};\sigma)\, .
\label{kzapp.12}
\eeq

\begin{figure}[ht]
	\begin{center}
\begin{tikzpicture}[line width=0.30mm]
\begin{scope}[yshift=-2cm, xshift=-8cm]
  \draw [arrows={-Stealth[width=1.6mm, length=1.8mm]}, color=gray] (-1.5,0) -- (3.5,0) node[above]{$\quad{\rm Re}(z)$};
  \draw [arrows={-Stealth[width=1.6mm, length=1.8mm]}, color=gray] (0,-1.8) -- (0,1.8)node[left]{${\rm Im}(z)$};
  \draw(0.8,1.2)node{$\bullet$}node[above]{$\tau/2$};
\draw(3.3,1.2)node{$\bullet$}node[above]{$\tau/2+1$};
\draw(-0.8,-1.2)node{$\bullet$}node[below]{$-\tau/2$};
\draw(1.7,-1.2)node{$\bullet$}node[below]{$-\tau/2+1$};
\draw(0.8,1.2)--(3.3,1.2);
\draw(-0.8,-1.2)-- (1.7,-1.2);
\draw(1.7,-1.2)--(3.3,1.2);
\draw(-0.8,-1.2)--(0.8,1.2);
\draw(0.45,-1.2)node[rotate=-20]{$|\!|$};
\draw(2.05,1.2)node[rotate=-20]{$|\!|$};
\draw(-0.2,-0.3)node[rotate=70]{$|\!|$};
\draw(2.3,-0.3)node[rotate=70]{$|\!|$};
\draw[arrows={-Stealth[width=1.6mm, length=1.8mm]}](0,0)--(0.3,0);
\draw(0,0)--(0.2,0);
\draw[arrows={-Stealth[width=1.6mm, length=1.8mm]}](2.3,0)--(2.4,0);
\draw(2.5,0)--(2.3,0);
\draw(0.2,0) .. controls (0.35,0) and (0.35,0.1) .. (0.5,0.1);
\draw(2.3,0) .. controls (2.15,0) and (2.15,0.1) .. (2.0,0.1);
\draw(0.5,0.1) -- (2.0,0.1);
  \draw[arrows={-Stealth[width=1.6mm, length=1.8mm]}](1.35,0.1)--(1.351,0.1);
\draw(1.57,0.5)node{A-cycle};  
\draw(2.5,0)node{$\bullet$};
\draw(3.2,-0.4)node{$z=1$};
\end{scope}
  \draw[fill=white] (0,0) circle (0.5cm);
  \draw[white, fill=white] (0.5,0) circle (0.1cm);
  \draw(2.0,0.1) -- (0.48,0.1);
  \draw(2.0,-0.1) -- (0.48,-0.1);
\draw[arrows={-Stealth[width=1.6mm, length=1.8mm]}](2.5,0)--(2.2,0);
\draw(2.5,0)--(2.3,0);
\draw(2.3,0) .. controls (2.15,0) and (2.15,0.1) .. (2,0.1);
\draw(2.3,0) .. controls (2.15,0) and (2.15,-0.1) .. (2,-0.1);
%
  \draw[arrows={-Stealth[width=1.6mm, length=1.8mm]}](1.07,0.1)--(1.06,0.1);
  \draw[arrows={-Stealth[width=1.6mm, length=1.8mm]}](1.24,-0.1)--(1.25,-0.1);
  \draw[arrows={-Stealth[width=1.6mm, length=1.8mm]}](0.25,0.433)--(0.24,0.44);
  \draw[arrows={-Stealth[width=1.6mm, length=1.8mm]}](-0.25,-0.433)--(-0.24,-0.44);
  \draw (2.5,0)node{$\bullet$} ;
  \draw(3.2,0.3)node{$\sigma=1$};
  \draw [arrows={-Stealth[width=1.6mm, length=1.8mm]}, color=gray] (-1.5,0) -- (3.5,0) node[below]{$\quad{\rm Re}(\sigma)$};
  \draw [arrows={-Stealth[width=1.6mm, length=1.8mm]}, color=gray] (0,-1.5) -- (0,1.8)node[right]{${\rm Im}(\sigma)$};
\draw(-6,0) .. controls (-6,1) and (-4,1) .. (-1.8,1); 
  \draw[arrows={-Stealth[width=1.6mm, length=1.8mm]}](-1.8,1)--(-1.79,1);
  \draw(-3.2,1.3)node{$\tau \rightarrow i \infty$};
   \draw(-3.2,0.7)node{$\sigma = e^{2\pi i z}$};
\draw(2,-1) .. controls (2.5,-1.7) and (2.5,-2.3) .. (2,-3);
 \draw[arrows={-Stealth[width=1.6mm, length=1.8mm]}](2,-3)--(1.95,-3.07);
\draw(3.4,-2)node{$\eta=1-\sigma$};
  \begin{scope}[yshift=-4cm]
 \draw[fill=white] (2.5,0) circle (0.5cm);
  \draw[white, fill=white] (2,0) circle (0.1cm);
  \draw [arrows={-Stealth[width=1.6mm, length=1.8mm]}, color=gray] (-1.5,0) -- (3.5,0) node[below]{$\quad{\rm Re}(\eta)$};
  \draw [arrows={-Stealth[width=1.6mm, length=1.8mm]}, color=gray] (0,-1.5) -- (0,1.8) node[right]{${\rm Im}(\eta)$};  
%
\draw[arrows={-Stealth[width=1.6mm, length=1.8mm]}](0,0)--(0.3,0);
\draw(0,0)--(0.2,0);
\draw(0.2,0) .. controls (0.35,0) and (0.35,0.1) .. (0.5,0.1);
\draw(0.2,0) .. controls (0.35,0) and (0.35,-0.1) .. (0.5,-0.1);
 \draw(2.5-0.48,0.1) -- (0.5,0.1);
  \draw(2.5-0.48,-0.1) -- (0.5,-0.1); 
  \draw[arrows={-Stealth[width=1.6mm, length=1.8mm]}](1.35,0.1)--(1.34,0.1);
  \draw[arrows={-Stealth[width=1.6mm, length=1.8mm]}](1.38,-0.1)--(1.39,-0.1);
  \draw[arrows={-Stealth[width=1.6mm, length=1.8mm]}]
  (2.25,-0.433)--(2.26,-0.44);
  \draw[arrows={-Stealth[width=1.6mm, length=1.8mm]}]
  (2.75,0.433)--(2.74,0.44);
  \draw (2.5,0)node{$\bullet$} ;
  \draw(2.5,-0.8)node{$\eta=1$};
  \end{scope}
\end{tikzpicture}
\end{center}
\caption{\textit{The degeneration $\tau \rightarrow i \infty$ of the torus with coordinate $z$ (left panel) yields a nodal sphere, where the image of the A-cycle connecting $z=0$ with $z=1$ is drawn in two different coordinates $\sigma$ and $\eta$ (right panel). The image of the A-cycle in the $\eta$ coordinate (lower-right panel) matches the deformation of the
loop ${\cal C}_y$ around $z=1$ in Figure \ref{figsphere}. Similar to Figure~\ref{figsphere}, the straight line portions of all the paths should be viewed as running along the real axis between 0 and 1; they have been slightly separated for visual convenience.}}
\label{figtorus}
\end{figure}
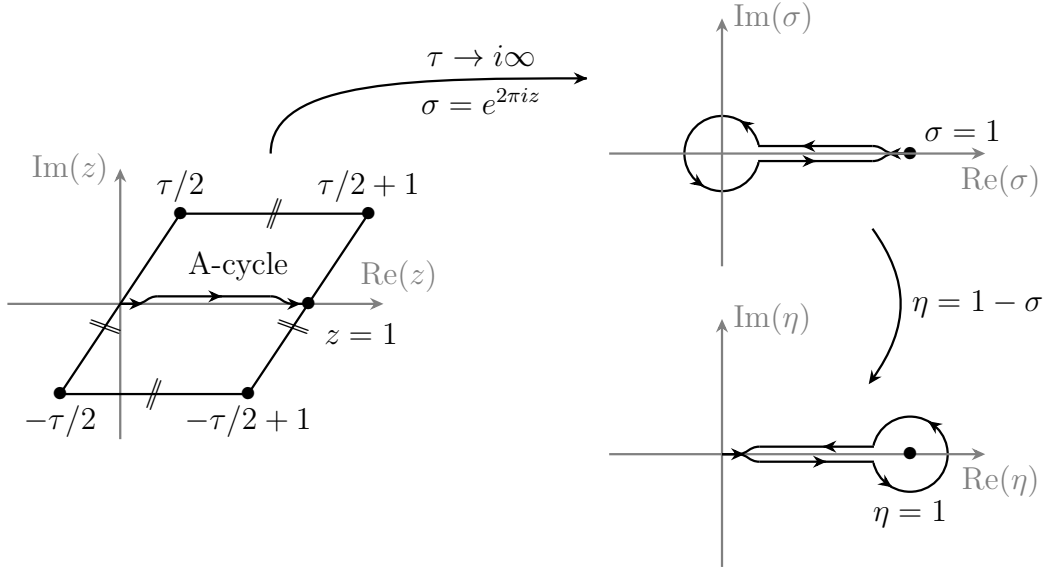

\subsection{Link between genus zero and genus one}
\label{app:deg.3}

The non-commutative arguments $t_{01},t_{12}$ obtained in the comparison (\ref{kzapp.12}) of KZ and KZB connections do not yet line up with (\ref{sec4eq.2}) and differ by a swap of $x$ and $y$. This can be fixed by an additional change of coordinates to $\eta=1-\sigma$ in the
degeneration of the KZB connection which is in fact necessary to map the origin $z=0$ of the torus to the origin $\eta=0$ of the nodal sphere (as opposed to $\sigma=1$). In this way, the homotopy deformation of the contour ${\cal C}_y$ of Figure \ref{figsphere} producing the action of zeta generators in genus zero is the image of the $A$-cycle of the torus $z \in (0,1)$ under the change of variables from $z$ via $\sigma=e^{2\pi i z}$ to $\eta=1-\sigma$, see Figure~\ref{figtorus}. Similar homotopy deformations of paths together with the degeneration (\ref{kzapp.12}) of the KZB connection were used by Enriquez to express the limit $\tau \rightarrow i \infty$ of elliptic associators in terms of $\Phi_{\rm KZ}$ \cite{Enriquez:Emzv}.

With the degenerate KZB connection 
in the coordinate $\eta=1-\sigma$ 
\beq
\lim_{\tau \rightarrow i\infty} J_{\rm KZB}(A,B;z|\tau) =
\bigg( \frac{t_{12}}{\eta}+ \frac{-t_{01}}{1-\eta} \bigg) d\eta= 
J(t_{12},-t_{01};\eta)\, ,
\label{kzapp.120}
\eeq
we obtain the factor 
\beq
g_w(t_{12},-t_{01}) = {\rm Pexp} \bigg( \int^1_0  \bigg[ \frac{t_{12}}{\eta} + \frac{-t_{01}}{1-\eta} \bigg] d\eta \bigg) \, \big|_{\zeta_w}
\label{kzapp.13}
\eeq
in the action (\ref{g01.11tau}) of genus-one zeta generators on $t_{01}$, in direct analogy with (\ref{kzapp.01}) in genus zero. Moreover, the realization  (\ref{kzapp.04}) of $D_{g_w}(y)$ in
genus zero generalizes to
\begin{align}
\tau_w(t_{01}) &= - {\rm Pexp} \bigg( \int_{{\cal C}_y}  \bigg[ \frac{t_{12}}{\eta} + \frac{-t_{01}}{1-\eta} \bigg]d\eta \bigg) \, \big|_{2\pi i \zeta_w} \, ,\label{kzapp.14} \\
&= - \lim_{\tau \rightarrow i \infty} {\rm Pexp} \bigg( \int_{0}^1 J_{\rm KZB}(A,B;z|\tau)
\bigg) \, \big|_{2\pi i \zeta_w} \, ,
\notag
\end{align}
i.e.\ the interpretation as a linearized monodromy passes through from genus zero to genus one. The loop ${\cal C}_y$ anchored at the 
tangential base point from $\eta=0$ to 1 on
the sphere around the point $\eta=1$ descends from the $A$-cycle $z \in (0,1)$ of the torus. The other part $\tau_w(t_{12})=0$ of the action (\ref{g01.11tau}) of genus-one zeta generators in turn follows from a loop around the origin of both the sphere ($\eta=0$) and the torus ($z=0$) which can be contracted to an infinitesimal circle and does not produce any odd zeta values through its periods, see (\ref{kzapp.05}).

In summary, this appendix derived the close analogy between the actions (\ref{kzapp.03}) and (\ref{g01.11tau}) of zeta generators in genus zero and one and justified the morphism \eqref{sec4eq.2} by comparing (i) the underlying connections of KZ- and KZB-type in the degeneration  of the torus to a nodal sphere and (ii) integration contours on the respective surfaces (loops around marked points and the pinched homology cycle of the degenerate torus).

\section{Brief recap of elliptic MZVs}
\label{app:eMZV}

We shall here review the definition and basic properties of elliptic MZVs, with particular emphasis on their interplay with the algebra $\mathfrak{u}$ of derivations of ${\rm Lie}[a,b]$ in Definition \ref{defmfu}.

Elliptic MZVs were introduced in the work of Enriquez \cite{EnriquezEllAss, Enriquez:Emzv} as coefficients of iterated integrals over KZB connections. In particular, Enriquez' $A$- $B$-elliptic MZVs are expansion coefficients of elliptic associators (obtained from the regularized $A$- and $B$-cycle holonomies of the KZB connection (\ref{kzapp.08}) \cite{KZB, EnriquezEllAss, Hain}) in the same way as multizetas at genus zero are coefficients of the KZ associator \cite{LeMura}. We will restrict the review of this appendix to the discussion of $A$-elliptic MZVs which generate $B$-elliptic MZVs through their modular S transformation $\tau \rightarrow -\frac{1}{\tau}$, see \cite{Broedel:2018izr, Zerbini:2018sox, Zerbini:2018hgs} for properties specific to $B$-elliptic MZVs. For further background on elliptic MZVs, the reader is for instance referred to \cite{Matthes:Thesis} for a comprehensive overview, to \cite{Broedel:2015hia, LMSonEMZV} for their algebraic and differential relations as well as the structure of the ring they generate, and to \cite{Broedel:2014vla, Broedel:2017jdo} for first applications to string amplitudes.

\subsection{Definition and basic properties}
\label{app:eMZV.1}

A convenient definition of individual $A$-elliptic MZVs involves the expansion coefficients $g^{(n)}(z| \tau)$ of the meromorphic Kronecker--Eisenstein series $F(z,\alpha | \tau)$ in (\ref{kzapp.08}) in its second argument $\alpha$,
\beq
F(z,\alpha | \tau) = \sum_{n=0}^{\infty} \alpha^{n-1} g^{(n)}(z| \tau)\, ,
\label{emzv.1}
\eeq
where $g^{(0)}(z| \tau) = 1$. The non-constant $g^{(n)}(z| \tau)$ at $n\geq 1$ are meromorphic in both $z$ in the universal cover of the torus and $\tau$ in the upper half-plane; their generating series in (\ref{kzapp.08}) yields explicit theta-function representations of $g^{(n)}(z| \tau)$ at arbitrary $n\in \mathbb N$ such as $g^{(1)}(z| \tau) = \partial_z \log \theta_1(z|\tau)$.

\begin{dfn}
\label{dfn:emzv} 
$A$-elliptic MZVs of length $r\geq 0$ are defined as iterated integrals of the Kronecker--Eisenstein coefficients in (\ref{emzv.1}) over points $z_i$ on the $A$-cycle \cite{Enriquez:Emzv},
\beq
\omega(n_1,\ldots,n_r|\tau) = \int^1_0
dz_{1}\, g^{(n_{1})}(z_{1}| \tau)  \int_{z_1}^1 dz_{2}\, g^{(n_{2})}(z_{2}| \tau) \ldots \int^1_{z_{r-1}}
dz_r\, g^{(n_r)}(z_r| \tau)\, .
\label{emzv.2}
\eeq
Following the enumeration conventions of  \cite{Broedel:2014vla}, they are specified by integers $n_1,\ldots,n_r \geq 0$ that specify their integration kernels. The endpoint divergences of the integrals (\ref{emzv.2}) with $n_1 = 1$ and/or $n_r=1$ are regularized according to the prescription in \cite{Broedel:2014vla} that is tailored to preserve the properties (i) and (ii) below of the convergent cases.
\end{dfn}

\begin{prop}[see \cite{Enriquez:Emzv, Broedel:2014vla, Broedel:2015hia}]
\label{prop.emzv}
The $A$-elliptic MZVs of Definition \ref{dfn:emzv} exhibit the following basic properties, where $u,v \in \{ 0,1,2,\ldots\}^\times$ are words in the integer entries $n_1,n_2,\ldots$ of (\ref{emzv.2}):
\begin{itemize}
\item[(i)] shuffle relations identical to those of multizetas in (\ref{appMZV.04}), 
\beq
\omega(u |\tau) \omega(v|\tau)
= \omega(u \shuffle v|\tau)\, ,
\label{emzv.3}
\eeq
see (\ref{appMZV.05}) for the shuffle product~$\shuffle$,
\item[(ii)] reflection relations
\beq
\omega(n_1,n_2,\ldots,n_r |\tau) 
= (-1)^{n_1+n_2+\ldots+n_r}\omega(n_r,\ldots,n_2,n_1 |\tau)\, ,
\label{emzv.4}
\eeq
\item[(iii)] Fourier expansions in non-negative powers of $q = e^{2\pi i \tau}$ with $\mathbb Q[(2\pi i)^{-1}]$-linear combinations $c_n(u)$ of multizetas as coefficients
\beq
\omega(u |\tau)
= \sum_{n=0}^\infty c_n(u) q^n\, ,
\label{emzv.5}
\eeq
\item[(iv)] The simplest non-constant $A$-elliptic MZVs occur at length $r=2$ with odd $n_1{+}n_2$ since (for integer $k,\ell \geq 0$)
\begin{align}
\omega(2k{+}1 |\tau) &= 0 \, , & \omega(2k{+}1 ,2\ell{+}1|\tau) &= 0 \, , \label{emzv.6} \\
\omega(2k |\tau) &= -2 \zeta_{2k} \, , &\omega(2k,2\ell|\tau) &= 2 \zeta_{2k}  \zeta_{2\ell} \, .
\notag
\end{align}
\end{itemize}
\end{prop}

\noindent {\bf Proof.} For $A$-elliptic MZVs with $n_1,n_r \neq 1$, i.e.\ convergent integrals (\ref{emzv.2}), (i) and (ii) readily follow from general properties of iterated integrals and the alternating parity $g^{(n)}(-z|\tau)=(-1)^n g^{(n)}(z|\tau)$ of the integration kernels. The extension of (i) and (ii) to divergent cases with $n_1=1$ and/or $n_r=1$ is ensured by the regularization of endpoint divergences according to \cite{Broedel:2014vla}. The Fourier expansion of the Kronecker--Eisenstein series \cite{ZagierF} and its coefficients $g^{(n)}(z| \tau)$ \cite{Broedel:2014vla} implies the statement of (iii) with coefficients $c_n(u) \in \mathbb C$; the stronger statement that $c_n(u)$ are $\mathbb Q[(2\pi i)^{-1}]$-linear combinations of multizetas can for instance be established from the degeneration formulae (\ref{kzapp.09}), (\ref{kzapp.10}) for the integration kernels or the degeneration formula of \cite{Enriquez:Emzv} for elliptic associators. Finally, the expressions of (iv) for $\omega(n_1|\tau)$ readily follow from direct integration of (\ref{emzv.2}) for $n_1\neq 1$ together with our regularization prescription for $n_1=1$; the results of (\ref{emzv.6}) for $\omega(n_1,n_2|\tau)$ with $n_1 {+}n_2$ even in turn are simple consequences of (i), (ii) and the expressions for $\omega(n_1|\tau)$.
\qed

\subsection{Elliptic MZVs and geometric derivations}
\label{app:eMZV.2}

As will be reviewed in this section, the $\tau$-dependence of $A$-elliptic MZVs is governed by the algebra $\mathfrak{u}$ of Tsunogai derivations $\epsilon_k$ and their action on the free Lie-algebra in two generators $a,b$, see section \ref{sec:sigman.0}. This link arises from the generating series of $A$-elliptic MZVs 
\beq
\mathfrak{A}_{a,b}(\tau) = e^{i\pi [a,b]} \sum_{r=0}^\infty (-1)^r \sum_{n_1,\ldots,n_r=0}^\infty
\omega(n_1,n_2,\ldots,n_r|\tau) \, {\rm ad}_a^{n_r}(b)\ldots {\rm ad}_a^{n_2}(b) {\rm ad}_a^{n_1}(b)\, ,
\label{emzv.8}
\eeq
known as the $A$-elliptic associator $\mathfrak{A}_{a,b}(\tau)$, and the appearance of holomorphic Eisenstein series ${\rm G}_{k}(\tau)$ in its $\tau$-derivative, with normalization conventions\footnote{The $q$-series in (\ref{emzv.7}) at  $k\geq 2$ line up with the definition
\[
{\rm G}_{k}(\tau) = \sum_{m,n \in \mathbb Z \atop{(m,n)\neq (0,0)}} \frac{1}{(m\tau+n)^k}
\]
of holomorphic Eisenstein series when employing the Eisenstein summation prescription \cite{DHoker:2022dxx} in the conditionally convergent case with $k=2$.}
\beq
{\rm G}_{k}(\tau) = \left\{
\begin{array}{cl}
\displaystyle 2 \zeta_k + \frac{2(2\pi i)^k}{(k{-}1)!} \sum_{m,n=1}m^{k-1} q^{mn}&: \ k\geq 2 \ {\rm even}\, , \\
 -1 &: \ k=0 \, .
\end{array}
\right.
\label{emzv.7}
\eeq

\begin{thm} [Enriquez \cite{EnriquezEllAss, Enriquez:Emzv}]
\label{thm:emzv}
 The $A$-elliptic associator $\mathfrak{A}_{a,b}(\tau)$ obeys the differential equation,
\beq
2\pi i \frac{\partial}{\partial \tau} \mathfrak{A}_{a,b}(\tau) =  \sum_{\ell=0}^\infty (1{-}2\ell) {\rm G}_{2\ell}(\tau) \epsilon_{2\ell}
\mathfrak{A}_{a,b}(\tau)\, ,
\label{emzv.9}
\eeq
with Tsunogai derivations acting
on the associator via
(\ref{epsaction}), (\ref{ep2ionb}) and their Leibniz property.
\end{thm}

The factor of $ e^{i\pi [a,b]} $ in the expansion (\ref{emzv.8}) of $\mathfrak{A}_{a,b}(\tau)$ ensures that regularized $A$-elliptic MZVs obey the reflection property (\ref{emzv.4}) and does not alter the differential equation (\ref{emzv.9}) since Tsunogai derivations $\epsilon_{2\ell}$ all annihilate the commutator $[a,b]$.

\begin{cor}[Broedel, Matthes, Schlotterer \cite{Broedel:2015hia}]
\label{cor:emzv}
$A$-elliptic MZVs are given by (tangential-base-point regularized \cite{Brown:mmv}) linear combinations of iterated integrals over holomorphic Eisenstein series (\ref{emzv.7}) with $\mathbb Q[(2\pi i)^{-1}]$ linear combinations of multizetas as coefficients. 
\end{cor}
\noindent {\bf Proof.} 
The discussion in section 4.3 of \cite{Broedel:2015hia} (with $x,y$ renamed to $a,b$) is equivalent to writing the $A$-elliptic associator as a path-ordered exponential in the conventions of (\ref{kzapp.02}),
\beq
\mathfrak{A}_{a,b}(\tau) =  {\rm Pexp}\bigg( \int^\tau_{i\infty} \frac{d \tau'}{2\pi i}\sum_{\ell=0}^\infty (1{-}2\ell) {\rm G}_{2\ell}(\tau') \epsilon_{2\ell} \bigg)
\mathfrak{A}_{a,b}(i\infty) \, .
\label{emzv.12}
\eeq
The initial value $\mathfrak{A}_{a,b}(i\infty)$ is obtained from the Drinfeld associator as spelt out in section 4.5 of \cite{KZB}, section 1.2 of \cite{Enriquez:Emzv} and section 2.3 of \cite{Broedel:2015hia}.
\qed

\begin{rmk}
$A$-elliptic MZVs only realize the proper subset of iterated Eisenstein integrals admitted by the differential equation (\ref{emzv.9}) and the relations in $\mathfrak{u}$, namely
\beq
{\rm ad}^{2k-1}_{\epsilon_0}(\epsilon_{2k}) = 0 \ \forall \ k \geq 1 \, , \ \ \ \ 
[\ep_2, \ep_{2\ell}] = 0 \ \forall \ \ell \geq 0
\label{emzv.11}
\eeq
and Pollack's relations in Remark \ref{Pollackrel}.
Each relation in the Lie-algebra $\mathfrak{u}$ of $\epsilon_{2\ell}$ leads to a dropout of a (shuffle-indecomposable) iterated Eisenstein integral in the expansion of (\ref{emzv.12}).

More precisely, the dropouts of iterated Eisenstein integrals from (\ref{emzv.12}) due to the closed-form relations among $\epsilon_{2\ell}$ in (\ref{emzv.11}) can be described as follows: 
\begin{itemize}
\item[(i)] 
The fact that $\epsilon_2$ is central in $\mathfrak{u}$, i.e.\ $[\ep_2, \ep_{2\ell}] = 0$, implies that all appearances of the quasi-modular holomorphic Eisenstein series ${\rm G}_2$ in the iterated Eisenstein integrals of $A$-elliptic MZVs can be reduced to polynomials in the depth-one integral $\int^\tau_{i\infty} d\tau' \, {\rm G}_2(\tau')$. The latter is related to the only shuffle-indecomposable divergent $A$-elliptic MZV (besides the vanishing $\omega(1|\tau)=0$) via $\omega(0,1|\tau) = \frac{i\pi}{2} -\frac{1}{2\pi i}  \int^\tau_{i\infty} d\tau' ( {\rm G}_2(\tau') - 2\zeta_2)$ \cite{Broedel:2015hia}.
\item[(ii)] 
The leftover iterated Eisenstein integrals over ${\rm G}_0 = -1$ and modular ${\rm G}_k$ at $k\geq 4$ in (\ref{emzv.12}) are expressible as Brown's iterated integrals over kernels $d\tau \, \tau^j {\rm G}_k(\tau)$ \cite{Brown:mmv}, where $k\geq 4$ and $0\leq j \leq k{-}2$ by virtue of ${\rm ad}^{k-1}_{\epsilon_0}(\epsilon_{k}) = 0$. The coefficients in these decompositions are $\mathbb Q[(2\pi i)^{\pm 1}]$ polynomials in $\tau$ that can for instance be determined from the generating functions in \cite{Broedel:2018izr}.
\end{itemize}
The Pollack relations of Remark \ref{Pollackrel} in turn induce further dropouts among Brown's iterated Eisenstein integrals over the above
$d\tau \, \tau^j {\rm G}_k(\tau)$
which are realized among elliptic MZVs. For instance, $[\ep_4,\ep_{10}] = 3 [\ep_6,\ep_8] $ interlocks the double integrals over ${\rm G}_{k_1}(\tau_1)  {\rm G}_{k_2}(\tau_2)$ with $(k_1,k_2) = (4,10), (6,8),(8,6), (10,4)$ such that only three $\mathbb Q$-linear combinations occur among $A$-elliptic MZVs. Similarly, the relation (\ref{tsurels2}) interlocks triple integrals (e.g.\ over ${\rm G}_4(\tau_1) {\rm G}_4(\tau_2) {\rm G}_8(\tau_3)$) with double integrals (e.g.\ over $\tau_1 {\rm G}_8(\tau_1)  {\rm G}_8(\tau_2)$).
\label{rmk:emzv}
\end{rmk}

\subsection{Elliptic MZVs versus modular graph forms}
\label{app:eMZV.3}

We conclude this section with comments on links between elliptic MZVs and the non-holomorphic modular graph forms \cite{DHoker:2015wxz, DHoker:2016mwo} of closed-string amplitudes. The main statement is that the one-to-one correspondence between relations in $\mathfrak{u}$ and dropouts of iterated Eisenstein integral from $A$-elliptic MZVs extends to modular graph forms.

An alternative way of constructing generating series of all $A$-elliptic MZVs akin to the associator is to expand the configuration-space integrals in genus-one open-string amplitudes in certain physical parameters \cite{Mafra:2019ddf, Mafra:2019xms}. The generating series of open-string integrals in the references over different numbers $n$ of marked points obey the differential equation (\ref{emzv.9}) of the $A$-elliptic associator with conjectural $(n{-}1)!\times (n{-}1)!$ matrix representations of the $\epsilon_{2\ell}$. In particular, the expansion of open-string integrals in terms of $A$-elliptic MZVs arises from the analogous matrix representations of the path-ordered exponential in (\ref{emzv.12}) and therefore features a dropout among the iterated Eisenstein integrals for each relation among the matrix representations of $\epsilon_{2\ell}$.

Similar generating series of configuration-space integrals adapted to closed-string (rather than open-string) amplitudes at genus one generate all modular graph forms \cite{Gerken:2019cxz, Gerken:2020yii}. The holomorphic differential equations of these closed-string integrals in $\tau$ again take the form of (\ref{emzv.9}) with almost identical conjectural matrix representations of $\epsilon_{2\ell}$ as seen in the open-string case. Hence, the matrix representation of the path-ordered exponential in (\ref{emzv.12}) -- along with the dropouts of iterated Eisenstein integrals for each relation in $\mathfrak{u}$ -- also features in generating series of modular graph forms.

More precisely, modular graph forms are modular combinations of iterated Eisenstein integrals, their complex conjugates, multizetas and rational functions of $\tau,\bar \tau$ \cite{DHoker:2015wxz, Broedel:2018izr, Gerken:2020yii, Dorigoni:2022npe}. Accordingly, their generating series combine the holomorphic path-ordered exponential (\ref{emzv.12}) with its complex conjugate (reversed in the concatenation order of the $\epsilon_{2\ell}$) and series in multizetas as in Brown's construction of equivariant iterated Eisenstein integrals \cite{Brown:2017qwo, Brown:2017qwo2}. The identification of Brown's construction with generating series of modular graph forms is described in \cite{Dorigoni:2022npe}, and the description of the series in multizetas via zeta generators can be found in the companion paper \cite{Dorigoni:2024oft} to this work.

In conclusion, the dropouts of elliptic MZVs due to relations in $\mathfrak{u}$ directly carry over to modular graph forms under the following assumption: the matrices in the generating series of open- and closed-string genus-one integrals in \cite{Mafra:2019ddf, Mafra:2019xms, Gerken:2019cxz, Gerken:2020yii} need to obey the entirety of all relations among geometric derivations $\epsilon_{2\ell}$ as supported by the arguments in section 4.5 of \cite{Mafra:2019xms}. Independently on this assumption, the generating series of equivariant iterated Eisenstein integrals of \cite{Brown:2017qwo, Brown:2017qwo2} exhibits one dropout of non-holomorphic modular form per relation in $\mathfrak{u}$.


\providecommand{\href}[2]{#2}\begingroup\raggedright\endgroup

\end{document}